\pdfoutput=1
\RequirePackage{ifpdf}
\ifpdf 
\documentclass[pdftex]{sigma}
\else
\documentclass{sigma}
\fi

\usepackage[all]{xy}
\usepackage{multirow}

\numberwithin{equation}{section}

\newtheorem{Theorem}{Theorem}[section]
\newtheorem{Corollary}[Theorem]{Corollary}
\newtheorem{Lemma}[Theorem]{Lemma}
\newtheorem{Proposition}[Theorem]{Proposition}
\newtheorem{Conjecture}[Theorem]{Conjecture}

\usepackage{mathrsfs}
\DeclareMathAlphabet{\mathscrbf}{OMS}{mdugm}{b}{n}

\usepackage[x11names]{xcolor}

\definecolor{burgundy}{rgb}{0.5, 0.0, 0.13}
\definecolor{carmine}{rgb}{0.59, 0.0, 0.09}
\definecolor{bostonuniversityred}{rgb}{0.8, 0.0, 0.0}
\definecolor{alizarin}{rgb}{0.82, 0.1, 0.26}

\def\BZ{\mathbb Z}
\def\BP{\mathbb P}
\def\BQ{\mathbb Q}
\def\BR{\mathbb R}
\def\BC{\mathbb C}

\def\BJ{\mathbb J}

\def\calA{\mathcal A}

\def\calD{\mathcal D}
\def\calL{\mathcal L}

\def\calT{\mathcal T}

\def\calP{\mathcal P}

\def\calT{\mathcal T}

\def\calH{\mathcal H}
\def\calP{\mathcal P}
\def\calO{\mathcal O}

\def\D{\Delta}

\def\l{\lambda}
\def\k{\kappa}

\def\la{\langle}
\def\ra{\rangle}

\def\ep{\epsilon}

\def\wt{\widetilde}
\def\wh{\widehat}

\def\SL{\mathrm{SL}}
\def\PSL{\mathrm{PSL}}
\def\GL{\mathrm{GL}}
\def\Sp{\mathrm{Sp}}
\def\pt{\partial}

\def\z{\zeta}

\def\Vol{\mathrm{Vol}}

\def\Phih{\widehat\Phi}
\def\Psih{\widehat\Psi}

\def\Psih{\widehat\Psi}

\def\e{\mathbf{e}}
\def\a{\alpha}
\def\p{\mathfrak p}
\def\l{\lambda}
\def\h{\hbar}
\def\={\;=\;}
\def\+{ + }

\def\R{\mathbb R}
\def\C{\mathbb C}
\def\Z{\mathbb Z}
\def\N{\mathbb N}
\def\Q{\mathbb Q}
\def\={\;=\;}
\def\+{ + }
\def\m{ - }
\def\si{\;\sim\;}
\def\i{^{-1}}
\def\e{\mathbf{e}}
\def\hb{\hbar}
\def\a{\alpha}
\def\b{\beta}
\def\l{\lambda}
\def\z{\zeta}
\def\G{\Gamma}
\def\s{\sigma}

\def\v{\varepsilon}
\def\W{\kappa} 

\def\sma#1#2#3#4{\bigl(\smallmatrix#1&#2\\#3&#4\endsmallmatrix\bigr)} 
\def\vsma#1#2#3#4{(\smallmatrix#1&#2\\#3&#4\endsmallmatrix)} 

\def\L{\text{{\rm Li}}}
\def\den{\mathrm{den}}
\def\num{\text{num}}
\def\O{O}
\def\o{o}

\def\G{\Gamma}
\def\J{\mathbf J} \def\J{J} 

\def\calJ{\mathcal J} 
\def\BJ{{\mathbf Q}} 
\def\Qbar{\overline{\Q}}
\def\V{\mathbf V}
\def\tV{v} 
\def\be{\begin{equation}}
\def\ee{\end{equation}}
\def\ve{\varepsilon}
\def\ssm{\smallsetminus}

\def\K#1{{K_{#1}}}

\def\Li{\mathrm{Li}}
\newcommand{\mb}{\mathbf}
\def\tcdot{\!\cdot\!}
\def\tiA{\wt A}
\def\tiB{\tilde{B}}
\def\opt{{\mathrm{opt}}}
\def\smo{{\mathrm{smooth}}}

\def\ul{\underline}
\def\g{\gamma}

\def\bH{\mathbf H}
\def\B{\mathcal B} \def\hB{\widehat\B}
\def\bK{\mathbf K}
\def\bx{\mathbf x}

\def\sQ{\mathcal{Q}}
\def\diag{\operatorname{\mathbf{diag}}} 
\def\bJ{J} \def\bJ{\mathbf{J}}
\def\sJ{\mathcal{J}}

\def\hol{\mathrm{hol}}
\def\conn{^\mathrm{conn}}
\def\sH{\mathcal H} 
\def\th{\theta}

\def\bPhi{\mathbf\Phi}
\def\bPhih{\mathbf{\Phih}}
\def\jp{\mathbf{j}}
\def\WW{\Lambda^2}

\def\vve{E}
\def\jt{\jp_\l} 
\def\jt{\jp_{\text{tw}}} 
\def\jt{\jp_{\textsf T}} 
\def\jt{\jp^{(\l)}}
\def\jt{\jp^{(\l)}}
\def\jt{\wt{\jp}}
\def\jhol{{\mathbf{j}}^{\mathrm{hol}}}
\def\EXACT{^{\text{\rm exact}}}
\def\RED{^{\text{red}}}

\def\FF{\mathfrak F}

\def\bJx{\bJ\RED} 
 
\def\bPhix{\bPhi\RED} \def\Phix{\bPhix}
\def\Wx{W\RED}
\def\calPx{\mathcal{P}\RED} 

\def\bPsi{\BJ}

\begin{document}


\newcommand{\arXivNumber}{2111.06645}

\renewcommand{\PaperNumber}{055}

\FirstPageHeading

\ShortArticleName{Knots, Perturbative Series and Quantum Modularity}

\ArticleName{Knots, Perturbative Series and Quantum Modularity}

\Author{Stavros GAROUFALIDIS~$^{\rm a}$ and Don ZAGIER~$^{\rm bc}$}

\AuthorNameForHeading{S. Garoufalidis and D. Zagier}

\Address{$^{\rm a)}$~International Center for Mathematics, Department of Mathematics,\\
\hphantom{$^{\rm a)}$}~Southern University of Science and Technology,
 Shenzhen, P.R.~China}
\EmailD{\href{mailto:stavros@mpim-bonn.mpg.de}{stavros@mpim-bonn.mpg.de}}
\URLaddressD{\url{http://people.mpim-bonn.mpg.de/stavros}}

\Address{$^{\rm b)}$~Max Planck Institute for Mathematics,
 Bonn, Germany}
\EmailD{\href{mailto:dbz@mpim-bonn.mpg.de}{dbz@mpim-bonn.mpg.de}}
\URLaddressD{\url{http://people.mpim-bonn.mpg.de/zagier}}

\Address{$^{\rm c)}$~International Centre for Theoretical Physics, Trieste, Italy}

\ArticleDates{Received April 25, 2023, in final form May 26, 2024; Published online June 24, 2024}

\Abstract{We introduce an invariant of a hyperbolic knot which is a map $\alpha\mapsto \boldsymbol{\Phi}_\alpha(h)$ from $\mathbb{Q}/\mathbb{Z}$ to matrices with entries in~$\overline{\mathbb{Q}}[[h]]$ and with rows and columns indexed by the boundary parabolic ${\rm SL}_2(\mathbb{C})$ representations of the fundamental group of the knot. These matrix invariants have a rich structure: (a)~their $(\sigma_0,\sigma_1)$ entry, where $\sigma_0$ is the trivial and~$\sigma_1$ the geometric representation, is the power series expansion of the Kashaev invariant of the knot around the root of unity ${\rm e}^{2\pi{\rm i} \alpha}$ as an element of the Habiro ring, and the remaining entries belong to generalized Habiro rings of number fields; (b)~the first column is given by the perturbative power series of Dimofte--Garoufalidis; (c)~the columns of $\boldsymbol{\Phi}$ are fundamental solutions of a linear $q$-difference equation; (d)~the matrix defines an ${\rm SL}_2(\mathbb{Z})$-cocycle $W_\gamma$ in matrix-valued functions on~$\mathbb{Q}$ that conjecturally extends to a smooth function on~$\mathbb{R}$ and even to holomorphic functions on suitable complex cut planes, lifting the factorially divergent series $\boldsymbol{\Phi}(h)$ to actual functions. The two invariants $\boldsymbol{\Phi}$ and $W_\gamma$ are related by a~refined quantum modularity conjecture which we illustrate in detail for the three simplest hyperbolic knots, the $4_1$, $5_2$ and $(-2,3,7)$ pretzel knots. This paper has two sequels, one giving a different realization of our invariant as a matrix of convergent $q$-series with integer coefficients and the other studying its Habiro-like arithmetic properties in more depth.}

\Keywords{quantum topology; knots; 3-manifolds; Jones polynomial; Kashaev invariant; volume conjecture; Chern--Simons theory; asymptotics; quantum modularity conjecture; quantum modular forms; hyperbolic 3-manifolds; dilogarithm; cocycles; $\mathrm{SL}_2(\mathbb{Z})$; denominators; Habiro-like functions; functions near $\mathbb{Q}$; Neumann--Zagier matrices; Nahm sums; $q$-holonomic modules}

\Classification{57N10; 57K16; 57K14; 57K10}

\part*{Part~0: Introduction and overview}
\pdfbookmark[1]{Part~0: Introduction and overview}{part0}

In this paper and the companion paper~\cite{GZ:qseries}, we will define and study
three different types of objects that can be associated to a hyperbolic knot:
\begin{itemize}\itemsep=0pt
\item periodic functions on $\Q$ with values in~$\Qbar$ with striking
arithmetic properties and belonging to a generalization of the Habiro ring;

\item divergent formal series in an infinitesimal variable~$h$, or more
precisely infinite collections of such power series, indexed by a rational number~$\a$
(here ``$h$'' is meant to remind one of Planck's constant and the perturbative expansions
of quantum field theory); and

\item $q$-series with integer coefficients, convergent in the unit disk
and also thought of via ${q={\rm e}^{2\pi{\rm i}\tau}}$ as holomorphic functions of a variable~$\tau$
in the upper half-plane.
\end{itemize}

The first of these generalizes the Kashaev invariant of the knot, while the
second and third correspond roughly to the two partition functions $Z(h)$ and
$\widehat{Z}(q)$ that are being studied in the ongoing program of Gukov et
al.\ \cite{Gukov:largeN, Gukov:BPS} for general 3-manifolds. We will study the first two
types of invariants in the present paper, and the functions of~$q$ or~$\tau$
in~\cite{GZ:qseries}. In all three cases, we will actually define a whole \emph{matrix}
of functions of the type described above, and in all three cases one of the central
questions will be the behavior of these functions under the action of the modular group
on the rational numbers or on the upper half-plane. Another key aspect is that each of
the three types of matrices constructed encodes the same information as the other two
and that all three can be interpreted as different realizations of the same abstract
object, a square matrix of ``functions-near-$\Q$'' that we believe is associated to
every hyperbolic knot, just as the different types of cohomology groups associated to
an algebraic variety over a number field, despite their very different properties,
are seen as different realizations of the same underlying ``motive''.

The starting point for our entire investigation is the Kashaev invariant of a knot and
the ``quantum modularity'' property for its Galois-equivariant extension that was
conjectured in~\cite{Za:QMF}. We will review these topics in detail in
Section~\ref{sec.QMCK}, but remind the reader briefly of the basic ingredients here. The
Kashaev invariant of a hyperbolic knot~$K$ is an element $\langle K\rangle_N$ of~$\Z\big[{\rm e}^{2\pi{\rm i}/N}\big]$ for every $N\in\N$ whose absolute value is conjectured to grow
exponentially like~${\rm e}^{cN}$, where~$c$ is~$1/2\pi$ times the hyperbolic volume of the
knot complement~$\mathbb S^3\ssm K$. This invariant can be extended to a function
$\J=J^{(K)}$ (we will omit the knot from the notation when it is fixed) from~$\Q/\Z$ to
$\Qbar$ by Galois equivariance. \big(This means that we write $\langle K\rangle_N$ as a
polynomial in~${\rm e}^{2\pi{\rm i}/N}$ with rational coefficients and define $J(a/N)$ for all $a$
prime to~$N$ as the same polynomial evaluated at~${\rm e}^{-2\pi{\rm i}a/N}$.\big) The quantum
modularity conjecture gives a formula for the ratio of the values of $J(X)$ and~$J(\g X)$ as an asymptotic series in $1/X$ as $X$ tends to infinity through integers,
or even through rational numbers with bounded denominator, where~\smash{$\g X=\frac{aX+b}{cX+d}$}
with~\smash{$\g=\sma abcd\in\SL_2(\Z)$}. The quantitative version of this conjecture is given
in equation~\eqref{QMCb} below, via a collection of well-defined formal power series
$\{\Phi_\a(h)\}_{\a\in\Q}$ with algebraic coefficients, but the conjecture can also be
visualized in a weaker qualitative form by comparing the graphs of $J(x)$ and of
$J(x)/J(\g x)$ as functions, as is done in the following figure (taken
from~\cite{Za:QMF}), which shows the plots of $ \log(J(x)) $
\begin{figure}[b!]\centering
 \includegraphics[height=5.3cm]{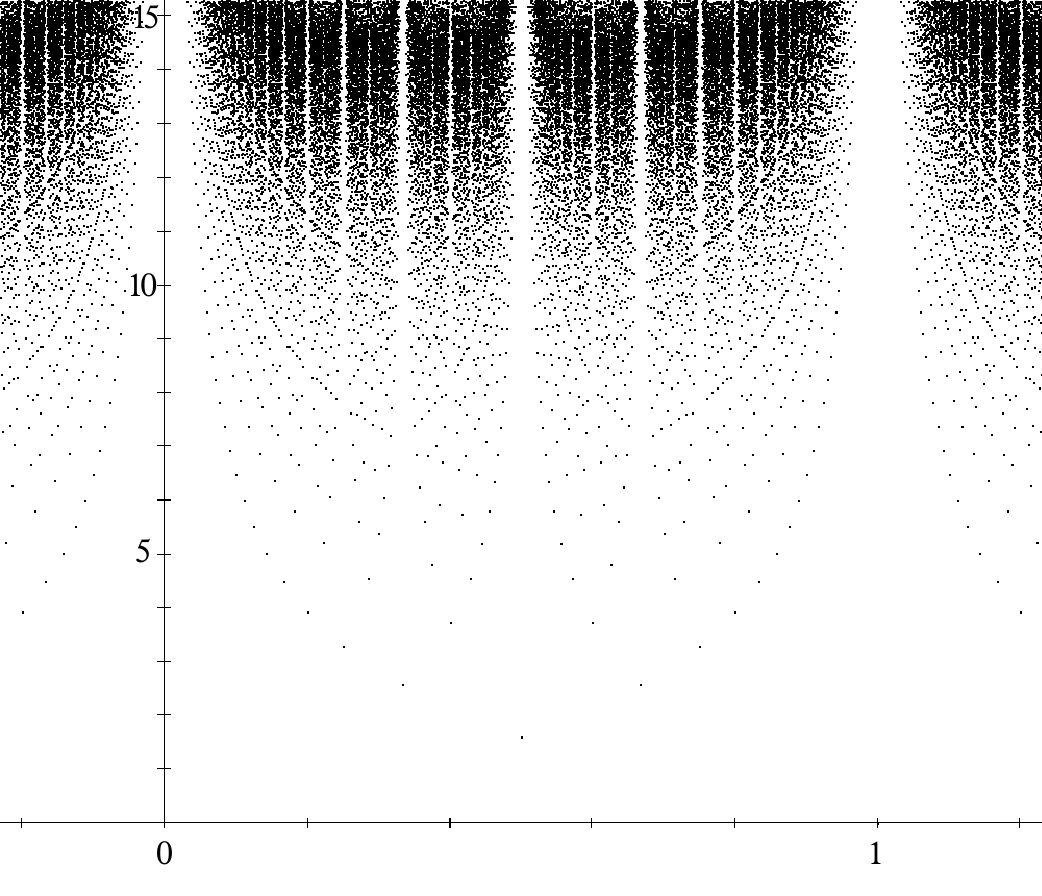}\hspace{13mm}
 \includegraphics[height=5.3cm]{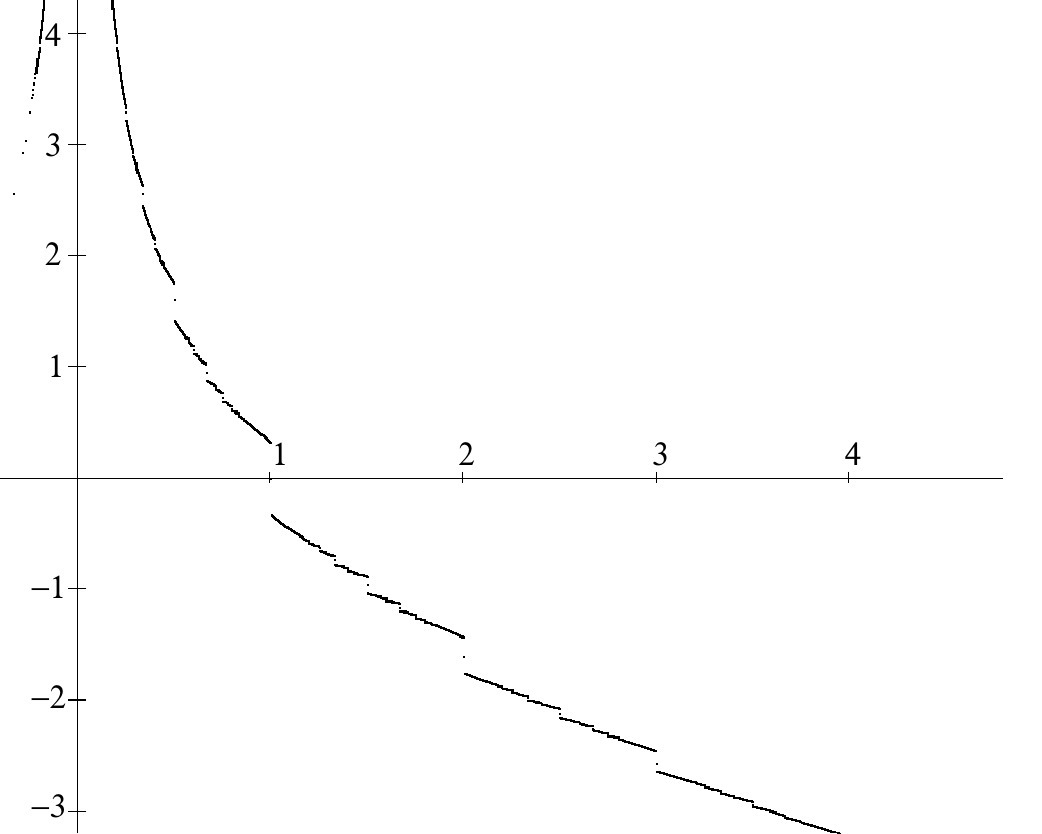}
\caption{The functions $\log(J(x))$ and $\log(J(x)/J(-1/x))$ for the $4_1$ knot.}
\label{OldJplot}
\end{figure}
and $\log(J(x)/J(-1/x))$ for $K=4_1$ (``figure~8 knot''), the simplest hyperbolic
knot. The former consists of a whole ``cloud'' of points and has no reasonable extension
to the real numbers, whereas the latter does extend to a well-defined function on~$\R$,
albeit one with infinitely many discontinuities.

\begin{figure}[t]\centering
\includegraphics[height=0.20\textheight]{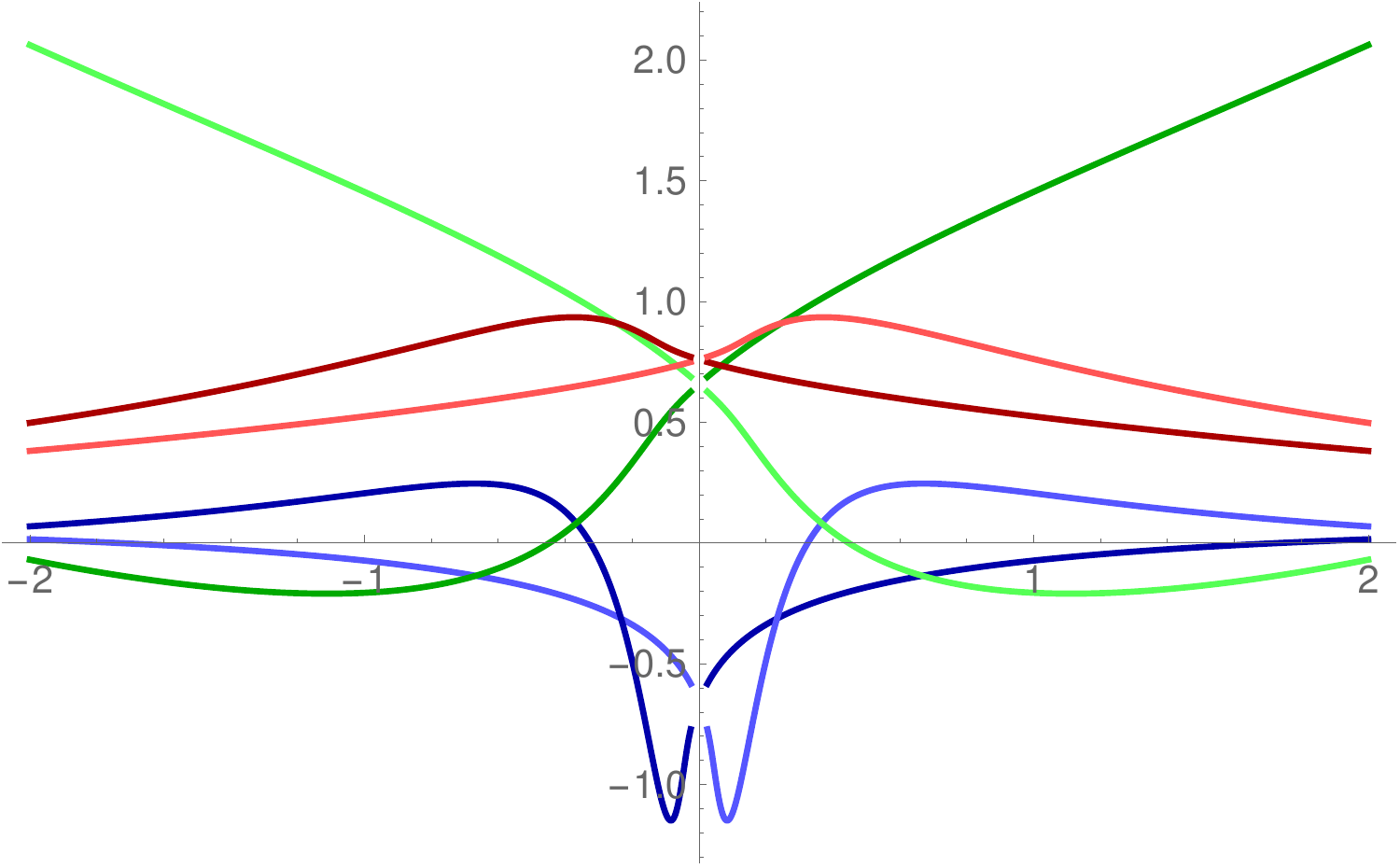}
\caption{Plots of the six nontrivial entries (rescaled) of the matrix $W^{(4_1)}_S(x) $.}
\label{LiftedPhi}
\end{figure}

The marked improvement of the graph on the right of Figure~\ref{OldJplot} as opposed to the one on the
left was already very striking and led to the introduction in~\cite{Za:QMF} of a notion
of ``quantum modular forms'' that has proved quite useful and has been exploited and
extended by several subsequent authors. But it was also somehow unsatisfactory, because
the function~${\log J(x)-\log J(-1/x)}$ still is far from smooth or even continuous, whereas
in all the other examples in~\cite{Za:QMF} the difference~${f(x)-f(\g x)}$ for a quantum
modular form $f\colon \BQ\to\BC$ extended to an analytic function on~$\BR$ minus a finite set.
This problem was~``solved'' in~\cite{Za:QMF} by defining quantum modular forms by the
weak requirement that the difference $f(x)-f(\g x)$ was ``analytically better behaved''
than $f$ itself, rather than demanding that it be analytic on the complement of a finite
set. But now it turns out that this cop-out is not needed, since the riddle of the missing
smoothness is solved completely by upgrading $J$ to a matrix of which it is only one entry.
Specifically, in the new picture, $J(X)$ is replaced by a certain matrix-valued invariant
$\bJ(X)=\bJ^{(K)}(X)$ (we use boldface letters to indicate matrices) which is defined and studied
in the course of this paper (Sections~\ref{sub.Formal}, \ref{sub.GQMC}, and \ref{sub.4.1}--\ref{sub.4.5}).
For the figure~8 knot this matrix
has the form
\[
\left(\begin{matrix} 1&J(X)&*\\0&*&*\\0&*&*\end{matrix}\right),
\]
 where
each of the six nontrivial components has a ``cloudlike'' graph like the first plot
in Figure~\ref{OldJplot}. But now instead of dividing the {\it scalar} invariant $J(X)$
by~$J(-1/X)$ as before, we look at the {\it matrix} product
$\bJ(-1/X)^{-1}\jt_S(X)\bJ(X)$,
where $\jt_\g(X)$ is the matrix-valued automorphy factor defined in~\eqref{autofactor}.
Then the graphs of the six nontrivial entries of this product matrix, multiplied by
suitable elementary factors to make them real and finite at the origin, look as
in Figure~\ref{LiftedPhi}
and are now smooth functions on the real line!\footnote{To
 avoid misconceptions, we note that the other striking property of this picture, its
 symmetry with respect to the vertical axis, is due to the accidental fact that the
 $4_1$ knot is amphicheiral. But already for the next simplest case of the $5_2$ knot,
 which will be the second of our three standard examples throughout the paper, the
 matrix $\bJ$ would be have size $4\times4$ with 12 non-trivial entries, all complex,
 and a graph of their 24 real and imaginary parts would be visually unintelligible.}
The same graph also beautifully illustrates the interrelationship of the different
matrix invariants of knots which we spoke of above, because the six curves that are
plotted are at the same time canonical lifts to $C^\infty(\R)$ of the six components of
the matrix of formal power series $\bPhi_\a(h)$ that is associated to the knot and to
every rational number~$\a$ (here for $\a=0$). In this way the matrix $\bJ$ of
Habiro-like functions determines the formal power series~$\bPhi_\a(h)$ (to get values
of~$\a$ other than~0, one would replace $-1/X$ by~$\g(X)$ for any $\g\in\SL_2(\Z)$ with
$\g(\infty)=\a$), and conversely the matrices~$\bPhi_\a(h)$ determines the matrix~$\bJ$ simply by $\bJ(\a)=\bPhi_\a(0)$.

This preliminary discussion and these pictures already give a first impression of the
content of this paper. However, giving the exact definitions of the various objects
being studied is not at all a straightforward business, because some of these are
based only on numerical data and cannot yet be justified by any theoretical
considerations, so that they can only reasonably be introduced after the numerical
investigations have been presented, while in other cases there are different candidate
definitions whose equality is only conjectural. In the body of the paper we
will therefore present the material in two stages. Part~I introduces the main players:
the perturbative power series, their conjectural analytic and number-theoretical
properties, and the emergence of an $\SL_2(\Z)$-cocycle through their quantum modularity
properties. The final (conjectural) statements are given in Section~\ref{sec.Matrix},
so that a reader who wants to see just the short version of the story right away can
skip directly to that section. Part~II then contains more detailed information about
the definitions and properties of the objects appearing in Part~I, including a
discussion of the numerical methods used, some of which are quite subtle. The paper ends
with an appendix containing tables of some of the functions studied for a few simple
hyperbolic knots.

Since the paper contains so many different types of objects, with rather intricate
inter\-connections and taking shape only gradually in the course of the exposition, it
seemed useful to end this introduction by giving a detailed overview of the main
ingredients. A further reason to include this rather long list here is that it contains
a number of items (the unexpected appearance of algebraic units, a description of the
Bloch group and extended Bloch group in terms of ``half-symplectic matrices'',
the notions of ``asymptotic functions near~$\Q$'' and of ``holomorphic quantum modular
forms'', a generalization of the Habiro ring to Habiro-like rings associated to number
fields other than~$\Q$, or a procedure to ``evaluate'' divergent power series numerically)
that are applicable or potentially applicable in domains quite separate from
that of quantum knot invariants and that therefore may be of independent interest.

$\bullet$ {\bf Indexing set.}
Both the rows and the columns of the matrices
associated to a knot~$K$ are indexed by a finite set~$\calP_K$ that can be described
either in terms of boundary parabolic representations of the fundamental group of the
knot complement~$\mathbb S^3\ssm K$ or in terms of flat connections, as explained in
detail in Section~\ref{sec.Phi}.

$\bullet$ {\bf Lift of the complex volume.}
The leading asymptotic exponent of the new matrices is a complex-valued function on
the set $\calP_K$ which agrees with the complexified volume of the boundary
parabolic representation, except that the latter is only well-defined modulo $4 \pi^2\BZ$.
Thus, a consequence of the refined quantum modularity conjecture is that the
complexified volume of a hyperbolic knot now has a canonical lift from $\BC^2/4\pi^2\BZ$ to~$\BC$.

$\bullet$ {\bf Level.}
As already mentioned, the matrices that we study also
depend on a rational number~$\a$. In all cases they are periodic in~$\a$, with the
period~$N=N_K$ however not always being the same: it is~1 for the $4_1$ and $5_2$ knots,
but~2 for the $(-2,3,7)$-pretzel knot. (These three knots will serve as our standard
illustrations throughout the paper.) Similarly, the modular invariance properties are
not always under the full modular group~$\SL_2(\Z)$, but sometimes under the
subgroup~$\G(N)$. We do not know what this ``level'' $N$ is in general, although we have
a~guess (in terms of the quasi-periodicity of the degrees of the colored
Jones polynomials), but its appearance in the numerical investigations was striking.

$\bullet$ {\bf Perturbatively and non-perturbatively defined power series.}
In the original version of the quantum modularity conjecture, the main statement was the
\emph{existence} of a collection of formal power series $\Phi_\a(h)$ describing the
relationship between $J(X)$ and $J(\g X)$ for large~$X$ and fixed~$\g\in\SL_2(\Z)$,
with no prediction of what this power series was. However, in two papers~\cite{DG,DG2}
by Tudor Dimofte and the first author explicit candidates for these power series as
perturbative series in~$h$ defined by Gaussian integration of a function explicitly
given in terms of a triangulation of~$M^3$, rendering the original conjecture much more
precise. These series will now form all but the top entry of the second column of
our matrix. \big(The first column is simply $(1 0 \dots 0)^t$.\big) The top entry
of the second column is defined in a completely different, non-perturbative way in
terms of the expansion near $q={\rm e}^{2\pi{\rm i}\a}$ of the Kashaev invariant of~$K$ seen as an
element of the Habiro ring. All of this will be explained in detail in
Section~\ref{sec.Phi}, while the definitions of the other entries of the matrix
$\bPhi_\a(h)$, which are again given by perturbatively defined power series in all but
the top row and (conjecturally) by elements of the Habiro ring for their top entries,
will emerge via the Refined quantum modularity conjecture discussed in
Sections~\ref{sec.RQMC}.\looseness=-1

$\bullet$ {\bf Arithmetic aspects.}
One of the main themes of this paper is
that the topology of a~knot involves a large amount of surprisingly complicated
algebraic number theory. This is valid both for the values of the Kashaev invariant
itself and of its generalizations as given by the matrix~$\bJ(x)$ at rational
arguments~$x$ and more generally for the coefficients of the entries~\smash{$\bPhi_\a^{(\s,\s')}(h)$} \mbox{($\a\in\Q$,~$\s,\s'\in\calP_K$)}, which conjecturally belong
to~$\Qbar(h)$. Among the most striking things that we found were the occurrence of
certain algebraic \emph{units}, which led to the paper~\cite{CGZ} with Frank Calegari
associating units (modulo~$n$th powers) in the $n$th cyclotomic extension of an
arbitrary number field to elements in the Bloch group of this field, congruence
properties of Ohtsuki type (which will be touched on only briefly here but will become
a main theme in the planned paper with Peter Scholze and Campbell Wheeler
on the construction of Habiro
rings associated to any number field), and universal bounds, independent of the knot,
for the denominators of the coefficients of the perturbative power series occurring
(see Section~\ref{sub.denominators}).

$\bullet$ {\bf Half-symplectic matrices and the extended Bloch group.}
The power series constructed in~\cite{DG,DG2} are given in terms of the
so-called Neumann--Zagier data describing the combinatorics of a triangulation of a knot
complement. This data takes the form of an~${N\times2N}$ integral matrix, where $N$ is
the number of simplices of the triangulation, together with a solution (corresponding
to the ``shape parameters'', i.e., the cross-ratios of the vertices of the ideal
tetrahedra) of a collection of algebraic equations defined by this matrix. The key
property here is
that the defining matrix is the upper half of a $2N\times2N$ symplectic matrix over~$\Z$.
This leads to a somewhat different description than the standard one of the Bloch
group and extended Bloch group. All of this is discussed in
Section~\ref{sec.perturbation}, together with the definition of the perturbative series
and a different appearance of the same construction in the context of Nahm sums.

$\bullet$ {\bf Unimodularity and inverse matrix.}
Experimentally, we find that the matrices that we construct are always unimodular, and also that there are explicit
formulas for their inverses as linear (or, for the top row, quadratic) rather than
higher-degree polynomials in the entries of the matrices themselves. Combined with the
behavior under complex conjugation, this leads to a kind of generalized unitarity
property for our matrices (of course, again only conjecturally, but we will not keep
repeating this since most of the properties we are discussing are only conjectural,
though based on such extensive data that they are very unlikely not to hold). All of
this will be discussed in Section~\ref{sec.Matrix}. The formula for the inverses of our
matrices (apart from the top row) can be interpreted as giving quadratic relations for
our power series, a special case of which appears in a recent paper of Gang, Kim and
Yoon~\cite{torsion} and which will be described in Section~\ref{sub.quad}.
These quadratic relations will take on a life of their own in the companion
paper~\cite{GZ:qseries} in terms of expressions for the ``state integrals'' defined by
Kashaev and others as bilinear expressions in power series in~$q={\rm e}^{2\pi{\rm i}\tau}$ and
$\wt q={\rm e}^{-2\pi{\rm i}/\tau}$.

$\bullet$ {\bf Extension property.}
As already stated, the rows and columns
of our matrices are indexed by the set~$\calP_K$ of flat connections. This set has a
canonical element, the trivial connection, which we put at the beginning of the list,
and the first row and column of each of our matrices then has a completely different
nature from the other entries. In particular, as we already saw, the first column
always consists of a~1 followed by~0s, so that the entire matrix is in
$(1+r)\times(1+r)$ block triangular form, where $1+r=|\calP_K|$. This means that
these matrices are describing structures which are $r$-dimensional extensions of
1-dimensional substructures. This can be seen clearly in the $q$-holonomy discussed
below, where the recursions satisfied by elements giving the top row contain a
constant term~1 and those of the other rows are the corresponding homogeneous recursions.\looseness=-1

$\bullet$ {\bf $\boldsymbol{q}$-holonomy.}
As already explained, the second column (the
first column being trivial) of our matrices of formal power series has a direct
definition in terms of the triangulation of the knot complement and the corresponding
Neumann--Zagier data. The other columns are defined in a more complicated way that we
still have not understood completely. Roughly speaking, each of the entries of the
original column belongs to a ``$q$-holonomic system'', meaning that it is part of a
sequence of functions of $q={\rm e}^{2\pi{\rm i}\a}$ that satisfy linear recursions over $\Q[q]$
and hence span a finite-dimensional space, and the other columns belong to, and in
fact span, the same space. This applies not only to the matrices of formal power
series in~$h$, but also to the Habiro-like matrices~$\bJ$ and to the matrices of
$q$-series studied in~\cite{GZ:qseries}, and will be discussed in detail in
Section~\ref{sec.NewHolonomy}. The mysterious point here is that the columns of our
matrices, which are completely and uniquely defined by the various properties embodied
in the refined modularity conjecture, give a canonical basis for these $q$-holonomic
modules, but that even in those situations where we know what the module
is or should be, we do not have an a priori description of this basis.

$\bullet$ {\bf Refined quantum modularity.}
As stated at the beginning of
this introduction, our whole story arises from the quantum modularity conjecture (QMC)
made in~\cite{Za:QMF}. In its original form, the QMC says that for all
$\gamma=\sma abcd$ in $\SL_2(\BZ)$ we have
\[
J(\g X) \approx (cX+d)^{3/2} J(X)
\Phih_{a/c}\left(\frac{2\pi {\rm i}}{c(cX+d)}\right), \qquad X\to\infty
\]
as $X$ tends to infinity with a bounded denominator, where $\Phih_\a(h)$ for $\a\in\Q$
is a ``completed'' version of $\Phi_\a(h)$ obtained by multiplying it by a suitable pure
exponential in~$1/h$, and where~``$\approx$'' denotes asymptotic equality to all orders
in~$1/X$ (or~$h$). This statement already refines the volume conjecture of
Kashaev~\cite{K95} and its arithmetic properties and extension to all orders as
described in~\cite{DGLZ}. In the refined quantum modularity conjecture (RQMC), which
will be developed step by step in the course of Sections~\ref{sec.inter}
and~\ref{sec.RQMC}, it is extended in two different ways. First of all, the above
asymptotic statement will be generalized by replacing the function~${J\colon\Q\to\Q}$ by the
other entries, initially of the first column and then by a kind of ``bootstrapping''
process (see below) of the whole matrix. More importantly, however, the right-hand
side will be sharpened by the addition of lower-order terms which are completed
versions of other entries of the matrix~$\bPhih_\a(X)$. Since the addition of an
exponentially smaller expression to a divergent power series does not make sense
a priori, this requires a process of numerical evaluation by ``optimal truncation'' and
then ``smoothed optimal truncation'' as listed in the bullet ``Numerical aspects'' below
and discussed in detail in Sections~\ref{sub.4.1} and~\ref{sub.optimal}. The final
result gives an asymptotic development to much higher precision of each of the
generalized Habiro-functions $\bJ^{(\s,\s')}$ evaluated at $\g X$ with $X$ tending to
infinity as a linear combination of $(r+1)$ of the power series $\bPhih_\a(h)$, with
$\a=a/c$ and~$h=2\pi {\rm i}/c(cX+d)$. This can then be written compactly in matrix form as
\[
\bJ^{(K)}(\g X)
 \approx \jt_\g(X) \bJ^{(K)}(X)
 \bPhih^{(K)}_{a/c}\left(\frac{2\pi {\rm i}}{c(cX+d)}\right)
\]
(= equation~\eqref{dQPhi}), which is the final version of the RMC. Here the
``automorphy factor'' $\jt_\g(X)$ is a diagonal matrix whose first entry is
$(cX+d)^{3/2}$ and whose other entries are pure exponentials in~$X+d/c$. Note that this
property relates the matrices $\bJ$ and $\Phih$, and allows in particular to compute
the second one from the first one.

$\bullet$ {\bf A matrix-valued cocycle.}
The matrix-valued form of RQMC as
just stated leads immediately to the definition of an $\SL_2(Z)$-cocycle with
coefficients in the space of matrix-valued functions on~$\Q$ (or more precisely---and
necessary in order to have an $\SL_2(Z)$-module structure---of almost-everywhere-defined
functions on~$\BP^1(\Q)$), defined by
\[
W_\g(x) = \bJ(\g x)^{-1} \jt_\g(x)  \bJ(x),
\]
where this time the ``automorphy factor'' $\jt_g(x)$ is a slightly different diagonal
matrix, again with elementary entries, as defined explicitly in
equations~\eqref{autofactor} and~\eqref{tweakdef}. The cocycle properties of this
automorphy factor imply that $W_\g$ is a multiplicative cocycle, meaning that
$W_{\g\g'}(X)$ is equal to $W_\g(\g'X)$ times~$W_{\g'}(X)$. Its remarkable properties are
summarized in the next two bullets.

$\bullet$ {\bf Analyticity and holomorphic quantum modular forms.}
The most important single discovery of this paper is that the cocycle
$W_\g(X)$, originally defined on rational numbers by the formula just given, extends
to a smooth function on the real numbers. This fact, which might have been found
purely experimentally by looking at the graphs of the components of~$W_\g(X)$, as
illustrated by Figure~\ref{LiftedPhi} above, and which can also be checked purely
experimentally, as explained in Section~\ref{sub.holomorphic}, was actually predicted
in advance on the basis of the occurrence of the same cocycle $\g\mapsto W_\g$
with a completely different construction in the companion paper~\cite{GZ:qseries}
to this one. Specifically, in that paper we construct a matrix $\bPsi^\hol(\tau)$ of
holomorphic functions of a~complex variable $\tau\in\C\ssm\R$, whose entries are power
series with integer coefficients in~$q={\rm e}^{2\pi{\rm i}\tau}$, and such that the coboundary
$\bPsi^\hol(\g \tau)\i\bPsi^\hol(\tau)$ extends holomorphically across
both the half-lines~${\big(\g\i(\infty),\infty\big)}$ and $\bigl(-\infty,\g\i(\infty)\bigr)$, with the
restrictions to these two half-lines coinciding with the function $W_\g$ there.
This extendability of $\bPsi^\hol(\g \tau)\i\bPsi^\hol(\tau)$ across subintervals
of the real line means that $\BJ^\hol$ is an example of a ``holomorphic quantum modular
form'', a new type of object that turns out to occur in many other contexts and that
will be described briefly in Section~\ref{sub.holomorphic} and in detail in the
papers~\cite{GZ:qseries,Z:HQMF}.

$\bullet$ {\bf ``Functions near $\boldsymbol{\BQ}$''.}
Each component $\Phi_\a^{\s,\s'}(h)$
of the matrix $\bPhi_\a$ has a natural completion, as explained in Section~\ref{sec.Phi}, defined
as its product with a certain exponential in~$1/h$ and (in the case of the top row)
a half-integral power of~$h$, and it is these completions that appear in the original
quantum modularity conjecture and its various extensions. It turns out that the
``right way'' to think of these collections of series is that they represent one single
``asymptotic function near~$\Q$'' defined by
\smash{$\BJ^{(\s,\s')}(\a-\hbar)=\Phih_\a^{\s,\s'}(2\pi {\rm i}\hbar)$}, where $\a$ varies over~$\Q$ and
$\hbar$ is infinitesimal. This notion of asymptotic functions near~$\Q$ (or simply
``functions near~$\Q$'' for short), which will be defined and explained more carefully
in Section~\ref{sec.nearQ}, sheds light on several properties of our knot invariants
(and also turns out to occur also in other contexts). In particular, the cocycle~$W_\g$,
which was initially defined (almost everywhere) on~$\Q$ by the formula~\smash{$W_\g(x)=\bJ(\g x)\i\jt_\g(X)\bJ(x)$}, is not a~coboundary in the space of functions
on~$\Q$, but \emph{is} one in the larger space of functions
near~$\Q$: $W_\g(x)=\BJ(\g x)\i \BJ(x)$. The occurrence of the \emph{same} cocycle with
two different representations as a coboundary in appropriate matrix-valued~$\SL_2(\Z)$-modules provides the link between the two papers and the reason for our
belief that both the matrix $\BJ$ of generalized Habiro functions and the matrix
$\BJ^\hol$ of $q$-series are realizations of the same underlying motive-like object.\looseness=1

$\bullet$ {\bf Numerical aspects.}
Everything in the paper is based on
numerical computations, and these have several non-obvious aspects, as discussed in
Section~\ref{sec.compKash}. In particular, we explain there how Kashaev invariants can
be computed rapidly and how one can then use extrapolation techniques to evaluate many
coefficients of the power series~\smash{$\Phi^{(\s,\s')}_a(h)$} numerically and recognize them
as real numbers. The calculations also have a ``bootstrapping'' aspect in which the
successively discovered relations among the series as described by the final refined
quantum modularity conjecture permit one to evaluate these series to increasing levels
of precision in a~recursive way. Finally, in order to identify the correct series in
the RQMC, it is crucial to be able to evaluate the divergent series in~$h$ occurring,
not only up to order~$h^N$ for any fixed integer~$N$, but up to exponentially small
error terms, where the constant occurring in the exponential can also be successively
improved in several steps. This is done by a process of ``smoothed optimal truncation''
which was originally a second appendix to this paper, but has now been relegated to a~separate publication~\cite{GZ:optimal} and is also briefly described in
Section~\ref{sub.optimal}.

We end with a few miscellaneous remarks on different aspects of the above constructions.

$\bullet$ {\bf Resurgence aspects.}
An important aspect of our paper are the matrices $\bPhi_\a(h)$ of factorially divergent
series and their distinguished lift $W_\g(x)$ to a matrix of analytic functions.
In this connection we find several properties that can be classified under the general
heading of ``resurgence.'' On the one hand, we find experimentally that the
coefficients of each entry of~$\bPhi_\a(h)$ are given asymptotically as \emph{integer}
linear combinations of certain divergent expansions involving the coefficients themselves
multiplied by gamma factors. The integrality of the so-called Stokes constants is a
phenomenon observed in the current paper and further studied in~\cite{GGM,GGM:peacock}.
A different connection of our results with the usual resurgence properties of the
perturbative series $\bPhi_\a(h)$ is the method of smoothed optimal truncation
mentioned just above, which can be seen as an alternative approach to lifting these
power series to actual functions than the standard method via Borel resummation and
Pad\'e approximation. Of course, the final emergence of a canonical lift coming from
the analyticity properties of the cocycle~$W_\g(x)$ eventually makes both numerical
procedures obsolete in our case, but this cocycle could not be found without having
them first.

$\bullet$ {\bf Equivalence of the various invariants.}
We observe that all of our invariants, assuming their conjectured properties, determine
each other and in particular all are determined by the colored Jones polynomials
and perhaps even by the Kashaev invariant alone. The point of the paper is therefore
not that our new invariants can distinguish knots, but that they reveal new properties
and make connection with $\SL_2(\BC)$-representations of the
fundamental group, hyperbolic geometry and non-trivial number theory.

$\bullet$ {\bf Eighth roots of unity and the Dedekind eta function.}
A further minor surprise was the appearance of Dedekind sums
at several places in the calculations, which we had not expected.
Notably, this occurred in the construction of state sums
for the non-Habiro-like elements of our matrix~$\bJ$ (see Section~\ref{sub.52statesum})
and in the ubiquitous but mysterious 8th root of unity that enters
all of our asymptotic formulas and that is related to the 8th root
of unity occurring in the modular transformation behavior of $\eta^3$.

$\bullet$ {\bf $\boldsymbol{3+1} $: a possible alternative interpretation.}
Finally, we mention a point that will not be discussed in the paper at all, but may
give a different way of looking at the objects studied here.
Namely, the $(\s,\s')$ entry of the matrix $\bJ(X)$ that we have
associated to a knot~$K$, and hence also to its complement~$M=\mathbb S^3\ssm K$, can be
thought of as numbers associated to the ``Witten cylinder'', which is the 4-manifold
$M\times \BR$ equipped with a pair of boundary-parabolic~$\SL_2(\BC)$ connections $\s$
and $\s'$ on its two ends, together with an integer~$k$, called the ``level'' in complex
Chern--Simons theory~\cite{Witten:JK,Witten:NZ} and related to the rational number $X$
by $k+2=\den(X)$. This suggests a possible interpretation of the entries of $\bJ(X)$
as expectation values of some yet-to-be-defined $(3+1)$-dimensional theory on
$M\times\BR$.

\begin{figure}[h!]\centering
\includegraphics[height=0.12\textheight]{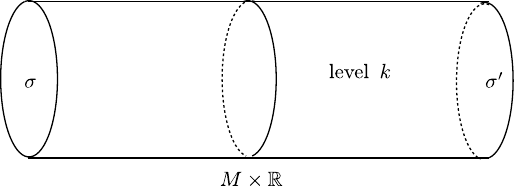}

\caption{The Witten cylinder.}\label{MtimesR}
\end{figure}

$\bullet$ {\bf Relation with the Kauffman bracket skein module.}
The objects that we have studied have another connection with a $3+1$ dimensional
theory that goes though the Kauffman bracket skein module and a conjecture of Witten
(see, for instance,~\cite{Witten5brane}), now a theorem of
Gunningham--Jordan--Safronov~\cite{Gunningham-Jordan}.
Recently, L\^e and the first author defined an explicit map from the Kauffman bracket
skein module of an integer homology 3-sphere to the Habiro ring~\cite{GL:skein}.
The image of this map (with can be effectively computed, and is naturally a topological
invariant) generates a finite-rank $\BZ\big[q^\pm\big]$-module which is a subset of
the Habiro ring, and a basis of this module conjecturally coincides with the
top row of the $\bJ$-matrix.

We end this introduction with a disclaimer. In this paper, we are trying only
to present interesting new phenomena and are \emph{not} striving for maximum
generality. Not only are we restricting our attention to knot complements
rather than arbitrary 3-manifolds, but we will usually assume that the knots
being considered have whatever properties (e.g., that their parabolic
character varieties are 0-dimensional) are convenient for our exposition.
In any case most of the material presented is empirical, based on extensive
experiments with a few very simple knots having all of these special
properties, and we prefer not to speculate on what modifications might be
needed in more general situations. Of course, we expect that the whole story is
quite general, and ongoing calculations by Campbell Wheeler that appear in his
thesis~\cite{Wheeler:thesis} indeed already
confirm that very similar types of statements will hold, for instance, for the
Witten--Reshitikhin--Turaev invariant of certain closed hyperbolic 3-manifolds.

\part*{Part I: The main story}
\pdfbookmark[1]{Part I: The main story}{part1}

\section{The original quantum modularity conjecture}
\label{sec.QMCK}

The starting point for this paper is the quantum modularity conjecture
of~\cite{Za:QMF}, which itself is a refinement of the famous volume conjecture
of Kashaev~\cite{K95}. Let us recall them briefly here.
Kashaev defined an invariant $\la K \ra_N \in \BZ[{\rm e}^{2 \pi {\rm i}/N}]$ for every
knot~$K$ and positive integer~$N$ and conjectured that if $K$ is hyperbolic
knot the absolute values of these numbers grow exponentially like ${\rm e}^{c_KN}$,
where $2\pi c_K$ is the hyperbolic volume of $S^3\ssm K$. A more precise version
of the conjecture~\cite{Gu} says that there is a full asymptotic expansion
\be
\label{K1}
\la K \ra_N  \sim
N^{3/2} {\rm e}^{v(K)N} \Phi^{(K)} \left(\frac{2\pi {\rm i}}N\right)
\ee
valid to all orders in $1/N$ where $v(K)$ is the suitably normalized complexified
volume of~$K$ and where $\Phi^{(K)}(h)$ is a (divergent) power series. \big(Here and
from now on we use the abbreviation~${\z_n=\e(1/n)}$ for $n\in\N$, where
$\e(x):={\rm e}^{2\pi{\rm i}x}$.\big) It was further conjectured in~\cite{DGLZ} and in~\cite{Ga:CS}
that~$\Phi^{(K)}(h)$ has algebraic coefficients, and more precisely that
it belongs to $\z_8 \delta^{-1/2}F_K[[h]]$, where $F_K$ is the trace
field of~$K$ and $\delta$ some nonzero element of $F_K$.

For instance, for the simplest hyperbolic knot $4_1$ (figure eight knot),
the Kashaev invariant is given explicitly by
\be
\label{J0for41}
\la 4_1\ra_N = \sum_{n=0}^{N-1} \bigl|(\z_N;\z_N)_n\bigr|^2
\ee
(see~\cite[equation~(2.2)]{K97})
with $(q;q)_n:=\prod_{j=1}^n\big(1-q^j\big)$, with values for $N=1,\dots,6$ and $100$
given numerically by
\begin{center}
 \def\arraystretch{1.5}
 \begin{tabular}{c|cccccccc}
$N$ & $ 1 $ & 2 & 3 & 4 & 5 & 6 & $ \cdots $ & 100 \\ \hline
$\la 4_1\ra_N$ & $ 1 $ & 5 & 13 & 27 & $46+2\sqrt5\approx50.47$ & 89
& $ \cdots $ & $8.2\times 10^{16}$ \\
 \end{tabular}
\end{center}
Here the trace field $F_K$ is $\BQ\big(\sqrt{-3}\big)$ and the series $\Phi^{(K)}(h)$ begins
\be
\label{as41}
\Phi^{(4_1)}(h)= \frac1{\sqrt[4]{3}}
\left(1 + \frac{11}{72\sqrt{-3}} h + \frac{697}{2 \big(72\sqrt{-3}\big)^2} h^2
 + \frac{724351}{30 \big(72\sqrt{-3}\big)^3} h^3 +\cdots\right) .
\ee
Similarly, for the knot $5_2$, where the Kashaev invariant is given by
formula~\eqref{K52} below, the trace field is $\BQ(\xi)$, where $\xi$ is the
root of $\xi^3-\xi^2+1=0$ with negative imaginary part, and~$\Phi^{(K)}(\hbar)$
has an expansion starting
\begin{align}
\Phi^{(5_2)}(h) ={}& \frac{\z_8}{\sqrt{3\xi-2}}
\left(1 + \frac{117\xi^2-222\xi+203}{24(3\xi-2)^3}h\right.\nonumber\\
&\left. {}  + \frac{117279\xi^2 - 209229\xi + 157228}{2\big(24(3\xi-2)^3\big)^2}h^2
 + \cdots\right).\label{as52}
\end{align}

The passage from the volume conjecture to the quantum modularity conjecture
(QMC for short) begins with the observation that the Kashaev invariant
$\la K\ra_N$ is the value at $x=-1/N$ of a 1-periodic function
$\J(x)=\J^{(K)}(x)$ on the rational numbers (i.e., it satisfies $\J(x+1)=\J(x)$ for
all $x$), defined uniquely by the further
requirement that it is a Galois-invariant function of~$\e(x)$. (The uniqueness
holds because every primitive $N$-th root of unity is a Galois conjugate of
$\z_N$.) As shown by Murakami and Murakami~\cite{MM}, this function
can also be identified with an evaluation of the colored Jones polynomial
$J_{K,N}(q)$~\cite{Jones, Tu1} by $\J(x) = J_{K,N}(\e(-x))$ for any~${N \in
\BZ}$ with $N x \in \BZ$.
In~\cite{Za:QMF}, it was found that~\eqref{K1} is just the special case
$\sma abcd=\sma0{-1}10$ of the more general (and of course still conjectural)
statement that
\be\label{QMCa}
\J^{(K)}\left(\frac{aN+b}{cN+d}\right) \sim (cN+d)^{3/2} {\rm e}^{v(K)(N+d/c)}
\Phi^{(K)}_{a/c}\left(\frac{2 \pi {\rm i}}{c(cN+d)}\right)
\ee
to all orders in $N$ as $N \to \infty$ for any matrix $\sma abcd \in \SL_2(\Z)$
with $c>0$, where $\Phi^{(K)}_\a(\hbar)$ is a~power series with
algebraic coefficients depending on $\a \in \BQ/\BZ$, with
\smash{$\Phi^{(K)}_0= \Phi^{(K)}$}. This does not yet look like a modularity statement,
but in~\cite{Za:QMF} it was further observed that~\eqref{QMCa} holds also
for non-integral values of $N$ (which we then denote by $X$ for clarity)
but with one crucial modification, namely that we have
\be
\label{QMCb}
\J^{(K)}\left(\frac{aX+b}{cX+d}\right) \sim (cX+d)^{3/2}
{\rm e}^{v(K)(X+d/c)} \Phi^{(K)}_{a/c}\left(\frac{2 \pi {\rm i}}{c(cX+d)}\right) \J^{(K)}(X)
\ee
to all orders in $1/X$ as $X\to\infty$ in~$\Q$ with bounded denominator,
with the same series $\Phi^{(K)}_\a(\hbar)$ as before but now with the
additional factor $\J^{(K)}(X)$ depending only on~$X$ modulo~1. (Here the condition
of $X$ having bounded denominator was included in the original conjecture, and will
be retained for its refinements in this paper, because all of our experiments were
done for $X$ with simple fractional part. However, it will be a consequence of
the smoothness statements cited above and discussed in Section~\ref{sub.lift} that in
fact~\eqref{QMCb} remains true for any sequence of rational numbers~$X$ tending to
infinity, and we have checked this experimentally for many cases.)

Formula~\eqref{QMCb} expresses the QMC in a quantitative form, in terms
of specific power series with algebraic coefficients, while the
plots of $J^{(K)}(X)$ and $J^{(K)}(X)/J^{(K)}(\g X)$ for~$K=4_1$ and~${\g=S}$ that
were shown (on a logarithmic scale because the functions grow exponentially) in
the introduction presents the same conjecture in a more qualitative visual form.
Both ways of looking at the conjecture will be refined greatly during the course of
this paper.

The QMC~\eqref{QMCb}, or even its special case~\eqref{QMCa}, give us not just
one, but a whole collection of power series \smash{$\Phi_\a^{(K)}(\hbar)$} associated to
any knot. These series \smash{$\Phi^{(K)}_\a(\hbar)$} have a striking arithmetical structure.
For example, in~\cite{Za:QMF} we found that for $K=4_1$ and $\a$ with
denominator~5 that
\be
\label{as41at5}
\Phi^{(4_1)}_{\a}(\hbar)= \pm  \sqrt[4]{3} \ve_\a^{1/5}
\bigl(A^{(\a)}_0 + A^{(\a)}_1 h +\cdots\bigr)
\qquad\text{if\quad$\a \in \frac15 \Z\ssm\Z,$}
\ee
where $\ve_\a$ is a unit (whose fifth root must be chosen appropriately) of the
cyclotomic field $\Q(\z_{15})=F_{4_1}(\z_5)$ ($\z_m:=\e(1/m)$) \big(explicitly
\smash{$\ve_\a= \frac{Z_\a^4-1}{Z_\a(Z_\a+1)^2}$}, with $Z_\a:=\e\big(\a-\frac13\big)$\big)
and the coefficients~\smash{$A_n^{(\a)}$} belong to the same field,
e.g., \smash{$A^{(a)}_0= 1+Z_\a^2+Z_\a^{-2}-Z_\a^4-Z_\a^{-4}$}, a prime of norm~29.
More generally, in the appendix to this paper extensive evidence is given for
the conjecture that the power series \smash{$\Phi^{(K)}_\a(\hbar)$} for any knot~$K$ and
rational number~$\a$ always belongs to
$\mu \delta^{-1/2}\sqrt[c]{\ve}F_{K,c}[[h]]$ with the same $\delta$ as before
and some $8c$-th root of unity~$\mu$, where~$c$ is the denominator of~$\a$, $F_{K,c}$
the field obtained by adjoining the $c$-th root of unity~$\e(\a)$ to~$F_K$, and
$\ve$ is an $S$-unit of $F_{K,c}$ for some finite set of primes~$S$ of $F_K$
independent of~$c$. This experimentally discovered property of quantum
invariants of knots in turn suggested the purely number-theoretical conjecture,
which was then proved in~\cite{CGZ}, that to an \emph{arbitrary} number
field~$F$ and element of the Bloch group of~$F$ one can canonically associate
a sequence of $S$-units, well defined up to $c$-th powers, in the~$c$ cyclotomic
extension $F(\z_c)$ for all~$c\ge1$, with~$S$ independent of~$c$.

We will give more details about this and other arithmetic properties of the
series \smash{$\Phi^{(K)}_\a$} (such as estimates of the denominators of their
coefficients) in Section~\ref{sec.arithmetic}, and will give a complete proof
in the case of the $4_1$ knot in Section~\ref{sec.QMC41}.
(This case and a few others were proven independently by Bettin and
Drappeau~\cite{BD}.) We will also give detailed numerical evidence for
several other knots, for several values of~$\a\in\BQ/\BZ$, and to a relatively
high degree in the power series \smash{$\Phi^{(K)}_\a(\hbar)$}, in the appendix to this
paper. The calculations required to obtain these values are not at all trivial,
since one has to be able to calculate the Kashaev invariants for (many)
arguments with large denominators and then use very precise extrapolation
methods to be able to find the coefficients to high enough accuracy to
recognize them numerically as algebraic numbers.

Presenting the numerical evidence for the QMC was the initial motivation
for this paper, and this already led to interesting numerical observations,
such as the appearance of the near unit~$\ve$ or the denominator estimates
mentioned above.
But in the course of doing these calculations we discovered that the QMC
is only part of a much larger story involving a whole collection of~power
series \smash{$\big\{\Phi^{(K,\s)}_\a(\hbar)\big\}_{\a \in \BQ/\BZ, \s\in\calP_K}$} indexed by a~certain finite set~$\calP_K$ (defined below) as well as by an index
in~$\Q/\Z$ as before. In the rest of Part~I, we explain what
these power series are, how they are related to each other, and how they
lead to new invariants and to a whole series of successive refinements of
the original quantum modularity conjecture.

\section{A collection of formal power series}
\label{sec.Phi}

\subsection[The indexing set P\_K]{The indexing set $\boldsymbol{\calP_K}$}
\label{sub.p}

The power series mentioned above are labeled by a finite set $\calP_K$
that coincides with the set of boundary parabolic $\SL_2(\C)$-representations
of $\G:=\pi_1\big(S^3\ssm K\big)$ (or equivalently, of flat connections
on~$S^3\ssm K$ whose restriction to the peripheral subgroup of $\G$ are
parabolic) whenever the latter is finite. For a hyperbolic knot, this
set has three distinguished elements, denoted $\s_0$, $\s_1$ and $\s_2$,
corresponding respectively to the trivial representation, the geometric
representation (given by the natural embedding of~$\G$ into the isometry group
of~$\mathbb H^3$) and the antigeometric representation, which is its
complex conjugate and corresponds to the geometric representation of the
orientation-reversed hyperbolic knot. We denote by $\calP\RED_K=\calP_K\ssm\{\s_0\}$
the reduced set of non-trivial representations (or connections), and often
number the elements of $\calP_K$ as $\s_0,\s_1,\dots,\s_r$, where
$r:=\big|\calP\RED_K\big|$ will be called the \emph{rank} of the knot. (We hope that the
superscript ``red'' will not confuse the reader into thinking that the representations
in $\calP\RED_K$ are reducible; in fact, quite the opposite is true.)
The points of~$\calP\RED_K$ correspond to the solutions (in~$\C$) of a set of polynomial
equations \big(the so-called Neumann--Zagier equations coming from a triangulation
of~$S^3\ssm K$\big) as explained in detail in Section~\ref{sec.perturbation}.
In particular, $\calP\RED_K$ comes equipped with an action of the absolute
Galois group $\text{Gal}\big(\overline{\BQ}/\BQ\big)$, so to every element
$\s \in \calP\RED_K$ is associated a number field~$F_\s$ \big(given either as the
field generated by the coordinates of the solution of the NZ equations or
as the fixed field of the stabilizer of $\s$ in $\text{Gal}\big(\Qbar/\BQ\big)$\big),
called its trace field, together with an embedding, also denoted~$\s$,
of~$F_\s$ into $\Qbar\subset\C$. The field $F_{\s_1}$
coincides with the trace field~$F_K$ of~$K$ as introduced above and
$F_{\s_2}$ is the same field with the complex conjugate embedding into~$\C$.
Two more important invariants of $\s\in\calP_K\RED$ are an element $\xi_\s$
of the Bloch group (or third $K$-group) of~$F_\s$ and a complexified
volume $\V(\s)=\V(K,\s)$ (= ${\rm i}$ times the usual volume plus the Chern--Simons
invariant), obtained as the image of~$\xi_\s$ under the Borel regulator
map or via the dilogarithm, or alternatively its renormalized version
$v(\s)=v(K,\s)=\V(\s)/2\pi {\rm i}$, which for the geometric representation is
the same as the number $v(K)$ occurring in~\eqref{K1}.
(We will usually omit the letter~$K$ in this and all similar notations when
the knot is not varying.) We extend all of these invariants to $\calP_K$
by setting~$F_{\s_0}=\Q$, $\V(s_0)=0$, $\xi_{\s_0}=0$.
In Section~\ref{sec.perturbation} of Part~II,
we will give more details about the set $\calP_K$ and its invariants, and
also say something about the situation when the variety of parabolic
representations contains positive-dimensional components. In
Sections~\ref{sec.Matrix} and~\ref{sec.perturbation}, we will also describe
a large (matrix- rather than vector-valued) collection of power series
associated to~$K$.

Before explaining how to associate a formal power series to
each~$\s\in\calP_K$ and $\a\in\Q/\Z$, we first would like to make the above
definitions more tangible by describing~$\calP_K$ explicitly for three
simple examples, the knots~$4_1$ and $5_2$ already used above and the
$(-2,3,7)$ pretzel knot (henceforth simply $(-2,3,7)$), which
will be our basic examples throughout the~paper. They have ranks 2, 3, and 6,
respectively. For $K=4_1$, the only elements
of~$\calP_K$ are the three universal ones~$\s_0$, $\s_1$ and $\s_2$, the
first corresponding to the trivial $\SL_2(\BC)$-representation with trace
field~$\Q$ and the other two both with trace
field~\smash{$\Q\big(\sqrt{-3}\big)$}, but with the complex embedding~${\sqrt{-3}\mapsto -{\rm i}\sqrt3}$ in the second case. The corresponding volumes
are $\V(\s_0)=0$, $\V(\s_1)={\rm i}V$, and $\V(\s_2)=-{\rm i}V$, where~${V=2.02\dots}$
is the usual volume, and are all real because the knot~$4_1$ is amphicheiral.
\big(In general, the mirror knot~$\overline K$ of a knot has trace fields
\smash{$F_{\overline K,\s}=\overline{F_{K,\s}}$} and~$\V(\overline K,\s)=\overline{\V(K,\s)}$.\big) For~${K=5_2}$ we again have
only two essentially different fields $\Q$ and $\Q(\xi)$ with
${\xi^3-\xi^2+1=0}$ (the cubic field with discriminant~$-23$), the latter with
three embeddings~$\s_1$, $\s_2$, and $\s_3$ corresponding to choosing the
root~$\xi\in\BC$ with negative, positive, or zero imaginary part,
respectively. But for the third knot~$K=(-2,3,7)$, $\calP_K$~consists of
seven elements, the trivial representation~$\s_0$, the three representations
$\s_1$, $\s_2$, $\s_3$ corresponding to the trace field of~$K$ (which is the
same as that of $5_2$, with its embeddings numbered the same way), and three
further elements~$\s_4$,~$\s_5$ and~$\s_6$ corresponding to the field
$\Q(\eta)$ with $\eta^3+\eta^2-2\eta-1=0$ (the abelian cubic field with
discriminant~$49$) together with the three embeddings into~$\BC$ given by
sending $\eta$ to~$2\cos(2\pi/7)$,~$2\cos(4\pi/7)$ and $2\cos(6\pi/7)$, respectively.
In general, to each knot~$K$ we associate the algebra $\calA_K=\Q\times\calA\RED_K$
defined as the product of the abstract fields~$F_\s$ with~$\s$ ranging over
representatives of the Galois orbits of~$\calP_K$, so that $\calP_K$ \big(resp.\ $\calP\RED_K$\big)
can be identified with the set of all algebra maps from~$\calA_K$ \big(resp.\ $\calA\RED_K$\big)
to~$\BC$; then for our three basic examples, we have
\be
\label{AK}
\calA\RED_{4_1} = \BQ\big(\sqrt{-3}\big) , \qquad \calA\RED_{5_2} = \BQ(\xi) , \qquad
\calA\RED_{(-2,3,7)} = \BQ(\xi)\times\BQ(\eta) .
\ee

\subsection[Four constructions of the power series Phi\^{}\{(K,s)\}\_a(h)]{Four constructions of the power series $\boldsymbol{\Phi^{(K,\s)}_\a(h)}$}
\label{sub.Formal}

We will now describe several different approaches to obtaining the formal
power series \smash{$\Phi^{(K,\s)}_\a(h)$} associated to an element~$\s$ of~$\calP_K$
and a number~$\a\in\BQ/\BZ$.

If $\s=\s_1$ is the geometric representation, then \smash{$\Phi^{(K,\s)}_\a(h)$} is
by definition just the power series~$\Phi^K_\a(h)$ whose existence is
asserted by the quantum modularity conjecture, and for $\s$ in the Galois
orbit of~$\s_1$ we simply apply Galois conjugation to this series (at the
level of its $n$-th power if~$\a$ has denominator~$n$), with some special
consideration for the roots of unity occurring. For example, for the
knot~$K=4_1$ the series \smash{$\Phi^{(K,\s)}_\a(h)$} for~$\a=0$ and $\a=a/5$ are the
ones given by~\eqref{as41} and~\eqref{as41at5}, respectively, for~$\s=\s_1$.
We then get \smash{$\Phi^{(K,\s_2)}_0$} simply by replacing $\sqrt{-3}$ by~$-\sqrt{-3}$
(or, in this case, replacing $h$ by~$-h$ and multiplying by~${\rm i}$) in~\eqref{as41},
and \smash{$\Phi^{(4_1,\s_2)}_{a/5}$} is obtained from~\eqref{as41at5} by performing the
same operation on both $\ve_{a/5}$ and the coefficients~\smash{$A_n^{(a/5)}$}. Similarly,
if $K$ is the $5_2$~knot then the value of \smash{$\Phi^{(K,\s)}_\a(h)$} at~$\a=0$ is given
by~\eqref{as52} if~$\s=\s_1$, and the values for $\s=\s_2$ or $\s_3$ are obtained
simply by replacing $\xi$ by its Galois conjugates. In general, the coefficients
of these power series lie in the product of a certain root of unity, the $c$-th
root (where $c$ is the denominator of~$\a$) of a unit in~$F_\s$, and a
conjugate of the same factor $\delta^{-1/2}$ as in the original QMC as described
in Section~\ref{sec.QMCK}. More details about the arithmetic of these numbers
will be given in Section~\ref{sec.arithmetic}.

The reader may have noticed that the QMC asserts the existence of the power
series \smash{$\Phi^{(K,\s)}_\a(h)$} for $\s=\s_1$ but gives no clue about how to define
them (and especially how to define those for~$\s$ not a Galois conjugate of
$\s_1$) given a hyperbolic knot $K$. A definition of the power series~\smash{$\Phi^{(K,\s)}_\a(h)$} for all $\s \neq \s_0$ was given by Tudor Dimofte and
the first author in the two papers~\cite{DG} (for~$\a=0$) and~\cite{DG2}
(for general~$\a$). What is more, the definition of the series uses as
input the gluing equation matrices of an ideal triangulation of the
knot complement, along with a solution to the Neumann--Zagier equations.
Roughly
speaking, one associates to an ideal triangulation of the knot complement
a collection of polynomial equations (the Neumann--Zagier equations) whose
solutions correspond to the elements of~$\calPx_K$, the solution for
each~$\s\in\calPx_K$ being a collection of algebraic numbers (the shape
parameters) belonging to the field~$F_\s$. One then associates to each
solution of these equations and for each $\a\in\Q/\Z$ a certain integral
that is evaluated perturbatively by the standard method of Gaussian
integration and Feynman diagrams (with a possible ambiguity of multiplication
by a power of $\e(\a)$). This process, whose details will be reviewed in
Section~\ref{sec.perturbation}, is completely effective and gives, for
instance, the three power series
\begin{align}
\Phi^{((-2,3,7),\s_j)}_0(h)
 ={}& \frac{\xi_j}{\sqrt{6\xi_j-4}}
\left(1 \m \frac{33\xi_j^2 - 123\xi_j + 128}{24(3\xi_j-2)^3} h
\right.\nonumber\\
&\left.- \frac{104172\xi_j^2 - 183417\xi_j + 130189}{2 \big(24(3\xi_j-2)^3\big)^2} h^2+ \cdots\right)\label{Phi237.123}
\end{align}
for the elements $\s_1$, $\s_2$, and $\s_3$ of $\calP_{(-2,3,7)}$, where
$\xi_1$, $\xi_2$, $\xi_3$ are the Galois conjugates of~$\xi$ as numbered above,
and the three totally different power series
\be
\label{Phi237.456}
\Phi^{((-2,3,7),\s_{j+3})}_0(h) =
 \sqrt{\frac{\eta_j-2}{14}} \left(1 \m \frac{43\eta_j^2 -21}{168} h
- \frac{3928\eta_j^2 + 63\eta_j -1491}{2\cdot 168^2} h^2 + \cdots\right)
\ee
for the elements $\s_4$, $\s_5$, and $\s_6$, where $\eta_j=2\cos(2\pi j/7)$
are the Galois conjugates of~$\eta$ in the ordering given above. The
coefficients of the power series \smash{$\Phi^{(K,\s)}(h)$} for all $\s\in\calP_K$
have similar arithmetic properties to the special case when $\s$ is Galois
conjugate to~$\s_1$.

As well as the ``straight'' power series \smash{$\Phi^{(K,\s)}_\a(h)$}, we will also
need the \emph{completed functions}
\be
\label{DefPhih}
\Phih^{(K,\s)}_\a(h) = {\rm e}^{\V(\s)/c^2h} \Phi^{(K,\s)}_\a(h),
\qquad c=\den(\a),\quad \s\ne\s_0 ,
\ee
which for the moment we think of as a purely formal expression (the
exponential of a Laurent series in~$h$ with a simple pole) but which will
be given a more precise sense later (cf.\ Section~\ref{sub.optimal}). It
is this combination that appear in all of our asymptotic formulas, e.g., the
right-hand side of~\eqref{QMCa} would become
\[
(cN+d)^{-3/2} \Phih_{a/c}^{(K)}\left(\frac{2\pi {\rm i}}{cN+d}\right)
\]
in this
notation. We should also mention here that in~\cite{DG2} the series
\smash{$\Phi_\a^{(K,\s)}(h)$} is defined only up to an $2n$-th root of unity, where
$n$ is the denominator of~$\a$. The generalized QMC that we will present in
the next section eliminates this ambiguity (at least up to a net sign depending
on~$\s$ but not on~$\a$).

An idea that will be crucial for this paper is that we have to associate
power series $\Phi_\a^{(\s)}(h)\in\Qbar[[h]]$ to the trivial representation
$\s=\s_0$ as well as to the non-trivial ones to get a coherent total picture.
Here there is no Neumann--Zagier data and we use instead a completely different
construction based on the Habiro ring. Recall that this ring is defined~by
\be
\label{HabDef}
\calH := \varprojlim \Z[q]/((q;q)_n) ,
\ee
where $(x;q)_n=\prod_{i=0}^{n-1}\big(1-q^ix\big)$ denotes the $q$-Pochhammer symbol
or ``shifted quantum factorial''. As mentioned in the introduction, Habiro
showed in~\cite{Ha} that the Galois-equivariantly extended Kashaev
invariant $\J^{(K)}(\a)$ ($\a\in\Q$) is the evaluation at~$q=\e(\a)$ of a
uniquely defined element, which we will denote by \smash{$\calJ^{(K)}(q)$}, of this
ring. We then \emph{define} the power series \smash{$\Phi_\a^{(K,\s)}(h)$} for~$\s=\s_0$ by
\[
\Phi_\a^{(K,\s_0)}(h) = \calJ^{(K)}\bigl({\rm e}(\a) {\rm e}^{-h}\bigr) \in \Q[\e(\a)][[h]].
\]
For example, for $K=4_1$ we have the explicit representation (equivalent
to~\eqref{J0for41} for $q=\z_N$)
\be
\label{Hab41}
\calJ^{(4_1)}(q) = \sum_{n=0}^\infty \big(q\i;q\i\big)_n (q;q)_n
 = \sum_{n=0}^\infty (-1)^n q^{-n(n+1)/2} (q;q)_n^2
\ee
of $\calJ^{(K)}(q)$ as an element of the Habiro ring, and setting
$q={\rm e}^{-h}\in\Q[[h]]$ we find
\be
\label{41HabAt0}
\Phi^{(4_1,\s_0)}_0(h) = 1 \m h^2 + \frac{47}{12} h^4
\m\frac{12361}{360} h^6 + \frac{10771487}{20160} h^8 \m \cdots
\ee
(which happens to be even because the knot~$4_1$ is amphicheiral),
while the Kashaev invariant for our second standard example~$5_2$
is given by formula~\eqref{K52} below and we find
\[
\Phi^{(5_2,\s_0)}_0(h) = 1 + h + \frac52 h^2 +\frac{49}6 h^3
+ \frac{797}{24} h^4 + \frac{19921}{120} h^5 + \cdots .
\]
For the $(-2,3,7)$ knot, we have no convenient Habiro-like formula for the
Kashaev invariant, but there is still a method (explained in Part~II) to
obtain its expansion to any order in~$h$ at any root of unity just from
the values at roots of unity, the expansion at~$q=1$ beginning
\[
\Phi^{((-2,3,7),\s_0)}_0(h) = 1 \m 12 h + 129 h^2 \m \frac{7275}4 h^3
\m \frac{384983}8 h^4 + \cdots .
\]
Note that the complexified volume vanishes for the trivial representation,
so that~\eqref{DefPhih} would suggest that we should define the completion
\smash{$\Phih_\a^{(K,\s_0)}(h)$} to be \smash{$\Phi_\a^{(K,\s_0)}(h)$}. But in fact, for reasons
that will appear clearly in Section~\ref{sec.inter}, it turns out to be better to
define \smash{$\Phih_\a^{(K,\s_0)}(h)$} in this case~by
\begin{gather}
\label{DefPhi0h}
\Phih_\a^{(K,\s_0)}(h)
= \left(\frac{ch}{2\pi {\rm i}}\right)^{\!3/2} \Phi_\a^{(K,\s_0)}(h) .
\end{gather}

We have now described constructions of the power series \smash{$\Phi_\a^{(K,\s)}(h)$}
for every~$\s\in\calP_K$, but based on very disparate ideas: if $\s$ is the
geometric representation or is Galois conjugate to it, we use the quantum modularity conjecture and Galois covariance, for other representations~$\s$
different from the trivial one we use a perturbative approach (which is
given in~\cite{DG} and~\cite{DG2} and conjectured there to agree with the
first definition when~$\s=\s_1$), and for the trivial representation we
define \smash{$\Phi_\a^{(K,\s)}(h)$} by a completely different formula based on the
Habiro ring. In fact, as already mentioned in the introduction, there
is even a fourth approach in which the series~\smash{$\Phi_\a^{(K,\s)}$} are obtained
from the asymptotics as $q$ tends radially to $\e(\a)$ of certain $q$-series
with integral coefficients. (This connection will not be discussed
further here but will be the main theme of~\cite{GZ:qseries}.) It is then
natural to ask why we consider these different series as being similar at
all and why we denote them in the same way. In the next two sections, we will
present a whole series of properties that justify this.

\section[Interrelations among the power series
 Phi\_a\^{}\{(s)\}(h)]{Interrelations among the power series
 $\boldsymbol{\Phi_\a^{(\s)}(h)}$}
\label{sec.inter}

In this section, we describe four empirically found properties, of very different
natures, that link and motivate the formal power series introduced above.

\subsection{The generalized quantum modularity conjecture}\label{sub.GQMC}

The function $\J^{(K)}$ from $\Q/\Z$ to $\Qbar$, which was
originally defined as the Galois-equivariant extension of the Kashaev invariant
$\langle K\rangle_N$, has now re-appeared as the constant term
\smash{$\Phi_\a^{(K,\s_0)}(0)$} of one of a collection of formal power series
\smash{$\Phi_\a^{(K,\s)}(h)\in\Qbar[[h]]$} indexed by the elements~$\s$ of a~finite
set~$\calP_K$ associated to the knot. This suggests that we should look
also at the constant terms of the other series as well, i.e., that we should
study the functions (\emph{generalized Kashaev invariants})
\begin{gather}
\label{cherries}
\J^{(K,\s)}\colon\ \Q/\Z \to \Qbar,\qquad \J^{(K,\s)}(\a) := \Phi_\a^{(K,\s)}(0)
\end{gather}
for all~$\s\in\calP_K$. These functions turn out to have beautiful
arithmetic properties generalizing in a non-obvious way the Habiro-ring
property of the original functions~$\J^{(K)}= \J^{(K,\s_0)}$. These
will be the subject of the subsequent paper~\cite{GSZ:habiro} and, apart from
a few numerical examples, will not be discussed further here.
Instead, we will concentrate on the asymptotic properties of the new
functions~\eqref{cherries}.
In particular, we can ask whether these functions satisfy an analogue of the quantum modularity conjecture for~$\J^{(K)}$. The answer turns out to be positive,
but to involve a~number of successive refinements arising from the numerical
data. We will present the simplest version here and the strongest versions,
which require more preparation, in Sections~\ref{sec.RQMC} and~\ref{sec.Matrix}.

We start once again with the simplest knot~$K=4_1$. Here the function
$\J^{(K,\s_0)}(\a)=\J^{(K)}(\a)$ is the one given by~\eqref{J0for41}
(with $\z_N$ replaced by~$\e(\a)$) whereas the new functions
$\J^{(K,\s_1)}(\a)$ and~$\J^{(K,\s_2)}(\a)$ are given explicitly by
\be
\label{Def41J1}
\J^{(K,\s_1)}(\a) = \frac1{\sqrt{c}\sqrt[4]{3}}
\sum_{Z^c = \z_6} \prod_{j=1}^c \bigl|1\m q^jZ\bigr|^{2j/c},
\qquad c=\den(\a),\quad q=\e(\a)
\ee
and $\J^{(K,\s_2)}(\a) ={\rm i} \J^{(K,\s_1)}(-\a)=\overline{\J^{(K,\s_1)}(\a)}$.
The original QMC says that $\J^{(4_1)}\bigl(\frac{aX+b}{cX+d}\bigr)$ is
asymptotically equal to
\smash{$(cX+d)^{3/2} \Phih^{(4_1)}_{a/c}\bigl(\frac{2\pi {\rm i}}{c(cX+d)}\bigr) \J^{(4_1)}(X)$}
for any matrix \smash{$\sma abcd\in\SL_2(\Z)$} as $X$ tends to infinity with bounded
denominator, where the ``completion'' $\Phih$ is defined by~\eqref{DefPhih}.
When we look at the corresponding asymptotics for the two new functions and
for the two simple matrices \smash{$\sma0{-1}10$} and \smash{$\sma1021$} of $\SL_2(\Z)$, we
see a similar behavior, but with two major differences: the ``automorphy factor''
$(cX+d)^{3/2}$ is no longer there, and there is a new exponential factor
involving the complex volume. Explicitly, what we find experimentally is
\be
\label{nontweak1}
\J^{(4_1,\s_1)}(-1/X)  \sim  {\rm e}^{\tV(K)/(\num(X)\cdot\den(X))}
\J^{(4_1,\s_1)}(X)  \Phih^{(4_1)}_0\left(\frac{2\pi {\rm i}}X\right)
\ee
(here ``num'' and ``den'' denote the numerator and denominator) and
\be
\label{nontweak2}
\J^{(4_1,\s_1)}(X/(2X+1))  \sim {\rm e}^{\tV(K)/((X+\frac12)\cdot\den(X)^2)}
\J^{(4_1,\s_1)}(X) \Phih^{(4_1)}_{1/2}\left(\frac{2\pi {\rm i}}{2(2X+1)}\right)
\ee
and similarly for $\Phi^{(\s_2)}$ but with $\tV(K)$ replaced by
$\tV(K,\s_2)=-\tV(K)$. The two equations~\eqref{nontweak1}
and~\eqref{nontweak2} can be written uniformly in the form
\[
\J^{(4_1,\s_1)}(\g X)  \sim {\rm e}^{\tV(K) \l_\g(X)}
\J^{(4_1,\s_1)}(X) \Phih^{(4_1)}_{a/c}\left(\frac{2\pi {\rm i}}{c(cX+d)}\right)
\]
for $\gamma=\sma abcd \in \SL_2(\BZ)$, where $\l_\g(x)$ is defined for
$x = r/s \in \BQ$ with $r$ and~$s$ coprime by
\be
\label{tweakdef}
\l_\g(x) := \frac{1}{\den(x)^2\big(x-\g^{-1}(\infty)\big)} = \frac{c}{s(cr+ds)}
= \pm \frac c{\den(x)\den(\g x)} .
\ee

The experiments show that the same thing happens for other knots~$K$ and all
representations~$\s$, i.e., we can formulate the
\emph{generalized quantum modularity conjecture}~(GQMC)
\be
\label{GQMC}
(cX+d)^{-\k(\s)} {\rm e}^{-\tV(\s)  \l_\g(X)}
\J^{(K,\s)}(\g X)  \sim  \J^{(K,\s)}(X)
\Phih^{(K)}_{a/c}\left(\frac{2\pi {\rm i}}{c(cX+d)}\right) ,
\ee
as $X \to \infty$ with bounded denominator (as usual), and where
$\gamma=\sma abcd \in \SL_2(\BZ)$ with $c>0$, and
where the \emph{weight}~$\W(\s)$ of the representation~$\s\in\calP_K$
is defined by
\be
\label{Defwsigma}
\W(\s) = \begin{cases} 3/2 & \text{if $\s=\s_0$,} \\
  0 & \text{otherwise.} \end{cases}
\ee
Notice that~\eqref{GQMC} coincides with the original QMC when $\s=\s_0$ because in this
case the factor~${\rm e}^{-\tV(\s)\l_\g(X)}$ on the left-hand side of~\eqref{GQMC} is
identically~$1$. We also see that the two different definitions~\eqref{DefPhih}
and~\eqref{DefPhi0h} of $\Phih^{(K,\s)}$ for~$\s\ne\s_0$ and $\s=\s_0$ can now be
written in a uniform way~as
\be
\label{UnifPhih}
\Phih_\a^{(K,\s)}(h) =
 |c\hbar|^{\W(\s)} {\rm e}^{v(\s)/c^2\hbar} \Phi_\a^{(K,\s)}(h),
 \qquad c=\den(\a),\quad \hbar:=h/2\pi {\rm i},
\ee
which will also be convenient at many other points. Notice that the convention
$\hbar=h/2\pi {\rm i}$ is almost, but not quite, the same as the one used in ordinary
quantum mechanics, and also that the factor $2\pi {\rm i}$ relating $h$ and~$\hbar$ is
the same as that used in our two different normalizations $\V(\s)$ and
$v(\s)=\V(\s)/2\pi {\rm i}$ of the volume, so that ${\rm e}^{v(\s)/c^2\hbar}={\rm e}^{\V(\s)/c^2h}$.

We end this subsection by proving a cocycle property of the arithmetic function
$\l_\g(X)$ that will be needed in Section~\ref{sec.Matrix}.

\begin{Lemma}
\label{lem.lambda}
For all $\g, \g'\in \PSL_2(\Z)$ and
$x\in\Q\ssm\big\{{\g'}^{-1}(\infty),{(\g\g')}^{-1}(\infty)\big\}$, we have
\[
\l_{\g \g'}(x) = \l_\g(\g' x) + \l_{\g'}(x) .
\]
\end{Lemma}

\begin{proof}
Let $\gamma=\sma abcd$, $\gamma'=\sma {a'}{b'}{c'}{d'}$ and
$\gamma \gamma'=\sma {a''}{b''}{c''}{d''}$. Then $c''d'-c'd''=c$, and hence
\[
\l_{\g \g'}(x) - \l_{\g'}(x) = \frac{c''}{s(c''r+d''s)}+\frac{c'}{s(c'r+d's)}
 = \frac c{s(c'r+d's)(c''r+d''s)} = \l_\g(\g' x)
\]
as required. A more enlightening way to say this is that $\l_\g(\g'(\infty))=C(\g\g',\g')$,
where
\[
C(\g_1,\g_2):=\frac{c\big(\g_1\g_2\i\big)}{c(\g_1)c(\g_2)}=\g_2\i(\infty)-\g_1\i(\infty),
\]
which is a coboundary and hence a cocycle.
\end{proof}

\subsection{Lifting the QMC from constant terms to power series}
\label{sub.QMCwithh}

In the previous subsection, we generalized the original QMC by replacing the Kashaev
invariant~$J^{(K)}$ by the generalized Kashaev invariants $J^{(K,\s)}$ for any
$\s\in\calP_K$. This in turn will be further refined in Section~\ref{sec.RQMC} by
adding terms of exponentially lower order to the right-hand side of the asymptotic
formula. Here we discuss instead a different refinement.

Our starting point, just as in Section~\ref{sub.GQMC}, is that the
Kashaev invariant $J^{(K)}(\a)$ is equal to the constant term \smash{$\Phi_\a^{(\s_0)}(0)$} of
the power series \smash{$\Phi_\a^{(\s_0)}(h)$} as defined in Section~\ref{sec.Phi}, so that the
original QMC~\eqref{QMCb} can be rewritten as
\[
\Phi_{\g X} ^{(\s_0)}(0)
 \sim  (cX+d)^{3/2} \Phi_{X} ^{(\s_0)}(0)
 \Phih_{a/c}\left(\frac{2\pi {\rm i}}{c(cX+d)}\right)
\]
for $X$ tending to infinity with fixed fractional part or with bounded denominator.
(Here we again omit the knot $K$ from the superscripts when it is not varying to
avoid cluttering up the notations. Recall also that~$c>0$.) It is then natural to ask
whether this asymptotic formula can be lifted to a corresponding statement
for the full series~\smash{$\Phi_\a^{(\s_0)}(h)$} rather than just its constant term.
The answer is affirmative, but with a little twist,
\be
\label{GQMCh}
\Phi_{\g X} ^{(\s_0)}(h^*)  \sim (cX+d)^{3/2}
\Phi_{X} ^{(\s_0)}(h) \Phih_{a/c}\left(\frac{2\pi {\rm i}}{c(cx+d)}\right) ,
\ee
where $x=X-\hbar$ with $\hbar$ as in~\eqref{UnifPhih} and \smash{$h^*=\frac h{(cx+d)(cX+d)}$}.

Let us explain what the asymptotic expansion~\eqref{GQMCh} means in the
simplest case of the figure~8 knot. Recall that \smash{$\Phi_X^{(\s_0)}(h)=\calJ\big(\e(X){\rm e}^{-h}\big)$}
where $\calJ(q)=\calJ^{(4_1)}(q)$ is the element of the Habiro ring given
by~\eqref{Hab41}, related to $\J(X)=\J^{(4_1)}(X)$ by $\J(X)=\calJ(\e(X))$.
Since the Habiro ring is closed under the operator $q {\rm d}/{\rm d}q$, it contains
the function $\calJ'$ defined by
\be \label{41Jderiv}
\calJ'(q) := q \frac {\rm d}{{\rm d}q} \calJ(q)
= \sum_{n=1}^\infty \big(q\i;q\i\big)_n (q,q)_n
\sum_{k=1}^nk \frac{1+q^k}{1-q^k} .
\ee
We then define the formal derivative $\J'\colon \Q/\Z\to2\pi {\rm i} \Qbar$ by
$\J'(X)=2\pi {\rm i}\calJ'(\e(X))$. Then the statement of~\eqref{GQMCh} in this case is
\be\label{CheckDeriv}
\frac1{(cX+d)^2}
\frac{\J'\bigl(\frac{aX+b}{cX+d}\bigr)}{\J\bigl(\frac{aX+b}{cX+d}\bigr)}
\m \frac{\J'(X)}{\J(X)} \approx 
\m \frac{2\pi {\rm i}}{(cX+d)^2}
\frac{\Phih'_{a/c}
 \bigl(\frac{2\pi {\rm i}}{c(cX+d)}\bigr)}{\Phih_{a/c}\bigl(\frac{2\pi {\rm i}}{c(cX+d)}\bigr)}
 ,
\ee
interpreted in the following sense. The left-hand
side of~\eqref{CheckDeriv} is $2\pi {\rm i}$ times an algebraic number belonging to some fixed
cyclotomic field for each fixed element $\g=\sma abcd\in\SL_2(\Z)$ and bound on the
denominator of~$X$, while the right-hand side is defined only as a divergent power
series in~$(cX+d)^{-1}$. The claim is then that when we compute both sides
of~\eqref{CheckDeriv} for fixed~$\g$ and for~$X$ tending to infinity
with bounded denominator, using~\eqref{41Jderiv} to compute the terms~$J'(X)$ and~$J'(\g X)$ as exact algebraic numbers, the two expressions agree numerically to
all orders in~$1/X$, and this is the statement that we verified numerically for several
elements~$\g$ and sequences of large rational numbers~$X$. Note that~\eqref{CheckDeriv}
is almost what we would get if we differentiated the original QMC formula~\eqref{QMCb}
logarithmically (which of course we are not allowed to do since it is only an
asymptotic statement valid for large rational numbers~$X$ with fixed denominator and
hence is rigid), except that then we would have an extra term \smash{$\frac32 \frac c{cX+d}$}
which is in fact not present because equation~\eqref{GQMCh} contains $(cX+d)^{3/2}$
rather than $(cx+d)^{3/2}$.

All of this was for the trivial connection~$\s_0$. If we consider instead an arbitrary
element $\s$ of $\calP_K$, then what we find is the obvious combination of~\eqref{GQMC}
(which was only for the constant terms $\Phi(0)$) and~\eqref{GQMCh} (which gave the
``twist'' needed to include~$h$), namely
\be
\label{GQMChh}
(cX+d)^{-\k(\s)} {\rm e}^{-\tV(\s)  \l_\g(X)}  \Phi_{\g X} ^{(\s)}(h^*)
  \sim  \Phi_X^{(\s)}(h) \Phih_{a/c}\left(\frac{2\pi {\rm i}}{c(cx+d)}\right) ,
\ee
with $x=X-\hbar$ and $h^*=h/(cx+d)(cX+d)$ as in~\eqref{GQMCh}.

Equation~\eqref{GQMChh} differs in two notable ways from the original~QMC~\eqref{QMCb}:
the appearance of the ``tweaking factor'' ${\rm e}^{-v(\s)\l_\g(X)}$ and the change of
infinitesimal variable from $h$ to~$h^*$. In fact, the first is explained very simply
by replacing the two power series $\Phi^{(\s)}$ in~\eqref{GQMChh} by their completions
as defined in~\eqref{UnifPhih}, because a short calculation shows that the number~$\l_\g(X)$ defined in~\eqref{tweakdef} is equal to the difference between
$1/\den(X)^2\hbar$ and $1/\den(\g X)^2\hbar^*$ with~${\hbar^*:=\h^*/2\pi {\rm i}=\hbar/(cx+d)(cX+d)}$, so that~\eqref{GQMChh} becomes simply
\be
\label{GQMChhh}
\Phih_{\g X} ^{(K,\s)}(h^*)  \sim (cx+d)^{-\k(\s)} \Phih_X^{(K,\s)}(h)
 \Phih^{(K,\s_1)}_{a/c}\left(\frac{2\pi {\rm i}}{c(cx+d)}\right) ,
\ee
where we have now again included the complete labels of the $\Phih$ series for clarity.
In this version both the tweaking factor ${\rm e}^{-v(\s)\l_\g(X)}$ and the automorphy
factor~$(cX+d)^{3/2}$ have been absorbed into the completed power series, but then
producing a new automorphy factor~${(cx+d)^{-3/2}}$. Finally, the ``twisting'' from $h$
to~$h^*$ is partly motivated by the calculation just given and the simplifications
in~\eqref{GQMChhh}, but more conceptually by observing that $x=X-\hbar$ implies
${\g x=\g X-\hbar^*}$. Equation~\eqref{GQMChhh} will then take on an even more natural
form in terms of the notion of ``functions near~$\Q$'' that will be introduced in
Section~\ref{sec.Matrix}.

\subsection{Quadratic relations}
\label{sub.quad}

The next interconnection among the power series \smash{$\Phi_\a^{(K,\s)}(h)$}
associated to a given knot~$K$ that we discover (experimentally, as always)
from the examples is that they satisfy an unexpected quadratic relation, namely
\be
\label{quadrel}
\sum_{\s \in \calPx_K} \Phi^{(K,\s)}_\a(h) \Phi^{(K,\s)}_{-\a}(-h) = 0 .
\ee
Notice that this relation is non-trivial even at the level of its constant term,
where it says, for example, that the value of the generalized Kashaev invariant
$J^{(5_2,\sigma)}(\alpha)$ defined in the last subsection belongs to the kernel
of the trace map from $\Q(\xi,\z_\a)$ to the trace field $\Q(\xi)$ of~$5_2$
for every rational number~$\a$. The special case of this when~$\a=0$
was observed independently by Gang, Kim and Yoon~\cite{torsion}.

The relation~\eqref{quadrel} is practically vacuous for the figure~8 knot,
since in that case it follows immediately from the identity
\smash{$\Phi_\a^{(4_1,\s_2)}(h)={\rm i}\Phi_{-\a}^{(4_1,\s_1)}(-h)$} mentioned at the beginning
of Section~\ref{sub.Formal}. (Stated differently, if we multiply the
series~\eqref{as41} by its value at $-h$, we obtain an element of
$\sqrt{-3}\Q\big[\big[h^2\big]\big]$, so that the trace down to~$\Q$ vanishes, and similarly
for~\eqref{as41at5}.) But for the $5_2$ knot
the identity is non-trivial even at~$\a=0$, where~\eqref{as52} gives
\begin{align*}
\Phi^{5_2}(h) \Phi^{5_2}(-h)
 ={}& \frac1{3\xi-2} + \frac{102\xi^2 - 183\xi + 135}{(3\xi-2)^7} h^2
\\
&\m \frac{143543 \xi^2 - 252029 \xi + 190269}{4 (3\xi-2)^{13}} h^4+ \cdots
\end{align*}
in which one can check that the three coefficients given, and in fact all
coefficients up to order~$h^{108}$, lie in the kernel of the trace map from~$\Q(\xi)$
to~$\Q$. Notice, by the way, that the series here has much simpler coefficients
(specifically, much smaller denominators) than the individual factors as given
by~\eqref{as52}. This is a special case of a more general phenomenon that will be
discussed in~\cite{GSZ:habiro}. When we look at~\eqref{quadrel} for this knot but
other values of~$\a$, the same thing happens: the $m$-th root of a unit in
$\Q(\xi,\z_m)$ that is a common factor of each \smash{$\Phi^{(5_2)}_\a(h)$} when
$\a$ has denominator~$m$ cancels when we multiply the series at~$\a$
and~$-\a$, and the series in~$\Q(\xi,\z_m)[[h]]$ that we find, although
it is no longer even when $\a$ is different from~0 or~$1/2$, always has coefficients
lying in the kernel of the trace map from $\Q(\xi,\z_m)$ to~$\Q(\z_m)$.

The above illustrates the relation~\eqref{quadrel} for our second simplest knot~$5_2$.
For our third standard example $K=(-2,3,7)$, this relation is even more surprising
because now $\calP_K$ has two Galois orbits, as discussed in Section~\ref{sec.Phi},
and the quadratic relation relates them to one another. Specifically, if we
consider separately the contributions from $\s_i$ for $1\le i\le 3$ and
for $4\le i\le6$, then equation~\eqref{Phi237.123} gives
\begin{gather*}
 \sum_{j=1}^3 \Phi^{(K,\s_j)}(h) \Phi^{(K,\s_j)}(-h)\\
\qquad= \operatorname{Tr}_{\Q(\xi)/\Q}\left(\frac{\xi^2}{2 (3\xi-2)}
+ \frac{605\xi^2 - 1217 \xi +878}{2^4 (3\xi-2)^7} h^2 + \cdots\right) \\
 \qquad =\frac12+ 0 h^2 \m \frac{13}{2^6} h^4
+ \frac{2987}{2^{11}\cdot3} h^6 + \frac{3517753}{2^{16}\cdot5} h^8 \m
 \frac{110362454561}{2^{19}\cdot3^3\cdot5\cdot7} h^{10} - \cdots
\end{gather*}
and equation~\eqref{Phi237.456} gives
\begin{align*}
&\sum_{j=4}^6 \Phi^{(K,\s_j)}(h) \Phi^{(K,\s_j)}(-h)\\
&\qquad= \operatorname{Tr}_{\Q(\eta)/\Q}\left(\frac{\eta-2}{2\cdot7}
+ \frac{18811 \eta^2 -78046 \eta+ 67485}{2^8 \cdot 3 \cdot 7^4} h^2+ \cdots\right)\\
& \qquad =- \frac12 + 0 h^2 +\frac{13}{2^6} h^4
 \m \frac{2987}{2^{11}\cdot3} h^6 \m \frac{3517753}{2^{16}\cdot5} h^8 +
 \frac{110362454561}{2^{19}\cdot3^3\cdot5\cdot7} h^{10} + \cdots .
\end{align*}
Each of these two series belongs to~$\Q\big[\big[h^2\big]\big]$. Computing many more terms
\big(we went up to~$ O\big(h^{38}\big) $\big), we find that their sum vanishes,
confirming the quadratic relation in a very striking way and at the same time showing
a subtle interdependence between the two cubic number fields associated
to this knot. We note, however, that these are only two of the three
number fields making up the algebra
$\calA_{(-2,3,7)} = \Q\times\BQ(\xi)\times\BQ(\eta)$ as defined in~\eqref{AK}.
We have not found \emph{any} relation between the power series
\smash{$\Phi^{(K,\s_0)}_\a(h)$ or $\Phi^{(K,\s_0)}_\a(h)\Phi^{(K,\s_0)}_{-\a}(-h)$}
and the power series~\smash{$\Phi^{(K,\s)}_\a(h)$} for $\s\ne\s_0$. This is reflected
in the fact that the summation in~\eqref{quadrel} is over~$\calPx_K$ and not over all
of~$\calP_K$.

We end this subsection by mentioning that, as well as the quadratic
relation~\eqref{quadrel}, there are also \emph{bilinear} expressions in
the~$\Phi^{(K,\s)}$ that are not zero, but (experimentally, and
in some cases provably) are \emph{convergent} rather than factorially divergent
power series. This will be discussed briefly in Section~\ref{sub.holomorphic}
and in detail in the companion paper~\cite{GZ:qseries}. Here we give only a
numerical example. In Proposition~\ref{prop.W1} below, we will give certain explicit
bilinear combinations of the $\Phi$-series which we believe are the Taylor expansions
of analytic functions and hence have a~positive (and known) radius of convergence.
In the simplest case (corresponding in the notation of Proposition~\ref{prop.W1}
to the $(\s_1,\s_1)$ component of the matrix \smash{$W_S^{(4_1)}(1+x)$}, where
$S=\sma0{-1}10$ as usual, combined with~\eqref{Orthogonality2}), this power series
is given by
\be
\label{DefR}
{\rm e}^{-v(4_1)} \Phi(2 \pi {\rm i} x)\Phi\biggl(-\frac{2 \pi {\rm i} x}{1+x}\biggr)
 - {\rm e}^{v(4_1)} \Phi\biggl(\frac{2 \pi {\rm i} x}{1+x}\biggr)\Phi(-2 \pi {\rm i} x) ,
\ee
with $\Phi=\Phi^{(4_1)}_0$ as given in~\eqref{as41}. The power series $\Phi$
has coefficients growing like $n!$ times an exponential
function (the precise asymptotics will be described in the next subsection)
and has~100th coefficient of the order of~$10^{94}$, but the combination~\eqref{DefR}
has radius of convergence~1 and, for instance, 100th coefficient of
order~$10^{-3}$. Notice that if we replace all $\Phi$'s in~\eqref{DefR}
by the corresponding $\Phih$'s, then the prefactors $ {\rm e}^{\pm v(4_1)}$
disappear.

\subsection{Asymptotics of the coefficients}
\label{sub.CoeffA}

The third interrelationship between the series $\Phi^{(\s)}_\a$ for different
elements $\s$ of~$\calP_K$ arises via the asymptotics of their coefficients.

For both theoretical and numerical purposes, we need to be able to compute
the ``values'' of the divergent series \smash{$\Phi^{(\s)}_\a(h)$} for very small~$h$,
and for this we need to know how their coefficients grow. We will write
\smash{$A^{(\s)}_\a(n)=A^{(K,\s)}_\a(n)$} for the coefficient of~$h^n$
in~\smash{$\Phi^{(K,\s)}_\a(h)$}.

As usual, we start with the simplest example $K=4_1$, $\s=\s_1$ (geometric
representation), and ~$\a=0$, where \smash{$\Phi^{(K,\s)}_\a(h)$} is just the
series~\eqref{as41}. Let us write just~$A(n)$ for its $n$-th coefficient
\big(so $A(0)=3^{-1/4}$, $A(1)=11 A_0/72\sqrt{-3}$\big). Calculating many
coefficients and using a standard numerical extrapolation method that
is recalled in Part~II, we find that $A(n)$ grows factorially like
$(n-1)! \l^{-n} \big(c_0+c_1n^{-1}+c_2n^{-2}+\cdots\big)$ for some constants
$\l$ and $c_i$. The numbers $\l$ and~$c_0$ are easily recognized to be
$2\V(K)=2{\rm i} \Vol(K)$ and~$3A(0)/2\pi$, respectively, but the further
coefficients $c_i$ have more and more complicated expressions. It turns
out that a much more convenient representation for the asymptotics is
as a sum of shifted factorials $(n-1-\ell)!$ rather than of terms
$n!/n^\ell$, because in this version we find the expansion
\be
\label{AnFirst}
A(n) \sim \frac3{2\pi}
\sum_{\ell\ge0} (-1)^\ell A(\ell) \frac{(n-\ell-1)!}{(2\V(4_1))^{n-\ell}}
\ee
with easily recognizable coefficients to all orders. If we now recall
that \smash{$\Phi_0^{(4_1,\s_2)}(h)$} equals \smash{${\rm i}\Phi_0^{(4_1,\s_1)}(-h)$} and hence
\smash{$A^{(4_1,\s_2)}_0(n) = (-1)^nA_n{\rm i}$}, then we can recognize~\eqref{AnFirst}
as one of a pair of coupled asymptotic expansions
\[
A^{(\s_1)}_0(n) \sim \frac3{2\pi {\rm i}}
\sum_{\ell\ge0} A^{(\s_2)}_0(\ell) \frac{(n-1-\ell)!}{(2\V(4_1))^{n-\ell}} ,
\qquad A^{(\s_2)}_0(n) \sim
\frac{-3}{2\pi {\rm i}}
\sum_{\ell\ge0} A^{(\s_1)}_0(\ell) \frac{(n-\ell-1)!}{(-2\V(4_1))^{n-\ell}} .
\]
This already looks quite nice, but the picture becomes even clearer when
we consider also the coefficients $B(0)=1$, $B(1)=0$, $B(2)=-1$, $\dots$ of the
third series $\Phi_0^{\s_0}$ as given in~\eqref{41HabAt0}. Since the $B(n)$
vanish for $n$ odd, it would first seem that one has to give separate
asymptotic formulas according to the parity of~$n$, but a better way is
to write \smash{$B(n)=A_0^{(\s_0)}(n)$} as a sum of \emph{two} asymptotic
expansions labelled by the two other elements $\s_1$ and~$\s_2$
of~$\calP_K $:
\be
\label{BnFirst}
B(n) \sim \sqrt{2\pi} \sum_{\ell\ge0}
A^{(\s_1)}_0(\ell) \frac{\G\big(n-\ell+\frac32\big)}{(-\V(4_1))^{n-\ell+3/2}}
\m \sqrt{2\pi} \sum_{\ell\ge0}
A^{(\s_2)}_0(\ell) \frac{\G\big(n-\ell+\frac32\big)}{\V(4_1)^{n-\ell+3/2}} .
\ee
Here we observe that the expressions $2\V(4_1)$, $-2\V(4_1)$, $-\V(4_1)$
and $\V(4_1)$ occurring in the denominators of the last two formulas can be
written in a uniform way as $\V(\s_1)-\V(\s_2)$, $\V(\s_2)-\V(\s_1)$,
$\V(\s_0)-\V(\s_1)$ and $\V(\s_0)-\V(\s_2)$, respectively. Exactly
analogous asymptotic statements turn out to hold for the coefficients of the
series \smash{$\Phi^{(4_1,\s)}_\a$} also for~$\a\ne0$, with the same coefficients,
leading for this knot to the uniform conjectural statement
\be
\label{CoeffAsymp}
A_\a^{(K,\s)}(n) \sim (2\pi)^{\W_\s-1} 
\sum_{\s'\ne\s}M_K(\s,\s') \sum_{\ell\ge0} A^{(\s')}_\a(\ell)
\frac{\G(n-\ell+\W_\s)}{\bigl(\V(\s)-\V(\s')\bigr)^{n-\ell+\W_\s}} ,
\ee
for all elements $\s\in\calP_K$ and all~$\a\in\Q$, where~$\W_\s$ is defined
as in~\eqref{Defwsigma} and where the coefficients~$M_K(\s,\s')$ are integers
independent on~$\a$, given for $K=4_1$ by
\[
M_{4_1} =
\begin{pmatrix} 0 & 1 & -1 \\ 0 & 0 & -3 \\ 0 & 3 & \hphantom{-}0 \end{pmatrix} .
\]
Experiments with our other two standard sample knots $5_2$ and $(-2,3,7)$
reveal the same asymptotic behavior~\eqref{CoeffAsymp}, with the matrix
$M_K$ given in these two cases by
\[
M_{5_2} = \begin{pmatrix}
 0 & \hphantom{-}1 & 1 & -1 \\
 0 & \hphantom{-}0 & 4 & -3 \\
 0 & -4 & 0 & -3 \\
 0 & \hphantom{-}3 & 3 & \hphantom{-}0 \end{pmatrix} ,
\qquad
M_{(-2,3,7)} =
\begin{pmatrix} 0 & \hphantom{-}0 & \hphantom{-}0 & 1 & 0 & 0 & 0 \\
 0 & \hphantom{-}0 & \hphantom{-}0 & 1 & 2 & 2 & 2 \\
 0 & \hphantom{-}0 & \hphantom{-}0 & 1 & 2 & 2 & 2 \\
 0 & -1 & -1 & 0 & 0 & 0 & 0 \\
 0 & -2 & -2 & 0 & 0 & 0 & 0 \\
 0 & -2 & -2 & 0 & 0 & 0 & 0 \\
 0 & -2 & -2 & 0 & 0 & 0 & 0 \end{pmatrix} .
\]

We end this subsection by making a number of remarks about the asymptotic
formula~\eqref{CoeffAsymp} and the matrices~$M_K$.

{\bf 1.}
The coefficients of the matrices~$M_K$ are much simpler invariants of~$K$ than
the coefficients of the power series \smash{$\Phi^{(\s)}_\a$}, because they are rational
integers rather than algebraic numbers and also do not depend on~$\a$.
It would be of considerable interest to have an direct topological definition
of these numbers rather than just an indirect one in terms of the (still
conjectural) asymptotic formula~\eqref{CoeffAsymp}. One possibility
in Section~\ref{sec.Phi} is that
they are related to the counting of flow lines in Floer homology. They are also
presumably the same as the skew-symmetric matrices of ``Stokes indices'' as
recently introduced by Kontsevich~\cite{MaximTalks}.

{\bf 2.}
The different forms of the asymptotics of the coefficients of $\Phi^{(\s)}_\a$
for $\s=\s_0$ and~$\s\ne\s_0$ are directly related to the different weights
and different completions of these series as given in equation~\eqref{UnifPhih}.

{\bf 3.}
A different asymmetry between the trivial and non-trivial
representations is seen in the fact~that $M_K(\s,\s_0)$ always vanishes but
$M_K(\s_0,\s)$ does not, meaning that the large-index coefficients of
the~$\Phi^{(\s_0)}$ series ``see'' the small-index coefficients of the
$\Phi^\s$ series for $\s\ne\s_0$ but not vice versa. It is interesting to note
that similar ``one-way phenomenon'' regarding the matrices appearing in
~\cite{Witten:continuation}, see also Gukov et al.~\cite{Gukov:resurgence,Gukov:BPS}.

{\bf 4.}
In all three examples given above, we further observe that apart from their
first column, which vanishes, and first row, which does not, the matrices
$M_K$ are skew-symmetric, i.e., ${M_K(\s,\s')=-M_K(\s',\s)}$ for $\s, \s'\ne\s_0$.
This phenomenon, which we expect to hold for all knots, will be shown below
to be a formal consequence of the quadratic relation~\eqref{quadrel}.

{\bf 5.}
We also observe that the lower $4\times4$ block of the matrix~$M_{(-2,3,7)}$
vanishes identically. In view of the numbering of the indices, this
means that $M_K(\s,\s')$ vanishes whenever~$\s$ and~$\s'$ are both real and distinct
from~$\s_0$. This in fact holds for all knots and is a special case of the more
general identity $M_K\big(\overline{\s\vphantom I},\overline{\s'}\big)=-M_K(\s,\s')$
for all $\s, \s'\ne\s_0$, which we can prove easily (assuming that the
expansion~\eqref{CoeffAsymp} is correct) simply by taking the complex conjugate
of~\eqref{CoeffAsymp} and noting that $\V(\overline\s)$ and \smash{$A_\a^{(\overline\s)}(n)$}
are the complex conjugates of $\V(\s)$ and \smash{$A_{-\a}^{(\s)}(n)$}, respectively (and, of
course, that the coefficients of~$M_K$ are real). The minus sign arises from the pure
imaginary prefactor~$(2\pi {\rm i})^{-1}$ in~\eqref{CoeffAsymp}.

{\bf 6.}
A corollary of~\eqref{CoeffAsymp} is the growth estimate
\[
A_\s^{(\s)}(n)= O\bigl(n^{\W_\s-1}n! \Delta(\s)^{-n}\bigr) ,
\]
where
\[
\Delta(\s)=\Delta(K,\s) = \min_{M_K(\s,\s')\ne0}\bigl|\V(\s)-\V(\s')\bigr| .
\]
This estimate will be important for the optimal truncation that is used in the
next section and discussed in more detail in Section~\ref{sec.compKash} and
in~\cite{GZ:optimal}.

{\bf 7.}
We should also mention that there is still some sign ambiguity in the definition
of the matrix~$M_K$. At the moment, even assuming the validity of the various
conjectures presented in the next two sections, we can only normalize the power
series \smash{$\Phi_\s^{(\s)}(h)$} up to the ambiguity of a~sign $\ve_\s\in\{\pm1\}$ independent
of~$\a$ but depending on~$\s$, and making this change would multiply~$M_K(\s,\s')$ by
$\ve_\s\ve_{\s'}$ (which would not affect either of the properties mentioned in~{\bf 3}~and~{\bf 4}~above). Similarly, when $\s=\s_0$ the formula defining $M_K(\s,\s')$ has
an inherent ambiguity coming from the choice of sign of square-root of
$V_\s-V_{\s'}=-V_{\s'}$ in~\eqref{CoeffAsymp} (only in the first term; the choices
for the other terms are then determined in the obvious way), so that each of the
matrix entries $M_K(\s_0,\s')$ is actually only well defined up to sign. Of course,
it is possible that there is some canonical way to normalize everything to eliminate
these ambiguities, but we do not yet know how to do this.

{\bf 8.}
Actually, however, there is a problem with all of these statements that we
have glossed over so far but that does need to be addressed. This is that the
right-hand side of~\eqref{CoeffAsymp} does not really make sense as it stands,
since the terms on the right-hand side are given by divergent series and hence
can be computed only up to some level of precision, but at the same time have
exponentially different orders of growth, so that it is not a priori clear what
it means to add them. In the case of~$4_1$, we did not see this problem, because
there is only one term in~\eqref{CoeffAsymp}. This point will be discussed
briefly in Section~\ref{sub.optimal} and in detail in~\cite{GZ:optimal}.

\section{Refining the quantum modularity conjecture}
\label{sec.RQMC}


In this section, we will show how one can go beyond the original QMC
as described in~Section~\ref{sec.QMCK} or its generalization as described
in~Section~\ref{sub.GQMC}. We will present this via a series of successive
refinements, each one found experimentally and building on its predecessors.
This will culminate in the complete, though of course still conjectural, definition
(in Sections~\ref{sub.4.1}--\ref{sub.4.4}) of the matrices~$\bJ$ and $\bPhi$
discussed in the introduction and of the final refinement (in Section~\ref{sub.4.5})
of the original quantum modularity conjecture.

\subsection{Improving the quantum modularity conjecture:
 optimal truncation}
\label{sub.4.1}

The QMC in its original form says that $\J^{(K)}(-1/X)$ agrees with
\smash{$X^{3/2} \J^{(K)}(X) \Phih_0^{(K)}(2\pi {\rm i}/X)$} to all orders in~$1/X$ as
$X$~tends to infinity with fixed denominator, with a similar statement when~$-1/X$ is replaced by $\frac{aX+b}{cX+d}$ for any~$\sma abcd\in\SL_2(\Z)$.
A natural question is whether we can do better than this and obtain an
asymptotic estimate, or even a precise asymptotic formula, for the \emph{difference}
of these two expressions. At first sight this seems to
makes no sense, since~\smash{$\Phih_0^{(K)}(h)$} \big(or more generally~\smash{$\Phih_{a/c}^{(K)}(h)$}\big)
is given in terms of a~divergent power series that a~priori
does not have a numerical value but rather gives only an approximation
up to any given order in~$h$. But we can remedy this by replacing the series
\smash{$\Phi(h) =\Phi_0^{(K)}(h)$} or \smash{$\Phi_{a/c}^{(K)}(h)$} by its ``optimal truncation''~$\Phi(h)^{{\rm opt}}$ obtained by truncating the divergent infinite series at
the value of~$N$ (depending on~$h$) where the terms of this series become smallest
in absolute value, a little like what is done in physics when for instance
the magnetic moment of the electron is computed to high accuracy by truncating a
divergent sum of Feynman integrals at a suitably small term. If $\Phi(h)=\sum A_nh^n$
with $A_n$ growing like $n!/B^n$ for some complex number~$B$, then this
``naive optimal truncation'' is given by $\sum_{n=0}^NA_nh^n$ with $N$ chosen
near to~$|B/h|$. Then the first term neglected, and hence also the expected
error, is of the order of magnitude of~${\rm e}^{-N}$, so we have a~way to define
$\Phi(h)$ up to an exponentially small error rather than only up to fixed
powers of~$h$. Of course, to get a completely well-defined function
$\Phi(h)^{{\rm opt}}$ we would have to fix a prescription for choosing~$N$, say
as the floor or ceiling or nearest integer to~$|B/h|$ (and perhaps also
dividing by~2 the last term retained), but since the terms with $n\approx|B/h|$ are all
very small the specific choice is not important and we will do better later anyway.

Using the description of the asymptotics of the coefficients of the
series~\smash{$\Phi^{(K,\s)}_\a(h)$} given in Section~\ref{sub.CoeffA} above, we can
compute their optimal truncations explicitly. Starting as usual with the
simplest example $K=4_1$ and $\a=0$ \big(and also $\Phi=\Phi^{(\s_1)}=\Phi^{(4_1,\s_1)}$,
the series occurring in the~QMC\big), we have from~\eqref{AnFirst} the estimate
$A_n=\O(n!/(2V)^n)$ with ${V=\Vol(4_1)=2.02988\dots}$, so the optimal truncation
occurs for $N$ near~$2V/|h|$. The expected error in $\Phi(h)$ for ${h=2\pi {\rm i}/X}$
is therefore of the order of ${\rm e}^{-2v(4_1)X}$, with $v(4_1)=V/2\pi=0.32306\dots$,
and since the completed function $\Phih(h)={\rm e}^{V/h}\Phi(h)$ grows like
${\rm e}^{v(4_1)X}$, this means that not only the relative but even the absolute expected
error in $\Phih(h)^{{\rm opt}}$ is exponentially small in this case. As a~numerical example, we consider the value~$X=100$. The Kashaev invariant
$\langle 4_1\rangle_{100}$, which we can compute to arbitrary precision
from~\eqref{Hab41} with $q\!=\!\z_{100}$, has the approximate value~$81985188380512462.9310054954341$, while the corresponding value
\smash{$100^{3/2} \Phih\bigl(\frac{2\pi {\rm i}}{100}\bigr)^{{\rm opt}}$} (obtained in this
case by retaining the first 66 coefficients of the divergent series) has the
numerical value $81985188380512461.9269943535808$ with an expected error of the
order of~$10^{-12}$. We see immediately that these two numbers are not equal
within the accuracy of the computation, so that the most obvious first guess for a~more precise version of the QMC is not true. But when we look at the difference
of these two numbers we find the numerical value
\[
\big\langle 4_1\big\rangle_{100} \m
100^{3/2} \Phih\left(\frac{2\pi {\rm i}}{100}\right)^{{\rm opt}} \approx 1.00401114185 ,
\]
which is very close to~1. Repeating the experiment for other large integral
values of~$X$, we find that this difference has the asymptotic expansion
$1-h^2+\frac{47}{12}h^4+\cdots$
(with $h=2\pi {\rm i}/X$ as before), which we recognize easily as the power series
$\Phi^{(4_1,\s_0)}(h)$ as given in~\eqref{41HabAt0}, and a~numerical computation
shows that indeed the optimal truncation of that series at $h=\frac{2\pi {\rm i}}{100}$
has precisely the same value~$1.00401114185$, to the same precision. Repeating the
calculations with other integral and non-integral values of~$X$ and also for
$\J^{(4_1)}(\g X)$ for matrices $\g\in\SL_2(\Z)$ other than~$\sma0{-1}1{\hphantom{-}0}$,
we find the same behavior in all cases, leading to the conjectural asymptotic formula
\begin{align}
& (cX+d)^{-3/2} \J^{(4_1)}\left(\frac{aX+b}{cX+d}\right) \nonumber\\
 &\qquad \overset?\approx{}
 \J^{(4_1)}(X) \Phih^{(4_1,\s_1)}\left(\frac{2\pi {\rm i}}{c(cX+d)}\right)
  + \Phih^{(4_1,\s_0)}\left(\frac{2\pi {\rm i}}{c(cX+d)}\right)\label{FirstRQMC}
\end{align}
for any matrix $\g=\sma abcd\in\SL_2(\Z)$ and for $X\to\infty$ with fixed
(or bounded) denominator, with the coefficient of the second series
\smash{$\Phih^{(4_1,\s_0)}$} being~1 for all~$X$ and~$\g$.

Here the natural question arises whether one can improve the precision
of~\eqref{FirstRQMC} even further by adding to the right-hand side a \emph{third}
term involving $\Phih^{(4_1,\s_2)}$, the last of the three completed series for the
$4_1$~knot. But for the moment we can't even make sense of this since the intrinsic
error in the optimal-truncation values of both $\Phih^{(4_1,\s_1)}(h)$ and
$\Phih^{(4_1,\s_0)}(h)$ has exponential decay of the order of ${\rm e}^{-v(K)X}$ (for the
first function because it grows like ${\rm e}^{+v(K)X}$ and has a relative error
${\rm e}^{-2v(K)X}$, as we have already seen, and for the second because it grows only
like a power of~$X$ but has a larger relative error ${\rm e}^{-v(K)X}$ by virtue
of~\eqref{BnFirst}). This is the same as the order of growth of the third function
$\Phih^{(4_1,\s_2)}(h)$, so that dividing the difference of the left- and
right-hand sides of~\eqref{FirstRQMC} by $\Phih^{(4_1,\s_2)}(h)$, with all
$\Phih$-series defined by optimal truncation, would give meaningless values. We will
return to this problem in Section~\ref{sub.4.3} below. Before doing
that, however, we first look at two other knots for which a new phenomenon appears
that is not visible for~$4_1 $.

\subsection{New elements of the Habiro ring}
\label{sub.4.2}

For the knot $K=5_2$ the set $\calP_K$ has four elements: the Habiro one, the
geometric and antigeometric ones, and the one corresponding to the real embedding
of the cubic field~$F_K=\Q(\xi)$. However, it has only three distinct real volumes:
the geometric volume $\operatorname{Im} \V(\s_1)=\Vol(K)$ (with the numerical value
$2.82812\dots$), the anti-geometric volume $\operatorname{Im} \V(\s_2)=-\Vol(K)$,
and~0 for both~${\s=\s_0}$ and $\s=\s_3$, and consequently only three distinct orders
of growth (one exponentially large, one exponentially small, and two of polynomial
growth) of the corresponding~$\Phih$-functions $\Phih^{(K,\s)}(h)$. (For simplicity
we concentrate for the moment only on $\a=0$ and omit it from the notations.)
This means that in the analogue of~\eqref{FirstRQMC} there is only one term that
is too small to be seen numerically when we replace the $\Phi$-series by their
optimal truncation, so that here one can hope to see \emph{three} distinct terms
on the right. To test this, we take the same values~${N=100}$,~$h=2\pi {\rm i}/N$ as before. Then
$\J^{(K)}\bigl(-\frac1N\bigr)$ is of the order of magnitude of~$10^{22}$ and
the difference \smash{$N^{-3/2}\J^{(K)}\bigl(-\frac1N\bigr)-\Phih^{(K,\s_1)}
\bigl(h\bigr)^{{\rm opt}} -\Phih^{(K,\s_0)}\bigl(h\bigr)^{{\rm opt}}$}
is of the order of~1 just as before, but now when we divide this difference by
$\Phih^{(K,\s_3)}\bigl(h\bigr)^{{\rm opt}}$ we obtain $2 + (1.22-5.23{\rm i})\cdot10^{-9}$,
strongly suggesting that the limiting value of this difference as~$X$ tends to
infinity through integers exists and is equal to~2. Further experiments for
non-integral values of~$X$ and for other matrices $\g=\sma abcd$ then lead to
the new conjectural asymptotic statement
\begin{align}
 (cX+d)^{-3/2} \J^{(5_2)}\left(\frac{aX+b}{cX+d}\right) \overset?\approx{}&  \J^{(5_2)}(X) \Phih^{(5_2,\s_1)}_{a/c}
 \left(\frac{2\pi {\rm i}}{c(cX+d)}\right)
 + \Phih^{(5_2,\s_0)}_{a/c}\left(\frac{2\pi {\rm i}}{c(cX+d)}\right)\nonumber\\
& + Q^{(5_2)}(X) \Phih^{(5_2,\s_3)}_{a/c}\left(\frac{2\pi {\rm i}}{c(cX+d)}\right)
 \label{52WeakRQMC}\end{align}
for all~$\g$ and all $X$ tending to infinity with bounded denominator, where
$Q^{(5_2)}(x)$ (which is a~temporary notation, only for this knot) is a function
that is independent of~$\g$ but that, unlike the constant
coefficient~1 of the Habiro term $\Phih^{(K,\s_0)}(h)$, is not independent of~$x$.
Instead, $Q^{(5_2)}(x)$ is numerically found to be a periodic function of period~1
taking on simple algebraic values, the first few being
\begin{gather*}
Q^{(5_2)}(0)=2,\qquad Q^{(5_2)}\big(\tfrac12\big)=8, \qquad
Q^{(5_2)}\bigl(\pm\tfrac13\bigr)=\frac{37\pm3\sqrt{-3}}2,
\qquad Q^{(5_2)}\bigl(\pm\tfrac14\bigr)=29\pm13{\rm i} .
\end{gather*}
(These values were found experimentally, using a Chinese-remainder-type interpolation,
and the existence of such functions for all knots is not known.)
Looking at more values (specifically, for all~$x$ with denominator up to~200),
we find that $Q^{(5_2)}(x)$ belongs to $\Z[\e(x)]$ and is Galois-invariant, so
we can write it as $\sQ^{(5_2)}(\e(x))$ where $\sQ^{(5_2)}(q)$ is an element
of $\Z[q]$ for every root of unity~$q$, the first values being given by
\begin{center}
\def\arraystretch{1.5}
\begin{tabular}{c|ccccccc}
Ord($q$) & $ 1 $ & 2 & 3 & 4 & 5 & 6\\ \hline
 $\sQ^{(5_2)}(q)$ & $ 2 $ & $ 8 $ & $ 20+3q $ & $ 29+13q $
 & $ 69+27q+37q^2+2q^3 $ & $ -46+69q $
\\
\end{tabular}
\end{center}
This suggests that $q\mapsto\sQ^{(5_2)}(q)$ might be an element of the Habiro
ring~$\calH$ defined in~\eqref{HabDef}, just as we know is the case for the
coefficient $\J^{(5_2)}(X)$ of the first
$\Phih$-term in~\eqref{52WeakRQMC}. This hypothesis can be tested numerically,
because a well-known property of any element $\sQ\in\calH$ (originally observed
by Ohtsuki~\cite{Oh} in the context of the WRT-invariants of integer homology
spheres even before the Habiro ring had been formalized) is that it satisfies
an infinite number of congruences, the simplest of which is that $\sQ(\z_p)$ for
every prime number~$p$ should be congruent modulo~$p$ to~\smash{$c_0+c_1\pi_p+c_2\pi_p^2+\cdots+c_{p-2}\pi_p^{p-2}$}, where $\pi_p=\z_p-1$ is the
prime dividing~$p$ in~$\Q(\z_p)$ and where the $c_i$ are rational integers
independent of~$p$. This means in our case that the coefficient of $x^i$ in
the polynomial $Q_p(1+x)\in\Z[x]$ should be congruent modulo~$p$ to a \emph{fixed}
integer~$c_i\in\Z$ for all primes $p>i+1$, and testing this for the numerically
obtained polynomials $Q_n$, we find that it is indeed true, with
$Q_p(\z_p)\equiv\text{Oh}(\pi_p)\!\pmod p$ for a power series
$\text{Oh}(x)\in\Z[[x]]$ beginning
\[
\text{Oh}(x) = 2 - 3x + 8x^2 - 28x^3 + 120x^4 - 614x^5 + 3669x^6 - 25125x^7 + \O\big(x^8\big) .
\]
In fact, later we were able to guess an explicit formula, given below in
Section~\ref{sub.qholo}, that is manifestly in the Habiro ring and that reproduces
the values of the polynomials $Q_n(q)$ and power series~$\text{Oh}(x)$ as given
above. But in many other cases, including the $(-2,3,7)$ knot discussed
below, we cannot give even conjectural explicit formulas
of the required kind, and in such cases it is important to be able to have a numerical
test of the Habiro-ness of a periodic function.

Notice that the right-hand side of~\eqref{52WeakRQMC} contains only three of the four
completed power series~$\Phih_{a/c}^{5_2,\s}$. Just as for the $4_1$ knot, this is not
because the last one isn't really there, but because our approximate evaluations are
not accurate enough at this point to detect the remaining term, which is exponentially
small. We will correct this in Section~\ref{sub.4.3}.

We were able to carry out similar calculations for the $(-2,3,7)$ pretzel knot,
though the numerical analysis required here was much more arduous due to the larger
number of series involved. Recall that this knot has rank~6, so that $\calP_K$
contains seven elements. What makes the
calculation feasible at all is that five of these seven elements are real (the Habiro
one and the ones corresponding to the real embedding of $\Q(\xi)$ and to all three
embeddings of $\Q(\eta)$), so that only one of the seven $\Phih$-functions is
exponentially small and hence invisible with optimal truncation. (Actually, the fact
that the other terms apart from the geometric one are of the order of~1 is not quite
enough: one also has to verify by using the formulas of Section~\ref{sub.CoeffA} and the
numerical values of the complex volumes~$\tV(K,\s_i)$ that the absolute error made
in calculating the exponentially large dominant term $\Phih^{(K,\s_1)}(h)$
using optimal truncation is exponentially small.) We find a formula exactly analogous
to~\eqref{52WeakRQMC}, but now with six terms on the right, namely
\be
\label{237WeakRQMC}
(cX+d)^{-3/2} \J^{(-2,3,7)}\left(\frac{aX+b}{cX+d}\right)  \overset?\approx
\sum_{\substack{0\le j\le 6\\j\ne2}}
Q^{(-2,3,7)}_j(X)
\Phih^{((-2,3,7),\s_j)}_{a/c}\left(\frac{2\pi {\rm i}}{c(cX+d)}\right) ,
\ee
where $j=2$ is omitted for the same reason as in~\eqref{52WeakRQMC}
(viz., that the corresponding term is too small to see at this stage) and where
\smash{$Q^{(-2,3,7)}_1(x)=\J^{(-2,3,7)}(x)$}, and \smash{$Q^{(-2,3,7)}_0(x)\equiv1$}, and the four
new periodic functions \smash{$Q^{(-2,3,7)}_j(x)=\sQ^{(-2,3,7)}_j(\e(x))$} take values in
$\Z[\e(x)]$ just as before, the first values (for $j\ne0,2$) being
\begin{center}
\def\arraystretch{1.2}
\begin{tabular}{c|ccccccc}
 Ord($q$) & $1 $ & 2 & 3 & 4 & 5 & 6\\ \hline
 $\sQ^{(-2,3,7)}_1(q)$ & $\hphantom{-}1 $ & $\hphantom{-1}1 $ & $ -5+6 q $ & $ 17-8 q $ &
 $ -21-27 q-5 q^2+4 q^3 $ & $ -107+108 q $
 \\
 $\sQ^{(-2,3,7)}_3(q)$ & $-4 $ & $ -12 $ & $ -15 - 10q $ & $ -16 - 2q $
 & $ -36 - 20q - 29q^2 - 24q^3 $ & $ 23 + 14q $ \\
 $\sQ^{(-2,3,7)}_4(q)$ & $\hphantom{-}2 $ & $ -10 $ & $ -16 - 12q $ & $ -46q $
 & $ -8 - 44q - 38q^2 - 48q^3 $ & $ 116 - 24q $ \\
 $\sQ^{(-2,3,7)}_5(q)$ & $-2 $ & $ \hphantom{1}{-}6 $ & $ -14 - 6q $ & $ 8 - 10q $
 & $ 32q - 4q^2 - 10q^3 $ & $ -82 + 122q $ \\
 $\sQ^{(-2,3,7)}_6(q)$ & $\hphantom{-}2 $ & $\hphantom{-1}2 $ & $ \hphantom{-1}4 - 8q $ & $ 10 - 12q $
 & $ -4 - 36q - 44q^2 - 34q^3 $ & $ 136 - 16q $
\end{tabular}
\end{center}
Just as with the $5_2$ knot, we can verify the Ohtsuki property for these functions
to a large number of terms and thus convince ourselves that each one belongs to the
Habiro ring, even though in this case we do not know an explicit formula that makes
this property manifest.

\subsection{Smoothed optimal truncation}
\label{sub.4.3}

We already mentioned at the end of Section~\ref{sub.4.1} that it would be natural
to expect a third term in~\eqref{FirstRQMC} involving the missing $\Phih$-function
$\Phih^{(K,\s_2)}(h)$, and the same applies even more strikingly to the two knots
discussed in Section~\ref{sub.4.2}, where we were obliged to omit the $\s_2$-term
in both~\eqref{52WeakRQMC} and~\eqref{237WeakRQMC} because it would have been absorbed
in the error terms of the other $\Phih$'s and hence could not be seen numerically if
these values were defined by naive optimal truncation. However, there is a more
precise way of turning the divergent series $\Phi(h)=\Phi_\a^{\s}(h)$ into functions
that are defined up to exponentially rather than merely polynomially small errors, but
with a much better exponent than before, by replacing the naive optimal truncation
$\Phi(h)^{{\rm opt}}$ by a \emph{smoothed} version $\Phi(h)^{{\rm sm}}$. The details of
the procedure to do this are somewhat complicated and play no role for the story
we are telling here, so will be given in detail in a separate
publication~\cite{GZ:optimal} and described briefly in Section~\ref{sub.optimal},
the only important point here being that the improvement is sufficiently
good, at least for our three standard knots, that we can unambiguously identify
the periodic coefficients of the missing $\Phih$-terms.

We start as usual with the knot~$K=4_1$ and the series \smash{$\Phi(h)=\Phi^{(\s_1)}_0(h)$}
whose initial terms are given in~\eqref{as41}. In Section~\ref{sub.4.1}, we considered
$X=100$, $h=\frac{2\pi {\rm i}}X$ and saw that the number~\smash{$100^{3/2} \Phih(h)^{{\rm opt}}\approx 8.195\times 10^{16}$} had an error
of the order of~$10^{-12}$, which was more than sufficient to identify
its difference with \smash{$\big\langle 4_1\big\rangle_{100}$} unambiguously as
$\Phih^{(\s_0)}(h)$ but not enough to see a possible contribution from the much
smaller~$\Phih^{(\s_2)}(h)$. If we replace optimal by smooth truncation, then the
error in $\Phih(h)$ decreases from (approximately) $10^{-15}$ to $10^{-44}$ and the
error in $\Phih^{(\s_0)}(h)$ from $10^{-15}$ to $10^{-42}$. We can
therefore compute the difference of the left- and right-hand sides
of~\eqref{FirstRQMC} (for $X=100$, $\g=S$) to 42 digits, finding that it vanishes,
and since the remaining $\Phih$-value $\Phih^{(\s_2)}(h)$ has the much larger order
of~$10^{-14}$, we see that this quantity, if it occurs at all, must have coefficient~0.
But when we replace $X=100$ by $100\frac13$, we find that the difference
\smash{$X^{-3/2} \J(-1/X)-\J\big(\frac13\big)\Phih^{(\s_1)}(h)^{{\rm sm}}
\m \Phih^{(\s_0)}\bigl(h\bigr)^{{\rm sm}}$}
no longer vanishes but instead is equal to $\Phih^{(\s_2)}\bigl(h\bigr)^{{\rm sm}}$
times $-1.732050807568877293527446341{\rm i}$, which coincides to this accuracy with
$-\sqrt{-3}$.
Doing the same for other large values of~$X$ with small denominators and other~$\g$,
we find that~\eqref{FirstRQMC} with all $\Phi$-values interpreted by smooth rather
than optimal truncation can be improved to
\be
\label{41RQMC}
(cX+d)^{-3/2} \J^{(4_1)}\left(\frac{aX+b}{cX+d}\right)
 \overset?\approx  \sum_{j=0}^2 Q^{(4_1)}_j(X)
\Phih^{(4_1,\s_j)}\left(\frac{2\pi {\rm i}}{c(cX+d)}\right),
\ee
where, just as for the $(-2,3,7)$ pretzel knot, \smash{$Q^{(4_1)}_0(x)=1$},
\smash{$Q^{(4_1)}_1(X)=\J^{(4_1)}(x)$} and \smash{$Q^{(4_1)}_2$} is a~1-periodic functions, the
notations in each case being a shorthand for \smash{$Q^{(4_1)}_{\s_j}$}. The first few
values of the periodic functions \smash{$Q^{(4_1)}_i(x)=\sQ^{(4_1)}_i(\e(x))$}
for $i=1$ and $2$ are given by
\begin{center}
\def\arraystretch{1.3}
\begin{tabular}{c|ccccccc}
Ord($q$) & $ 1 $ & 2 & 3 & 4 & 5 & 6\\ \hline
 $\sQ^{(4_1)}_1(q)$ & $1 $ & $ 5 $ & $ 13 $ & $ 27 $
 & $ 44-4 q^2-4q^3 $ & $ 89 $ \\
 $2\sQ^{(4_1)}_2(q)$ & $0 $ & $ 0 $ & $ -2 - 4q $ & $ -14q $ &
 $ -15 - 30 q - 22 q^2 - 8 q^3 $ & $ 46 - 92 q $
\end{tabular}
\end{center}
Just as with the functions~\smash{$Q^{(5_2)}(x)$} and \smash{$Q^{(-2,3,7)}_i(x)$} ($i=3,4,5,6$)
found for the $5_2$ and $(-2,3,7)$ knots in the previous subsection, the
function \smash{$Q^{(4_1)}_2$} (whose values we found by the numerical procedure just
outlined for all $x$ with denominators up to~200), multiplied by~2, turned out to
always belong to $\Z[\e(x)]$ and to satisfy all of the necessary Ohtsuki-type
congruences near~0 and~$1/2$ required for it to be an element of the Habiro ring.
In this case, following a tip by Campbell Wheeler, we were actually able to guess
a formula that reproduced all of the numerically found values and (after
multiplication by~2) was visibly in the Habiro ring, namely the following simple
variant of equation~\eqref{Hab41}:
\be
\label{Q0for41}
\sQ^{(4_1)}_2(q)
= \frac12
\sum_{n=0}^\infty \big(q^{n+1}-q^{-n-1}\big) \big(q^{-1};q^{-1}\big)_n (q;q)_n .
\ee

When we recompute the examples of Section~\ref{sub.4.2} with smooth rather than
optimal truncation, the situation is exactly similar and we are able to add a
\smash{$\Phih_\a^{(K,\s_2)}(h)$} term to the right-hand sides of both~\eqref{52WeakRQMC}
and~\eqref{237WeakRQMC}, obtaining for both knots a conjectural approximate formula
of the form
\be
\label{strongRQMC}
(cx+d)^{-3/2} \J^{(K)}\left(\frac{aX+b}{cX+d}\right)  \overset?\approx
\sum_{\s\in\calP_K}
Q^{(K)}_{\s}(X)
 \Phih^{(K,\s)}_{a/c}
\left(\frac{2\pi {\rm i}}{c(cX+d)}\right) ,
\ee
where \smash{$Q^{(K)}_{\s_0}(x)=1$}
and \smash{$Q^{(K)}_{\s_1}(x) = \J^{(K)}(x)$}.
In Sections~\ref{sub.qholo} and~\ref{sec.arithmetic}, we give more information about
these numbers including a formula for \smash{$Q^{(5_2)}_{\s_2}(x)$}
as an element of the Habiro ring.

\subsection{Strengthening the generalized quantum modularity conjecture}
\label{sub.4.4}

So far we have generalized the original quantum modularity conjecture~\eqref{QMCb}
in two very different ways: in Section~\ref{sub.GQMC}, we extended it from the
Kashaev invariant $\J=\J^{(\s_0)}$ to the functions defined in~\eqref{cherries},
and in the last three subsections we refined it by adding additional terms of lower
order to the right-hand side. Not surprisingly, these two can be combined, but
with some new aspects.

If we repeat the calculations described in the previous subsection (using smooth
truncation for all the $\Phi$-series occurring) but with the function
$\J^{(K)}=\J^{(K,\s_0)}$ replaced by the function defined in~\eqref{cherries} with
$\s\ne\s_0$, then instead of~\eqref{strongRQMC} we find
\be
\label{RQMCa}
{\rm e}^{-\tV(\s)\l_\g(X)}
\J^{(K,\s)}\left(\frac{aX+b}{cX+d}\right)
 \overset?\approx
\sum_{\substack{\s'\in\calP_K\RED }}
\J^{(K,\s,\s')}(X)
 \Phih^{(K,\s')}_{a/c}\left(\frac{2\pi {\rm i}}{c(cX+d)}\right) ,
\ee
where $\J^{(K,\s,\s')}$ are 1-periodic functions on $\BQ$ with
$\J^{(K,\s,\s_1)}(x) =\J^{(K,\s)}(x)$ (cf.\ \eqref{GQMC}).
There are, however, three main differences with~\eqref{strongRQMC}.
The first is that the automorphy factor ${(cX+d)^{3/2}}$ is replaced for $\s\ne\s_0$
by the factor
${\rm e}^{-\tV(\s)\l_\g(X)}$ involving the $\s$-th volume~$v(\s)$ (which is zero for~${\s=\s_0}$).
The second is that the Habiro power series \smash{$\Phih^{(K,\s_0)}_{a/c}$}, which
in~\eqref{strongRQMC} had the coefficient~1, now does not occur at all.
The third is that the new functions $\J^{(K,\s,\s')}(x)$ are now no longer elements
of the Habiro ring
when considered as functions of $q=\e(x)$, as was the case for the functions
\smash{$\J^{(K,\s_0,\s')}(x)=Q_{\s'}^{(K)}(x)$}. But they are still $\Qbar$-valued and have various
``Habiro-like'' properties, including the following:
\begin{itemize}\itemsep=0pt
\item
$\J^{(K,\s,\s')}(x)$ for $x\in\Q/\Z$ is the constant term of a power
series~\smash{$\Phi^{(K,\s,\s')}_x(h)$}
lying in the same space as the power series~\smash{$\Phi^{(K,\s)}_x(h)$}, as discussed briefly
after~\eqref{as41at5} and in more detail in Section~\ref{sec.arithmetic}, i.e.,
it belongs to $\mu \delta^{-1/2}\sqrt[m]{\ve} F_\s(\z_m)[[h]]$ with the same root
of unity $\mu$, the same element $\delta_\s$ of~$F_\s^\times$, the same set~$S$ of
primes of $F_\s$ (independent of~$x$) and the same $S$-unit $\ve=\ve_x$ of~$F_\s(\z_m)$
as for \smash{$\Phi^{(K,\s)}_x(h)$}.
\item
one can interpret~$\Phi^{(K,\s,\s')}_x(h)$ as $\sQ\big(\e(x){\rm e}^{-h}\big)$ where $\sQ$ is an
element of a Habiro ring $\calH_{F_\s}$ generalizing the ordinary Habiro ring
$\calH=\calH_\Q$ whose definition and arithmetic properties will be discussed in a
planned joint paper with Peter Scholze and Campbell Wheeler~\cite{GSZ:habiro}.
In particular, for primes~$p$
that split completely in~$F_\s$, there are congruence properties modulo~$p$ relating,
for instance, the first~$p$ coefficients of $\sQ\big({\rm e}^{-h}\big)$ to the value of~$\sQ(\z_p)$.
\end{itemize}

We can write equation~\eqref{RQMCa} more uniformly by allowing the
case of $\s=\s_0$, but remembering that there is then an automorphy factor
$(cX+d)^{-3/2}$ that is not present for $\s\in\calP_K$. Then all of the formulas
found so far can be collected into a single conjectural formula
\begin{gather}
\label{RQMC}
{\rm e}^{-\tV(\s)\l_\g(X)} (cX+d)^{-\k(\s)}
\J^{(K,\s)}\left(\frac{aX+b}{cX+d}\right)
 \overset?\approx
\sum_{\s'\in\calP_K}
\J^{(K,\s,\s')}(X)
\Phih^{(K,\s')}_{a/c}\left(\frac{2\pi {\rm i}}{c(cX+d)}\right)\!\!\!
\end{gather}
valid for all $\g=\sma abcd\in\SL_2(\Z)$ and all $\s\in\calP_K$ for $X\to\infty$
with bounded denominator, where \smash{$\J^{(K,\s,\s')}$} are periodic functions belonging
to a generalized Habiro ring and satisfying~\smash{$\J^{(K,\s,\s_0)}=\delta_{\s,\s_0}$} and
\smash{$\J^{(K,\s_0,\s)}=Q^{(K)}_\s$} (as introduced in the previous section), with the weight
$\k_\s$ and the multiplier~$\l_\g(X)$ defined as in~\eqref{tweakdef}
and~\eqref{Defwsigma}.

As a concrete illustration of the refined QMC~\eqref{RQMC} we take once again the
figure~8 knot with $\s=\s_1$. Here,~\eqref{RQMC} involves three terms
$\s'=\s_0,\s_1,\s_2$, with two of the three coefficients already known
(the first vanishes and the second is $\J^{(4_1,\s_1)}(X)$) but with the third one being
a~new periodic function on $\BQ$ given explicitly by
\begin{gather}
\J^{(4_1,\s_1,\s_2)}(x)
= \frac {\rm i}{2\sqrt[4]{3}\sqrt{c}}
\sum_{Z^c = \z_6} \big(Z q -Z^{-1} q^{-1}
\big) \prod_{j=1}^c \bigl|1\m q^jZ\bigr|^{2j/c},\nonumber
\\
c=\den(x),\quad q=\e(x) ,\label{Q1for41}
\end{gather}
which is related to the function given in~\eqref{Def41J1} in exactly
the same way as $Q^{(4_1)}_{\s_2}(x)$ and $\J^{(4_1)}(x)$ are related
by~\eqref{Q0for41} and~\eqref{Hab41}. Likewise, the refined QMC~\eqref{RQMC}
for $4_1$ and for $\s=\s_2$ leads to the periodic function
$\J^{(4_1,\s_2,\s_2)}(x) =-{\rm i} \J^{(4_1,\s_1,\s_2)}(-x)$. We also find
the \emph{bilinear identity}
\be
\label{FirstBilin}
\J^{(4_1,\s_1,\s_1)}(x)  \J^{(4_1,\s_2,\s_2)}(x) \m
\J^{(4_1,\s_1,\s_2)}(x)  \J^{(4_1,\s_2,\s_1)}(x)  = 1
\ee
for all $x\in\Q/\Z$. \big(Note also that $\J^{(\s_1,\s_1)}(x)=\J^{(\s_1)}(x)$
and $\J^{(\s_2,\s_1)}(x)=\J^{(\s_2)}(x)$.\big)
This identity will be generalized to all knots in Section~\ref{sec.Matrix}.

\subsection{The refined quantum modularity conjecture}
\label{sub.4.5}

The refinement of the quantum modularity conjecture that we have obtained so far,
equation~\eqref{RQMC}, has two noteworthy aspects. One is that, although we find
new collections of ``Habiro-like'' functions $\J^{(\s,\s')}$ for the asymptotic
expansion as $X\to\infty$ of the functions $\J^{(\s_0,\s)}(\g X)$ for different
$\s\in\calP_K$ (here we continue the practice of omitting the knot from all
notations when it is not varying), these arise as coefficients of the \emph{same}
completed formal power series $\Phih^{(\s')}(h)$ as we found for the initial
Galois-extended Kashaev invariant $\J^{(\s_0)}$. The other is that among the new
coefficients $\J^{(K,\s,\s')}$, the subset corresponding to $\s'=\s_1$ coincides
precisely with the set of functions $\J^{(\s)}$ whose asymptotic behavior near fixed
rational points is being studied. It is therefore natural to ask whether the
functions $\J^{(K,\s,\s')}$ for $\s'$ different from $\s_1$ also have a quantum
modularity property, and if so, what new power series are involved. In this final
subsection, we will study both questions and give our (nearly) final version of the~QMC.

As usual, we look first at the $4_1$ knot. Here, as well as the three periodic functions
$\J^{(\s_0,\s_1)}(x)\allowbreak:=\J(x)$ and $\J^{(\s_0,\s_1)}(x) := \J^{(\s_1)}$ and
$\J^{(\s_0,\s_2)}:=\J^{(\s_2)}$ we had studied earlier, we found two new periodic functions
$\J^{(\s_1,\s_2)}(x)$ and $\J^{(\s_2,\s_2)}(x)$, given explicitly by
formulas~\eqref{Q1for41} and related to the others by~\eqref{FirstBilin}. If we look
numerically at the asymptotics
of both functions with $x=-1/X$ for $X\to\infty$ with bounded denominator, we find
\begin{gather*}
X^{-3/2}\J^{(\s_0,\s_2)}\left(-\frac1X\right)  \sim  \J^{(\s_0,\s_1)}(X)
\Psih^{(1)}\left(\frac{2\pi {\rm i}}{X}\right) , \\
{\rm e}^{-\tV(\s_1)\l_S(X)}  \J^{(\s_1,\s_2)}\left(-\frac1X\right)  \sim
\J^{(\s_1,\s_1)}(X) \Psih^{(1)}\left(\frac{2\pi {\rm i}}{X}\right)
\end{gather*}
with the \emph{same} completed power series
\[
\Psih^{(1)}(h) = {\rm e}^{\V(4_1)/h} \Psi^{(1)}(h),\qquad
\Psi^{(1)}(h)= \frac{{\rm i} \sqrt[4]3}{2}
\left(1 \m \frac{37}{72\sqrt{-3}} h \m \frac{1511}{2(72\sqrt{-3})^2} h^2+\cdots\right)
\]
in both cases. Based on the analogy with the asymptotics of the functions
$\J^{(\s_1,\s_1)}(-1/X)$ as given in~\eqref{41RQMC}, we would now expect
the more accurate
approximations
\begin{align}
 & X^{-3/2}\J^{(\s_0,\s_2)}\left(-\frac{1}{X}\right) \approx \Psih^{(0)}\left(\frac{2\pi {\rm i}}{X}\right) +
 \J^{(\s_0,\s_1)}(X) \Psih^{(1)}\left(\frac{2\pi {\rm i}}{X}\right)
 + \J^{(\s_0,\s_2)}(X) \Psih^{(2)}\left(\frac{2\pi {\rm i}}{X}\right) ,
 \nonumber\\
 & {\rm e}^{-\tV(\s_1) \l_S(X)} \J^{(\s_1,\s_2)}\left(-\frac{1}{X}\right) \approx
 \J^{(\s_1,\s_1)}(X) \Psih^{(1)}\left(\frac{2\pi {\rm i}}{X}\right) +\J^{(\s_1,\s_2)}(X) \Psih^{(2)}\left(\frac{2\pi {\rm i}}{X}\right) ,\nonumber
 \\
 & {\rm e}^{-\tV(\s_2) \l_S(X)} \J^{(\s_2,\s_2)}\left(-\frac{1}{X}\right) \approx
 \J^{(\s_2,\s_1)}(X) \Psih^{(1)}\left(\frac{2\pi {\rm i}}{X}\right) +\J^{(\s_2,\s_2)}(X) \Psih^{(2)}\left(\frac{2\pi {\rm i}}{X}\right) ,
 \label{FirstMatrix}
\end{align}
where $\Psih^{(2)}(h)$ is the completed series ${\rm e}^{-\V(4_1)/h} \Psi^{(2)}(h)$ with
$\Psi^{(2)}(h)=-{\rm i}\Psi^{(1)}(-h)$ and where $\Psih^{(0)}(h)$ is the completed series
$(h/2\pi {\rm i})^{3/2} \Phi^{(0)}(h)$ (cf.\ \eqref{DefPhi0h}) with
\[
\Phi^{(0)}(h) = -h + \frac{11}6 h^3 \m \frac{1261}{120} h^5
+ \frac{611771}{5040} h^7 \m \cdots
\]
the power series in $h \Q\big[\big[h^2\big]\big]$ obtained by replacing $q=\e(X)$ by $q={\rm e}^{-h}$
in formula~\eqref{Q0for41}, as well of course as similar formulas for
\smash{$\J^{(\s,\s')}(\g X)$} for other matrices $\g=\sma abcd$ in $\SL_2(\Z)$ with the
completed power series \smash{$\Psih^{(j)}(h)$} replaced by suitable new completed power
series
\[\Psih_{a/c}^{(j)}(h)={\rm e}^{\V(\s_j)/ch}(h/2\pi {\rm i})^{\k(\s_j)}\Psi_{a/c}^{(j)}(h)
\]
but with the same periodic coefficients.

To test~\eqref{FirstMatrix} or its generalizations to other $\g\in\SL_2(\Z)$
directly we would need to find many terms of the power series \smash{$\Psi_\a^{(j)}(h)$}
and carry out the smoothed optimal truncation as described earlier in this section,
because the different exponential growths of their completions would mean that the
contributions with $j\ne1$ would not be numerically visible at the level of mere
formal power series. This could be done, but an easier test of the prediction is to
take linear combinations of the first two or last two lines in~\eqref{FirstMatrix}
to eliminate the dominant $\Psih^{(1)}$-term. This (together with~\eqref{FirstBilin})
produces the two new asymptotic predictions
\begin{align*}
 &X^{-3/2} \J^{(\s_1,\s_1)}(X) \J^{(\s_0,\s_2)}\left(-\frac{1}{X}\right) \m
 {\rm e}^{-\tV(\s_1) \l_S(X)}\J^{(\s_0,\s_1)}(X)
 \J^{(\s_1,\s_2)}\left(-\frac{1}{X}\right)\\
 &\qquad\approx \J^{(\s_1,\s_1)}(X) \Psih^{(0)}\left(\frac{2\pi {\rm i}}{X}\right) , \\
 &{\rm e}^{-\tV(\s_1) \l_S(X)} \J^{(\s_2,\s_1)}(X) \J^{(\s_1,\s_2)}\left(-\frac{1}{X}\right)\m
 {\rm e}^{-\tV(\s_2) \l_S(X)} \J^{(\s_1,\s_1)}(X) \J^{(\s_2,\s_2)}\left(-\frac{1}{X}\right)\\
 &\qquad \approx \Psih^{(2)}\left(\frac{2\pi {\rm i}}{X}\right) ,
\end{align*}
both of which can be tested directly since they do not involve functions of
different orders of growth on the right, and both of which we confirmed numerically
to very high precision. We omit the details, having given more than enough
descriptions of analogous numerical calculations in this section already.

Generalizing the above discussion to other knots, we find as our nearly-final version
of the QMC the asymptotic statement
\begin{gather*}
(cX+d)^{-\k(\s)} {\rm e}^{-\tV(\s)\l_\g(X)}
\J^{(K,\s,\s')}\left(\frac{aX+b}{cX+d}\right)\\
\qquad \overset?\approx
\sum_{\s'' \in \calP_K}\J^{(K,\s,\s'')}(X)
\Phih^{(K,\s'',\s')}_{a/c}\left(\frac{2\pi {\rm i}}{c(cX+d)}\right)
\end{gather*}
for $X\in\Q$ tending to infinity with bounded denominator and for every
$\sma abcd\in\SL_2(\Z)$ with~${c>0}$, where the functions \smash{$\J^{(K,\s,\s')}$} are the
``Habiro-like'' functions that we found in Section~\ref{sub.4.4}, given
as the constant terms of certain power series \smash{$\Phi^{(K,\s,\s')}_\a(h)\in\Qbar[[h]]$},
and where \smash{$\Phih^{(K,\s,\s')}_\a(h)$} are the completions defined by
\[
\Phih^{(K,\s,\s')}_\a(h) = (ch/2\pi {\rm i})^{\k(\s)} {\rm e}^{\V(\s)/c^2h}
\Phi^{(K,\s,\s')}_\a(h), \qquad \s,\s' \in \calP_K .
\]
To get the final version, we upgrade this statement about constant terms to a statement
about the full (completed) power series in the same way as we did in \ref{sub.QMCwithh},
obtaining:

{\bf Refined quantum modularity conjecture (RQMC):}
\emph{For fixed $\g=\sma abcd\in\SL_2(\Z)$ with~${c>0}$, we have
\[
\Phih^{(K,\s,\s')}_{\g X}(h^*)  \overset?\approx
(cx+d)^{\k(\s)} \sum_{\s'' \in \calP_K}\Phih^{(K,\s,\s'')}_X(h)
\Phih^{(K,\s'',\s')}_{a/c}\left(\frac{2\pi {\rm i}}{c(cx+d)}\right)
\]
to all orders in~$1/X$ as $X\in\Q$ tending to~$+\infty$ with bounded denominator, where}
\[x=X-\hbar$ and~${h^*=h/(cx+d)(cX+d)}.\]

We end this section by observing that the two versions of the refined quantum modularity
conjecture that we just stated can both be written more succinctly in matrix form as
\be
\label{dQPhi}
\bJ^{(K)}(\g X)
 \approx \jt_\g(X) \bJ^{(K)}(X) \bPhih^{(K)}_{a/c}\left(\frac{2\pi {\rm i}}{c(cX+d)}\right)
\ee
and
\be
\label{MatQMC}
\bPhih^{(K)}_{\g X}(h^*)  \overset?\approx  \jp_\g(x)
\bPhih^{(K)}_X(h) \bPhih^{(K)}_{a/c}\left(\frac{2\pi {\rm i}}{c(cx+d)}\right)
\ee
as $X\to\infty$ with bounded denominator for a fixed knot~$K$ and element
$\g=\sma abcd\in\SL_2(\Z)$, where~\smash{$\bJ^{(K)}$} and \smash{$\bPhih^{(K)}$} denote the matrices
of Habiro-like functions and completed formal power series, with columns and rows
indexed by $\calP_K$, with entries \smash{$\J^{(K,\s,\s')}(x)$} and \smash{$\Phi^{(K,\s,\s')}_{a/c}(h)$},
respectively, and where $\jp$ and $\jt$ are the matrix-valued automorphy factors
defined by
\be
\label{autofactor}
\jp_\g(x) = \diag\bigl(|cx+d|^{\kappa(\s)}\bigr) ,\qquad
\jt_\g(x) = \diag\bigl({\rm e}^{\tV(\s) \l_\g(x)} |cx+d|^{\kappa(\s)}\bigr) ,
\ee
the second of which is the ``tweaked'' version of the first.
Note that both of these factors are unchanged if we replace $\g$ by~$-\g$, and
hence are actually automorphy factors on~$\PSL_2(\Z)$. Also, from the fact
that $\l$ is an additive cocycle (see Lemma~\ref{lem.lambda}) we deduce that both
$\jp$ and~$\jt$ are matrix-valued cocycles on~$\PSL_2(\Z)$, meaning that they satisfy
\be
\label{jcoc}
 \jp_{\g \g'}(x) = \jp_{\g'}(x)  \jp_\g(\g' x) , \qquad
\jt_{\g \g'}(x) = \jt_{\g'}(x)  \jt_\g(\g' x) ,
\ee
for all $\g, \g'\in \PSL_2(\BZ)$. This will be important in the next section.


\section{The matrix-valued cocycle associated to a knot}
\label{sec.Matrix}

Let us define, for a fixed knot $K$ (suppressed from the notation as usual),
matrix~$\g\in\SL_2(\Z)$ and number $x\in\Q\ssm\big\{\g^{-1}(\infty)\big\}$,
an $(r+1)\times (r+1)$ matrix~$W_\g(x)$ by
\be
\label{Wdef}
W_\g(x) = \bJ(\g x)^{-1} \jt_\g(x)  \bJ(x)  ,
\ee
where $\jt$ is the automorphy factor defined in~\eqref{autofactor}.
(This formula makes sense because the matrix~$\bJ$ is conjecturally invertible, and
even unimodular, as discussed in~\eqref{Orthogonality1} below.) This
function has remarkable properties. On the one hand, the refined quantum modularity
conjecture as presented above can now be rewritten as the asymptotic statement
\be
\label{WRQMC}
W_\g(X) \approx \bPhih_{a/c}{\left(\frac{2\pi {\rm i}}{c(cX+d)}\right)}^{\i}
\ee
for $X\in\Q$ tending to infinity with bounded denominator. In particular, unlike the
completely discontinuous function~$\bJ(x)$ in terms of which it is defined, $W_\g(X)$
has an asymptotic behavior at infinity that depends only on $X$ as a real number and
not on its numerator and denominator separately, and in Section~\ref{sub.lift}
we will present very strong evidence that this is true not only asymptotically at
infinity, but also for finite values of the argument, so that $W_\g(x)$ becomes a~smooth (and in fact real analytic) function of its argument away from the singularity
at~${x=\g\i(\infty)}$. On the other hand, the cocycle property~\eqref{jcoc} of~$\jt$
immediately implies that the function $\g\mapsto W_\g(\cdot)$ is a cocycle
for the group $\PSL_2(\Z)$ acting on the multiplicative group of
almost-everywhere-defined invertible matrix-valued functions on~$\BP^1(\Q)$, meaning
that it satisfies
\begin{gather}
\label{cocycle}
W_{\g\g'}(x) = W_\g(\g' x)  W_{\g'}(x)
\end{gather}
for all $\g$ and $\g'$ in $\PSL_2(\BZ)$. But this cocycle property then immediately
extends by continuity to imply that $W_\g$ on~$\R$ is also a $\PSL_2(\Z)$-cocycle,
but now in the space of piecewise smooth matrix-valued
functions on~$\BP^1(\R)$. We can then use the smoothness to define a canonical lift
of each of the formal power series \smash{$\Phi_\a^{(\s,\s')}(h)$} to an actual function of~$h$,
simply by requiring~\eqref{WRQMC} to be an exact rather than merely an asymptotic
equality.

These various properties will be described in detail in this section. The first
subsection treats the elementary properties (behavior under complex conjugation,
determinant, and inverse) of the matrices $\bJ(x)$ and $W_\g(x)$. The discussion of
the smoothness properties and the lifting of the perturbative series
\smash{$\Phi_\a^{(\s,\s')}(h)$} to well-defined functions of~$h$ will be given in
Section~\ref{sub.lift}, while the brief final subsection treats the expected equality
between the cocycle~$W_\g(x)$ and the cocycle constructed in the companion
paper~\cite{GZ:qseries} using state integrals, which gives the real explanation for
its smoothness and even analyticity.

\subsection{The Habiro-like matrix and the perturbative matrix}
\label{sub.habmat}

In Section~\ref{sec.RQMC}, we saw how successive refinements of the original quantum modularity conjecture~\eqref{QMCb} led to a whole matrix $\bJ^{(K)}(x)$ of generalized
Kashaev invariants and to a collection of matrices \smash{$\bPhi_\a^{(K)}(h)$} of formal power
series having~\smash{$\bJ^{(K)}(\a)$} as their constant term. The existence of these new matrices
and the description of their properties is the main content of this paper. We emphasize
that, although the refined QMC which led to the definition of these matrices and to the
means of finding them numerically is still conjectural, the matrices themselves are
well-defined, at least in terms of a chosen triangulation: Their first columns are
trivial (a~one followed by $r$~zeros). Their second columns were defined in
Section~\ref{sec.Phi} in terms of the original Kashaev invariant and of the perturbative
series defined in~\cite{DG,DG2}. The further columns of the~$\bPhi$-matrix can also
be given by
Gauss-type integrals like those in~\cite{DG,DG2}, and in principle one could also always
find explicit formulas for the $\bJ$ matrix, as has been written out for the~$4_1$
knot in detail
in Section~\ref{sub.4.3} (equations~\eqref{Q0for41} and~\eqref{Q1for41}) and will be
discussed more generally in \hbox{Sections~\ref{sub.qholo}--\ref{sub.statesum}} in the
context of $q$-holonomic systems, with full details for the knot~$5_2$. In general, it
is not known that these quantities are topological invariants, since their definitions
depend a priori on the choice of an
ideal triangulation and are believed but not proven to be invariant under Pachner moves.
But we expect this invariance to be true, and in any case the new matrices are
completely computable in practice, as we seen, and have extremely interesting
properties. In this subsection, we look at the properties that are directly visible,
and in the following one at the deeper properties of the associated cocycle~$W$.

{\bf Extension property.}
From their very definitions, the matrices $\bJ$ and $\bPhi$ (now omitting the knot from
the notation) both have a $(1+r)\times(1+r)$ block triangular form
\be
\label{Jblock}
\bJ(x) = \begin{pmatrix} 1 & Q(x) \\
 0 & \bJx(x) \end{pmatrix} , \qquad
\bPhi_\a(h) = \begin{pmatrix} 1 & \sQ\bigl(\e(\a){\rm e}^{-h}\bigr) \\
0 & \Phix_\a(h) \end{pmatrix} ,
\ee
where \smash{$Q(x)=\bigl(Q^{(\s)}(x)\bigr)_{\s\in\calPx_K}$} is the vector of length~$r$
whose entries are given by the periodic functions found in Section~\ref{sec.RQMC},
\smash{$\sQ(q)=\bigl(\sQ^{(\s)}(q)\bigr)_{\s\in\calPx_K}$} is the corresponding function
in terms of~${q=\e(x)}$ (which we believe to be elements of the Habiro ring
$\calH\otimes\Q$ and therefore to be defined not just at roots of unity, but also at
infinitesimal deformations of roots of unity), and~$\bJx$ and $\Phix_\a$ are $r\times r$
matrices with rows and columns indexed by the elements of~$\calPx_K$. We are mainly
interested in the larger matrices, but we will sometimes want to consider the ``reduced''
matrices separately because they sometimes occur separately (notably in
Section~\ref{sub.quad}, where only $\calPx_K$ occurs in~\eqref{quadrel}, and in the
statements below about the inverse matrices of~$\bJ\RED$ and $\bPhi_\a\RED$). This
block triangular property, trivial though it is, should have a deeper meaning as the
statement that the $(r+1)$-dimensional objects associated to a knot~$K$ (specifically,
the $q$-holonomic modules that will be the subject of Section~\ref{sec.NewHolonomy}),
with a basis parametrized by the set of representations~$\calP_K$, are in fact
extensions of $r$-dimensional objects with a basis parametrized by~$\calPx_K$
by something one-dimensional.

{\bf Behavior under complex conjugation.}
The next point is the following compatibility with complex conjugation, namely
\be
\label{Jsym1}
\overline{\bJ^{\s,\s'}(-x)}
= \begin{cases} \hphantom{i}\bJ^{\s,\s'}(x)& \text{if $\s$ is real,} \\
 -{\rm i}\bJ^{\overline\s,\s'}(x)&\text{if $\s$ is not real,}\end{cases}
\ee
where
``$\s$ real'' means~$\s=\overline\s$. In matrix form this becomes
\be
\label{Jsym}
\overline{\bJ(-x)} = B \bJ(x) ,
\ee
where $B$ is the unimodular symmetric unitary matrix with $B^{(\s,\s')}$
equal to~1 if ${\s'=\s=\overline\s}$, to~$-{\rm i}$ if $\s'=\overline\s\ne\s$, and to~0
if $\s'\ne\overline\s$.
The symmetry~\eqref{Jsym1} has as the special case ${(\s,\s')=(\s_0,\s_1)}$ the behavior
$\overline{\J(-x)} = \J(x)$ of the Kashaev invariant for rational numbers $x$ under
complex conjugation, which holds because the colored Jones polynomials have real (even
integer) coefficients or alternatively because $\J$ is an element of the Habiro ring.
The same argument applies conjecturally to all $J^{(\s_0,\s)}$, since they also belong
to the Habiro ring, and if we use the full RQMC it also suffices to establish
the general case. Actually, equation~\eqref{Jsym} can be strengthened to
\be
\label{Jsymb}
\overline{\bPhi_{-\a}\big(\overline h\big)} = B  \bPhi_{\a}(h),
\ee
which specializes at $h=0$ to~\eqref{Jsym}. We remark that equation~\eqref{Jsymb}
remains true if we replace both $\bPhi$'s by their completions $\bPhih$ (except for
the top rows, which differ by a factor of~${\rm i}$), because~${\overline{V(\s)}=V(\overline\s)}$
for all $\s\in\calPx_K$.

{\bf Unimodularity.}
The next statement, generalizing equation~\eqref{FirstBilin}, is that the matrices $\bJ$
and even $\Phix$, are experimentally found to be unimodular. More precisely, this is
definitely true for the $4_1$ and $5_2$ knots, for which we have closed formulas for
all of the entries of the Habiro-like matrices and can compute numerically; for other
knots, we are convinced, and willing to conjecture, that the determinant is~$\pm1$,
but we have no really convincing reason except aesthetics that it should always be~$+1$.
The unimodularity implies in particular that the $\bJ$-matrices are always invertible,
a fact that is of course crucial even to define the cocycle~$W$ in~\eqref{Wdef}.
Notice that it is compatible with~\eqref{Jsym} and~\eqref{Jsymb}, since
the matrix~$B$ is unimodular.

{\bf Inverse/Unitarity.}
The final property that we want to mention, again only conjectural, is more surprising.
This is that the inverse of $\bJx$ (but not of the full matrix~$\bJ$, for which we
have no corresponding guess) can be given explicitly by the formula
\be
\label{Orthogonality1}
\bJx(x)^t \bJx(-x) = \overline{B\RED} ,
\ee
where we have set $B=\sma 100{B\RED}$. In fact, even this statement can be strengthened,
namely to
\be
\label{Orthogonality2}
\Phix_\a(h)^t \Phix_{-\a}(-h) = \overline{B\RED} ,
\ee
which specializes to~\eqref{Orthogonality1} when we set~$h=0$. In view of~\eqref{Jsym},
the first of these equations can be rewritten as
\[
 \overline{\bJx(x)^t} B\RED \bJx(x) = B\RED ,
\]
which we see as a kind of unitarity or rather sesqui-unitarity, since if~$B\RED$ were
the identity matrix they would simply say that the matrices $\bJx(x)$ and~$\Phix_\a(h)$
are unitary.
Note that equation~\eqref{Orthogonality2}
remains true also if we replace $\Phi$ by $\Phih$, since the volume factors cancel,
and also that the very special case~$\s=\s'=\s_1$, for which
the right-hand side of~\eqref{Orthogonality2} vanishes, is just the quadratic
relation~\eqref{quadrel} that was discussed in Section~\ref{sub.quad}. It is
also worth remarking explicitly that the expression on the left of~\eqref{Orthogonality2}
is a priori an element of $\Q(\z_c)[[h]]$ ($c=\den(\a)$), at least if the predicted
algebraic properties of the power series \smash{$\Phi_\a^{(\s,\s')}$} as discussed in
Section~\ref{sub.field} are true, because the extra factors (root of unity and $c$-th
root of an $S$-unit in $F_{\s''}(\z_c)$) cancel in the products
$\Phi^{(\s'',\s)}(h)\Phi^{(\s'',\s')}(-h)$ and because the sum over~$\s''$ implicit
in the matrix multiplication gets us from $F_{\s''}(\z_c)$ down to~$\Q(\z_c)$.
Finally, we should mention that equation~\eqref{Orthogonality1}
also gives us a formula for the inverse of the full matrix $\bJ(x)$, because of the block
triangular form of the latter as given in~\eqref{Jblock}, namely
\[ \bJ(x)\i = \begin{pmatrix} 1 & -Q(x)B\RED \bJx(-x)^t \\
 0 & B\RED \bJx(-x)^t \end{pmatrix} ,
\]
in which the elements of the top row are bilinear in the entries of $\bJ(x)$ and
$\bJ(-x)$ rather than merely linear as is the case for the other rows.

However, the real interest to us of the final point above is not just that there
are explicit formulas for the inverses of the matrices~$\bJx$ (or even~$\bJ$)
and $\bPhi$, but above all that the inverse of $\bJx$ is expressed \emph{linearly} (more
correctly, sesquilinearly) in terms of the entries of the matrix itself.
This means in particular that the entries of the cocycle $W\RED$ (the bottom right
$r\times r$ block of~$W$) are expressed \emph{bilinearly} in terms of those of~$\bJ\RED$.
This remark will come into its own in the sequel~\cite{GZ:qseries}, where this reduced
cocycle will arise in a completely different way as a bilinear combination in functions
of $q={\rm e}^{2\pi\tau}$ and $\wt q={\rm e}^{-2\pi/\tau}$ as a consequence of the factorization of
state integrals.

We end this subsection by listing the properties of the function $W_\g$ that it
inherits by virtue of its definition~\eqref{Wdef} from the corresponding properties
of~$\bJ$ listed above. These will become important in the next subsection, when we
extend $W_\g$ from~$\Q$ to~$\R$.

The ``extension property'' is immediate: the matrix $W_\g(x)$ has the block triangular
form \smash{$\sma 100{W_\g\RED(x)}$} for an $r\times r$ ``reduced'' matrix $W_\g\RED(x)$ which
is again a cocycle. The complex conjugation property for $W$ takes the form
\be
\label{Wsym1}
\overline{W_\g(x)} = W_{\v\g\v}(-x)
\ee
for all $\g \in \SL_2(\BZ)$ and $x \in \BQ$, where $\v =\sma{-1}{0}{\hphantom{-}0}{1}$.
(Note that $\v \g \v \in \SL_2(\BZ)$.) This is a~consequence of the following
short calculation using~\eqref{Jsym} and the easy conjugation behavior of the automorphy
factor~\smash{$\jt_\g(x)$}:
\[
\overline{W_\g(x)} = \overline{\bJ(x)^{-1}\jt_\g(-x)\bJ(\g(-x))} =
\bJ(x)^{-1} B^{-1} \overline{\jt_\g(-x)} B \bJ(-\g(x)) = W_{\v\g\v}(-x) .
\]
A nice consequence of~\eqref{Wsym1} is that we can now extend equation~\eqref{WRQMC},
which described the asymptotics of~$W_\g(X)$ as $X$ tends to infinity with bounded
denominator on the assumption of the RQMC, to give the corresponding asymptotic behavior
of $W_\g(X)$ also as $X\to-\infty $:
\be
\label{NegativeX}
W_\g(X) \approx B \bPhih_{a/c}{\left(\frac{2\pi {\rm i}}{c(cX+d)}\right)}^{\i},
\qquad X\to-\infty .
\ee
The third property is that the determinant of~$W_\g(x)$ is given by
\[
\det W_\g(x) = |j(\g,x)|^{-3/2} ,
\]
where $j(\g,x)$ is defined as $cx+d$ for $\g=\sma abcd$. This follows from the
(conjectural) unimodularity of~$\bJ$, the definition of~\smash{$\jt_\g(x)$}, and the fact that
\smash{$\sum_{\s\in\calP_K}v(\s)$} vanishes (``Galois descent'').
Finally, from~\eqref{Orthogonality1} we immediately deduce the
corresponding formula for the inverse matrix of~$W\RED_\g(x)$:
\[
W\RED_\g(x)^{-1} = \overline{B\RED} W\RED_{\v\g\v}(-x)^t .
\]
In this connection, we note that \eqref{Wdef} and~\eqref{Orthogonality1} also imply that
\[
\Wx_\g(x) = \big(\overline{B\RED}\big)^{-1} \bJ\RED(-\g x)^t \jt\RED_\g(x) \bJ\RED(x)
  .
\]
In other words, $\Wx$ is bilinear in the entries of the matrices~$\bJ$, an important
property that is also shared by the functions defined by state integrals.

\subsection{Smoothness}
\label{sub.lift}

\begin{figure}[htpb!]\centering
\includegraphics[height=0.20\textheight]{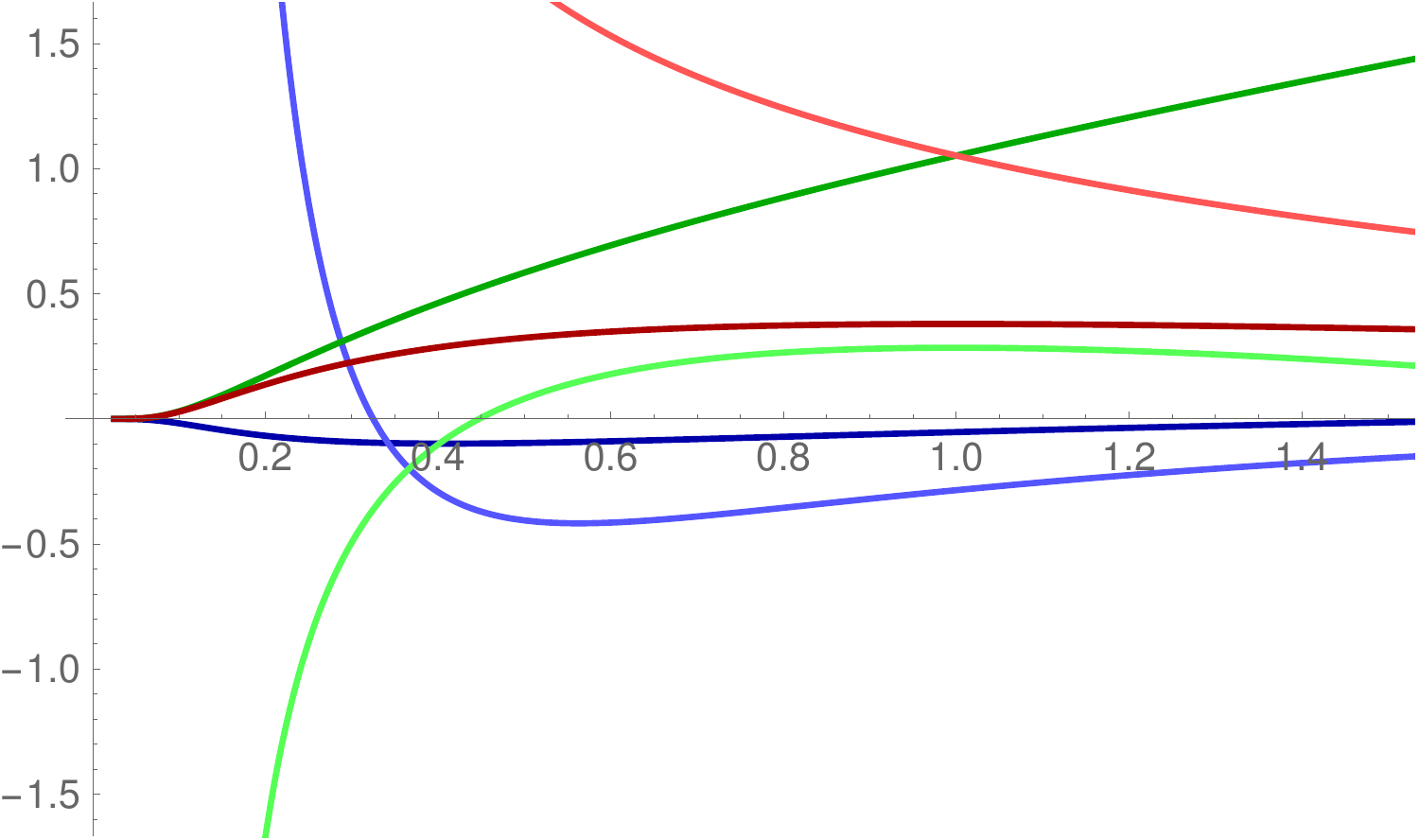}
\caption{Plots of the six nontrivial entries of the matrix $W_S(x)$ for the $4_1$ knot.}\label{WS41plot}
\end{figure}

We now come to the really exciting point. The cocycle~$W_\g(x)$ is defined in terms of
the ``Habiro-like'' matrix~$\bJ$ by~\eqref{Wdef}. The entries of $\bJ^{(4_1)}$, one of
which was shown in Figure~\ref{OldJplot} of the introduction, would all have a
``cloudlike'' structure like the one seen there. But when one graphs the entries of
the matrix $W_\g(x)$, they are all smooth! For instance, Figure~\ref{WS41plot} shows
the graphs of the six nontrivial components of the $3\times3$ matrix $W_S(x)$ for
the figure~8 knot (with three of them divided by~$i$ to make them real), where
$S=\sma0{-1}10$ as usual, plotted in each case for all rational numbers in $(0,2]$
with denominator at most~40 (so for roughly 1000 data points).

We formalize this by stating the following conjecture:
\begin{Conjecture}\label{conj.Wsmooth}
The function $W_\g$ defined on $\BQ\ssm\big\{ \g^{-1}(\infty) \big\}$
extends to a $C^\infty$ function on $\BR\ssm\big\{ \g^{-1}(\infty) \big\}$.
\end{Conjecture}

A first consequence of this is that the cocycle property~\eqref{cocycle},
which held for the restriction of~$W$ to~$\BP^1(\Q)$ by equation~\eqref{Wdef} and the
cocycle property of~$\jt$, is then automatically true for the extended function
on~$\BP^1(\R)$, even though there is no longer any ``coboundary-like formula'' of
type~\eqref{Wdef}. This new cocycle now takes values in the much smaller space of
almost-everywhere-defined matrix-valued functions on $\BP^1(\R)$.

The conjectural smoothness of the function~$W_\g$ has another important consequence
that was already mentioned in the introduction to this section, namely that we can
invert the asymptotic statement~\eqref{WRQMC} to get a definition of \emph{exact}
matrix-valued functions for all \hbox{$\a\in\Q$} which are smooth on all of~$\R$ and
whose Taylor expansions (after the ``Wick rotation'' \hbox{$h\mapsto\hbar=h/2\pi {\rm i}$})
agree with the divergent power series~$\bPhi_\a(h)$. To do this, we simply define a
new function $\bPhi_\a\EXACT$ by requiring~\eqref{WRQMC} to be an exact rather than
just an asymptotic equality, i.e., by defining
\be
\label{PhiExact}
\bPhi_\a(2\pi {\rm i}x)\EXACT := \diag\bigl(|cx|^{-\k(\s)}{\rm e}^{-v(\s)/c^2x}\bigl)
 W_\g\left(\frac1{c^2x}-\frac dc\right)\i
\ee
for any $\g=\sma abcd \in\SL_2(\Z)$ with $\a=\g(\infty)=a/c$.
This should then be an everywhere smooth and almost everywhere analytic function
on~$\R$ whose Taylor expansion at~0 agrees with the original divergent series. The
six functions obtained in this way from the non-trivial elements of $W_S(x)$ for the
figure~8 knot as plotted above (and multiplied by suitable powers of~${\rm i}$ to make them
all real) are the ones shown in Figure~\ref{LiftedPhi} in the introduction.

Apart from the numerical data, there are at least four reasons why we
should expect this smoothness property of the function~$W_\g$ to hold:
\begin{itemize}\itemsep=0pt
\item[{\bf 1.}] At the simplest level, equation~\eqref{WRQMC} tells us that the matrix
$W_\g(X)$ is at least asymptotically smooth in the limit as $X\to\infty$ through rational
numbers with bounded denominator, since it agrees to all orders in~$1/X$ with a
power series in~$1/X$ with coefficients that do not depend on the denominator or
other arithmetic properties of~$X$. Stated more visually, if we were to display
the components of $W_\g(X)$ by plotting their values, for instance, for
rational values of~$X$ between 1000 and~1001 and with denominators less than 100,
then these data points would have to lie on a very smooth curve to very high precision.
\item[{\bf 2.}] In fact this same argument can be pushed much further, since by using
the cocycle property~\eqref{cocycle} for $x=X$ tending to infinity with bounded
denominator we get a description of the asymptotic behavior of~$W_g(x)$ in the
neighborhood of any rational point, not merely at infinity, and hence an explicit
formula for its Taylor expansion at any rational point near which it has a smooth
expansion. This will be carried out in Proposition~\ref{prop.W1} below.
\item[{\bf 3.}] But the real reason that we expected the smoothness property is
much deeper and also predicts (and in some cases leads to a proof of) much more: the
entries of the matrix-valued function $W_\g(\cdot)$ for a fixed~$\g$ extend to functions
that are not merely smooth, but actually \emph{analytic}, on~$\R\ssm\{-d/c\}$. This
comes from the study of $q$-series associated to a~knot and their relation to state
integrals, as carried out in the companion paper~\cite{GZ:qseries} to this one,
and will be discussed in more detail in the final subsection of this section.
\item[{\bf 4.}] Finally, once one expects the real-analyticity, one can check it numerically
using only the matrices studied in this paper, without any reference to either $q$-series
or state integrals, by computing many Taylor coefficients of $W_\g$ at any rational point
using Proposition~\ref{prop.W1} and seeing that they now grow only polynomially rather
than factorially. This point too will be discussed in more detail
in Section~\ref{sub.holomorphic} below.
\end{itemize}

In the context just of this paper, where we are considering only functions
on~$\Q$ and formal power series in~$h$, but not holomorphic functions in the
upper or lower half-planes or on cut planes, we cannot justify the statement
about analyticity or even continuity of the entries of the matrix $W_\g(\cdot)$, i.e.,
we cannot show that the function $W_\g(\cdot)\colon \Q\ssm\{-d/c\}\to\C$ has
any natural extension to a matrix-valued function on~$\R\ssm\{-d/c\}$.
However, as indicated in point~{\bf 2} above, we can deduce a weaker statement
if we assume the RQMC. To explain what this means, we must first discuss the various
possible senses in which a function $f\colon\Q\to\C$ can be continuous or differentiable.
There are at least three different natural notions. Usually one considers the set
of rational numbers with either the discrete topology or else the topology inherited
from their embedding into the reals. In the first sense, of course every
function from $\Q$ to~$\C$ is continuous (i.e., $f(\a+\v_i)\to f(\a)$ for any
sequence of rational numbers $\a+\v_i$ converging to $\a\in\Q$, since any such sequence
is eventually constant) and in fact even $C^\infty$ (with the ``Taylor expansion''
of $f$ at an arbitrary rational point $\a$ being just the constant power series~$f(\a)$).
In the second sense, one means that $f(\a+\v_i)\to f(\a)$ or
$f(\a+\v_i)=P_{\a,d}(\v_i)+o\big(\v_i^d\big)$ as $i\to\infty$ for every~$\a$ and
every $d\in\N$, where $P_{\a,d}$ is a polynomial of degree~$d$ and the $\v_i$
are a sequence whose absolute values tend to~0 as $i$ tends to~$\infty$. Such
a function of course need not extend as a~$C^\infty$ or even continuous
function to~$\R$ \big(an obvious counterexample being $f(x)=1/\big(x-\sqrt2\big)$\big),
but if it does then this extension is unique, so that the space
$C^\infty_{\text{strong}}(\Q)$ of smooth functions in
this sense contains $C^\infty(\R)$ as a~subspace. But there is a~third, weaker, sense, in
which one requires~${f(\a+\v_i)=f(\a)+o(1)}$ or~${f(\a+\v_i)=P_{\a,d}(\v_i)
+o\big(\v_i^d\big)}$ only for sequences $\{\v_i\}$ of rational numbers that have bounded
\emph{numerators} but denominators tending to infinity (so that in particular they tend
to~0 in the usual sense). We then have the strict inclusions
\[
C^\infty(\R)\subsetneqq C^\infty_{\text{strong}}(\Q) \subsetneqq
C^\infty_{\text{weak}}(\Q)\subsetneqq C^\infty_{\text{discrete}}(\Q) = \C^\Q .
\]
An example (courtesy of Peter Scholze) of a function $f\colon\Q\to\R$ that is $C^\infty$ in
the weak sense but not in the strong sense is given by choosing a sequence of rational
numbers $\{x_n\}$ tending to~0 and disjoint intervals $I_n\ni x_n$ with $I_n$ containing
no rational numbers with numerator~$\le n$; then define $f$ to be 0 at~$x=0$ and to be
the restriction of a $C^\infty$ function on $\R^*$ supported on~$\bigcup_nI_n$ and with
$f(x_n)=1$, in which case $f$ is obviously smooth in the strong sense always from~0 and
in the weak sense at~0 (since the values of $f$ on any sequence of rational numbers
tending to~0 with bounded numerators stabilizes to~0), but is not even continuous at~0.

After this lengthy preliminary discussion, we can state the result on the smoothness
properties of the cocycle~$W_\g$, with an explicit formula for the power series of
$W_\g(\a+\v)$ near any~$\a\in\Q$. We will write $\v$ as $-\hbar$ to match our previous
conventions.

\begin{Proposition}[assuming RQMC]\label{prop.W1}
The function $W_\g$ belongs to $C^\infty_{\text{\rm weak}}\big(\Q\ssm\big\{\g\i(\infty)\big\}\big)$ for
every $\g\in\PSL_2(\Z)$. Explicitly, $W_\g(\a-\hbar)$ for $\a\ne\g\i(\infty)$ and
$\hbar$~tending to~$0$
with bounded numerator is given to all orders in~$\hbar$ by the power series
\be\label{WTaylor}
W_\g(\a-\hbar) \approx \bPhi_{\g\a}(2\pi {\rm i}\hbar^*)\i
 \diag\left( \left|\frac{\den(\g\a) \hbar^*}{\den(\a) \hbar}\right|^{\k(\s)}
 {\rm e}^{v(\s)\l_\g(\a)}\right) \bPhi_\a(2\pi {\rm i}\hbar) ,
\ee
with $\hbar^*=\hbar/((c\a+d)(c\a-c\hbar+d))$ if $\g=\sma abcd$.
\end{Proposition}
\begin{proof}
By the definition of weak smoothness on~$\Q$, we have to show that $W_\g(\a+\v)$
for fixed~${\g\in\SL_2(\Z)}$ and~$\a\in\Q$ is given to all orders by a power series in~$\v$
depending only on~$\a$ and~$\g$ as $\v$ tends to~0 through rational numbers with bounded
numerator. If we write $\a$~as~$a'/c'$ with $a'$ and~$c'$ coprime and extend
$\bigl(\!\smallmatrix a'\\c'\endsmallmatrix\!\bigr)$ to a matrix
$\g'=\sma{a'}{b'}{c'}{d'}\in\SL_2(\Z)$, then the condition of~$\v$ having bounded
numerator is easily seen to be equivalent to the condition that $\a+\v=\g'X$
with $X$ tending to~$\pm\infty$ with bounded denominator. We consider first the case
when~$X\to+\infty$ (meaning that $\a+\v$ tends to~$\a$ from the left). By the cocycle
property~\eqref{cocycle} and the basic asymptotic property~\eqref{WRQMC} of~$W_\g$, we
have
\begin{align*}
W_\g(\a+\v) &= W_\g(\g'X) = W_{\g\g'}(X) W_{\g'}(X)\i \notag \\
 & \approx \bPhih_{a''/c''}\left(\frac{2\pi {\rm i}}{c''(c''X+d'')}\right)\i
  \bPhih_{a'/c'}\left(\frac{2\pi {\rm i}}{c'(c'X+d')}\right) ,
\end{align*}
where we have written $\g=\sma abcd$ and $\g\g'=\g''=\sma{a''}{b''}{c''}{d''}$ and where
$\approx$ as usual means that the two expressions being compared are equal to
all orders in~$1/X$ as $X$ tends to infinity with bounded denominator or $\v$ to zero
with bounded numerator. If we now use our previous conventions, writing
\[
-\v = \hbar = \frac1{c'(c'X+d')} , \qquad
\hbar^* = \frac1{c''(c''X+d'')} = \frac{{c'}^2\hbar}{c''(c''-cc'\hbar)} ,
\]
and also use that the ``tweaking function'' $\l_\g$ satisfies
\be
\label{tweakCoboundary}
\frac1{{c'}^2\hbar} \m \frac1{{c''}^2\hbar^*}
= \left(X+\frac{d'}{c'}\right)\m \left(X+\frac{d''}{c''}\right)
=\frac{c}{c'c''} = \l_\g(\a) ,
\ee
then we get equation~\eqref{WTaylor} above. This equation expresses $W_\g(\a+\v)$ as
a product of three matrices of power
series in~\hbox{$\hbar=-\v$} and hence shows that it is itself such a matrix. This
proves the assertion in the first case~\hbox{$X\to+\infty$}.
To treat the case~$X\to-\infty$ we use equation~\eqref{NegativeX} instead
of~\eqref{WRQMC} and find
that $W_\g(\a-\v)$ is given by the \emph{same} formula as a product of three matrices
of power series for $\v>0$ as it was for~$\v<0$, because the prefactors $B$
in~\eqref{NegativeX} cancel. This completes the proof that $W_\g$ is a two-sided
smooth function on the rational numbers in the weak sense.
\end{proof}

\begin{Corollary}
\label{PhiExactSmooth}
The function $x\mapsto\bPhi_\a(2\pi {\rm i}x)\EXACT$ on~$\Q$ defined by~\eqref{PhiExact}
is differentiable in the weak sense for every $\a\in\Q$.
\end{Corollary}

\begin{proof}
This follows directly from Proposition~\ref{prop.W1} away from $x=0$, since the
diagonal pre\-factor in~\eqref{PhiExact} is smooth away from~0, and
$\bPhi_\a(2\pi {\rm i}x)\EXACT$ simply agrees with $\bPhi_\a(2\pi {\rm i}x)$ to all orders in~$x$
as~$x\to0$, so it is smooth there too.
\end{proof}

As a final remark in this subsection, we recall that in order to specify
the cocycle $\g\mapsto W_\g$ completely, it suffices to give its values for
the two special matrices~$T=\sma1101$ and $S=\sma0{-1}1{\hphantom{-}0}$, since these generate the
whole modular group. The function $W_T(x)$ is elementary (constant and conjugate
to a diagonal matrix of $N$th roots of unity, where $N$ is the level of the knot).
In the case where the level is~1, such as the $4_1$ or $5_2$ knots, $W_T$ is
simply the identity matrix and the whole cocycle is determined by the single
matrix-valued function $W_S(x)$. The fact that~$T\mapsto1$,~$S\mapsto W_S$ extends to
a cocycle on the whole group is then equivalent to the requirement that~$W_S(x)$
satisfies the symmetry property $W_S(x) = W_S(-1/x)^{-1}$,
together with the three-term Lewis functional equation
\[
W_S(x) = W_S(x+1) W_S(x/(x+1)) ,
\]
familiar from the theory of period polynomials of holomorphic modular forms on
$\SL_2(\Z)$ or of period functions in the sense of~\cite{Lewis-Zagier} of Maass forms
on~$\SL_2(\Z)$. Our cocycles thus belong in some sense to the same family as periods
of modular forms.

\subsection[``Functions near Q'']{``Functions near $\boldsymbol{\Q}$''}\label{sec.nearQ}

We now come to an important and somewhat subtle point. In the calculation that we gave
to prove Proposition~\ref{prop.W1}, we used only the refined quantum modularity
conjecture in its ``rational version''~\eqref{dQPhi}, since the statement of
Proposition~\ref{prop.W1} involves only the values of~$W_\g$ at rational arguments.
If we had used instead the full version~\eqref{MatQMC} of the RQMC, we would have
obtained a stronger version of the ``weakly smooth'' condition that applies to
approximating a rational number not just by rational numbers that differ from it by
a small rational number with bounded numerator, but also by infinitesimal variations
of such numbers. To make sense of a statement of this type, we now introduce a notion
that will shed more light on the two cocycles $\g\mapsto \l_\g$ and of~$\g\mapsto W_\g$
and that is also relevant in connection with the notion of ``holomorphic quantum
modular forms'' that will be touched on briefly in Section~\ref{sub.holomorphic} and
developed in more detail in~\cite{GZ:qseries} and in the survey paper~\cite{Z:HQMF}.
This is the notion of \emph{asymptotic functions near~$\Q$}. The basic idea here is to
specify a particular type of asymptotic behavior (such as a formal power series)
in an infinitesimal neighborhood of every rational point, where ``neighborhood'' can mean
that we approach the rational number from the right and the left on the real line,
or in other contexts from above and below in the upper and lower complex half-planes.
Since there are many types of behavior that may be of interest, and since it is hard
 to give a general definition that includes all of the examples and all of the properties
that one wishes to include, we will restrict here to the particular classes that arise
in the context of knot invariants.

The simplest version of this notion is just given by a collection
$\{f_\a(\ve)\}_{\a\in\Q}$ of formal power series with complex coefficients indexed by
the rational numbers. Here we want to think of the infinitesimal power series
variable~$\ve$ as the difference between the rational number~$\a$ and an
infinitesimally nearby real ``number'' $\a+\ve$, i.e., we want to think of the whole
collection of power series $\{f_\a\}$ as a single ``asymptotic function near~$\Q$'',
i.e., as a ``function''~$f$ defined in infinitesimal neighborhoods of all rational
points~$\a$ by $f(\a+\ve)=f_\a(\ve)$. Of course $f$ is not a function at all in the
traditional sense, since one cannot evaluate it at numerical values of its argument,
but as we will see in a moment, this point of view is nevertheless very fruitful. It
originally showed up in the paper~\cite{Za:QMF}, where ``quantum modular forms'' were
first defined simply as functions on~$\Q$ (more precisely, as almost-everywhere-defined
functions on~$\Q$) but then upgraded to a notion of ``strong quantum modular forms''
where the original values at rational numbers became the constant terms of a
collection of formal power series.

The set of asymptotic functions near~$\Q$ (from now on we omit the quotation marks,
trusting the reader to remember that these are not actually functions) of this special
type forms a ring via pointwise addition and multiplication if we think of its elements
as collections of formal power series, and by straight addition and multiplication if
we think of them as functions defined in infinitesimal neighbourhoods of all rational
points. To understand its elements, it is helpful to think of the following two
extreme cases.
\begin{itemize}\itemsep=0pt
\item[(i)]
 Each $f_\a(\ve)$ is the Taylor expansion $\sum f^{(n)}(\a) \ve^n/n!$ of a function
 $f\in\C^\infty(\R)$ at the point~$\a$. Here the various asymptotic expansions near
 rational points fit together nicely into a single smooth function on~$\R$.
\item[(ii)]
 Each $f_\a(\ve)$ is the formal power series expansion at $q={\rm e}^{2\pi{\rm i}(\a+\ve)}$ of an
 element $A(q)$ of the Habiro ring $\calH =\varprojlim\Z[q]/((q;q)_n)$. Here the
 different power series do not in general fit together smoothly at all, and even their
 constant terms jump around wildly, as illustrated in Figure~\ref{OldJplot} in the
 introduction.
\end{itemize}

At first sight this definition seems to be pointless because we are not requiring
any compatibility at all between the different power series $f_\a$ and therefore the
ring we have just introduced is canonically isomorphic to the direct product
$\C[[\v]]^\Q=\prod_\a\C[[\v]]$ of one copy of the power series ring $\C[[\v]]$
for every rational number~$\a$. The point, however, is that if we pass to the quotient
ring $\FF_{0,0}\approx\prod_\a\C[[\v]]/\bigoplus_\a\C[[\v]]$ of asymptotic functions
in the neighborhood of all but a finite set of rational points, then the modular
group $\G_1=\SL_2(\Z)$ acts by sending $f$ to $f\circ\g$ for~${\g\in\G_1}$, and this
action does not simply permute the different power series~$f_\a$ but twists them
as well. Specifically, if $f(x)$ is represented near~$\a$ by $f(\a+\v)=f_\a(\v)$
then $f(\g(x))$ is represented near $\a^*=\g(\a)$ by $f_{\a^*}(\ve^*)$ rather than
simply by $f_{\a^*}(\ve)$, where $\ve^*=\g(\a+\ve)-\g(\a)$, or more explicitly
$\v^*=\v/(c\a+d)(c\a+d+c\v)$ if~$\g=\sma abcd$. This is precisely the twist that
we already encountered in Section~\ref{sub.QMCwithh} (equations~\eqref{GQMCh}
and~\eqref{GQMChh}) and in Proposition~\ref{prop.W1} above, except that we have
changed the previous variable~$h$ to $-2\pi {\rm i}\ve$ here. (The rescaling of~$h$ by a
factor $2\pi {\rm i}$ was introduced for our knot invariants only to make the power series
coefficients algebraic and there is no reason to make this change of variable in
the general situation.)

We now generalize the above notion by introducing two complex parameters $v$ and~$\k$
and considering the vector space of asymptotic functions near~$\Q$ whose local
form $f_\a(\v)=f(\a+\v)$ in a real infinitesimal neighborhood of any $\a\in\Q$ is
given by
\be
\label{DefFkv}
f_\a(\v) = |\den(\a) \v|^\k {\rm e}^{-v/\den(\a)^2\v} \phi_\a(\v)
\ee
for some power series $\phi_\a(\ve)\in\C[[\v]]$. Again we pass to the quotient
$\FF_{\k,v}$ of almost-everywhere-defined asymptotic functions on~$\Q$, i.e., we
identify two collections of completed power series if they differ for only finitely
many~$\a$. The space $\FF_{\k,v}$ is a free module of rank~1 over the ring~$\FF_{0,0}$
introduced above, and is again isomorphic to $\prod_\a\C[[\v]]/\bigoplus_\a\C[[\v]]$
via $f\mapsto\{\phi_\a\}$, but with a~different action of $\SL_2(\Z)$ than before.
Specifically, $\g\in\SL_2(\Z)$ sends $f$ to the function near~$\Q$ that corresponds
via~\eqref{DefFkv} to the collection of power series
$\{\phi^*_\a(\v)={\rm e}^{v\l_\g(\a)}\phi_{\a^*}(\v^*)\}$ with $\a^*$ and~$\v^*$ as above
and with the ``tweaking cocycle'' $\l_\g(\a)$ introduced in~\eqref{tweakdef}.
Alternatively, in terms of the variable $x=\a+\v$ infinitesimally near~$\a\in\Q$ we
can write this action as the ``slash action'' (familiar from the theory of modular
forms if~$\k$ is an even integer, and from the theory of the principal series
representation of~$\SL_2(\R)$ if not) given by $(f|_\k\g)(x) = |cx+d|^{-\k}f(\g x)$
for~${\g=\sma abcd\in\SL_2(\Z)}$, where $|cx+d|^{-\k}:=|c\a+d|^{-\k}(1+c\v/(c\a+d))^{-\k}$.
This also explains the reason for including the perhaps strange-looking factors
$\den(\a)$ and $\den(\a)^{-2}$ in~\eqref{DefFkv}, since without them there would
be no action of the modular group. We also point out that, because of the ``tweaking''
factor ${\rm e}^{v\l_\g(\a)}$ in the definition of the action, $\FF_{\k,v}$ is a free module
of rank~1 over the ring $\FF_{0,0}$ as a vector space, but not as
an $\SL_2(\Z)$-module: one cannot choose a generator in an $\SL_2(\Z)$-invariant way.

Of course from the point of view of this paper the main reason for introducing the
parameters~$\k$ and $v$ and the definition~\eqref{DefFkv} is that this is exactly the
behavior that we found experimentally from the refined quantum modularity conjecture,
with $v=v(\s)$ and $\k=\k(\s)$ being the normalized volume and weight associated to
a parabolic flat connection~$\s$ and with~$\phi_a(\v)$ and~$f_\a(\v)$ being the power
series \smash{$\Phi_\a^{(K,\s)}(h)$} \big(or more generally \smash{$\Phi_\a^{(K,\s,\s')}(h)$}\big) and its
completion~\smash{$\Phih_\a^{(K,\s)}(h)$} \big(or \smash{$\Phih_\a^{(K,\s,\s')}(h)$}\big) as
in~\eqref{UnifPhih}, with $h=-2\pi\v$. But it is worth noting that the space
$\FF_{\k,v}$ also contains classical modular forms on the full modular group, since
a holomorphic modular form~$f(\tau)$ of (necessarily even) weight~$k$ on~$\SL_2(\Z)$
canonically defines a function near~$\Q$ of type~\eqref{DefFkv} with $\k$ equal
to~$-k$, with $v$ equal to $-2\pi {\rm i}$ times the valuation of~$f$ at infinity (the
smallest exponent of $q=\e(\tau)$ in the Fourier expansion of~$f(\tau)$), and with
each power series $\phi_\a(\v)$ reducing to its constant term~$\phi_\a(0)$, as one
sees easily by using
the modular transformation property of~$f$ to compute the asymptotic development
of $f(\a+\v)$ for $\a\in\Q$ and $\v$ tending to~0 with positive imaginary part. More
generally, mock modular forms (whose definition we omit) also define elements
of~$\FF_{\k,v}$, where $\k$ is again the negative of the weight, but in that case the
power series $\phi_\a(\v)$ are in general factorially divergent rather than constant.
We do not elaborate on any of this since it is far from the theme of this paper, but
it is nice to observe that classical modular and mock modular forms have properties
in common with the asymptotic functions occurring here.

There are two further points that we should mention in connection with the
definition~\eqref{DefFkv}. One is that the absolute value appearing there is only
appropriate for~$\v$ real, which is our original situation when we think of $\a+\v$
as being a deformation of the rational number~$\a$ on the real line or when we take
$\v=-1/c(cX+d)$ with $X$ a rational number tending to infinity as in the RQMC. But
when we consider functions near~$\Q$ in the complex as well as in the real domain, the
absolute value sign would destroy holomorphy. If $\k$ is an even integer, the problem
does not arise, since we can simply replace $|\v|^\k$ by $\v^\k$, which is
holomorphic. If this is not the case then if we consider only functions near
rational points in the upper or lower half-plane, we can still replace $|\v|^\k$
in the definition by~$\v^\k$, which makes sense because $\v$ has a well-defined
logarithm in either half-plane. (We will never encounter functions in~$\FF_{\k,v}$
for $\k\ne0$ that are defined in a 360$^o$ complete complex neighborhood of~$\a$;
our functions will either be defined for nearby real points or for nearby non-real
points, or sometimes in a cut plane $\C\ssm(-\infty,0]$ or $\C\ssm[0,\infty)$, in
which case we can extend $|\v|^\k$ holomorphically as $\v^\k$ or $(-\v)^\k$,
respectively.) However, when $\k$ is not an integer and we want to discuss
the $\SL_2(\Z)$ action on $\FF_{\k,v}$, then we have to include some kind
of multiplier system, as familiar from the theory of modular forms of arbitrary weight.
Again, we omit details.

The second minor comment is that one can further generalize $\FF_{0,0}$ by
introducing a level~$N$ as well as the parameters~$\k$ and~$v$. This generalization
is necessary if we want to include modular or mock modular forms of level~$N$ (say
on~$\G=\G_0(N)$ or~$\G(N)$) into our definition, but also for our knot invariants if
the knot has a level~$>1$, as we found to be the case for the $(-2,3,7)$-pretzel knot.
Here the power series $\phi_\a(\v)$ and their completions will have period~$N$ rather
than~1 with respect to~$\a$, and more importantly, the number $v$ in~\eqref{DefFkv} is
no longer constant but must be replaced by a number~$v_\a$ that depends on the
$\G$-equivalence class (``cusp'') of~$\a$. Again we omit details, since this is not
our main subject.

We now return to the functions studied in this paper and to the reason why we
introduced asymptotic functions near~$\Q$ in the first place. Consider first the
tweaking function defined by equation~\eqref{tweakdef}. We showed in
Lemma~\ref{lem.lambda} that the map $\g\mapsto\l_\g$ is a cocycle in the space of
almost-everywhere-defined functions on~$\BP^1(\Q)$. It is easily checked that it is
not a coboundary in that space. But if we extend~$\l_\g$ to a function near~$\Q$ by
setting $\l_\g(\a+\v)=\l_\g(\a)$ (constant power series), then
equation~\eqref{tweakCoboundary} says that it now is a coboundary:
$\l_\g(x)=\mu(x)-\mu(\g x)$ where $\mu$ is the function near~$\Q$ defined by
$\mu(\a+\v)=1/\den(\a)^2\v $. More interestingly, if for each~$\s$ and~$\s'$
in~$\calP_K$ we define a function near~$\Q$ by
\smash{$Q^{(K,\s,\s')}(\a-\hbar)=\Phih_\a^{(K,\s,\s')}(h)$} and then put them together as a
matrix-valued function~$\BJ$ near~$\Q$ given by \smash{$\BJ(\a-\hbar)=\bPhih^{(K)}_\a(h)$},
then using equation~\eqref{tweakCoboundary} again we see that the complicated
equation~\eqref{WTaylor} can be replaced by the much simpler equation
\be
\label{WAsCoboundary}
W_\g(x) = \BJ(\g x)\i \BJ(x) .
\ee
Notice that in this equation the $(\s,\s')$-entry on the left-hand side is the sum
over~$\s''\in\calP_K$ of the product of the $(\s,\s'')$-entry of $\BJ(\g x)\i$ and the
$(\s'',\s')$-entry of $\BJ(x)$, which belong to~$\FF_{-\k(\s''),-v(\s'')}$ and
$\FF_{\k(\s''),v(\s'')}$, respectively. Thus each of the terms of the sum belongs
to~$\FF_{0,0}$ and we never encounter the problem of having to make sense of sums of
asymptotic functions of different orders of growth. The fact that the entries
of $W_\g$ all belong to $\FF_{0,0}$ is, of course, a~necessary prerequisite for the
final statement that they actually belong to its subring~$C^\infty(\R)$.

Equation~\eqref{WAsCoboundary} tells us the cocycle~$\g\mapsto W_\g$, which was
not a coboundary in the space of almost-everywhere-defined matrix-valued functions
on~$\BP^1(\Q)$ or of piecewise smooth functions on~$\BP^1(\R)$, becomes one when we
pass to the space of matrix-valued functions near~$\Q$. Both of these can be seen
as manifestations of a general phenomenon that one finds in almost all mathematical
contexts where notions of homology or cohomology play a role: even though one is
really only interested in cocycles that are not coboundaries, the cocycles that one
studies are almost always constructed as coboundaries in some bigger space.

\subsection{Analyticity}
\label{sub.holomorphic}

In Section~\ref{sub.lift}, we discussed the surprising smoothness properties
of the function $W_\g$ on $\R\ssm\{\g\i(\infty)\}$. In this subsection, we
come to a point much deeper than the smoothness, namely analyticity properties
of functions defined in a cut plane. These functions are closely related to
\emph{state integrals}. Such integrals appeared originally in the work of Hikami,
Andersen, Kashaev and others (see, for example, \cite{AK, DGLZ, Hikami}) in relation to
the partition function of complex Chern--Simons theory and to quantum Teichm\"uller
theory, and reappear in our context in~\cite{GZ:qseries}, the companion paper to this one.
We refer to these papers for details and describe the main points here in qualitative form only.

State integrals are analytic functions with several key features:
\begin{itemize}\itemsep=0pt
\item
 They are holomorphic for all $\tau \in \BC'=\BC\ssm(-\infty,0] $.
\item
 Their restrictions to $\BC\ssm\BR$ factorize bilinearly as finite sums of products of
 a $q$-series and a $\tilde q$-series, where $q=\e(\tau)$ and $\tilde q=\e(-1/\tau)$;
 see~\cite{holomorphic-blocks,GK:qseries}.
\item
 Their evaluation at positive rational numbers also factorizes bilinearly as a~finite
 sum of a~product of a~periodic function of $\tau$ and a~periodic function of
 $-1/\tau$; see~\cite{GK:evaluation}.
\end{itemize}
They are defined as multidimensional integrals of a product of quantum dilogarithms
times the exponential of a quadratic form. The quantum dilogarithm, invented
by Faddeev~\cite{Faddeev,FK-QDL}, is a~remarkable meromorphic function of two variables.
The structure of its poles implies that the state integrals are holomorphic
functions of $\tau$ in the cut plane~$\BC'$. The quantum dilogarithm is also a
quasi-periodic function with two quasi-periods, and this has two consequences, one of
which is directly related to the third ``feature'' above, and the other to the second
``feature'' and to the paper~\cite{GZ:qseries}.

The first consequence is the fact that one can apply the
residue theorem to give an exact formula for the values of the state-integrals at
positive rational numbers. Such a formula was given explicitly for the one-dimensional
state integrals considered in \cite[equation~(1), Theorem~1.1]{GK:evaluation}, and
those one-dimensional state integrals cover the case of the three knots that we consider
here, namely the $4_1$, $5_2$ and $(-2,3,7)$ pretzel knot. It turns out that
equation~(15) of~\cite{GK:evaluation} applied to the case of $(A,B)=(1,2)$ gives a
function on $\BQ^+$ which is none other than one of our four entries of \smash{$W^{(4_1)}_\g$}
when $\g=S$. To get the other three entries of \smash{$W^{(4_1)}_S$} one can apply the proof
of~\cite{GK:evaluation} to a $2 \times 2$ matrix of state-integrals of the $4_1$ knot
introduced in \cite[Theorem~3]{GGM}. And finally, to get the full matrix \smash{$W^{(4_1)}_\g$}
for all $\g$, one can apply the proof of~\cite{GK:evaluation} to a $2 \times 2$ matrix
of state-integrals of the $4_1$ knot that depend on a modular version of Faddeev's
quantum dilogarithm~\cite{GKZ:cocycle}.

The second consequence is perhaps even more interesting. Not only is each component of
the state integral matrix a finite sum of products of a $q$-series and a
$\tilde q$-series, but this sum precisely corresponds to matrix multiplication and says
that the whole state integral matrix~$W_S(\tau)$, whose restriction to a real half-line
is our function~$W_S$, factors as the product of a matrix of~$\tilde q$-series multiplied
by a matrix of $q$-series. More explicitly, $W_S(\tau)$ factors in the upper and
lower complex half-planes as \smash{$\BJ^\hol(-1/\tau)  \jhol_S(\tau) \i\BJ^\hol(\tau)$},
where $\BJ^\hol(\tau)$ is an~${(r+1) \times (r+1)}$ matrix with holomorphic and
periodic entries and $\jhol_S(\tau)$ is a diagonal matrix of automorphy factors.
Furthermore, the
equivariant extension~$W_\g$ of the state integrals mentioned above is again a~holomorphic
function in the cut plane whose restriction to the real half-line is our cocycle $W_\g$
from Section~\ref{sec.nearQ} and whose restriction to $\C\ssm\R$ factors for every~$\g$
as \smash{$\BJ^\hol(\g(\tau))  \jhol_\g(\tau) \i\BJ^\hol(\tau)$}
with the same periodic function~$\BJ^\hol(\tau)$.
The fact that this quotient extends analytically across a half-line, even though the
matrix-valued holomorphic function $\BJ^\hol(\tau)$ does not, is an example of a
(matrix-valued) \emph{holomorphic quantum modular form}, a new and quite general
context that is discussed in much more detail in~\cite{GZ:qseries,Z:HQMF}, and of which the mock modular
forms mentioned in the previous subsection give another nice example.
The fact that $W_\g(\tau)$ factors as $\BJ^\hol(\g(\tau))  \jhol_\g(\tau)
\i\BJ^\hol(\tau)$ is another
instance of the general principle (``cocycles are constructed as coboundaries in some
larger space'') mentioned at the end of Section~\ref{sec.nearQ}. So we have now
represented the original cocycle $\g\mapsto W_\g$ on the real line as a coboundary in
two different worlds: functions defined in a small open neighborhood of
$\BP^1(\R)\ssm X$ in $\BP^1(\C)$ for some finite set~$X$, and asymptotic functions
near~$\BQ$. But in fact these two representations~\smash{${W_\g(\tau)=\BJ^\hol(\g(\tau))  \jhol_\g(\tau) \i\BJ^\hol(\tau)}$} and~\smash{$W_\g(x)=\BJ(\g(x))\jt_\g(x) \i\BJ(x)$} are not independent:
as we will see in~\cite{GZ:qseries},
the periodic holomoprhic function~$\BJ^\hol$ has an asymptotic development as one
approaches any rational number from above or below in~$\C\ssm\R$, and this is a~representative of the~\emph{same} asymptotic functions near~$\BQ$ that we obtained
from the Habiro-like functions on~$\BP^1(\Q)$. It is this manifestation of the same
abstract object in two completely different realizations that we referred to in the
opening paragraph of this paper as an analogue in our context of the notion of motives.

The above discussion explains why one can expect, and in a few cases even prove,
the analyticity of the cocycle function~$W_\g(x)$. But it also seems worth observing
that, once one has predicted this analyticity, one can check it numerically using
only the matrices studied in this paper, without any reference to either $q$-series
or state integrals. Specifically, Proposition~\ref{prop.W1} gives the Taylor expansion
of $W_\g$ at any rational point, and since the coefficients of this series are
effectively computable, we can calculate a large number of them and see experimentally
that the series has a non-zero radius of convergence, as was already done in
equation~\eqref{DefR} for the special case of the expansion of $W_S$ for the $4_1$
knot around $x=1$. In fact, the coefficients can be computed in two different ways,
either by using the refined quantum modularity conjecture numerically with the help
of optimal and smooth truncation of divergent series, as was done in
Section~\ref{sec.RQMC}, or else by using the \emph{exact} formulas (when they are
available, e.g., for the~$4_1$ and $5_2$ knots) for $W_\g$ on~$\Q$ to compute the values
of this function at many rational points near a given point and then interpolating
numerically by the method recalled in ``Step~3'' of Section~\ref{sub.comp2}.
In this way, we can
verify the predicted real-analyticity to high precision and in a very convincing way
using only the data coming from the Kashaev invariant and its associated functions.
The simplest example is equation~\eqref{DefR} given in Section~\ref{sec.Phi} for the
$4_1$ knot and $\g=S$. The improvement of convergence in this case is very dramatic:
the 150th coefficient of $\Phi(2\pi {\rm i}x)$ (the last one that we computed) is about
$10^{284}$, but the 150th coefficient of the bilinear combination of power series
occurring on the right-hand side of~\eqref{DefR} is only $0.002 $!

But here we can actually do even more; by changing the variables one gets a new series
that not only again (conjecturally and experimentally) has radius of convergence~1,
but that now also gives numerical confirmation of the prediction that $W_S(x)$ extends
holomorphically to the whole cut plane. Specifically, if we make the change of variables
\smash{$1+x=\bigl(\frac{1+w}{1-w}\bigr)^2$}, under which~${x=0}$ corresponds to $w=0$ and
the condition $1+x\in\BC'$ is equivalent to~$|w|<1$, then we get a~power series
$B(w)\in\BR\big[\big[w^2\big]\big]$ defined by
\begin{align}
\label{QRS3}
B(w) &=
{\rm e}^{-v(4_1)} \Phi\biggl(\frac{8 \pi{\rm i} w}{(1-w)^2}\biggr)
\Phi\biggl(-\frac{8 \pi{\rm i} w}{(1+w)^2}\biggr)
 -
{\rm e}^{v(4_1)} \Phi\biggl(-\frac{8 \pi{\rm i} w}{(1-w)^2}\biggr)
\Phi\biggl(\frac{8 \pi{\rm i} w}{(1+w)^2}\biggr)
\end{align}
(with $\Phi(x)\in\BR[[x]]$ again given by~\eqref{as41}) which should have radius
of convergence~1. In fact, the numerical calculation, described in~\cite{GZ:qseries},
show that the 150th coefficient of $B$ is about $-7.5 \cdot 10^{10}$, again far
smaller than the original~$10^{284}$. The fact that this number is much bigger than
the corresponding number~$0.002$ for the bilinear combination~\eqref{DefR} is not
because the series~$B(w)$ is worse than the one in~\eqref{DefR}, but precisely because
it is better: in order to get the full domain~$\BC'$ of holomorphy of $W_S(x)$ we have
had to produce a power series that has singularities on the entire unit circle rather
than at just one point, and the coefficients correspondingly have much larger growth
(namely exponential in the square-root of the index, just as in the Hardy--Ramanujan
partition formula, rather than being only of polynomial growth, or in this case even
of polynomial decay). But in any case, whether we use~\eqref{QRS3} or just~\eqref{DefR},
we see that the single divergent power series $\Phi(h)$, which describes the asymptotic
behavior of \smash{$W_S^{(4_1,\s_1,\s_1)}(x)$} near either~$\infty$ or~0, suffices in an explicit
manner to determine this function everywhere on all of~$\R^*$. For general knots,
the corresponding statement would only hold if we consider the entire matrix~$\bPhi$
rather than just one entry. In fact, as the whole discussion of Sections~\ref{sec.RQMC}
and~\ref{sec.Matrix} shows, if we assume the whole RQMC, then at least in favorable
cases it is probably true that the single power series \smash{$\Phi^{(K)}_0(h)$} coming from the
modularity of the original Kashaev invariant actually determines everything.

We can summarize this whole subsection as the following conjecture for the cocycle $W_\g$.

\begin{Conjecture}
\label{conj.Wanalytic}
The function $W_\g$ on $\BQ\ssm\big\{\g\i(\infty)\big\}$ extends to a real-analytic function on~$\BR\ssm\big\{\g\i(\infty)\big\}$, and its restriction to each component of
$\BR\ssm\big\{\g\i(\infty)\big\}$ extends to a holomorphic function on the cut plane consisting
of this half-line and~$\BC\ssm\BR$.
\end{Conjecture}

\subsection{The non-hyperbolic case}
The main thrust of this paper, and all of the examples which we have treated in detail,
concern the case of hyperbolic knots, for which the volume is positive. We expect that
matrix-valued cocycles exist for nonhyperbolic 3-manifolds, with or without boundary,
and know that this is so for the example of the complement of the trefoil~\cite{Za:strange}
(where the corresponding invariant is sometimes known by the name Kontsevich--Zagier series)
as well as for WRT invariant of the Poincar\'e homology sphere (a spherical 3-manifold),
which was studied by Lawrence and Zagier~\cite{LZ}. In these examples and many
others that have been treated since, the series that occur are Taylor series of mock modular
forms, and we think that this will always happen for manifolds for which all of the volumes
vanish modulo~$\pi^2\Q$ (e.g., torus knots, Seifert-fibered manifolds or, in the closed case, spherical
manifolds). When it happens, the entries in the~$\calP\RED_K$-part of the matrix are the product of
an elementary exponential term \big(a rational power of~\smash{${\rm e}^{\pi^2/h}$}\big) and a rational power of~$q$,
so that the corresponding $\Phi$-series is purely exponential in~$h$, while the entries in the top row
of the matrix (which are again elements of the Habiro ring) still have factorially divergent
$h$-series as in the hyperbolic case, are now elementary functions, with coefficients that
are special values of Dirichlet $L$-series and a Borel transform which is simply a~trigonometric
function. However, we should emphasize that this simple behavior is not expected for all
non-hyperbolic knots or manifolds, but only for those for which {\it all} solutions
of the Neumann--Zagier equations are torsion in the Bloch group, so that all volumes~$v(\s)$
are rational multiples of~$2\pi {\rm i}$. Some knots, the like $(2,1)$-cabling of the $4_1$ knot,
are non-hyperbolic, so have vanishing volume~$\V(\s_1)$ modulo~$4\pi^2$ but have some
$\V(\s)$ with non-zero imaginary part, and then one expects to find non-trivial $h$-series.\looseness=-1

\part*{Part II. Complements}
\pdfbookmark[1]{Part II. Complements}{part2}

\section{Half-symplectic matrices and their perturbative series}\label{sec.perturbation}

In Section~\ref{sec.Phi}, we introduced a finite set $\calP_K$ associated to a
knot~$K$ and the formal power series~\smash{$\Phi_\a^{(K,\s)}(h)$} for each $\a\in\Q$ and
$\s\in\calPx_K=\calP_K\ssm\{\s_0\}$, as defined by Dimofte and the first author
in~\cite{DG,DG2} in terms of the Neumann--Zagier data of a triangulation of
the knot complement. In this section, we provide details and also a somewhat more
general construction, depending on more general data consisting of a
``half-symplectic matrix'' (defined below), an integral vector, and a solution of the
associated Neumann--Zagier equations. This more general class has a $q$-holonomic
structure that will be studied in Section~\ref{sec.NewHolonomy} and will also include
the formal power series \smash{$\Phi_\a^{(K,\s,\s')}(h)$} ($\s,\s'\in\calP_K$) that we found in
Sections~\ref{sec.RQMC} and~\ref{sec.Matrix}, as well as the asymptotic series
of Nahm sums near roots of unity. These half-symplectic matrices give a new perspective
on the classical Bloch group and the extended Bloch group.

\subsection{Half-symplectic matrices and the Bloch group}
\label{sub.HalfSymplectic}

To each knot~$K$ and each element $\s\in\calP_K$ there is an associated element of the
Bloch group (or third algebraic $K$-group) of~$\Qbar$ that plays a central role for
many of the constructions and that can be described in terms of the Neumann--Zagier data
of a triangulation of the knot complement. In fact, this construction produces an
invariant lying in a set defined by ``half-symplectic matrices'' (= upper halves of
symplectic matrices over~$\Z$) which is a refinement of the usual Bloch group that has
several nice aspects and seems not to have been considered in the literature.
In this subsection and the following one, we will describe this set and how one obtains
elements in it from the data of a triangulation. In the final subsection, we will
explain how to associate a~formal power series in~$h$ to any such element, the two
cases of primary interest being the matrix of power series \smash{$\Phi_\a^{(K,\s,\s')}(h)$}
associated to a knot and the power series describing the asymptotics of Nahm sums
near rational points.

For each positive integer~$N$, we denote by $\bH_N$ the set of $N\times2N$
\emph{half-symplectic matrices}, by which we mean matrices
$H=(A B)\in M_{N\times2N}(\Z)$ satisfying the two conditions
\begin{itemize}\itemsep=0pt
\item[(i)]
the $2N$ columns of~$H$ span $\Z^N$ as a~$\Z$-module, and
\item[(ii)]
the matrix $AB^t$ is symmetric.
\end{itemize}
The name refers to the fact that such matrices arise as the upper half
of symplectic matrices, i.e., of matrices $M=\vsma ABCD\in\mathrm{GL}_{2N}(\Z)$ satisfying
\smash{${\sma ABCD}^{-1}=\sma{D^t}{-B^t}{-C^t}{A^t}$}. To each \hbox{$H\in\bH_N$}, we associate
the generically zero-dimensional variety~$V_H$ defined as the set of \hbox{$N$-tuples}
${z=(z_1,\dots,z_N)}$ in \smash{$\big(\mathbb A^1\ssm\{0,1\}\big)^N$} \big($\mathbb A^1 =$ affine line\big)
satisfying the equations
\be
\label{NZeq}
V_H \colon\ \prod_{j=1}^N z_j^{A_{ij}} =
(-1)^{(AB^t)_{ii}\vphantom{y_{y_y}}} \prod_{j=1}^N (1-z_j)^{B_{ij}},
\qquad i=1,\dots,N .
\ee
To define the associated power series, we will need both an element of $V_H(\C)$ and a
slightly stronger discrete datum than~$H$, namely a pair (or triple)
\be
\label{DefXi}
\Xi = (H, \nu) = ((A B), \nu) \qquad \text{with}
\quad \nu \in \text{diag}\bigl(AB^t\bigr)+ 2\Z^N .
\ee
Equation~\eqref{NZeq} can then be written in abbreviated form as $z^A=(-1)^\nu(1-z)^B$.

We observe that there is a second description of the variety $V_H$ as the set of
$N$-tuples ${x=(x_1,\dots,x_N)}$ satisfying the trinomial equations
\be
\label{Pteq}
1 = (-1)^{\sum_j C_{ji}\nu_j} \prod_{j=1}^N x_j^{A_{ji}} +
(-1)^{\sum_j D_{ji}\nu_j} \prod_{j=1}^N x_j^{B_{ji}}, \qquad i=1,\dots,N
\ee
or in abbreviated form $1=(-1)^{C^t\nu}x^{A^t}+(-1)^{D^t\nu}x^{B^t}$, which is
isomorphic to $V_H$ via the bijections~\smash{$x\mapsto z=(-1)^{D^t\nu}x^{B^t}=1-(-1)^{C^t\nu}x^{A^t}$} and
\smash{$z\mapsto x=z^{-C^t}(1-z)^{D^t}$}. The $x$ are the Ptolemy coordinates as discussed
in Section~\ref{sub.Ideal} below in the case of knots and their logarithms are the
vectors~$w$ used below. Note that the signs $(-1)^{C^t\nu}$ and $(-1)^{D^t\nu}$
in~\eqref{Pteq} formulas do not depend on the choice of~$\nu$, since its value
modulo~2 is fixed by~\eqref{DefXi}. They do depend on the choice of symplectic
completion $(C D)$ of the half-symplectic matrix $(A B)$, but only in a~trivial
way: any other choice $(C^* D^*)$ of $(C D)$ has the form $(C D)+S(A B)$ for
some symmetric integral~${N\times N}$ matrix~$S$ (this corresponds to multiplying
$M=\sma ABCD$ on the left by the symplectic matrix~$\sma10S1$), and this simply
replaces $x$ by \smash{$(-1)^{S^t\nu}x$}, i.e., it changes the signs of some of the~$x_i $.

To any complex solution~$z$ of the system of equations~\eqref{NZeq} one can
associate a complex volume~$\V(z)$ in $\C/4\pi^2\Z$ that is defined roughly as the
sum of the dilogarithms of the $z_i$ plus a~suitable logarithmic correction. More
concretely, the imaginary part of~$\V(z)$ is a well-defined real number given
as $\sum_j D(z_j)$, where $D(z)=\operatorname{Im}\bigl(\Li_2(z)+\log|z|\log(1-z)\bigr)$ is
the Bloch--Wigner dilogarithm, which is single-valued. To define the full
value of $\V(z)$ modulo~$4\pi^2$ requires more work, because the function
$\Li_2(z)$ itself, defined by analytic continuation from its value~$\sum_{n\ge1}z^n/n^2 $ for~${|z|<1}$, is multivalued on $\C\ssm\{0,1\}$.
However, the function ${F(v)=\Li_2(1-{\rm e}^v)}$ has the derivative $v/({\rm e}^{-v}-1)$,
which is meromorphic with residues in~$2\pi {\rm i}\Z$. Hence, $F$~is a
well-defined function from $\C\ssm2\pi {\rm i}\Z$ to $\C/4\pi^2\Z$, satisfying the
easily checked functional equation~${F(v+2\pi {\rm i}n)=F(v)-2\pi {\rm i}n\log(1-{\rm e}^v)}$ for $n\in\Z$.
(See~\cite{Zagier-Gangl}.) We can then define\looseness=-1
\[
\V(z)=\V_\Xi(z)\in\C/4\pi^2\Z
\]
by the formula
\be
\label{defV}
\V(z) = \sum_{j=1}^N\left(F(v_j)+\frac{u_jv_j}2
+ \frac{\pi {\rm i}\nu_j}2 (Cu-Dv)_j \m \frac{\pi^2}6 \right) ,
\ee
where $u_j$ and $v_j$ are any choice of logarithms of $z_j$ and $1-z_j$ satisfying
$Au-Bv=\pi {\rm i}\nu$ (which automatically exist as a consequence of the conditions on~$\nu$
in~\eqref{DefXi} and the condition~(i) on~$H$) and where $(C D)$ is the bottom
half of a completion of $H$ to a full symplectic matrix. To see that this number,
which is only well-defined modulo~$4\pi^2$, is independent of the choice of $u$ and~$v$,
we observe that any other choice $(u^* v^*)$ of logarithms of~$z$ and~$1-z$ satisfying
$Au^*-Bv^*=\pi{\rm i}\nu$ has the form $(u^* v^*)=(u v)+2\pi {\rm i} \big(B^t A^t\big) n$ for
some~$n\in\Z^N$ (this follows easily from the conditions (i) and~(ii)), and then
using the functional equation of~$F$, we find
\[
\V^*\m\V=\pi {\rm i}n^t \bigl(-2Au+Au+Bv+2\pi {\rm i} AB^t n\m\pi {\rm i} \nu\bigr)
= 2\pi^2\big(n^t\nu\m n^tAB^tn\big) ,
\]
which is 0~modulo~$4\pi^2$ because
$AB^t$ is symmetric and integral with diagonal congruent to $\nu$ modulo~2. On
the other hand, the expression~\eqref{defV} \emph{does} depend on the choice of the
${2N\times N}$ integral matrix $(C D)$, but only very mildly, by a multiple of~$\pi^2/2$,
since changing $(C D)$ to~${(C D)+S(A B)}$ for some symmetric integral $N\times N$
matrix~$S$ changes the right-hand side of~\eqref{defV} by \smash{$-\frac{\pi^2}2\nu^tS\nu$}.
We believe, but have not checked, that it should be possible to lift the
formula~\eqref{defV} to a formula giving $\V_\Xi(z)$ modulo~$4\pi^2$ rather than
just modulo~$\pi^2/2$ in terms of $H$ alone by adding to the right-hand side a term
$r_M\pi^2/2$ where $\e(r_M/8)$ is the 8th root of unity occurring in the
transformation law of Siegel theta series with characteristics under the action of
the symplectic matrix $M=\sma ABCD$ as given by Igusa \cite[Theorem~3, p.~182]{Igusa}.

We make two small remarks on the above formulas before proceeding. The first is that
the term $\frac\pi2 \nu^t(Cu-Dv)$ in~\eqref{defV} is needed, not only to make the
expression on the right independent of the choice of logarithms $u$ and~$v$
modulo~$4\pi^2$ (it is already independent of this choice modulo~$\pi^2$ even if
this term is omitted), but in order to get the right imaginary part: the imaginary
part of $F(v)$ is $D(z)+\operatorname{Im}(u\bar v)/2$ for ${\rm e}^u=1-{\rm e}^v=z$, where $D(z)$ is the
Bloch--Wigner dilogarithm as above, and it is only if we include the term with $Cu-Dv$
in~\eqref{defV} that its imaginary part has the correct value~$\sum_jD(z_j)$. The
other is that the vector $w:=Dv-Cu$ whose scalar product with~$\nu$ gave the
correction term in~\eqref{defV} also gives a parametrization of the $N\times2$
matrix~$(u v)$ as $\big(B^t A^t\big)w+{\rm i}\pi\big(D^t,C^t\big)\nu$.
\big(To see this, just write the relationship of~$(u v)$ to $\nu$ and~$w$ as~$M (\smallmatrix u\\-v\endsmallmatrix)=(\smallmatrix \pi {\rm i}\nu\\-w\endsmallmatrix)$
and use the formula for~$M\i$.\big) This is simply the logarithmic version of the
alternative characterization of the variety~$V_H$ given in~\eqref{Pteq}, with
$w=\log x$. Changing~$(u v)$ by $2\pi {\rm i}\big(B^t A^t\big)n$ with $n\in\Z^N$
corresponds to taking a different logarithm $w$ of the same~$x$.

We next turn to the relation between half-symplectic matrices and the Bloch group.
The latter is an abelian group $\B(F)$ which is defined for any field~$F$ of
characteristic zero as the quotient of the kernel of the map $d\colon\Z[F]\to\WW(F^\times)$
sending $[x]$ to $x\wedge(1-x)$ for $x\ne0,1$ by the subgroup generated by the 5-term
relation of the dilogarithm. But the precise definition varies slightly in the
literature because of delicate 2- and 3-torsion issues arising from the particular
definition of the exterior square (for instance, does one require $x\wedge x=0$ for
all $x$ or just~${x\wedge y=-y\wedge x }$?) and the particular choice of the 5-term
relation, which potentially comes in $5^6$ versions obtained from one another by
replacing each of the 5 arguments by its images under the group generated by
$x\mapsto1/x$ and $x\mapsto1-x$. In fact, we will need the extended Bloch group as
introduced by Neumann~\cite{Neumann:CS} and studied further by Zickert and others
in~\cite{Goette,Zickert:extended}, but here also there are several versions.
We recall the definition from~\cite{Zickert:extended} here, and then
describe a small refinement and the relation to half-symplectic matrices.

Denote by $\wh\C$ the set of pairs of complex numbers $(u,v)$ with ${\rm e}^u+{\rm e}^v=1$.
This is an abelian cover of $\C^\times\ssm\{0,1\}$ via $z={\rm e}^u=1-{\rm e}^v$, with Galois
group isomorphic to~$\Z^2$. The extended Bloch group $\hB(\C)$ as defined
in~\cite{Goette,Zickert:extended} is the kernel of the map
$\wh d\colon \Z\big[\wh\C\big]\to\WW(\C)$, where~$\WW(\C)$ is defined by requiring only
$x\wedge y+y\wedge x=0$ (rather than $x\wedge x=0$, which is stronger by 2-torsion)
and where $\wh d$ maps $[u,v]:=[(u,v)]\in\Z\big[\wh\C\big]$ to~$u\wedge v$, divided by the
lifted version of the 5-term relation, namely, the $\Z$-span of the set of elements
\smash{$\sum_{j=1}^5(-1)^j[u_j,v_j]$} of \smash{$\Z\big(\wh\C\big)$} satisfying $(u_2,u_4)=(u_1+u_3,u_3+u_5)$
and $(v_1,v_3,v_5)=(u_5+v_2,v_2+v_4,u_1+v_4)$. There is an extended
regulator map from $\hB(\C)$ to~$\C/4\pi^2\Z$ given by mapping $\sum[u_j,v_j]$ to
$\sum\calL(u_j,v_j)$, where \smash{$\calL(u,v)=F(v)+\frac12uv-\frac{\pi^2}6$}, which one
can check vanishes modulo~$4\pi^2$ on the lifted \hbox{5-term} relation. One can
also define $\hB(F)$ for any
subfield~$F$ of~$\C$, such as an embedded number field, by replacing $\wh\C$ by the
subset $\wh F$ consisting of pairs $(u,v)$ with ${\rm e}^u=1-{\rm e}^v\in F$. The relation of the
Bloch group and the extended Bloch group to algebraic $K$-theory is that~$B(F)$ for
any field~$F$
is isomorphic up to torsion to the algebraic $K$-group~$K_3(F)$ \cite{Suslin},
with the Borel regulator map from $K_3(\C)$ to $\C/\pi^2\Q$ being given at the level
of the Bloch group by dilogarithms, while the extended Bloch group of a number
field~$F\subset\C$ is isomorphic to $K_3^{\text{ind}}(F)$ \cite{Zickert:extended}, for
which the Borel regulator lifts to~$\C/4\pi^2\Z$.

We now extend this group slightly by replacing $\Z\big[\wh\C\big]$ by the larger group
$\Z\big[\wh\C\big] \oplus \C$
and $\wh d$ by a map from this group to $\WW(\C)/({\rm i}\pi\wedge {\rm i}\pi)$, still given on
\smash{$\Z\big[\wh\C\big]$} by
$[u,v]\mapsto u\wedge v$ but now also defined on~$\C$ by $\wh d(x)=x\wedge(x+\pi{\rm i})$,
which despite appearances
is a linear map because of the antisymmetry of~$\wedge$. We then divide the kernel
of this new~$\wh d$
by a larger set of relations, namely the same lifted 5-term relation as before (with
$\C$ component equal to~0) together with the relations
$([u,v]+[v,u]-[u',v']-[v',u'],0)$ and $([u,v]+[-u,v-u+\pi{\rm i}],u)$ for all $(u,v)$ and
$(u',v')$ in~$\wh\C$, corresponding to the elements $[z]+[1-z]$ and $[z]+[1/z]$.
The extended regulator map to~$\C/4\pi^2\Z$ is now defined by mapping
$\bigl(\sum[u_j,v_j],x\bigr)$ to $\sum\calL(u_j,v_j)-x\pi{\rm i}/2$, which agrees with
the previous definition when $x$ 
is~0 and which can be checked to vanish also on the new relations. The advantage
of this further
extension of the Bloch group is that the solutions~$(u,v)$ of the logarithmic
Neumann--Zagier equations
(i.e., the set of $(u,v)\in\C^{2N}$ with~${(u_j,v_j)\in\wh\C}$ for each~$j$ and
$Au-Bv=\pi{\rm i}\nu$
with~$\nu$ as in~\eqref{DefXi}) now give an element of~$\hB(\C)$, namely the
class~$\xi$ of
the pair \smash{$\bigl(\sum_{j=1}^N[u_j,v_j], w\nu^t \bigr)$}, where $w=Cu-Dv$ as before.
Using the parametrization $ (u,v)=\big(B^t A^t\big)w+\big(C^t D^t\big)\nu {\rm i}\pi $ discussed
above, we check easily
that the image of this in~$\WW(\C)$ under \smash{$\wh d$} is
$\big(\nu^tCD^t\nu\big) ({\rm i}\pi)\wedge({\rm i}\pi)$, and its
image under the regulator map is precisely the number~$\V(z)$ defined in~\eqref{defV}.
When $(u,v)$ comes from a triangulation of a 3-manifold, then the effect of the
extended 5-term relation is
precisely that of a $(2,3)$-Pachner move (changing one triangulation to another by
replacing two
tetrahedra with a common face by three tetrahedra with the same set of vertices),
so that the
element $\xi\in\hB(\C)$ is a topological invariant of the manifold.\looseness=-1

We end this subsection by explaining briefly how half-symplectic matrices actually
give a~new description of the extended Bloch group as a quotient of $\Sp_\infty$ by
suitable relations. Here for convenience we are writing $\Sp_N$ rather than
$\Sp_{2N}$ for the group of symplectic matrices of size~${2N\times2N}$ over~$\Z$,
and $\Sp_\infty$ for the direct limit of these groups with respect to the natural
inclusions $\Sp_N\hookrightarrow\Sp_{N+1}$. It also turns out to be more convenient
to define $\Sp_N$ as the space of matrices $M$ satisfying $MJ_N^*M^t=J_N^*$ instead
of $MJ_NM^t=J_N$ used above, where~\smash{$J_N=\sma0{-1_N}{1_N}0$} and $J_N^*$ is the block
diagonal matrix with $N$ copies of~$J_1$ on the diagonal, in which case the
inclusion just sends $M$ to $M^+=\sma M00{1_2}$, and similarly the lifted 5-term
relations become much simpler with this convention.
The relations that we divide by are roughly as follows. The first is stability
(identify
$[M]$ and $[M^+]$). A second is that we identify $[M]$ and~${\bigl[\sma10S1 M\bigr]}$
with~$S$ integral and symmetric are equivalent. (This corresponds to working with
half-symplectic rather than full symplectic matrices.) A third is that we identify~\smash{$M\in\Sp_N$} with~\smash{$\sma g00{{g^t}\i}M$} for any~${g\in\GL_N(\Z)}$.
(This corresponds to permuting the $N$ relations~\eqref{NZeq} or multiplying one
of them by a monomial in the others.) A fourth is to identify~$M$ with~$M\sma P00P$
for any $N\times N$ permutation matrix~$P$, corresponding in the geometric case to
changing the numbering of the $N$ simplices, and yet another (which maybe can be
omitted) corresponds to relabelling the edges so that the
shape parameter $z$ goes to $z'$ or $z''$. The main one, of course, is a
symplectic-matrix version
of the 5-term relation. This was first discovered in the special case
corresponding to the
Nahm sums~\eqref{Nahm1} by Sander Zwegers in an unpublished 2011 conference talk
and then given
in various versions for arbitrary symplectic matrices by Dimofte and the first
author in~\cite{DG}
and in unpublished work by Campbell Wheeler and Michael Ontiveros~(MPIM). The set
of equivalence classes
becomes an abelian group by setting $[M]+[M']$ equal to the class of $\sma M00{M'}$ and
$-[M]$ to the class of $M\i$. To get a map from this group to the extended Bloch group
of $\C$, we have to first enlarge it by looking at equivalence classes, not just of
half-symplectic matrices~$H$ (which is enough by the second of the equivalence relations
listed above), but of pairs consisting of a half-symplectic matrix $H=(A B)$ together
with a~solution $(u,v)\in\C^{2N}$ of the logarithmic NZ equations $Au-Bv=\pi\nu$ with
$\nu$ as in~\eqref{DefXi}, with corresponding lifts of the 5-term and of the various
other relations.
The map from this larger group to the extended Bloch group is then the one described
in the previous
paragraph. It is injective because the 5-terms relations defining the extended Bloch
group all lift to corresponding relations at the (half-) symplectic level.
It is also surjective, as one can show using elements of the set
$\Sp_{N,N'} = \big\{M\in M_{N\times N'}(\Z)\mid MJ_{2N}M^t=J_{2N'}\big\}$
of ``non-square symplectic matrices'' (note that this set is just $\Sp_N$ if $N=N'$
and reduces to~0 if $N'>N$) together with the obvious composition maps
\smash{$Sp_{N,N'}\times\Sp_{N',N''}\to\Sp_{N,N''}$}, in order to eliminate superfluous relations.
(Roughly speaking, if \smash{$\sum_{j=1}^N[u_j,v_j]$} is the $\Z\big[\wh\C\big]$-component of an
element of $\hB(\C)$ as defined above, then we define $N'\le N$ as
the rank of the group generated by all $u_j$ and~$v_j$ and obtain an element of
$Sp_{N,N'}$ by
writing the $u$'s and $v$'s in terms of these generators, which then always satisfy a
collection of NZ equations.) A more detailed discussion of this and of the whole
relationship
between half-symplectic matrices and Bloch groups, including our versions of the~5-term
relation lifted to symplectic and half-symplectic matrices, is also given
in~\cite{GZ:neumann}.

This concludes our discussion of half-symplectic matrices and the equations~\eqref{NZeq}.
These objects arise in (at least) two different contexts, in 3-dimensional topology
and in the
study of special $q$-hypergeometric series (Nahm sums). The former is of course the
one that
is of most relevance for this paper, and will be discussed in more detail in the next
subsection, but after that we will also say something about Nahm sums because
they will play a role in the sequel~\cite{GZ:qseries} to this paper and also
because they give the most elementary approach to defining the associated formal
power series that are our main subject of interest.

\subsection{Ideal triangulations and the Neumann--Zagier equations}
\label{sub.Ideal}

\begin{figure}[htpb!]
\centering
 \includegraphics[height=4cm]{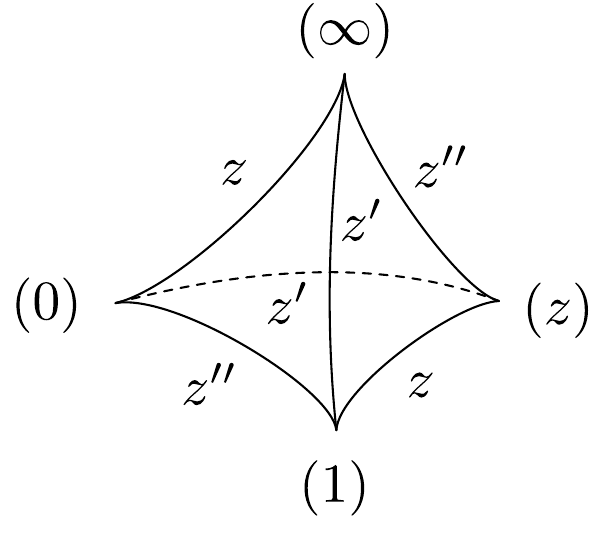}
\caption{A tetrahedron with shape parameters.}\label{fig.tetrahedron}
\end{figure}

In 3-dimensional geometry, the shape of an ideal tetrahedron in $\mathbb H^3$
is encoded by a complex number (``shape parameter'') $z \in \BC\ssm\{0,1\}$, the
tetrahedron being isometric to the convex hull of the four points
$0,1,\infty,z \in\mathbb P^1(\C)=\partial\big(\mathbb H^3\big)$. The shape $z$ has three
forms $z$, $z'=1/(1-z)$ and~${z''=1-1/z}$, each corresponding to the choice of a pair of
opposite edges of the tetrahedron as shown in Figure~\ref{fig.tetrahedron}.

An ideal triangulation of a
3-manifold with torus boundary components give rise to an equation~\eqref{NZeq}, where
the variables $z_i$ solving the equations~\eqref{NZeq} are the shape parameters
of the tetrahedra and the equations are the ``gluing conditions'' relating the shape
parameters of the tetrahedra incident on the various edges of the triangulation and on
the cusp. These gluing equations originated in the work of Thurston~\cite{Th}
and further studied in~\cite{NZ}, where the key symplectic property of the
matrices $(A B)$ was found. We explain very briefly how this works for 3-manifolds
whose boundary component is a torus, equipped with an ideal triangulation with $N$
tetrahedra. Each edge of the triangulation gives rise to a gluing equation asserting that
the product of the shape parameters of all tetrahedra incident to that edge equals to~1.
Every peripheral curve (i.e., a curve in the boundary torus of the 3-manifold) also
has an equation of this form (often called the holonomy equation, following Thurston),
obtained by setting the product of the shape parameters as the curve intersects the triangulated
boundary, equal to~1.
The product of the gluing equations corresponding to all edges is identically~1, so
one gluing equation is redundant and can be removed and replaced by the holonomy
equation of a nontrivial peripheral curve. Since there are $N$ edges, this gives
a collection of $N$ gluing equations. If one is interested in the geometric solution
that describes the complete hyperbolic structure, where all the shape parameters
have positive imaginary part, the above gluing equations are
replaced by their stronger logarithmic form, where the right-hand side is now $2\pi {\rm i}$
for each edge and $0$ for the peripheral curve. Using the fact that the three shape
parameters $z$, $z'=1/(1-z)$ and~${z''=1-1/z}$ satisfy the relation $z z' z''=-1$,
and in logarithmic form $\log z + \log z' + \log z''=\pi{\rm i}$, it follows that we can
eliminate one of the three variables at each tetrahedron (after choosing a pair of
opposite edges for each tetrahedron). Doing so, the logarithmic form of the gluing
equations now become linear equations for $\log z_i$ and $\log (1-z_i)$, whose
coefficients give rise to the Neumann--Zagier matrices $A$ and $B$, and where
right-hand side is a distinguished flattening~$\nu$ that should satisfy the mod~2
congruence given in~\eqref{DefXi}. (This congruence can presumably be deduced from
the ``parity condition'' for
ideal triangulations proved by Neumann~\cite{WN:comb}, but we have not checked this.)
Neumann--Zagier's theorem is that the above matrix $(A  | B)$ is the upper half of
a symplectic
matrix with integer entries. Note that the corresponding pairs $(H,z)$ and $(\Xi,z)$
are called
``NZ datum'' and ``enhanced NZ datum'' in~\cite{DG}. Note also that a different choice
of opposite edges in each tetrahedron cyclically permutes the triple \smash{$(z_j,z_j',z_j'')$}
and changes the corresponding Neumann--Zagier matrices, but does not change the
corresponding element of~\smash{$\hB(\C)$}.

The connection between gluing equations and symplectic matrices involves not only
the shapes of ideal tetrahedra, but also their Ptolemy variables. The latter is
an assignment of nonzero complex numbers $x_i$ at each edge of an ideal triangulation
that satisfy the Ptolemy equations, namely at each tetrahedron we have a quadratic
equation $x_1 x_2 \pm x_3 x_4 \pm x_5 x_6 =0$ (with suitable signs). The signs require
either ordered triangulations or a choice of a Ptolemy cocycle and a~detailed
description is given in~\cite[equation 12.2]{GTZ:complex} and also in~\cite[Section~3]{GGZ}.
(The equivalence between the shape and the Ptolemy description of character varieties
of surfaces is discussed in detail by Fock--Goncharov~\cite{Fock-Goncharov}.)
In the 3-dimensional case of a knot complement,
these $x_j$ are exactly the ones introduced
in~\eqref{Pteq}
(and related to $w=Du-Cv$ by $x_j={\rm e}^{w_j}$), which here becomes a system of quadratic
trinomial relations after rescaling because in each of the column of the gluing
equation matrices there are at most six non-zero entries, corresponding to the
six edges of the tetrahedron corresponding to that column.

We mention in passing that the variety defined by just the first $N-1$ edge gluing
equations is 1-dimensional (for a suitably chosen triangulation) and that this
curve maps to the $\PSL_2(\BC)$-character variety (via the developing map which
assigns a solution to the gluing equations a~$\PSL_2(\BC)$-representation of
the fundamental group of the manifold, well-defined up to conjugation). The
$\PSL_2(\BC)$-character variety maps to $\BC^* \times \BC^*$
(modulo a $\BZ/2^2$ quotient) and its image is described by the vanishing of the
A-polynomial $A(\ell,m)$
(where $\ell$ is the longitude) as introduced subsequently in~\cite{CCGLS}.
The variety obtained by adding taking the first $N-1$ relations together with the
relation $m^p\ell^q=1$ for coprime integers $p$ and~$q$ corresponds to the compact
3-manifold obtained by doing a $(p,q)$ Dehn surgery on the knot complement. A detailed
discussion of the choices involved to write down these matrices can be found in~\cite[Section~2, Appendix~A]{DG}. All of this data is standard in knot theory,
and is computed explicitly for any given knot complement (or more generally, an ideal
triangulation of a cusped hyperbolic 3-manifold) by the computer implementation of
\texttt{SnapPy}~\cite{snappy}.

Once an ideal triangulation $\D$ of a 3-manifold $M$ as above has been
fixed, a solution $z$ of its gluing equations gives rise via a developing
map to a representation $\rho_z$ (i.e., a group homomorphism) of $\pi_1(M)$
in $\mathrm{PSL}_2(\BC)$, well-defined up to conjugation. If we choose the
Neumann--Zagier equations as above, the representation $\rho_z$ is
boundary-parabolic and gives rise to an element of the extended Bloch
group~\cite{Zickert} and has a well-defined complex volume; see~\cite{Neumann:CS}
and also~\cite{GTZ:complex}. Thus, if $\D$ is an ideal
triangulation of a the complement of a knot $K$, we have a map
$z \mapsto \rho_z$ from $V_H(\BC)$ to $\calP_K$, and the complex volume of
$\rho_z \in \calP_K$ coincides with the complex volume of~$z$, as follows
from the work of Neumann~\cite{Neumann:CS} and Zickert~\cite{Zickert:extended}
on the extended Bloch group.

There are, however, several subtleties of the above construction which
we should point out. For instance, there exist triangulations of hyperbolic knots
for which the map $V_H(\BC)\to\calP\RED_K$ is not onto or even for which the complex
solutions set $V_H(\BC)$ is empty (this can happen even for triangulations of the
complement of the $4_1$ knot). In this paper, we will ignore these issues
and assume that we are dealing with ideal triangulations for which the map
is onto. We will further restrict our attentions to knots for which the
set $\calP_K$ is finite. (There are known to be knots for which the variety
$\calP_K$ has strictly positive dimension, but they are too complicated for
the calculations in this paper to be carried out. We believe that in such
cases the right indexing set of our formal power series would be
the set of components of~$\calP_K$ or of the variety~$V_H(\BC)$.)

An alternative approach to the definition of the set $\calP_K$ comes from the
branches of the $A$-polynomial curve above the point~$m=1$, where $m$
is the eigenvalue of the meridian. Even if the $\SL_2(\BC)$ character variety
of decorated representations of a knot complement has positive-dimensional
components, its image in $\BC^* \times \BC^*$, as given by the eigenvalue
of the meridian and the longitude, is one-dimensional, and (ignoring
any zero-dimensional components) is defined by the zeros of the
$A$-polynomial of the knot. The $A$-polynomial is discussed in detail
in the appendix of~\cite{BDRV}. We will focus on knots that
satisfy the property that the number of parabolic
representations $\s$ coincides with the degree of the $A$-polynomial
of a knot with respect to the longitude. Note that the Galois group
$\text{Gal}\big(\overline{\BQ}/\BQ\big)$ acts on the set $\calP_K$ of boundary
parabolic representations. What's more, in a boundary
parabolic representation, the longitude has eigenvalue $\pm 1$ and this
partitions the set $\calP_K$ into two subsets $\calP^{\pm}_K$,
each of which is stable under the action of $\text{Gal}\big(\overline{\BQ}/\BQ\big)$.
The geometric representation lies in $\calP^{-}_K$;
see~\cite[Lemma~2.2]{Ca:geometric}.

\subsection[Nahm sums and the perturbative definition of the Phi-series]{Nahm sums and the perturbative definition of the $\boldsymbol{\Phi}$-series}
\label{sub.PhiPerturbative}

In this final subsection, we describe how to attach to~$\Xi=((A B),\nu)$ as
in~\eqref{DefXi}, a solution~$z$ of the equation~\eqref{NZeq} and a number~$\a\in\Q$
a completed formal power series belonging to ${\rm e}^{\V(z)/\den(\a)^2h} \C[[h]]$.
As already stated in Section~\ref{sec.Phi}, this was done in~\cite{DG} (for $\a=0$)
and~\cite{DG2} (for general~$\a$) in the context of knot complements and
Neumann--Zagier data. However, there is a~completely different situation where the same
formal power series are attached to the same data $(\Xi,z,\a)$, namely the asymptotics
near roots of unity of special \hbox{$q$-hypergeometric} series called Nahm sums.
Since these are a little more elementary we will use them to explain the derivation of
the formal power series.

We begin by recalling what Nahm sums are. The simplest one is defined by
\be
\label{Nahm1}
F_{A,b}(q)= \sum_{n_1,\dots,n_N\ge0}
\frac{q^{\frac12n^tAn+b^tn}}{(q;q)_{n_1}\cdots(q;q)_{n_N}}\in \Z[[q]] ,
\ee
where $A$ is an even positive definite symmetric matrix in $M_N(\Z)$ and $b$ an
element of~$\Z^N$. Changing each $n_j$ by~1, we see that the stationary points of the
summand (i.e., the places where nearby terms are asymptotically equal, giving the
expected main contributions to the whole sum)
 are given in the limit $q\to1$ by $q^{n_j}=z_j+\o(1)$, where $z=(z_1,\dots,z_N)$
is a solution of Nahm's equation $1-z=z^A$ (which is the special case $B=\mathbf1_N$
of~\eqref{NZeq}, with $(A \mathbf1_N)$ being half-symplectic). Nahm observed
that for any solution $z$ of this equation the element $\sum_i[z_i]$ belongs
to the Bloch group $\B(\C)$ and conjectured that $F_{A,b}(q)$ (up to a rational
power of~$q$, and considered as a function of $\tau$ with $q=\e(\tau)$)
can only be a modular function if at least one solution of the Nahm equation
has a trivial class in the Bloch group, and conversely that $F_{A,b}$ (again up
to a power of~$q$ and as a function of~$\tau$) is a modular function of~$\tau$
for \emph{some} $b$ if all solutions of the Nahm equation have trivial class
in the Bloch group. The first assertion was proved in~\cite{CGZ}; the second
is still open.

If we now generalize the Nahm sum to
\be
\label{Nahm2}
F_\Xi(q)=\sum_{n\in\Z^N, B^tn\ge0} \frac{(-1)^{\nu^tn} q^{\frac12(n^tAB^tn+\nu^tn)}}
{(q;q)_{\vphantom{x^\l}(B^tn)_1}\cdots(q;q)_{\vphantom{x^\l}(B^tn)_N}}
=\sum_{\substack{m,n\in\Z^N_{\ge0}\times\Z^N\\ m=B^tn}}
\frac{(-1)^{\nu^tn} q^{\frac12(n^tAm+\nu^tn)}}{(q;q)_{m_1}\cdots(q;q)_{m_N}} ,
\ee
with $\Xi=((A B),\nu)$ as in~\eqref{DefXi}, which is still a power series in $q$
because of the congruence condition on~$\nu$, then the same consideration as
before shows that the stationary points of the sum correspond via $z=q^{B^tn}$ to the
solutions of the equation~\eqref{NZeq}. A formal computation of the contribution
of the summands near these stationary point will lead to the perturbative series
in~$h$ that we are looking for, where $q={\rm e}^{-h}$ with $h\to0$, and in some cases
one can show that these $h$-series actually do describe the radial asymptotics
of the Nahm sum (see~\cite{GZ:asymptotics}, where this is shown to be the case
for the original Nahm sum~\eqref{Nahm1} and the real solution of the Nahm
equation), but in general the calculation is purely formal because unless all
$z_i$ are between 0 and~1 the series corresponding to~$h$ will not correspond
to any subsum of~\eqref{Nahm1} or~\eqref{Nahm2}. We also mention that there
are even more general Nahm sums whose $n$-th summand $\big(n\in\Z^N\big)$ is the product
of a~root of unity, a power of~$q$ given by a quadratic function of~$n$,
and a~product of Pochhammer symbols (possibly to integer powers) with linear
forms in~$n$ as arguments, which occur in several places in quantum
topology, e.g., the 3D-index~\cite{DGG} and many of the $q$-series in~\cite{GZ:qseries}.\looseness=-1

We now explain how to associate to the datum $\Xi=(H,\nu)$ and point $z$
on $V_H$ a formal power series \smash{$\Phi_\a^{(\Xi,z)}(h)$} for each~$\a\in\Q/\Z$.
We will do this first for the easier case $\a=0$, and then discuss how the
formula changes in the general case. The calculations for $\a=0$ were done
first for the simplest Nahm sum~\eqref{Nahm1} in~\cite{Zagier} and~\cite{VZ} and
for general half-symplectic matrices $(A B)$ (though under the assumption that
$B$ is invertible over~$\Q$) in the context of knots in~\cite{DG}. The power series
that were obtained in these two different contexts were syntactically identical, and
this coincidence persisted for general $\a$, with the perturbative series of~\cite{DG2}
being syntactically equal to the asymptotics of Nahm sums at roots of
unity~\cite[Section 5]{GZ:asymptotics}, with the formulas in all cases being
given in terms of (sums of) formal Gaussian integrals. It is for this reason
that we can use the easier Nahm sums to motivate the precise form of the integral
to be studied. We only sketch the argument, referring to the papers above for more
details.

We begin by rewriting the first definition in~\eqref{Nahm2} in the form
\begin{gather}
\label{Nahm3}
F_{(A B),\nu}(q) =
\frac1{(q;q)_\infty^N} \sum_{n\in\Z^N} (-1)^{\nu^tn} q^{\frac12(n^tAB^tn+\nu^tn)}
\prod_{j=1}^N \bigl(q^{(B^tn)_j+1};q\bigr)_\infty ,
\end{gather}
where we no longer have to restrict to $n$ with $B^tn\ge0$ because
$\big(q^{m+1};q\big)_\infty$ vanishes for $m\in\Z_{<0}$. We must assume for now that
the symmetric matrix $AB^t$ is positive definite to ensure the convergence of the
series~\eqref{Nahm2} or~\eqref{Nahm3}, but this is not important at the end since
the final formulas will be purely algebraic and make sense without this assumption.
The key point is that if $q={\rm e}^{-h}$ with $h$ small then the sum will be approximated
to all orders by the corresponding integral, with~$\Z^N$ replaced by~$\R^N$ and the
summation sign by an integral sign. (This is a~consequence of the Poisson summation
formula, which represents the sum over~$\Z^N$ of a sufficiently smooth function of
sufficiently rapid decay as the sum over~$\Z^N$ of its Fourier coefficient, whose
constant term is the integral corresponding to the original sum and whose other
terms are of smaller order.) We then look at the expansion of the integrand around its
stationary points and approximate each by a Gaussian times a power series in a small
local variable, as is always done in perturbation theory. The stationary
points are indexed by the complex points~$z$ of~$V_H$, as already indicated,
the correspondence being given by $q^{(B^tn)_j} \sim z_j$.
On the other hand, for~${z\in\C^*}$,~$q={\rm e}^{-h}$ with~$h$ tending to~0, and
$t$ either fixed or growing more slowly than any power of~$1/h$, we have the asymptotic
formula
\begin{align}
\frac1{\big(z{\rm e}^{t\sqrt h};q\big)_\infty} & \sim
\exp\left(\sum_{m=0}^\infty \frac{B_m\bigl(t/\sqrt h\bigr)}{m!}
\Li_{2-m}(z) h^{m-1}\right)\nonumber \\
& = \exp\left(\frac{\Li_2(z)}h + \left(\frac t{\sqrt h}-\frac12\right)
\log\left(\frac1{1-z}\right)
+ \frac{t^2}2 \frac z{1-z} + \text{(small)}\right) ,\label{newpsixz}
\end{align}
where $B_m(t)$ denotes the $m$-th Bernoulli polynomial and ``(small)'' is an explicit
power series in $t$ and $\sqrt h$ with no constant term in~$\sqrt h$. (The first
statement is \cite[Lemma, p.~53]{Zagier} and the second follows because all
contributions from $B_m\big(t/\sqrt h\big)$ with $m\ge2$ except for the quadratic part of
the $B_2$-term are small.) Inserting this into the parts near the stationary points
of the integral corresponding to the sum~\eqref{Nahm3}, we find after some calculation
that the total contribution of the stationary part corresponding to a given solution
$z$ of~\eqref{NZeq} is ${\rm e}^{\V(z)/h}$ times an explicit power series in~$h$ (initially
in~$\sqrt h$, but then in~$h$ because of the parity properties of Bernoulli
polynomials) which is written out in~\cite{DG}. We only mention here that the power
series obtained has coefficients in $\Q(z)$ (and hence in $F_\s$ in our application
to knots) except for a~prefactor~${\det(A+B \diag(z_j/(1-z_j)))^{-1/2}}$
coming from the determinant of the quadratic part of the Gaussian.

When $\a=a/c$ is not integral, the calculations, done in~\cite{DG2} in the
general case (still with~$B$ invertible) and in~\cite{GZ:asymptotics} for the
special Nahm sums~\eqref{Nahm1}, are much more complicated and we refer to those
papers for the explicit formulas. A key point is that the stationary points of the
integral are now indexed by the $c$-th roots of the solutions~$z$ of~\eqref{NZeq},
but with the quadratic form appearing in the Gaussian depending only on~$z$ and
not on the choice of $c$-th root. This means that each of the formal power series
\smash{$\Phi_\a^{(\Xi,z)}(h)$} has the form of a sum over $(\Z/c\Z)^N$ \big(after choosing some
fixed $\sqrt[c]z$\big) of expressions similar to those occurring for the simpler
case~$c=1$. The reader can get a feeling for the nature of the formulas appearing
by looking at Section~\ref{sec.QMC41} of this paper, where they are carried out
in detail for the Kashaev invariant of the $4_1$ knot, this case however being
deceptively simple because of the positivity of all of the terms occurring.

We make one final remark. The specific formulas given in~\cite{DG,DG2}
gave only the series~\smash{$\Phi_\a^{(K,\s)}$} as discussed in Section~\ref{sec.Phi}, i.e., only
the first column of our matrix \smash{$\bPhi^{(K)}$}, because the vector~$\nu$ was always
assumed to be the one coming from the geometric ``flattening''. By varying~$\nu$,
one can get the other columns of~$\bPhi$. This variation produces a $q$-holonomic
system that turns out to be closely related to the ones for the generalized Kashaev
invariants that will be discussed in Section~\ref{sub.qholo}. This will be the theme of the next
section.

\section[Two q-holonomic modules]{Two $\boldsymbol{q}$-holonomic modules}\label{sec.NewHolonomy}

In Part~I we were led by the refined quantum modularity conjecture to find an entire
matrix~$\bJ^{(K)}$ of Habiro-like functions generalizing the Kashaev invariant, the
first column being the vector of formal power series found in~\cite{DG,DG2}.
In this section, we will study the structure of the other columns of this matrix and
will see that they have a natural ``$q$-holonomic structure'' in terms of an infinite
collection of functions that satisfy a recursion of finite length and hence all lie
in a finite-dimensional module. In Section~\ref{sub.qholo}, we explain this in detail
for the case of the~$4_1$ knot, where explicit formulas for the entries of~$\bJ^{(K)}$
were already given in Part~I. In the next subsection, we give the corresponding
formulas for the $5_2$ knot, where the matrix in question has size $4\times4$
instead of $3\times3$.
These are considerably more complicated than for the $4_1$ knot and have the interesting
new feature that Dedekind sums appear. In Section~\ref{sub.statesum},
we explain how these formulas could be guessed. This ansatz involves studying
the $q$-hypergeometric series defining the original Kashaev invariant via stationary
points and formal Gaussian summation, analogous to what was done in
Section~\ref{sub.PhiPerturbative} for Nahm sums and what will be done in
Section~\ref{sec.QMC41} for the Kashaev invariant of the $4_1$ knot. In the final
subsection, we discuss the point already alluded to at the end of
Section~\ref{sec.perturbation} that the power series in~$h$ studied there have a
$q$-holonomic structure with the \emph{same} coefficients as the one associated to
the matrix~$\bJ$. This is for the momenta purely experimental and is one of the many
mysteries associated with the subject. In Section~\ref{sub.qholonomicPhi}, we also
briefly mention two further conjectural objects associated to knots (or more generally
to half-symplectic matrices) that we believe share the same $q$-holonomic structure.

\subsection{Descendant Habiro-like functions}\label{sub.qholo}

It turns out that the first row of $\bJ$ and the first row of $\bPhih$ are basis
elements of the span of an inhomogeneous recursion, and the same holds (but now with
the corresponding homogeneous recursion) for each of the remaining rows of $\bJ$ and
of $\bPhih$ as well as for the matrix of $q$-series of~\cite{GZ:qseries}. To illustrate how this works, we give the complete formulas
for the matrix for the $4_1$ knot. The corresponding formulas for the $5_2$ knot, which
are considerably more complicated and illustrate several further refinements (like the
appearance of Dedekind sums), will be given in Section~\ref{sub.statesum}.

Collecting together our previous results for the $4_1$ knot for the reader's convenience,
we obtain that the matrix $\bJ=\bJ^{(4_1)}$ of periodic functions on $\BQ$ has the form
\[
\bJ(x) = \begin{pmatrix}
 1 & \J^{(0,1)}(x) & \J^{(0,2)}(x) \\
 0 & \J^{(1,1)}(x) & \J^{(1,2)}(x) \\
 0 & \J^{(2,1)}(x) & \J^{(2,2)}(x)
\end{pmatrix}
 =\sJ(q) = \begin{pmatrix}
 1 & \sJ^{(0,1)}(q) & \sJ^{(0,2)}(q) \\
 0 & \sJ^{(1,1)}(q) & \sJ^{(1,2)}(q) \\
 0 & \sJ^{(2,1)}(q) & \sJ^{(2,2)}(q)
\end{pmatrix}
\]
(with $q=\e(x)$ and omitting~$K$ as usual), where the elements of the first
row are given by
\begin{gather}
 \sJ^{(0,1)}(q) =\sQ^{(4_1)}_1(q)=\sJ^{(4_1)}(q)
 = \sum_{n=0}^\infty (q;q)_n \big(q^{-1};q^{-1}\big)_n ,\nonumber\\
 \sJ^{(0,2)}(q) =\sQ^{(4_1)}_2(q)
 =\frac{1}{2}\sum_{n=0}^\infty \big(q^{n+1}-q^{-n-1}\big) (q;q)_n \big(q^{-1};q^{-1}\big)_n\label{Q41J0Q0}
\end{gather}
(equations~\eqref{Hab41} and~\eqref{Q0for41}), with $\sQ^{(4_1)}_i(q)$ being
the elements of the Habiro ring defined and tabulated in Section~\ref{sub.4.3},
and that the elements of the other two rows are given by
\begin{align}
\sJ^{(1,1)}(q) &= \frac{1}{\sqrt{c}\sqrt[4]{3}} \sum_{Z^c = \z_6}
 \prod_{j=1}^c \bigl|1\m q^jZ\bigr|^{2j/c} ,\nonumber\\
\sJ^{(2,1)}(q) &= \frac{\rm i}{\sqrt{c}\sqrt[4]{3}}  \sum_{Z^c = \z_6^{-1}}
 \prod_{j=1}^c \bigl|1\m q^jZ\bigr|^{2j/c} ,\nonumber\\
\sJ^{(1,2)}(q) &= \frac{1}{2\sqrt{c} \sqrt[4]{3}}
\sum_{Z^c = \z_6} \bigl(Zq\m Z\i q\i\bigr) \prod_{j=1}^c \bigl|1\m q^jZ\bigr|^{2j/c},\nonumber
\\
\sJ^{(2,2)}(q) &= \frac{\rm i}{2\sqrt c \sqrt[4]{3}}
\sum_{Z^c = \z_6^{-1}} \bigl(Zq\m Z\i q\i\bigr)
\prod_{j=1}^c \bigl|1\m q^jZ\bigr|^{2j/c}\label{Q41abcd}
\end{align}
(equations \eqref{Def41J1}, \eqref{Q1for41} and the accompanying text).
The syntactical similarity between equations~\eqref{Q41J0Q0} and
equations~\eqref{Q41abcd} is striking, and leads directly to the $q$-holonomy.

To see this, we rewrite the two formulas in~\eqref{Def41J1} as
\be
\label{41holonomy0}
\sJ^{(0,1)}(q)=\sH^{(0)}_0(q) ,\qquad
\sJ^{(0,2)}(q)=\frac12\bigl(q\sH^{(0)}_1(q)-q\i\sH^{(0)}_{-1}(q)\bigr) ,
\ee
where \smash{$\big\{\sH^{(0)}_m(q)\big\}_{m\in\Z}$} is the sequence of elements of the Habiro
ring defined by
\be
\label{qHol41}
\sH^{(0)}_m(q) = \sum_{n=0}^\infty (q;q)_n \big(q^{-1};q^{-1}\big)_n q^{m n},
\qquad m\in\Z.
\ee
It is easy to see that this sequence satisfies the recursion relation
\be
\label{rec41}
q^{m+1} \sH^{(0)}_{m+1}(q) + (1-2 q^m) \sH^{(0)}_{m}(q)
+ q^{m-1} \sH^{(0)}_{m-1}(q) = 1,
\qquad m \in \BZ
\ee
(a similar, but homogeneous, recursion relation for the descendants of certain
$q$-series associated to the $4_1$ knot was given in~\cite[equation~(14)]{GGM} and
used in~\cite{GZ:qseries}) and also that the~$\BQ\big[q^\pm\big]$-module they span is
free of rank~3 with the top row of the matrix $\bJ$ as a basis. If we now introduce
two further sequences of functions of~$q$ (or periodic functions of~$x$, where
$q=\e(x)$) by
\begin{align}
\sH^{(1)}_m(q) &= \frac{1}{\sqrt c \sqrt[4]{3}}
\sum_{Z^c = \z_6} Z^m \prod_{j=1}^c \bigl|1\m q^jZ\bigr|^{2j/c} ,\nonumber
\\
\sH^{(2)}_m(q) &= \frac{\rm i}{\sqrt c \sqrt[4]{3}}
 \sum_{Z^c = \z_6^{-1}} Z^m \prod_{j=1}^c \bigl|1\m q^jZ\bigr|^{2j/c} ,\label{41holonomy12}
\end{align}
then \eqref{Q41abcd} says that the non-trivial elements of the
second and third rows of $\bJ$ are given by
\[
\sJ^{(i,1)}(q)=\sH_0^{(i)}(q) ,\qquad
\sJ^{(i,2)}(q)=\frac12\bigl(q\sH^{(i)}_1(q)-q\i\sH^{(i)}_{-1}(q)\bigr),
\qquad i=1,2 .
\]
Furthermore, we see that the first column of the matrix~$\bJ$, trivial
though it is, nevertheless belongs to the same $q$-holonomic module as the
other columns, since as well as equations~\eqref{41holonomy0}
and~\eqref{41holonomy12} we also have the relation
\be
\label{41holonomy1stcol}
\sJ^{(i,0)}(q)=
q \sH^{(i)}_1(q) + \sH^{(i)}_0(q) + q\i \sH^{(i)}_{-1}(q), \qquad i=0,1,2 ,
\ee
as we see by specializing the recursion~\eqref{rec41} and its counterparts for
$\sH^{(1)}_m$ and $\sH^{(2)}_m$ \hbox{to~$m=0$}. Then the quantitative version of the
``syntactical similarity'' noted above is that we can write the formulas~\eqref{Q41J0Q0}
and~\eqref{Q41abcd} or~\eqref{41holonomy12} and~\eqref{41holonomy1stcol} uniformly and
more compactly in matrix form as
\begin{gather}
\label{41matrixholonomy}
\sJ^{(4_1)}(q) =
\begin{pmatrix} \sH^{(0)}_{-1}(q) & \sH^{(0)}_0(q) &\sH^{(0)}_1(q) \\
 \sH^{(1)}_{-1}(q) & \sH^{(1)}_0(q) &\sH^{(1)}_1(q) \\
 \sH^{(2)}_{-1}(q) & \sH^{(2)}_0(q) &\sH^{(2)}_1(q) \end{pmatrix}
\begin{pmatrix} q^{-1} & 0 & \frac12 q \\ 1 & 1 & 0 \\ q & 0 & -\frac12 q^{-1}
\end{pmatrix} .
\end{gather}
Note that none of these equations are unique, since any one of them could be written
in infinitely many other ways by using the $q$-holonomy property, e.g., we could
specialize the recursions to any value of~$m$ other than~0 to get formulas for
$\sJ^{(i,0)}(q)$ different from~\eqref{41holonomy1stcol}. Similarly, we could
rewrite~\eqref{41matrixholonomy} by taking three other columns (or linear combinations
of columns) of the $\sH$-matrix for the first factor on the right, with the
corresponding new matrix of Laurent polynomials in the second factor.

More interesting is that there is also nothing sacred about the particular collection
$\sH^{(i)}(q)$ of functions of~$q$ that we chose to define our $q$-holonomic system,
and that there infinitely many other collections, even with completely different
indexing sets (e.g., $\Z^2$ instead of~$\Z$) that could be used instead and that might
have been found it we had given a different combinatorial description of knot. However,
the module over $\Q\big[q,q\i\big]$ that they generate is at least conjecturally intrinsic to
the knot and is simply the span of the columns of~\smash{$\bJ^{(K)}$}, which therefore
constitute a canonical basis indexed by~$\calP$. This is one of the most mysterious
aspects of our matrix invariants. We will return to it in at the end of this section
in connection with other possible representations of the same abstract $q$-holonomic
module.

We end the subsection with a final remark. Despite the apparent similarity in the
formulas for the elements of the first row and all other rows of the matrix $\bJ$,
there is a crucial difference between formulas like~\eqref{Q41J0Q0} or~\eqref{qHol41}
for the top rows of our matrix and formulas like~\eqref{Q41abcd}
or~\eqref{41holonomy12}
for the other rows: the former are sums over the lattice points of a cone and
hence satisfy an inhomogeneous linear $q$-difference equation, whereas the latter
are sums over periodic groups~$\BZ/c\BZ$ and hence have no boundary terms and
satisfy a homogeneous equation. Another difference, to which we hope to return
in~\cite{GSZ:habiro} in the context of Habiro rings for general number fields, is
that~\eqref{Q41J0Q0} and~\eqref{qHol41} obviously give algebraic integers when $q$~is
a root of unity, whereas~\eqref{Q41abcd} or~\eqref{41holonomy12} give algebraic
integers in some non-evident way, since it is not obvious (but in fact true)
that the sums in these formulas are divisible by~$\sqrt c$. We will find
exactly the same behavior for the elements of the $\bJ$-matrix for the $5_2$ knot
in the next subsection.

\subsection[The J-matrix for the 5\_2 knot]{The $\bJ$-matrix for the $\boldsymbol{5_2}$ knot}
\label{sub.52statesum}

In this subsection, we describe that analogues of the formulas just given for our second
standard knot~$5_2$, because as usual the figure~8 knot has such special properties that
some of the interesting features are obscured. 

The Kashaev invariant of the $5_2$ knot is given by
\be
\label{Kashaev.52}
 \calJ^{(5_2)}(q) =\sum_{m=0}^{\infty} \sum_{k=0}^m
q^{-(m+1)k} \frac{(q;q)_m^2}{\big(q^{-1};q^{-1}\big)_k} .
\ee
(See~\cite[equation~2.3]{K97}.) This is manifestly an element of the Habiro ring. We
generalize it to the two-parameter family of elements of the Habiro ring given by
\be
\label{52desc}
\sH^{(0)}_{a,b}(q) =\sum_{m=0}^{\infty} \sum_{k=0}^m
q^{-(m+1)k + am+bk} \frac{(q;q)_m^2}{\big(q^{-1};q^{-1}\big)_k}, \qquad a,b \in \BZ .
\ee
These again form a $q$-holonomic module in the sense of~\cite{WZ}, meaning
that they satisfy recursions like~\eqref{rec41} (though in this case more complicated,
and omitted here) and hence generate a~$\Q\big[q,q\i\big]$-module of finite rank. Here the rank
is~4 and the $q$-holonomic module is generated (as we expect to hold for every knot)
by the first row of the matrix $\bJ$ of the knot,
\begin{align*}
 \sJ^{(0,0)}(q) &
 =- \sH^{(0)}_{0,0}(q)+q^{-1} \sH^{(0)}_{-1,0}(q)+\sH^{(0)}_{0,-1}(q)=1 , \qquad
\sJ^{(0,1)}(q) =\sH^{(0)}_{0,0}(q) , \\
\sJ^{(0,2)}(q) & =\sH^{(0)}_{0,0}(q) \m q^{-1}\sH^{(0)}_{-1,0}(q) , \qquad
\sJ^{(0,3)}(q) = 2\sH^{(0)}_{0,0}(q) \m q^{-1} \sH^{(0)}_{-1,0}(q)
+ \sH^{(0)}_{-1,1}(q) .
\end{align*}

Just as in the case of the~$4_1$ knot, we find that the further three rows are given
by the \emph{same} linear combinations of three other two-parameter families
$\sH^{(i)}_{a,b}$ ($1\le i\le3$) of functions, i.e., we have
\begin{align*}
 \sJ^{(i,0)}(q) &
 = - \sH^{(i)}_{0,0}(q)+ q^{-1}\sH^{(i)}_{-1,0}(q)\m\sH^{(i)}_{0,-1}(q) ,\qquad
\sJ^{(i,1)}(q) =\sH^{(i)}_{0,0}(q) , \\
\sJ^{(i,2)}(q) & =\sH^{(i)}_{0,0}(q) \m q^{-1}\sH^{(i)}_{-1,0}(q) , \qquad
\sJ^{(i,3)}(q) = 2\sH^{(i)}_{0,0}(q) \m q^{-1} \sH^{(i)}_{-1,0}(q)
+ \sH^{(i)}_{-1,1}(q)
\end{align*}
for $i=0,1,2,3$. The formulas for the functions for $i\ne0$, whose origin will be
indicated in Section~\ref{sub.statesum}, are of the same type as the corresponding
ones for the~$4_1$ knot (equation~\eqref{Q41abcd}), though considerably more complicated,
but are completely different from~\eqref{52desc}, namely
\be
\label{52H}
\sH_{a,b}^{(i)}(x) = \frac 1{c \sqrt{3\xi_i-2}}  \frac{\theta_{1,i}^{c-1}
 \calD_{\z}(\z \theta_{1,i})^2}{\calD_{\z}\big(\z^{-1} \theta_{2,i}^{-1}\big)}
\sum_{k,m \bmod c} \zeta^{-(k+1)m} \theta_{1,i}^{-m} \theta_{2,i}^{-k}
\frac{(\zeta \theta_{1,i};\zeta)_k^2 }{\big(\zeta^{-1}\theta_{2,i}^{-1};\zeta^{-1}\big)_m} ,
\ee
where $c=\den(x)$, $\z=\e(-x)$ and
\smash{$\theta_{1,i}^c=-\xi_i^{-3}$} and \smash{$\theta_{2,i}=\xi_{i}^{-2}$} are any
choice of~$c$-th roots of~\smash{$-\xi_i^{-3}$} and $\xi_{i}^{-2}$ and $\xi_1$ (resp.,
$\xi_2$, $\xi_3$) the complex root of the equation $\xi^3-\xi^2+1=0$ (as in
Section~\ref{sub.p}) with negative (resp., positive, zero) imaginary part.
Here $\calD_\z(x)$ is the renormalized version of the cyclic quantum dilogarithm
$D_\z(x)$ defined for $q=\e(a/c)$ by
\be
\label{Dm}
\calD_q(x) = {\rm e}^{-2 \pi {\rm i} s(a,c)/2} D_q(x)
= {\rm e}^{-2 \pi {\rm i} s(a,c)/2}
\exp \Biggl(\sum_{j=1}^{c-1} \frac{j}{c} \log\big(1-q^j x\big) \Biggr) ,
\ee
where $s(a,c)$ is the Dedekind sum~(cf.~\cite{HZ,Rademacher}) and
where the logarithm is the principal one away from the cut at the negative real
axis and is defined on the cut as the average of the principal branches just above
and just below. The cyclic quantum dilogarithm
appears in the expansion of Faddeev's quantum dilogarithm at roots of unity (see,
for example,~\cite{GK:evaluation,Kashaev:star}) and plays a key role in the
definition of the near units associated to elements of the Bloch group~\cite{CGZ}.

It is worth mentioning that the formulas~\eqref{Q41abcd} and~\eqref{41holonomy12}
for the $4_1$ knot can also be written in terms of the modified cyclic quantum
dilogarithm~$\calD_q$, because \smash{$\prod_{j=1}^c \bigl|1\m q^jZ\bigr|^{2j/c}$} can be
rewritten as $\calD_q(Z) \calD_{q^{-1}}\big(Z^{-1}\big)$. In fact, we expect formulas of
this type, involving multiplicative combinations of the $\calD_q$'s corresponding
to the combinations defining the element of the Bloch group of~$F$ corresponding
to the knot, to exist for all knots.

\subsection{State-sums}
\label{sub.statesum}

In this subsection, we explain where the formulas just given come from. More precisely,
we discuss a heuristic method to discover a formula for the first column of the
matrix $\bJ^{(5_2)}$ given a formula for its top entry i.e., for the Kashaev invariant
of the knot. This method produces periodic functions similar to the constant term of
the formal power series \smash{$\Phi^{(\s)}_\a(h)$} discussed in Section~\ref{sec.perturbation}.
It also generalizes to the further columns, by replacing the Kashaev invariant by
the other in its top row (i.e., in the row of the matrix that is expected always to
have entries belonging to the rational Habiro ring), thus producing predictions for
the entire matrix~$\bJ$. This is useful in particular for the numerical confirmation
of the generalized quantum modularity conjecture.\looseness=-1

Our starting point is the formula~\eqref{Kashaev.52} for the Kashaev invariant of the
$5_2$ knot.
Let
\[
b_{k,\ell}(q) =q^{-(\ell+1)k} \frac{(q;q)_\ell^2}{\big(q^{-1};q^{-1}\big)_k}
\]
denote the summand of the Kashaev invariant of $5_2$ in
equation~\eqref{Kashaev.52}. The function $b_{k,\ell}(q)$ (which is proper
$q$-hypergeometric in the sense of~\cite{WZ})
satisfies the linear $q$-difference equations
\begin{gather}
\label{eq.b1}
\frac{b_{k+1,\ell}(q)}{b_{k,\ell}(q)} = q^{-\ell} \big(1-q^{k+1}\big)^2,
\qquad
\frac{b_{k,\ell+1}(q)}{b_{k,\ell}(q)} = q^{-(k+1)} \frac{1}{1-q^{-\ell-1}},
\end{gather}
whose right-hand sides are in $\BQ\big(q,q^k,q^\ell\big)$.
It follows that for natural numbers $r$, $s$ we have
\[
\frac{b_{k+r,\ell}(q)}{b_{k,\ell}(q)} = q^{-r\ell} \big(q^{k+1};q\big)_r^2,
\qquad
\frac{b_{k,\ell+s}(q)}{b_{k,\ell}(q)} = q^{-(k+1)s} \frac{1}{\big(q^{-\ell-1};q^{-1}\big)}
\]
and hence
\begin{align*}
\frac{b_{k+r,\ell+s}(q)}{b_{k,\ell}(q)} &= q^{-k s -r \ell-(r+1)s}
\frac{\big(q^{k+1};q\big)_r^2}{\big(q^{-\ell-1};q^{-1}\big)_s} .
\end{align*}
Setting $q^k=z_1$, $q^\ell=z_2$, $q=1$ and equating the ratios of
equations~\eqref{eq.b1} to~$1$, we get the gluing equations for $(z_1,z_2)$
\begin{equation}
\label{eq.ge}
z_2^{-1}(1-z_1)^2=1, \qquad z_1^{-1}\big(1-z_2^{-1}\big)^{-1}=1 .
\end{equation}
Although the summation for the Kashaev invariant when $q$ is a primitive
$N$-th root of unity is a subset of $[0,N-1]^2$ and when $\big(q^k,q^\ell\big)$
is near $(z_1,z_2)$ is outside the summation range, we will pretend
that we have performed analytic continuation. Choose $\zeta=\e(a/c)$
where $(a,c)=1$ and~${c>0}$ and \smash{$(\th_1,\th_2)=\big(z_1^{1/c},z_2^{1/c}\big)$}.
In other words, we choose $\th_i$ to be arbitrary $c$-th roots of~$z_i$~($i=1,2$).
Then we can define $a_{r,s}(\th_1,\th_2;\z)$ by
\begin{align*}
 a_{r,s}(\th_1,\th_2;\z)&=
\frac{b_{k+r,\ell+s}(q)}{b_{k,\ell}(q)} \big|_{q^k=\th_1, q^\ell=\th_2, q=\zeta} = \zeta^{-(r+1)s} \th_1^{-s} \th_2^{-r}
\frac{(\zeta \th_1;\zeta)_r^2 }{\big(\zeta^{-1}\th_2^{-1};\zeta^{-1}\big)_s} .
\end{align*}
The principle of equipeaked Gaussians in the asymptotics of $\J(\g X)$
with $\g=\sma abcs \in \SL_2(\BZ)$ (as used in Section~\ref{sub.proof41}
for the case of the $4_1$ knot) suggests the expression
\begin{align}
\label{eq.a1}
S(\th_1,\th_2;\zeta)=\sum_{r,s=0}^{c-1} a_{r,s}(\th_1,\th_2) .
\end{align}
The first observation is that the sum in equation~\eqref{eq.a1} is
$c$-periodic, i.e., that $r,s \in \BZ/c\BZ$.
This follows from the fact that $(z_1,z_2)$ satisfy the gluing
equations~\eqref{eq.ge}. A second observation, which we will not make use of,
if the fact that $b_{k,\ell}(q)$ determines $a_{r,s}(\th_1,\th_2)$ according to
the above definitions. Conversely, $a_{r,s}(\th_1,\th_2;\z)$ determines
$b_{k,\ell}(q)$ by $a_{r,s}(1,1;\z)=b_{r,s}(\zeta) $.
A curious consequence of this is that $
S(1,1;\z)=\J(\z)$
recovers the Kashaev invariant.
The gluing equations~\eqref{eq.ge} can be solved as follows: $z_1=-\xi^{-3}$, $z_2=\xi^{-2}$,
where $\xi^3-\xi^2+1=0$. The three solutions give rise to the three
embeddings of the trace field of $5_2$ into the complex numbers.
For~${\z=\e(a/c)}$, let $F_c=F(\z)$ and $F_{G,c}=F_c(\th_1,\th_2)$,
giving extensions $
F \subset F_c \subset F_{G,c}$,
where~$F_{G,c}/F_c$ is an abelian Galois (Kummer) extension with group $(\BZ/c\BZ)^2$
and $S(\th_1,\th_2) \in F_{G,c}$. To find how $S(\th_1,\th_2;\z)$
transform under the Galois group, we compute
\begin{align*}
 \frac{a_{r,s}(\zeta \th_1,\th_2;\z)}{a_{r+1,s}(\th_1,\th_2;\z)}
 &= \th_2 (1-\zeta \th_1)^{-2} = a_{1,0}(\th_1,\th_2;\z)^{-1},\\
 \frac{a_{r,s}(\th_1,\zeta \th_2;\z)}{a_{r,s+1}(\th_1,\th_2;\z)}
 &= \zeta \th_1\big(1-\zeta^{-1}\th_2^{-1}\big) = a_{0,1}(\th_1,\th_2;\z)^{-1}
\end{align*}
(where the left-hand side of the above equations is independent of
$r$ and $s$ hence it must equal to the right-hand side).
Since the sum in equation~\eqref{eq.a1} is $c$-periodic, it follows that
\begin{subequations}\label{eq.Sq1Sq2}
\begin{align}
\label{eq.Sq1}
 S(\zeta \th_1,\th_2;\z) &= S(\th_1,\th_2;\z) \th_2 (1-\zeta \th_1)^{-2}
 = S(\th_1,\th_2;\z)a_{1,0}(\th_1,\th_2;\z)^{-1},\\
\label{eq.Sq2}
 S(\th_1, \zeta \th_2;\z) &= S(\th_1,\th_2;\z) \zeta \th_1 \big(1-\zeta^{-1} \th_2^{-1}\big)
 =S(\th_1,\th_2)a_{0,1}(\th_1,\th_2;\z)^{-1}.
\end{align}
\end{subequations}
To fix the Galois invariance of $S(\th_1,\th_2;\z)$, we consider the product
\[
P(\th_1,\th_2;\z) = \prod_{r=0}^{c-1} \big(1-\zeta^{r+1} \th_1\big)^{2r}
\prod_{s=0}^{c-1} \big(1-\zeta^{-s-1}\th_2^{-1}\big)^{-s}.
\]
We can rewrite the above product using the cyclic quantum dilogarithm
function~\eqref{Dm} as follows%
\begin{gather}
\label{eq.PD}
P(\th_1,\th_2;\z)
= z_1^{-1} z_2^{-1}
\frac{D_\z(\th_1)^2}{D_{\z^{-1}}\big(\th_2^{-1}\big)} .
\end{gather}
From the transformation property for the cyclic quantum dilogarithm
\[
\frac{D_\z(x)}{D_\z\big(\z^{-1}x\big)} =\frac{(1-x)^c}{1-x^c}, \qquad
D_\z(x) D_{\z^{-1}}(x) = (1-x^c)^c(1-x)^c
\]
and the fact that $(z_1,z_2)$ solve the gluing equations~\eqref{eq.ge}, we obtain that
\begin{subequations}\label{eq.Ps1Ps2}
\begin{align}
\label{eq.Ps1}
P(\zeta \th_1, \th_2;\z) &=
P(\th_1,\th_2;\z) \big( \th_2^{-1} (1-\zeta \th_1)^{2} \big)^c
= P(\th_1,\th_2;\z)  a_{1,0}(\th_1,\th_2;\z)^c, \\
\label{eq.Ps2}
P(\th_1, \zeta \th_2;\z) &= P(\th_1,\th_2;\z) \big(
\zeta^{-1} \th_1^{-1} \big(1-\zeta^{-1}\th_2^{-1}\big)^{-1} \big)^c
= P(\th_1,\th_2;\z)  a_{0,1}(\th_1,\th_2;\z)^c .
\end{align}
\end{subequations}
Equations~\eqref{eq.Sq1Sq2} and~\eqref{eq.Ps1Ps2}
imply that
\[
P^{1/c}(\th_1,\th_2;\z) S(\th_1,\th_2;\z) \in \ve^{1/c} F_c,
\]
where $\ve$ is a unit, which in fact coincides with the one constructed
in~\cite{CGZ}.

The expression given in the above equation, after multiplication by
a prefactor, coincides with \smash{$\sH^{(i)}_{0,0}(x)$} of equation~\eqref{52H}
if we choose $\th_1$ and $\th_2$ corresponding to the root $\xi_i$ of
$\xi^3-\xi^2+1=0$.

In this way, we have succeeded in guessing the entries of the first column of
the mat\-rix~$\bJ^{(5_2)}$ of the $5_2$ knot starting from the formula~\eqref{Kashaev.52}
for its Kashaev invariant. All of this seems to reek a~little of ``black magic''.
But the same method applied to the case of the $4_1$ knot (whose Kashaev invariant
is given in~\eqref{Q41J0Q0}) reproduces the formulas given in~\eqref{Q41abcd}. In
fact, we believe that this will work for any knot, giving each entry of the first
column of the $\bJ$-matrix as a sum of products of cyclic quantum dilogarithms with
summands modelled on the solution of the Neumann--Zagier gluing equations of the knot
triangulation in the same way that the expression~\eqref{eq.PD} is modelled on the
gluing equations~\eqref{eq.ge}.

\subsection[The q-holonomic module of formal power series]{The $\boldsymbol{q}$-holonomic module of formal power series}
\label{sub.qholonomicPhi}

We now explain one of the most mysterious aspects of our story, the appearance
of two very different realizations of the same $q$-holonomic system in the
contexts of state sums and of perturbative formal power series. In fact, as
we will indicate briefly at the end, we believe that there are actually four
\hbox{$q$-holonomic} systems, of totally different origins, given by recursions
with the same Laurent polynomials as coefficients.

We begin by recalling the derivation of the perturbative series in~$h$ from Nahm sums,
as described in Section~\ref{sub.PhiPerturbative}.
The Nahm sums $F_\Xi(q)$ as defined in~\eqref{Nahm2}, with $H$ fixed and $\nu$
varying over $\text{diag}\big(AB^t\big)+2\Z^N$, form a module of finite rank over
the ring $R=\Z\big[q^{\pm1}\big]$ of Laurent polynomials in~$q$. For instance, for the
original Nahm sums as defined in~\eqref{Nahm1}, we have the recursions\looseness=-1
\[
F_{A,b}(q) \m F_{A,b+e_k}(q) = q^{\frac12e_k^tAe_k+b^te_k}
F_{A,b+Ae_k}(q), \qquad k=1,\dots,N
\]
(here $e_k$ denotes the $k$-th basis vector of~$\Z^N$), as one sees by noting
that $(1-q^{n_k})/(q;q)_{n_k}$ vanishes if~$n_k=0$ and equals $1/(q;q)_{n_k-1}$
if~$n_k\ge1$, so that the difference on the left corresponds simply to
shifting the multi-index~$n$ by~$e_k$. A similar but more complicated calculation
(again corresponding to the shift $n\mapsto n+e_k$ in the definition of the sum
and using the relationship between Pochhammer symbols with nearby indices) gives
a collection of $N$ recursion relations among the various $F_{((A B),\nu)}$ with
fixed~$(A B)$. This system is always $q$-holonomic \cite{WZ}, meaning in
particular that the solution space is finite-dimensional.

When we discussed the asymptotic behavior of the Nahm sum~\eqref{Nahm2} in
Section~\ref{sub.PhiPerturbative}, we first rewrote the sum as in~\eqref{Nahm3} and
then replace the sum over~$n\in\Z^N$ by an integral over~$x\in\R^N$. It is then
clear that exactly the same argument (replacing $x$ by $x+e_k$) shows that the formal
power series arising from Gaussian integrals near the various stationary points of
the integral satisfy the same system of recurrences, and hence also form a
\hbox{$q$-holonomic} module. In favorable cases, including all the ones we have looked
at, the rank of this system will be equal to the cardinality of~$\calP_K$, because
the characteristic variety of the system coincides with the variety~$V_H$ as defined
in Section~\ref{sub.HalfSymplectic}.

The surprising discovery is that the abstract \hbox{$q$-holonomic} module that we
obtain this way is the \emph{same} as the one that we found in the first three
subsections of this section from the Habiro-like functions, i.e., although the functions
of~$q$ are completely different and are even defined in different places, the modules in
question are spanned by sequences of elements indexed in the same way and satisfying the
same recursions over~$\Q\big[q^{\pm1}\big]$, and moreover that the special basis indexed
by~$\calP_K$ is given in both systems by the \emph{same} linear combination of these
elements. (Compare the discussion at the end of Section~\ref{sub.qholo}.)
We believe that this coincidence of two~$q$-holonomic structures will hold,
not only for the matrix invariants of knot complements, but more generally for
corresponding objects associated to any half-symplectic matrix in the sense of
Section~\ref{sec.perturbation}. This will be further studied in~\cite{GSZ:habiro}.

However, a big surprise of the sequel~\cite{GZ:qseries} to this paper is that the very
same $q$-holonomic structure actually occurs a third time in terms of the $q$-series
coming from state integrals that are studied there. We believe that this coincidence
holds because these three objects are simply different realizations of the same
``function-near-$\Q$'' belonging to a generalized Habiro ring, with the ``nearness''
being realized for the Habiro-like functions by approaching a given rational number
through nearby rational numbers of slowly growing height, and in the case of the
$q$-series by approaching a rational number from above in the upper half-plane
(or equivalently, approaching a root of unity radially in the unit disk). We have
checked the agreement of the relations over~$\Q\big[q,q\i\big]$ among the columns of the
$\bJ$- and $\bPhi$-matrices for both the~$4_1$ and~$5_2$ knots (although we do not
describe that verification in this paper because the specific formulas that were used
for~$\bJ$ for these two knots in Sections~\ref{sub.qholo}
and~\ref{sub.52statesum}--\ref{sub.statesum} are not the same as the ones coming from
ideal triangulations and Neumann--Zagier data and are rather complicated), while the
corresponding verification for the $q$-series for the same two knots is given
in~\cite{GGM} and discussed in~\cite{GZ:qseries}. We conjecture that these
recursive relations coincide with the ones defined in current work of Rinat
Kashaev and the first author~\cite{GK:desc} from the colored Jones polynomials.
But we should emphasize that we still have no idea
why any of these $q$-holonomic modules has a canonical basis indexed by~$\calP$.

\section[Proof of the modularity conjecture for the 4\_1 knot]{Proof of the modularity conjecture for the $\boldsymbol{4_1}$ knot}
\label{sec.QMC41}

In this section, we give our proof of the quantum modularity conjecture
for the figure 8 knot, announced several years ago. Another proof was given
by \cite{BD}, as well as proofs of the quantum modularity conjecture for a few
other knots, but we give our proof for completeness and because the
point of view here is somewhat different from the one there.

We denote by $\J(x)$ the $\J$-function for the~$4_1$ knot, as given explicitly by
equation~\eqref{Hab41} with~${q=\e(x)}$. We have to show that
\be
\label{tag2}
\J\left(\frac{aX+b}{cX+d}\right) \sim
(cX+d)^{3/2} \J(X) \Phih_{a/c}\left(\frac{2\pi {\rm i}}{c(cX+d)}\right)
\ee
to all orders in~$1/X$ as $X$ tends to infinity with bounded denominator with
${\g=\sma abcd\in\SL_2(\Z)}$ fixed and~$c>0$, where \smash{$\Phih_{a/c}(h)={\rm e}^{V/c^2h} \Phi_{a/c}(h)$}
is the completed version of a formal power series~$\Phi_{a/c}(h)$ with algebraic
coefficients and where
\[
V=\text{Vol}\big(S^3\ssm 4_1\big)=2.02988\dots
\] is the volume of
the complement of this knot. We give the proofs separately for the special case $\g=S$, $\a:=a/c=0$ and for the general case, since all the main ideas are already visible
for the former and the details are much simpler.

\subsection[The case of a=0]{The case of $\boldsymbol{\a=0}$}
\label{sub.proof41.0}

We begin with the case $\a=0$, which makes it clear where the
factor~$\J(X)$ in equation~\eqref{tag2} comes from. We use the notation
\[
P_n(x)=|(q;q)_n|^2 ,
\]
for $q=\e(x)$ with $x$ rational, so that $P_n(x)$ is the $n$-th summand in the
definition of~$\J(x)$, and denote by $c_r$ ($r\ge0$) the numbers defined by the
Taylor expansion
\[
\cot\left(\frac\pi6\m\frac x2\right) = \sum_{r=0}^\infty c_r \frac{x^r}{r!} ,
\]
the first values being given by the table
\begin{center} \def\arraystretch{1.2}
\begin{tabular}{|c|cccccccc|}\hline $r$ & 0 & 1 & 2 & 3 & 4 & 5 & 6 & 7 \\ \hline
 $c_r$ & $\sqrt3$ & 2 & $2\sqrt3$ & $10$ & $22\sqrt3$ & $182$ & $602\sqrt3$ & $6970$
 \\ \hline
\end{tabular}
\end{center}
Note that these numbers can also be written $c_r=2 \operatorname{Im}({\rm i}^{-r}\L_{-r}(\e(1/6)))$
and hence have a~natural extrapolation backwards by $c_{-1}=0$, $c_{-2}=-V$.
We want to study $\J(-1/X)=J(1/X)$ as $X$ tends to infinity in the fixed
residue class~$\b\pmod1$, with $\b$ rational. The summands in~\eqref{J0for41}
are all positive (that is what makes this case much easier to treat than the
general one), and it is easy to find their local peaks, which occur
near $n=\big(m+\frac56\big)X$ for $0\le m<\den(\b)$, where~$\den(b)$ is the denominator
of $b$. (Notice that the terms for larger values of~$m$ are 0 anyway, since
$P_n(x)$ vanishes for $n\ge\den(x)$.) The following proposition, which is valid at
a fixed peak even for $X$ real, gives the asymptotic value of the summand $P_n(x)$
for $n$ in each of these peaks. As usual, $B_r(x)$ denotes the $r$-th Bernoulli
polynomial.

\begin{Proposition}
\label{prop.a1}
Fix an integer $m\ge0$, and set $M=m+\frac56$. Then for $X$ tending to infinity
and $n$ an integer of the form $MX-\nu$ with $|\nu|\ll X$ we have the asymptotic
expansion
\begin{equation}
\label{tag3}
\log P_n\left(\frac1X\right) \si \frac V{2\pi} X + \log X + \log P_m(X)
+ \sum_{k=1}^\infty c_{k-1} \frac{B_{k+1}(\nu)}{(k+1)!} \left(-\frac{2\pi}X\right)^k .
\end{equation}
\end{Proposition}

\noindent{\bf Note:} By the above remark we can omit the first term and sum over
$k\ge-1$ instead.

\begin{proof}
We first note that for $q=\e(1/X)$ we have
\[
\log\left(\frac{P_n(1/X)}{P_{n-1}(1/X)}\right)
= \log\bigl|1\m q^n\bigr|^2 = \log\left|1\m \e\left(\frac16 +\frac\nu X \right)\right|^2
= -\sum_{k=1}^\infty \frac{c_{k-1}}{k!} \left(\frac{2\pi\nu}X\right)^k
\]
(here we have used that $\frac {\rm d}{{\rm d}x}\log|1-\e(x)|^2=2\pi\cot(\pi x) $),
in agreement with equation~\eqref{tag3} to all orders in~$1/X$ since
$B_{r+1}(\nu+1)-B_{r+1}(\nu)= (r+1)\nu^r$. This proves~\eqref{tag3} up to a power
series independent of~$n$ (but depending a priori on~$\a$ and~$m$). The
full assertion uses the shifted Euler--Maclaurin summation formula; we omit the details.
\end{proof}

Note that equation~\eqref{tag3} does not make sense if~$X$ is rational and $m\ge\den(X)$,
since then~$P_n(1/X)$ and $P_m(X)$ vanish, but we will use it only in the
exponentiated form
\begin{equation}
\label{tag4}
P_n\left(\frac1X\right) \si P_m(X) X {\rm e}^{V/h}
\exp\Biggl( \sum_{r\ge1} (-1)^rc_{r} \frac{B_{r+1}(\nu)}{(r+1)!} h^r\Biggr),
\qquad h=\frac{2\pi}X ,
\end{equation}
which holds also in this case. The key point here is that the only
dependence on~$m$ of the expression on the right-hand
side is the factor $P_m(X)$, which equals $P_m(\b)$ if $X$ goes to infinity
in the fixed class $\Z+\b$ modulo~1.
Moreover, since $B_2(\nu) = \nu^2+\O(\nu)$ and $B_{r+1}(\nu)=\O\big(\nu^{r+1}\big)$ for
$r>2$, the exponential factor in equation~\eqref{tag4} has
the form \smash{${\rm e}^{-\sqrt 3\nu^2h/2}\phi\big(\nu\sqrt h,\sqrt h\big)$}, where
\[
\phi\big(\v\sqrt h,\sqrt h\big)=
\exp\bigg(\frac{\sqrt{3}}{2} h \left(\v-\frac{1}{6} \right)+
\sum_{r\ge 2} (-1)^rc_{r} \frac{B_{r+1}(\v)}{(r+1)!} h^r\bigg)
 .
\]
The contribution to $\J(1/X)$ from the $m$-th peak is equal to $P_m(\b)X{\rm e}^{V/h}$ times
the sum of this exponential factor over all $\nu$ with $|\nu|\ll X$ in a fixed residue
class $\nu_0\pmod1$, where $\nu_0\equiv-MX\pmod1$. But by the Poisson summation
formula and the fact that the Fourier transform of a~Gaussian decays
exponentially, we have that
\begin{align}
\sum_{\substack{\nu\equiv\nu_0 \text{(mod~1)}\\|\nu|\ll X}}
{\rm e}^{-\sqrt 3\nu^2h/2}\phi\big(\nu\sqrt h,\sqrt h\big)
={}& \int_{-\infty}^\infty {\rm e}^{-\sqrt 3 \nu^2h/2}\phi\big(\nu\sqrt h,\sqrt h\big) {\rm d}\nu\nonumber\\
= {}&\sqrt X I_{\sqrt3}\bigl[\phi\big(t,\sqrt h\big)\bigr] ,\label{VZ}
\end{align}
to all orders in~$h$, where $I_\l$ for $\l>0$ denotes the linear map
from $\C[[t]]$ to~$\C$ defined by
\[
I_\l\bigl[\phi(t)\bigr] = \frac1{\sqrt{2\pi}} \int_{-\infty}^\infty
{\rm e}^{-\l t^2/2}\phi(t) {\rm d}t ,
\qquad I_\l\bigl[t^n\bigr] =
\begin{cases}
(n-1)!! \lambda^{-(n+1)/2} &\text{if $n$ is even,}\\ 0
&\text{if $n$ is odd.}
\end{cases}
\]
Here $(n-1)!!$ as usual denotes the ``double factorial''
$(n-1)\times(n-3)\times\cdots\times3\times1$. A~discussion of the estimates that prove
equation~\eqref{VZ} is given in~\cite{VZ,Zagier} and~\cite{GZ:asymptotics}. It
follows that the contribution to $\J(1/X)$ from the $m$-th peak is equal to
$P_m(\b)X^{3/2}{\rm e}^{V/h}\Phi_0(h)$ to all orders, where $\Phi_0(h)$ is the power series
given by equation~\eqref{Phi410}. Hence, $\J(1/X)$ itself equals
$\J(\b)X^{3/2}{\rm e}^{V/h}\Phi_0(h)$ to all orders, as claimed. Note that $\Phi_0(h)$ equals
$3^{-1/4}$ times a power series in~$h$ with coefficients in~$\Q\big(\sqrt3\big)$ and leading
coefficient~1, since $I_\l[\phi(t)]$ has coefficients in \smash{$\l^{-1/2}\Q(\l)$} for any power
series~$\phi(t)$ with rational coefficients. (That it is a power series in~$h$ rather
than merely~\smash{$\sqrt h$} follows from the fact that $I_\l$ annihilates odd functions.)

This concludes the proof of equation~\eqref{tag2} when $\a=0$, with $\Phi_0(h)$ given by
\be
\label{Phi410}
\Phi_0(h) = I_{\sqrt{3}}\big(\phi\big(t\sqrt{h},\sqrt{h}\big)\big) .
\ee

\subsection{The general case}
\label{sub.proof41}

We now apply the same analysis to the expansion of $\J(x)$ around an
arbitrary rational number~$\a$. The second part of the argument, replacing
sums by integrals and computing them by using the functional~$I_\l$, is
unchanged, but the analogue of Proposition~\ref{prop.a1} is now slightly more
complicated, since the asymptotic formula for $P_n(x)$ near
the $m$-th peak depends on both~$m$ and the residue class of~$n$ modulo the
denominator of~$\a$. We use the notations given above, i.e.,
\smash{$x=\frac{aX+b}{cX+d}$} with~$X$ tending to infinity in the class
$\b\!\pmod1$ and $M=m+\frac56$
with $0\le m<\den(\b)$, but now we also fix a residue class $r\!\pmod c$ and
consider $n$ satisfying
\[
n \equiv r+m d\pmod c , \qquad n = \frac M\hb+c \nu
\]
with $|\nu|\ll X$. (Notice that $\nu$ has the opposite sign to the one
used for~$c=1$.)

\begin{Proposition}
\label{prop.a2}
For fixed $m<\den(\b)$ and $r\in\Z/c\Z$ and for $n$ and $X$ tending to
infinity as in~$(5)$, we have the asymptotic formula
\begin{gather}\label{tag6}
\log P_n\bigl(x\bigr) \si \log\left(\frac1\hb\right) + \log P_m(\b)
+ \sum_{k=-1}^\infty C_k^{(r)}(\nu) h^k ,
\end{gather}
valid to all orders in~$h$, where $C_k^{(r)}$ is the polynomial of degree
$k-1$ defined by
\begin{gather}
\label{tag7}
C_k^{(r)}(\nu) = - 2 \operatorname{Re}\biggl[\frac{{\rm i}^{-k}}{(k+1)!}
 \sum_{j=1}^c \L_{1-k}\big(\z_\a^{r+j}Z\big) B_{k+1}\left(\nu+\frac jc\right)\biggr],
\qquad\! r\in\Z/c\Z, \quad\! k\ge-1\!\!\!
\end{gather}
with $\z_\a=\e(\a)$ and $Z=\e(-5/(6c))$.
\end{Proposition}
The proof of this proposition, which we omit, is similar to that of
Proposition~\ref{prop.a1}, the main point again being that the difference
of the right-hand sides of~\eqref{tag6} for $n$ and $n-1$ is given by
\begin{align*}
\sum_{k=-1}^\infty \big[C_k^{(r)}(\nu)\m C_k^{(r-1)}(\nu-1)\big] h^k
& = - 2 \operatorname{Re}\Biggl[\sum_{k=0}^\infty \L_{1-k}(\z_\a^rZ)
\frac{(-{\rm i}h\nu)^k}{k!}\Biggr] \\
 &
= \log\bigl|1\m\z_\a^rZ{\rm e}^{-{\rm i}\nu h}\bigr|^2 = \log|1-q^n|^2 \\ &
= \log P_n(x) \m \log P_{n-1}(x) ,
\end{align*}
because
\[
q^n = \e\left(n\left(\a\m\frac1{c(cX+d)}\right)\right) =
\e\left((r+md)\a\m\frac{m+5/6}c\m\nu\hb\right) = \z_\a^rZ{\rm e}^{-{\rm i}\nu h}
\]
and
\begin{align*}
C_k^{(r)}(\nu) \m C_k^{(r-1)}(\nu-1) &
= - 2 \operatorname{Re}\Biggl[\frac{{\rm i}^{-k}}{(k+1)!}
 \Biggl(\sum_{j=1}^c\m\sum_{j=0}^{c-1} \Biggr)
\L_{1-k}\big(\z_\a^{r+j}Z\big) B_{k+1}\left(\nu+\frac jc\right)\Biggr] \\
 & = - 2 \operatorname{Re}\biggl[{\rm i}^{-k}\L_{1-k}(\z_\a^rZ) \frac{B_{k+1}
(\nu+1)\m B_{k+1}(\nu)}{(k+1)!}\biggr] \\ &
 = - 2 \operatorname{Re}\biggl[  \L_{1-k}(\z_\a^rZ) \frac{(\nu/{\rm i})^k}{k!} \biggr] .
\end{align*}
Note that in the above calculation we used that \smash{$C_{-1}^{(r)}(\nu)$} is
independent of both~$\nu$ and~$r$. In fact, it is given by
\[
C_{-1}^{(r)}(\nu) = 2 \operatorname{Im}\Biggl[\sum_{j=1}^c \L_2\big(\z_\a^{r+j}Z\big)\Biggr]
 = \frac2c \operatorname{Im}[\L_2(Z^c)] = \frac V c ,
\]
where we have used the well-known ``distribution'' property of the
dilogarithm. The corresponding distribution property
of the 1-logarithm $\L_1(z)=-\log(1-z)$ shows that \smash{$C_0^{(r)}(\nu)$} is also
independent of~$\nu$ and is given by
\[
C_0^{(r)}(\nu) = \sum_{j=1}^c \left(\nu+\frac jc-\frac12\right)
\log\bigl|1\m\z_\a^{r+j}Z\bigr|^2 = \log E_r(\a) ,
\]
where $E_r(\a)$ is the real algebraic number defined by
\begin{equation}
\label{tag8}
E_r(\a) = \prod_{j=1}^c \bigl|1\m\z_\a^{r+j}Z\bigr|^{2j/c},
\qquad r\in\Z/c\Z .
\end{equation}
Hence, the exponentiated version of~(6) can be written
\[
P_n(x) \si P_m(\b) \frac{{\rm e}^{V/c\hb}}\hb E_r(\a)
\exp\left(\sum_{k=1}^\infty C_k^{(r)}(\nu) h^k\right) ,
\]
where again the exponential factor at the end has the form
${\rm e}^{-c\sqrt 3\nu^2h/2}\phi_{c,r}\big(\nu\sqrt h,\sqrt h\big)$ with
\[
\phi_{c,r}\big(\v\sqrt h,\sqrt h\big)=
\exp\biggl(
\widetilde C_1^{(r)}(\nu) +
\sum_{k \geq 2} C_k^{(r)}(\nu) h^k
\biggr),
\]
where \smash{$\widetilde C_1^{(r)}(\nu)$} is given by the same formula
as the right-hand side of~\eqref{tag7} (with $k=1$) and with
$B_2(x)=x^2-x+1/6$ replaced by $B_2(x)-x^2$. Note that $\phi_{c,r}(t,\v)$
is a power series in~$\v$ with coefficients in $\Q(\z_\a,Z)[t]$ and
leading coefficient~1.
The same reasoning as for~$c=1$ now shows that the sum of the values of
$P_n(x)$ for $n$ running over the $m$-th peak and in the residue
class~$r+md\!\pmod c$ is equal to \smash{$\hb^{-3/2} {\rm e}^{V/c\hb} P_m(\b)
\Phi_\a^{(r)}(h)$} for some power series \smash{$\Phi_a^{(r)}(h)$} with leading
coefficient $E_r(\a)$, and summing this over all~$r$ gives
equation~\eqref{QMCa} for the $4_1$ knot with
\[
\Phi^{(4_1,\s_1)}_\a(h)=3^{-1/4}c^{-1/2} \sum_{r\!\!\!\!\pmod c}\Phi_\a^{(r)}(h) .
\]
Note that the formal Gaussian integration formula for
the power series $\Phi_\a(h)$ requires to expand the integrand up to
order $O\big(h^{3k+1}\big)$ in order to obtain the coefficient of $h^k$ in the
series $\Phi_\a(h)$.

It remains to look at the units $E_r(\a)$. Write $F$ for $\Q(\e(1/6))$,
the trace field of the figure~8 knot, and $F_c$ for its cyclotomic
extension $F(\z_\a)=\Q(Z)$. We claim that both $E_r(\a)/E_0(\a)$ and~$\prod_{r\!\!\pmod c}E_r(\a)$ belong to~$F_c$. The second claim follows from
the first, since it is clear that~$E_0(\a)^c$ belongs to~$F_c$, and
the first claim follows from the calculation
\[
\frac{E_r(\a)}{E_{r-1}(\a)}
 = \frac{\prod_{j=0}^c \bigl|1-\z_\a^{r+j}Z\bigr|^{2j/c}}{\prod_{j=0}^{c-1}
\bigl|1-\z_\a^{r+j}Z\bigr|^{2(j+1)/c}}
 = \frac{\bigl|1-\z_\a^rZ\bigr|^2}{\prod_{n\!\!\!\!\pmod c} \bigl|1-\z_\a^nZ
\bigr|^{2/c}} = |1\m\z_\a^rZ\bigr|^2 ,
\]
from which we get by induction the formula
\[
E_r(\a) = E_0(\a) \bigl|(\z_\a Z,\z_\a)_r\bigr|^2
\]
for all $r$. In particular, we can write our asymptotic formula to
leading order as
\[
\J\left(\frac{aX+b}{cX+d}\right)\Bigr/ \J(X) \si
\frac{c E_0(\a) S(\a)}{3^{1/4}}
 X^{3/2} \exp\left(\frac V{2\pi} \left(X+\frac dc\right)\right)
\]
as $X\to\infty$ with bounded denominator, where
\begin{equation}
\label{tag10}
S(\a) = \sum_{r\!\!\!\!\pmod c}\frac{E_r(\a)}{E_0(\a)} =
\sum_{n=0}^{c-1} \bigl|(\z_\a Z,\z_\a)_n\bigr|^2 \in F_c .
\end{equation}

It is the factor $S(\a)$ which for $c=5$ contains the funny prime $\pi_{29}$
occurring in \cite[p.\ 14]{Za:QMF}, while~$E_0(a)$ is the unit
analyzed in~\cite{CGZ}. Note that the special properties~\eqref{cunit}
of this unit are clear from the definition~\eqref{tag8}
since if~$c$ is prime to~$6$ and we denote by $\s_k$ the Galois automorphism
of~$F_c$ over~$F$ sending a primitive $c$-th root of unity to its $k$-th
power, then it is easy to see from~(8) that $\s_k(E_r(a/c)^c)=E_r(ka/c)^c$
and that the quotient $E_r(ka/c)^k/E_r(a/c)$ belongs to~$F_c$.

For other knots~$K$, there is a similar story, but we can no longer
rigorously prove anything, since the terms in
the sum defining $\J^K(x)$ are no longer positive (or even real) and
there is cancellation. However, this sum still has the form of an
$N$-dimensional sum of products or quotients of Pochhammer symbols, where
$N$ is the dimension of some terminating $q$-hypergeometric series (related
to the number of simplices in a triangulation of $S^3\ssm K$),
and we can formally look at the parts of this sum where the summands are
locally constant (``stationary phase''), even if those ``parts'' now lie
outside of the original domain of summation. This leads to a conjectural,
but completely explicit, formula of the same general form as~(2), and in
particular to an asymptotic formula like~(9), but with $E_r(\a)=E^K_r(\a)$
now depending on an element $r$ of $(\Z/c\Z)^N$ rather than just $\Z/c\Z$ and
with the sum in~(10) replaced by one over~$(\Z/c\Z)^N$. For the $5_2$ knot
and its sister, the $(-2,3,7)$ pretzel knot, we worked out the
corresponding expressions and for small values of~$c$ obtained both the units
and the pre-factors $S^K(\a)\in F_c$ ($F$ = trace field of~$K$) that we had
previously found numerically. These are, however, much more complicated
than in the~$4_1$ case; for instance, the factor
\smash{$\p_{29}=2-\ve^{(a)}_1+\ve^{(a)}_2+2\ve^{(a)}_3$}, a prime of norm 29
that occurred for the $4_1$ knot and $c=5$ (see equation~\eqref{as41at5} and the
discussion in the next section) is replaced for the $5_2$ knot by $\p_7^2 \p_{43}$
if~$c=3$ and by $\p_{9491}\p_{1227271}$ if~$c=5$, where each $\p_p$ denotes a
prime of norm~$p$ in~$\Q(\xi)$.


\section{Arithmetic aspects}
\label{sec.arithmetic}

In this section, we discuss the arithmetic properties of the power series
\smash{$\Phi_\a^{(K,\s)}(h)$}, in particular the identification of the number fields in which
their coefficients lie and the integrality properties of these coefficients.

\subsection{Algebraic number theory aspects}
\label{sub.field}

A detailed study of the power series \smash{$\Phi_\a^{(K,\s)}(h)$} (or more generally
\smash{$\Phi_\a^{(K,\s,\s')}(h)$}) reveals several interesting algebraic number theoretical
aspects, especially concerning the field of definition, transformation under the
action of the Galois group, and above all the appearance of non-trivial algebraic units.

We begin by looking in more detail at the series \smash{$\Phi_{a/5}^{(4_1)}$} because this
example is quite illuminating. The first few terms of the series were given
in~\cite[p.\ 670]{Za:QMF}, as
\begin{align}
\Phi^{(4_1,\s_1)}_{a/5}(h)
={}&\sqrt[4]3\sqrt[10]{\vve^{(a)}}
\big(\bigl(2-\vve^{(a)}_1+\vve^{(a)}_2+2\vve^{(a)}_3\bigr)\nonumber \\ &
+ \frac{2678-943\vve^{(a)}_1+1831\vve^{(a)}_2+2990\vve^{(a)}_3}{2^3 3^2 5^2
\sqrt{-3}} h +\cdots\bigg) ,\label{at5}
\end{align}
where $\vve^{(a)}_k=2\cos\bigl(\frac{2\pi(6a-5)k}{15}\bigr)$ and
$\vve^{(a)}=\vve^{(a)}_2/\bigl(\vve^{(a)}_1\bigr)^3\vve^{(a)}_3$,
except that the formula was given there in terms of $\log\Phi$, which
introduced spurious denominators in all terms of the expansion. Actually,
this is one of the first insights from the numerical calculations:
earlier papers had always worked with the logarithm, which is what one sees
if one does a Feynman diagram expansion and looks at the contribution of connected
graphs only, but (as in many other combinatorial problems) one gets much
simpler numbers by looking at the exponentiated sum, corresponding to
summing over all graphs rather than just the connected ones. In the case at hand,
this meant that the coefficients of $h$ and~$h^2$ in~\cite{Za:QMF} contained
mysterious powers of the prime \[\pi_{29}^{(a)} = 2-\vve^{(a)}_1+\vve^{(a)}_2+2\vve^{(a)}_3\]
of $\Q\bigl(\cos\bigl(\frac{2\pi}{15}\bigr)\bigr)$, which simply disappear as
soon as one goes from the logarithm of the series to the series itself. But the
few terms of \smash{$\Phi^{(4_1,\s_1)}_{a/5}(h)$} given in~\eqref{at5} also suffice to illustrate
several other key arithmetic points:
\begin{itemize}\itemsep=0pt
\item[(a)]
 The most striking feature of~\eqref{at5} is the appearance of the 10th root of the
 algebraic unit~$\vve^{(a)}$ as a prefactor. From this and the corresponding numbers
 found for other values of~$\a$ and for other knots we were led to conjecture the
 existence of algebraic units in cyclotomic extensions of any number field determined
 by elements of the Bloch group of this field, a prediction that was then confirmed
 in the joint paper~\cite{CGZ} with Frank Calegari.
\item[(b)]
 The case of the $4_1$ knot has the somewhat misleading special property
 that the Kashaev invariant~\eqref{J0for41} is always positive, so that we seem to be
 seeing elements in the real part~${\Q\bigl(\cos\bigl(\frac{2\pi}{15}\bigr)\bigr)}$ of
 the cyclotomic extension $F(\z_5)=\Q(\z_{15})$ of the trace field
 $F=F_{4_1}=\Q\big(\sqrt{-3}\big)$ rather than in this cyclotomic extension itself. In
 particular, as was not observed in~\cite{Za:QMF}, the unit \smash{$\vve^{(a)}$} is, up to
 sign, the \emph{square} of an element of $F(\z_5)$, so that its 10 root is, up to
 a root of unity, in fact a fifth root of a unit in this larger field. Specifically,
 we~have\looseness=-1
 \[\sqrt{-\vve^{(a)}}=\big(\z_{15}^{(a)}- \big(\z_{15}^{(a)}\big)^{-1}\big) \vve^{(a)}_2/\vve^{(a)}_1,\]
 permitting us to rewrite~\eqref{at5} in the form~\eqref{as41at5} given in
 Section~\ref{sec.QMCK}. This, too, turned out to be true for the general case studied
 in~\cite{CGZ}, where one associates to a number field~$F$ and an element of its
 Bloch group the $c$-th root of a unit (or at least $S$-unit for a finite set of
 primes~$S$ independent of~$c$) in $F(\z)$ for every primitive $c$-th root of
 unity~$\z$, and not a $(2c)$-th root. This unit, for $F=F_\s$ and $\z=\e(\a)$, is
 expected to appear as a prefactor in~\smash{$\Phi_\a^{(K,\s)}(h)$} for every~$K$, $\s$, and~$\a$.
\item[(c)]
 Apart from the unit prefactor \smash{$\sqrt{-\vve^{(a)}}$} (which equals \smash{$\ve_{a/5}$} in the
 notation of~\eqref{as41at5}), there is a further prefactor $3^{1/4}$ that coincides
 with the torsion \smash{$\delta(4_1)^{-1/2}=\big(\sqrt{-3}\big)^{-1/2}$} up to a~root of unity and
 an element of $F_{4_1}$.
\item[(d)]
 After we remove these factors, the remaining power series has coefficients in the
 cyclotomic extension $F(\z_5)$ of the trace field.
\item[(e)]
 The denominators of this remaining power series, when we calculate it to many more
 terms using the methods described in Section~\ref{sub.comp2}, have powers of $3$
 (the ramified prime already occurring in (c)) and $D_n$, where
\be
\label{Dn}
D_n = 2^{3n + v_2(n!)} \prod_{\substack{\text{$p$ prime} \\ p>2}}
 p^{ \sum_{i\ge0}[n/p^i(p-2)]} .
\ee
(Note that the exponent $v_p(D_n)$ of $p >2$ in $D_n$ can be written
as $r+v_p(r!) = v_p((p r)!)$ where $r=[n/(p-2)]$.) We will
return to this point in the next subsection.
\item[(f)]
 The unit $\ve_{a/5}$ occurring in~(a) and~(b), the term under the square-root sign
 in~(c), and the coefficients of the ``remaining power series'' as defined in~(d) are
 not only in $F(\z_5)$, but transform under the Galois group
 \smash{$\{\z\mapsto\z^r\}_{5\nmid r}$} of $F(\z_5)/F$ in the ``obvious'' way, i.e., each of
 these numbers is a polynomial in $\z={\rm e}^{2\pi{\rm i}a/5}$ whose coefficients lie in $F_{4_1}$
 and are independent of~$a$.
\item[(g)]
 Finally, the unit $\ve_{a/5}$ of $F(\z_5)$, considered in the quotient
 $F(\z_5)^\times/F(\z_5)^{\times 5}$, transforms under the action of the Galois
 group~$ \text{Gal}(F(\z_5)/F)=(\Z/5\Z)^\times$ in two different ways:
\be
\label{cunit}
\ul\s_r(\ve_{a/5}) = \ve_{ar/5} = (\ve_{a/5})^{1/r}, \qquad r\in(\Z/5\Z)^\times,
\ee
where $\ul\s_r$ is the Galois automorphism defined by $\ul\s_r(\z_5)=\z_5^{ r}$.
\end{itemize}

We conjecture that these properties (a)--(g) hold for all hyperbolic knots~$K$, all
representations~$\s$ in~$\calP_K$ and all roots of unity $\z_\a$, with $F$ replaced
by the trace field~$F_\s$ and 5~by the denominator of~$\a$, as well as a~few other
small modifications (in particular, that instead of a~unit one may get an $S$-unit
for small finite set~$S$ of primes, essentially the ones occurring in the shape
parameters of a triangulation of $S^3\ssm K$, which was empty for $4_1$). In other
words, the power series \smash{$\Phi^{(K,\s)}_{\a}(h)$} can be written in the form
\be
\label{Phiprec2}
\Phi_\a^{(K,\s)}(h) =\mu_{\s,\a}\cdot(\ve_{\s,\a})^{1/c}
\cdot \delta_\s^{-1/2} \sum_{n=0}^\infty \tiA^{(K,\s)}_{\a,n} h^n ,
\qquad \tiA^{(K,\s)}_{\a,n} \in F_\s(\z_\a)
\ee
(so that \smash{$\tiA^{(K,\s)}_{\a,n}$} is the product of an algebraic number independent of~$n$
and the coefficient denoted \smash{$A^{(K,\s)}_\a(n)$} in Section~\ref{sub.CoeffA}),
where $F_\s$ is defined as in Section~\ref{sec.Phi}, $\mu_{\s,\a}$ is an $(8c)$-th root
of unity, and $\ve_{\s,\a} \in F_\s(\z_\a)^\times$ is a near-unit, canonically defined only
up to $c$-th powers, that transforms up to $c$-th powers as in~\eqref{cunit} (with 5
replaced everywhere by~$c$) and that conjecturally depends only on the element of the
Bloch group $B(F_\s)$ determined by $\s$ and in fact coincides with the near-unit that
was constructed in~\cite{CGZ}, and with the same denominator bound~$D_n$ as
in~\eqref{Dn}, independent of~$K$, $\s$ and~$\a$.

\subsection{Denominators and integrality properties}
\label{sub.denominators}

The universal denominator statement given in formula~\eqref{Dn} above was found empirically
on the basis of the extensive numerical data for the $4_1$, $5_2$ and the $(-2,3,7)$ pretzel
knots presented in the appendix to this paper. In this section, we prove it for the denominators
of the power series defined in terms of Gaussian-type integrals in~\cite{DG2}.
This proof only applies to $\s\in\calP\RED_K$, since there is no such integral representation
for the trivial representation, but the corresponding denominator statement is true here also
and can in fact be strengthened because the power series in that case come
from the Habiro ring, as explained at the end of this section.

\begin{Theorem}
\label{thm.denom}
For each knot $K$, representation $\s\in\calP_K$, and number $\a \in \BQ/\BZ$, we have
\[
D_n  \tiA^{(K,\s)}_{\a,n} \in \calO_S\big[\z_\a,c^{-1}\big],
\]
where \smash{$\Phi_\a^{(K,\s)}(h)$} is as in~{\rm \cite{DG2}}, $D_n$ is as in~\eqref{Dn},
\smash{$\tiA^{(K,\s)}_{\a,n}$} is as in~\eqref{Phiprec2}, $c$ is the denominator of~$\a$,~$\calO$ is the ring of integers of $F_\s$ and $S$ is a finite set of primes
of $F_\s$ that depends on $K$ but not on $n$ or on~$\a$.
\end{Theorem}

The first few values of $D_n$ are given by
\begin{align*}
&   1, \ 24, \ 1152, \ 414720, \ 39813120,  \ 6688604160, \ 4815794995200, \  115579079884800, \\
&  22191183337881600, \ 263631258054033408000, \ 88580102706155225088000, \\
&   27636992044320430227456000, \ 39797268543821419527536640000, \ \dots.
\end{align*}
Campbell Wheeler pointed out to us that the above sequence appears (with no proof)
to equal to the sequence~\texttt{A144618} of the online-encyclopedia of integer
sequences~\cite{OEIS}, the latter related to Stirling's formula with half-shift
$D_n=\text{den}(a_n)$, where
\[
z! \sim \sqrt{2\pi} (z+1/2)^{z+1/2} {\rm e}^{-z-1/2}
\sum_{n=0}^\infty \frac{a_n}{(z+1/2)^n},
\qquad z \to \infty .
\]

The numbers $D_n$ grow rapidly, for example
\[
D_{50} = 2^{197} 3^{72} 5^{19} 7^{11} 11^5 13^4 17^3 19^2 23^2 29^1
 31^1 37^1 41^1 43^1 47^1
\]
or $D_{50}/24^{50}50!=5^7 7^3 11^1 13^1 17^1$. \big(In general, $D_n/24^nn!$ is an
integer whose $n$th root tends to~\smash{$\prod_{p\ge5}p^{2/(p-1)(p-2)}=1.8592481285\dots$}.\big)
We give the first 50 values of~$D_n$ by tabulating the ratio $\delta_n=D_n/3D_{n-1}$
(after removing the power of~2, and omitting the values equal to~1):

\begin{center}
\small
\begin{tabular}{|ll|ll|ll|ll|ll|ll|}
\hline
$n$ & $\delta_n$ & $n$ & $\delta_n$ & $n$ & $\delta_n$ & $n$ & $\delta_n$ &
$n$ & $\delta_n$ & $n$ & $\delta_n$
\\ \hline
3 & $3 \tcdot 5$ &
11 & $13$ &
20 & $7$ &
27 & $3^3 \tcdot 5 \tcdot 11 \tcdot 29\!$ &
35 & $7^2 \tcdot 37$ &
42 & $3 \tcdot 5 \tcdot 23$
\\ 
5 & $7$ &
12 & $3 \tcdot 5$ &
21 & $3 \tcdot 5 \tcdot 23\!$ &
29 & $31$ &
36 & $3^2 \tcdot 5 \tcdot 11$ &
44 & $13$
\\ 
6 & $3 \tcdot 5$ &
15 & $3 \tcdot 5^2 \tcdot 7 \tcdot 17\!$ &
22 & $13$ &
30 & $3 \tcdot 5^2 \tcdot 7 \tcdot 17$ &
39 & $3 \tcdot 5 \tcdot 41$ &
45 & $3^2 \tcdot 5^2 \tcdot 7 \tcdot 11 \tcdot 17 \tcdot 47\!$
\\ 
9 & $3^2 \tcdot 5 \tcdot 11\!$ &
17 & $19$ &
24 & $3 \tcdot 5$ &
33 & $3 \tcdot 5 \tcdot 13$ &
40 & $7$ &
48 & $3 \tcdot 5$
\\ 
10 & $7$ &
18 & $3^2 \tcdot 5 \tcdot 11$ &
25 & $7$ &
34 & $19$ &
41 & $43$ &
50 & $7$
\\ \hline
\end{tabular}

\end{center}

\noindent
We also remark that $D_{n_1}D_{n_2}|D_{n_1+n_2}$ for all $n_1,n_2\ge0$ and hence
that the subgroup
\[
R_D[[h]] = \left\{ \sum_{n=0}^\infty \frac{a_n}{D_n} h^n  \,\big|\,
a_n \in R \  \text{for all} \   n \right\}
\]
of $K[[h]]$ is a subring for every subring $R$ of a field~$K$ of characteristic zero.

\begin{proof}
We give the proof only for the case~$\a=0$, $c=1$, using the formulas in~\cite{DG}.
The general case can be proved along the same lines using the more complicated
formulas in~\cite{DG2}, in which the Bernoulli numbers are replaced by Bernoulli
polynomials, but we do not give the details here.
We will also ignore the prime~2 in our proof since it behaves somewhat differently
and in any case can be added to the finite set of excluded primes~$S$ in the statement
of the theorem.

The power series \smash{$\Phi_0^{(K,\s)}(h)$} attached to a triangulation $\calT$
of~$\mathbb S^3\ssm K$ were defined in~\cite{DG} as formal Gaussian integrals
$\langle f_{\calT}\rangle$
of the formal power series
\be
\label{psixz}
f_{\calT}(x;z) = \exp\Bigg(\sum_{j=1}^N \sum_{\substack{r, k\ge0 \\ 2r+k-2 > 0}}
\frac{B_r}{r!} \frac{(-x_j)^k}{k!} \Li_{2-r-k}(z_j)  h^{r+\tfrac k2-1} \Bigg)
\ee
in a multi-variable $x=(x_1,\dots,x_N)$, where $z_1,\dots,z_N$ are the shape parameters
of~$\calT$ and where~${\langle f\rangle=\langle f(x)\rangle_Q}$ is the mean value
defined by Gaussian integration with respect to a certain quadratic form~$Q$ with
coefficients in the field~$F_\s$. This form is essentially the one given by the
symmetric matrix $A+B \,\diag(z_j/(1-z_j))$ that occurred in the discussion
of~\eqref{Nahm3}
in Section~\ref{sub.PhiPerturbative}, and the function~\eqref{psixz} is essentially
the product of the functions occurring in~\eqref{newpsixz}, except that the terms
with $r+k=m$ fixed were combined there into a single Bernoulli polynomial $B_m(x)$
for~$m\ge3$ or $B_2(x)-x^2$ for~$m=2$, and that the normalizations used in~\cite{DG}
were slightly different from the ones used in Section~\ref{sub.PhiPerturbative}.

We now expand the exponential in~\eqref{psixz} as the product of the exponentials
of the monomials in the sum, and recall that $\Li_{2-m}(z)\in\Z[1/(1-z)]$ for
every~$m\ge2$, to deduce that \smash{$\Phi_\a^{(K,\s)}(h)=\langle f_{\calT}\rangle$} is an
$R$-linear combination of Gaussian averages $\langle T\rangle$ of terms~$T$ of the
form
\be
\label{term}
T = \prod_{j=1}^N\prod_{\substack{r, k\ge0 \\ 2r+k\ge3\substack}}
\frac1{\l_j(r,k)!} \biggl(\frac{B_r}{r!}
\frac{x_j^k}{k!} h^{r+\tfrac k2-1}\biggr)^{\l_j(r,k)}
\ee
with non-negative multiplicities $\l_j(r,k)$, where $R=R_{\calT,\s}$ denotes the
ring generated over~$\Z$ by the numbers $(1-z_j)\i$. We write the monomial $T$ as
$c(T) h^n \bx^\bK/\bK!$ with
\be
\label{nAndK}
n = \sum_{j,r,k}\l_j(r,k)\left(r+\frac k2-1\right),\qquad K_j=\sum_{r,k}\l_j(r,k)k
\ee
and where for notational convenience have written $\bx$ and $\bK$ for the
$N$-tuples $(x_1,\dots,x_N)$ and~${(K_1,\dots,K_N)}$ (thus deviating from the
convention in the rest of the paper where boldface denotes matrices) and
$\bx^\bK/\bK!$ for the divided power \smash{$\prod x_j^{K_j}/K_j! $}. To prove the theorem,
we have to bound both the denominators of the numerical coefficient~$c(T)$ and the
further denominators coming from the Gaussian averaging
$ \bx^\bK/\bK! \mapsto\big\langle\bx^\bK/\bK! \big\rangle $.

We begin with the latter question. For this, we recall first that the Gaussian
average $\langle f\rangle_Q$ is given by \smash{${\rm e}^{\D_Q}(f)\bigl|_{x=0}$} for any power
series~$f$, where $\D_Q$ is the Laplacian associated to~$Q$, and hence is equal to
$\D_Q^\ell(f)/\ell!$ if $f$ is a homogeneous polynomial of degree $2\ell$. (For
polynomials of odd degree it of course vanishes trivially.) The Laplacian $\D_Q$
is a quadratic polynomial in the derivatives $\pt_i=\pt/\pt x_i$, and we can enlarge
the ring~$R$ by adjoining to it the coefficients of this polynomial, so since the
image of any divided factorial under any product $\pt_1^{\ell_1}\cdots\pt_N^{\ell_N}$
is an integer, we then certainly have that the Gaussian integral
$\big\langle\bx^\bK/\bK! \big\rangle$ is $1/\ell!$ times an element of~$R$. Unfortunately,
it turns out that this estimate is not good enough for our purposes, and we have
to work a little harder.

Writing $\D_Q$ as an $R$-linear combination of binomials $\pt_i\pt_j$ with
$1\le i\le j\le N$, and applying the multinomial theorem, we see that
$\D_Q^\ell/\ell!$ is an $R$-linear combination of terms
\smash{$\prod_{i\le j}(\pt_i\pt_j)^{m_{ij}}/m_{ij}!$} with $m_{ij}\in\Z_{\ge0}$. Define an even
symmetric $N\times N$ matrix $M$ by setting $M_{ij}=M_{ji}=m_{ij}$ for $i<j$ and
$M_{ii}=2 m_{ii}$. Then \smash{$\prod_{i\le j}(\pt_i\pt_j)^{m_{ij}}=\prod_j\pt_j^{K_j}$}
with $\bK=M\mb1$, where $\mb1$ is the vector consisting of~$N$~1's, and this
sends $\bx^\bK/\bK!$ to~1 and all other monomials to~0.
It follows that a universal denominator of $\big\langle\bx^\bK/\bK! \big\rangle$ is
the number
\begin{align*}
\D(\bK) := \text{l.c.m.}\Biggl\{ \prod_{1\le i\le N}(M_{ii}/2)!\prod_{1\le i<j\le N}M_{ij}!
 \,\Big|\, \text{$M=M^t\in M_{N,N}(\Z_{\ge0})$ even, $M\mb1=\bK$}\Biggr\} .
\end{align*}
Notice that this does divide~$\ell!$, because $\ell$ is the sum of the diagonal
$M_{ii}/2$ and of the $M_{ij}$ with~${i<j}$, so that this statement refines the bound
given above, but $\D(\bK)$ is in general much smaller than~$\ell!$ and this
improvement will be needed for the proof. A usually sharper multiplicative upper
bound for $\D(\bK)$ is the largest integer~$S$ whose square divides the product of
the~$K_j!$ (the proof of this also uses only the integrality of multinomial
coefficients), and then of course~$\D(\bK)$ also divides the g.c.d.\ of $\ell!$ and
$S$, which in general is smaller than either one. (For instance, for $\bK=(6,9,9,10)$
we have $\ell!=355687428096000$, $S=3135283200$, and~${(\ell!,S)=S/3}$.)
Either of these two latter upper bounds would be sufficient for our proof, but in
fact there is an easy upper bound that is stronger than either one of them and is
extremely sharp (in particular, it
is equal to $\D(\bK)$ for all $\bK$ with $N\le4$ and $\max(K_j)\le 30$), namely
\[
\D^*(\bK) := \prod_{\text{$p$ prime}} p^{\delta_p(\bK)} \quad \text{with} \quad
\delta_p(K_1,\dots,K_N) := \sum_{s\ge1}
\left[\frac12\sum_{j=1}^N\left[\frac{K_j}{p^s}\right]\right] .
\]
To show that $\D(\bK)$ divides~$\D^*(\bK)$, we need
$\sum_iV_p(M_{ii}/2)+\sum_{i<j}V_p(M_{ij})\le\delta_p(\bK)$ for every prime~$p$ and
every even symmetric matrix~$M$ with non-negative entries and row sums~$\bK$, where~${V_p(m):=v_p(m!)}$ for any~$m\ge0$ denotes the largest power of~$p$ dividing~$m!$.
In view of the standard formula $V_p(m)=\sum_{s\ge1}[m/p^s]$, it suffices for this
to show that
\[\sum_i[M_{ii}/2q]+\sum_{i<j}[M_{ij}/q]\le\frac12\sum_j[K_j/q]\]
for every prime power~$q$, and this follows immediately from the obvious facts
$[x/2q]=[[x/2]/q]$ and $[x]+[y]+\cdots \le[x+y+\cdots]$ valid for arbitrary real
numbers $x, y,\dots$.

Now going back to our main problem, we now see that it suffices to show that the
product~$\D^*(\bK) c(T)$ has denominator dividing~$D_n$ for all terms~$T$ as above,
with $\bK=(K_1,\dots,K_N)$ and~$n$ defined by~\eqref{nAndK}. We will prove this
one prime at a time (ignoring the prime~2), which is natural in view of the fact
that the upper bound~$\D^*(\bK)$ is defined by its prime power decomposition. To do
this, we will split both the term~$T$ and the corresponding numerical
coefficient~$c(T)$, and also each of the $N$ factors $T_j$ and $c(T_j)$ of which
they are comprised, as the product of four factors in a way depending on the
prime~$p$ being studied, labelled ``$s$'' (terms with $r=0$ and~$k$ smaller than~$p$),
``$b$'' (terms with $r=0$ and $k$ bigger than or equal to~$p$),
``1''~(terms with~${r=1}$), and ``$\ge2$'' (terms with $r\ge2$), with a similar splitting
of the individual weights~$K_j$ into the sum of four pieces
\smash{$K_j^{(s)}=\sum_{3\le k<p}\l_j(0,k) k$}, \smash{$K_j^{(b)}=\sum_{k\ge p}\l_j(0,k) k$}, \smash{$
K_j^{(1)}=\sum_{k\ge1}\l_j(1,k) k$}, and~\smash{$K_j^{(\ge2)}=\sum_{r\ge2, k\ge0}\l_j(r,k) k$}.
We also decompose the number~$n$ (the exponent of~$h$ in~$T$) in~\eqref{nAndK} as
$\ell+n'-t$ with
\[
\ell := \frac12 \sum_{j=1}^NK_j , \qquad
n' :=\sum_{\substack{1\le j\le N \\ r\ge2, k\ge0\substack}}\l_j(r,k) (r-1) , \qquad
t :=\sum_{\substack{1\le j\le N \\k\ge3\substack}}\l_j(0,k)
\]
and also split $t$ as \smash{$t^{(s)}+t^{(b)}$} according as $3\le k\le p-1$ or $k\ge p$ in
the summation, with each of \smash{$t^{(s)}$} and \smash{$t^{(b)}$} splitting into the sum over
$1\le j\le N$ of pieces \smash{$t_j^{(s)}$} and \smash{$t_j^{(b)}$} in the obvious way.

The numerical coefficient $c(T)$ can be decomposed as
\be
\label{cTdecomp}
c(T) = \prod_{1\le j\le N}\frac{K_j!}{P_j(T)} \cdot
\prod_{\substack{1\le j\le N \\ r\ge2, k\ge0\substack}}
\frac1{\l_j(r,k)!} \biggl(\frac{B_r}{r!}\biggr)^{\l_j(r,k)}
\ee
with
\[
P_j(T) = \prod_{0\le r\le1, k\ge1}\l_j(r,k)! k!^{\l_j(r,k)} \cdot
\prod_{r\ge2, k\ge0}k!^{\l_j(r,k)} ,
\]
which we can split up further as \smash{$P_j^{(s)}(T)P_j^{(b)}(T)P_j^{(1)}(T)P_j^{(\ge2)}(T)$}.
The reason that we have included the factor $\l_j(r,k)$ into the definition of
$P_j(T)$ for $0\le r\le1$ but not for~$r\ge2$ is that~$\l!k!^\l$ divides $(k\l)!$
for all $\l\ge0$ if $k$ is strictly positive but not if $k=0$, and the terms with~${r\ge2}$ can have~$k=0$. Then the product \smash{$P_j^{(b)}(T)P_j^{(1)}(T)P_j^{(\ge2)}(T)$}
divides \smash{$\big(K_j-K_j^{(s)}\big)!$}, while the first factor \smash{$P_j^{(s)}(T)$} divides \smash{$t_j^{(s)}!$}
up to a $p$-adic unit because here $k$ is always less than~$p$ and therefore~$k!$
is not divisible by~$p$. (Here we have made repeated use of the integrality of
multinomial coefficients.) On the other hand, by Lemma~\ref{LemmaB} below
and the submultiplicativity of~$D_n$, the second factor in~\eqref{cTdecomp} has
denominator dividing~$D_{n'}$. Putting this all together, we deduce that~$c(T)$
is $G(T)/D_{n'}$ times a $p$-adic integer for every~$p$ (always different from~2
and not dividing the denominators of the elements of~$R$), where
\smash{$G(T)=\prod_{j=1}^N\bigl(K_j!/t_j^{(s)}! \big(K_j-K_j^{(s)}\big)!\bigr)$}. Using the
submultiplicativity of $D_n$ again, this reduces the problem to showing that
${\delta_p(\bK)\le v_p(D_{\ell-t} G(T))}$ for each~$p$, and in view of
the definitions of $\delta_p(\bK)$ and $D_n$ and of the above-mentioned formula~${V_p(m)=\sum_{s\ge1}[m/p^s]}$ for the $p$-adic valuation of factorials, this in turn
will follow if we can show that
\be
\label{Final}
\Biggl[\frac12\sum_{j=1}^N\biggl[\frac{K_j}q\biggr]\Biggr] \le
\biggl[\frac{\ell-t}{q^*}\biggr]
+ \sum_{j=1}^N \Bigg(\biggl[\frac{K_j}q\biggr] \m \Bigg[\frac{t_j^{(s)}}q\Bigg]
\m \Bigg[\frac{K_j-K_j^{(s)}}q\Bigg] \Bigg)
\ee
for each prime power $q=p^s$ with $s\ge1$, where $q^*:=p^{s-1}(p-2)$. 

For this final step, we first note that
\[
\frac{\ell-t}{q^*} = \sum_{j=1}^N
\frac{K_j^{(s)}+K_j^{(b)}+K_j^{(1)}+K_j^{(\ge2)}-2t_j^{(s)}-2t_j^{(b)}}{2q^*}
 \ge \sum_{j=1}^N\frac{K_j-2t_j^{(s)}}{2q}
\]
since
\[
\frac{\big(K_j^{(b)}-2t_j^{(b)}\big)}{q^*}\ge \left(1-\frac{2}{p}\right) \frac{K_j^{(b)}}{q^*}=\frac{K_j^{(b)}}{q}
\]
 \big(because $k\ge p$ in the terms defining \smash{$K_j^{(b)}$} and \smash{$t_j^{(b)}$}\big) and~$q^*<q$.
Using that $\big[\frac x{2q}\big]=\big[\frac12\big[\frac xq\big]\big]$, we deduce that
\[
\biggl[\frac{\ell-t}{q^*}\biggr] \ge
\Biggl[\frac12\sum_{j=1}^N\Bigg[\frac{K_j-2t_j^{(s)}}{q}\Bigg]\Biggr]
\]
and hence (since $x\le y$ certainly implies $[x/2]\le[y/2]$) the
inequality~\eqref{Final} will follow if we have the inequality
\[
\biggl[\frac{K_j}q\biggr] \le \Bigg[\frac{K_j-2t_j^{(s)}}{q}\Bigg]
+ 2 \Bigg(\biggl[\frac{K_j}q\biggr] \m \Bigg[\frac{t_j^{(s)}}q\Bigg]
\m \Bigg[\frac{K_j-K_j^{(s)}}q\Bigg] \Bigg)
\]
for every $1\le j\le N$. But this inequality is trivial, since
\smash{$K_j-K_j^{(s)}\le K_j-3t_j^{(s)}\le K_j-2t_j^{(s)}$} (because every $k$ in the
definition of \smash{$K_j^{(s)}$} is $\ge3$) and
\[
\biggl[\frac{K_j}q\biggr]\ge\Bigg[\frac{K_j-2t_j^{(s)}}q\Bigg]
+2\Bigg[\frac{t_j^{(s)}}q\Bigg].
\]
 This completes the proof of Theorem~\ref{thm.denom} modulo that of the following lemma.
\end{proof}

\begin{Lemma}\label{LemmaB}
For any integers $r\ge2$ and $\l\ge0$, we have
\be
\label{BernEstimate}
\frac1{\l!} \left(\frac{B_r}{r!}\right)^\l \in \frac1{D_{\l(r-1)}} \Z .
\ee
\end{Lemma}

\begin{proof}
We prove this one prime at time. By well-known results of
von Staudt and Clausen, the Bernoulli number $B_r$ ($r>0$) has $p$-adic valuation
$-1$ if~${(p-1)|r}$ and $B_r/r$ is $p$-integral if~${p-1}$ does not divide~$r$. From this
we deduce that the $p$-adic valuation of the denominator of $B_r/r!$ is bounded
above by $\bigl[\frac r{p-1}\bigr]$, which is $\le\bigl[\frac {r-1}{p-2}\bigr]$
since $\frac{r-1}{p-2}\ge\frac r{p-1}$ if $r\ge p-1$ and both expressions
vanish otherwise. The $p$-adic valuation of the denominator of
\smash{$\frac1{\l!}\bigl(\frac{B_r}{r!}\bigr)^\l$} is therefore bounded above by
\smash{$V_p(\l)+\l\bigl[\frac{r-1}{p-2}\bigr]$}. On the other hand, from the definition
of~$D_n$ we have $v_p(D_n)=m+V_p(m)=V_p(pm)$, where $m=\bigl[\frac n{p-2}\bigr] $.
We must therefore show that
\[
V_p(\l)+\l \biggl[\frac{r-1}{p-2}\biggr] \le
V_p\biggl(p \biggr[\frac{\l(r-1)}{p-2}\biggr]\biggl)
\]
for all $\l\ge0$ and $r\ge2$. For this, we make a case distinction: if $ 2\le r<p-1$, then
\[
\text{l.h.s.} = V_p(\l) = V_p\biggl(p \biggl[\frac\l p\biggr]\biggr) \le
V_p\biggl(p \biggl[\frac\l{p-2}\biggr]\biggr) \le \text{r.h.s.} ,
\]
while if $r\ge p-1$, then we set $h=\bigl[\frac{r-1}{p-2}\bigr]\ge1$ and have
instead
\[
\text{l.h.s.} = V_p(\l)+\l h \le \l h+V_p(\l h) = V_p(p\l h) \le \text{r.h.s.}
\]
because $[\l x]\ge \l [x]$ for any positive real number~$x$.
\end{proof}

We end this subsection with several further observations concerning the denominators
and integrality properties of the coefficients of our divergent power series. The
first is that the bound in Theorem~\ref{thm.denom} is not only sharp in the strong
sense that it is best possible for \emph{every} integer~$n\ge0$ and not merely that
it is attained for some~$n$, but that this optimality is reached in two very
different extreme ways: in the above lemma if $r=p-1$ and $\l\ge0$ is arbitrary
(in which case both sides of~\eqref{BernEstimate} have the same value~$V_p(p\l)$)
and again in the calculation~\eqref{Final} in the case when only \smash{$K_j^{(b)}$} occurs
and all $k_i$ are~equal to $p$, so that \smash{$K_j^{(b)}$} is exactly \smash{$p t_j^{(b)}$} (in
other words whenever the dominating contribution in~\eqref{term} comes from the
terms with $(r,k)=(p-1,0)$ or~$(0,p)$). The fact that two completely different
mechanisms lead to the same function $n\mapsto D_n$ suggests that this function
may be a more fundamental one than appears at first sight and may have a broader
domain of applicability.

The second observation is that the universal denominator statement given by
Theorem~\ref{thm.denom} can be sharpened by considering the series at the logarithmic
level, or equivalently, by studying the denominators of the contributions from
connected rather than from all Feynman diagram. This was motivated by the observation
of Peter Scholze that the logarithm of the series \smash{$\Phi_0^{(4_1)}(h)$} in~\eqref{as41},
which we had calculated up to order $ \O\big(h^{150}\big)$, had coefficients with smaller
denominators than those of the series itself. Specifically, he found experimentally
that the first occurrence of~$p^k$ for small primes $p$~($\ne2, 3$) and $k\ge1$ in
\smash{$\log\Phi_0^{(4_1)}(h)$} occurred for the coefficient of $h^n$ with~${n=k(p-1)-1}$
rather than $n=k(p-2)$ as for the power series \smash{$\Phi_0^{(4_1)}(h)$} itself.
At first sight this statement seems to contradict the intuition mentioned at the
beginning of the section that the arithmetic of the series~$\Phi_\a(h)$ is much
simpler if one does not take their logarithms. But in fact both statements are true!
The reason is that in general $\Phi_\a(h)$ is a linear combination of finitely
many power series corresponding to the stationary points of the function being
integrated (specifically, there are~$c^M$ such series, where $c$ is the denominator
of~$\a$ and $M$ can be taken to be the number of tetrahedra in a triangulation of
the knot complement), and it is not reasonable to take logarithms of sums. But
for~$\a=0$ there is only one summand, so here it is reasonable to take the logarithm,
and for general~$\a$ the logarithm of each of the finitely many summands of~$\Phi_a(h)$ coming from the contribution to the state integral near an individual
stationary point is indeed simpler than this summand itself, because it corresponds
to a sum over only connected rather than over all Feynman diagrams. The final
statement is given in the following theorem.

\begin{Theorem}
\label{thm.denomc}
For each integer $n\ge1$, define
\[
D_n\conn = \prod_{\text{$p$ {\rm prime}}} p^{[(n+1)/(p-1)]}  .
\]
Then the coefficient of $h^n$ in \smash{$\log\Phi_0^{(K,\s)}(h)$} for any knot $K$ and
any $\s\in\calP_K$ has denominator dividing $D_n\conn$ apart from a finite set
of primes depending only on $K$ and~$\s$. More generally, for any $\a\in\Q$ we
have \smash{$\Phi_\a^{(K,\s)}(h)\in R\otimes\exp\bigl(\sum_{n\ge1} R h^n/D_n\conn\bigr)
\subset R_D[[h]]$}.
\end{Theorem}

The first few values of $D_n\conn$ are given by
\begin{align*}\label{Dnconnvalues}
&   2, \ 12, \ 24, \ 720, \ 1440, \ 60480, \ 120960, \ 3628800, \ 7257600, \ 479001600, \
958003200, \\ &   2615348736000, \ 5230697472000, \ \dots.
\end{align*}
This sequence too appears in the online-encyclopedia of integer
sequences~\cite{OEIS} under the name \texttt{A091137} and with the formula given
above, and coincides with the denominator of the Todd polynomials given in
Hirzebruch's book~\cite[Lemmas~1.5.2 and~1.7.3]{Hirzebruch}
without proof and quoted from the paper~\cite{Atiyah-Hirzebruch}.

The proof of the above theorem (which actually implies Theorem~\ref{thm.denom})
is similar to the proof of Theorem~\ref{thm.denom} and is omitted.

Note that the ``connected denominators'' $D_n\conn$ are considerably smaller than
the ``additive denominators''~$D_n $: $ D_n/n!$ is an integer growing exponentially
like $(44.621{\dots}+o(1))^n$, as already mentioned, while $D_n\conn/(n+2)!$
is an integer of subexponential growth.

The final remarks concern the relation of the above results with the known integrality
properties of elements of the Habiro ring. The proof of Theorem~\ref{thm.denom} as
given above only works for the power series \smash{$\Phi^{(\s)}_\a$} with $\s\ne\s_0$, because
the perturbative expansion does not apply to the case~$\s=\s_0$. However, as we
know, this remaining case is actually simpler because it belongs to the Habiro ring
and therefore satisfies \smash{$\Phi_0^{(\s_0)}(h)\in\Z\big[\big[{\rm e}^h-1\big]\big]$}, and more generally~\smash{$\Phi_\a^{(\s_0)}(h)\in\Z(\e(\a))\big[\big[\e(\a){\rm e}^{-h}-1\big]\big]$}. The corresponding property no
longer holds for $\s\ne\s_0$, even for $\a=0$ and the figure~8 knot, and in some sense
should not even be expected, because the ``natural'' invariant here is the completed
power series \smash{$\Phih_0^{(4_1)}$}, which contains a transcendental factor ${\rm e}^{\V(4_1)/h}$.
However, if we consider the product
$\sqrt{3}\Phi(h)\Phi(-h)=-\sqrt{-3}\Phi^{(1)}(h)\Phi^{(2)}(h)$, which would be
unchanged if we replaced the power series by their completions, then we \emph{do}
find experimentally that it belongs to the ring~$\Z[1/3]\big[\big[{\rm e}^{-h}-1\big]\big]$, i.e., after the
change of variables from $q={\rm e}^{-h}$ to $q=1+x$ it becomes a 3-integral power series
in~$x$. We expect, and have checked numerically, that the same is also true for
$5_2$ if one multiplies all three series \smash{$\Phi_0^{(5_2,\s_i)}$} ($i=1,2,3$), and for
$(-2,3,7)$ for both products of three series corresponding to the two number fields
$\Q(\xi)$ and $\Q(\eta)$ corresponding to this knot. These properties are explained
by the properties of Habiro rings for general number fields as being developed
in~\cite{GSZ:habiro}.

Finally, we found experimentally that we can obtain a power series that is already
integral (away from 2 and~3) in~${\rm e}^{-h}-1$ from \smash{$\Phi(h)=\Phi^{(4_1,\s_1)}(h)$}
\emph{without} multiplying it by $\Phi(-h)=-{\rm i}\Phi^{(4_1,\s_2)}(h)$ if we multiply
instead by
\smash{$\mathcal E^{(4_1,\s_1)}(h)
:=\exp\bigl(-\sum_{r\ge1}^\infty \frac{|B_{r+1}|}{(r+1)!}C_rh^r\bigr)$} with $C_r$ as
in Section~\ref{sub.proof41.0}. (Notice that this implies the previous statement
becomes $\mathcal E(h)\mathcal E(-h)=1$.) We expect that there will be a similar
correction factor $\mathcal E^{(K,\s)}(h)$ for any \smash{$\Phi_\a^{(K,\s)}(h)$} and that
the corrected $\Phi$-series can be seen as the $h$-deformed versions of the units
constructed in~\cite{CGZ}, and hope to study this too in~\cite{GSZ:habiro}.


\section{Numerical aspects}
\label{sec.compKash}

In this section, we describe how the power series whose existence is predicted
by the modularity conjecture can be computed effectively via a
numerical computation of the Kashaev invariant, extrapolation, and
recognition of algebraic numbers in a known number field. We also
describe other methods that are applicable to the power series \smash{$\Phi_\a^{(K,\s)}(h)$}
for $\s$ different from~$\s_1$, as well as the smooth truncation
methods of ``evaluating'' factorially divergent power series at non-zero arguments
that were discussed in Section~\ref{sub.4.3}.
The actual numerical data for several knots will then be presented in the appendix.

\subsection[Computing the power series Phi\_a\^{}\{(K,s)\}(h)]{Computing the power series $\boldsymbol{\Phi_\a^{(K,\s)}(h)}$}
\label{sub.comp2}
In this subsection, we explain the various methods that can be used to compute the
coefficients of the power series \smash{$\Phi_\a^{(K,\s)}(h)$} for all $\a\in\Q$ and
$\s\in\calP_K$ numerically and then exactly as algebraic numbers.

{\bf Step 1: Computing the colored Jones polynomials.}
To compute the Kashaev invariant~$\la K \ra_N$ of a knot~$K$, we use
the Murakami--Murakami formula $\la K \ra_N=J_N\big({\rm e}^{2 \pi {\rm i}/N}\big)$, where~${J_n(q)\!=\!J_{K,n}(q)}$
is the $n$-th colored Jones polynomial, together with a theorem of T.T.Q.~L\^e and the
first author~\cite{GL} that asserts the existence of a recursion relation for
the polynomials~$J_n(q)$. This reduces the problem to that of
computing/guessing this recursion relation concretely for a given knot.
This in turn has been solved for several knots in joint work of the first author,
X.~Sun and C.~Koutschan~\cite{GK,GK2,GS}. The solution required a modulo $p$
computation of the $N$-th colored Jones polynomial (for several primes $p$
and several thousand values of $N$), together with a careful guess of the
supporting coefficients of such a recursion. In particular, the recursion
was computed in~\cite{GS} for the twist knots $K_p$ with $|p|\le15$, was
guessed in~\cite{GK} for the $(-2,3,3+2p)$ pretzel knots with $|p|\le5$,
and was computed (when $p=2$) or guessed (when~${p=3,4,5}$) in~\cite{GK2}
for the double twist knots $K_{p,p}$ with $2\le p\le5$.

{\bf Step 2: Computing the Kashaev invariant for large~$\boldsymbol{N}$.}
In order to get numerical information about the asymptotics
of the Kashaev invariant $\la K \ra_N$, we need to be able to compute it
numerically to high precision for large values of $N$, say of the order
of~$N=5000$. Although both the Kashaev invariant and the colored Jones
polynomial are given by finite-dimensional terminating $q$-hypergeometric sums,
and the latter have been programmed in \texttt{Mathematica}~\cite{B-N},
these programs can only give the value of the colored Jones polynomial
and of the Kashaev invariant for modest values of $N$, say, up to~$N=20$,
which is far less than we need for the numerical extrapolation. By using the
recursion, we can compute $J_n(\z_N)$ numerically to high precision for
$N$ large and $0<n<N$. (This is far faster than computing the colored Jones
polynomials $J_n(q)$ for these values of~$n$ and substituting $q=\z_N$ at
the end.) However, this does not give the Kashaev invariant
$\la K \ra_N=J_{K,N}(\z_N)$ because the recursion gives $P(q,q^n)J_n(q)$
as a linear combination of a bounded number of previous values $J_{n-i}(q)$,
where $P(q,x)$ is a~fixed polynomial that is always divisible by~$1-x$.
To overcome this, we use the recursion relation and its first $r$ derivatives,
where $(1-x)^r\|P(q,x)$, to get a recursion for the length-$r$ vector~\smash{$\big(J_n(q),J_n'(q),\dots,J_n^{(r)}(q)\big)$}. We can the use the recursion to compute
the whole vector numerically for~$q=\z_N$ and~$0<n<N$, and the single value
$J_n(q)$ for~$n=N$. It follows that

\begin{Proposition}
 \label{prop.ckashaev}
The Kashaev invariant $\la K \ra_N$ of a knot $K$ can be computed
in $O(N)$ steps.
\end{Proposition}

This linear-time algorithm, which seems to be new even for the Kashaev invariant,
can be used equally well to compute $J(\g N)$ for~$\g=\sma abcd\in\SL_2(\Z)$ fixed and
$N$ tending to infinity, or even $J(\g X)$ with~$\g=\sma abcd$ fixed and $X$~tending to
infinity with fixed fractional part, since this simply the value of $J_n(\g X)$ with $n$
equal to the denominator of $\g X$ and can be calculated by the same trick. Note
that the computation takes time $O(N)$ numerically and $O\big(N^3\big)$ if we work over
$\BQ[\z_N]$.

{\bf Step 3: Computing $\boldsymbol{\Phi_\a^{(K,\s_1)}(h)}$.}
Once we know how to compute $J(\g N)$ for large integers~$N$ (or even $J(\g X)$
for large~$X$ with bounded denominator), we can we can obtain the first few
coefficients of the power series \smash{$\Phi^{(K,\s_1)}_{\a}(h)$} numerically for a fixed
rational number~$\a=\g(\infty)$ by combining the quantum modularity conjecture~\eqref{QMCa} (or~\eqref{QMCb})
together with the extrapolation method of the second author (as described in detail
in~\cite{GrMor} or the appendix of~\cite{GV}) or the closely related
Richardson transform~\cite[Chapter~8]{Bender}. This is quite effective and gives, for
instance, the first hundred coefficients of the series~\eqref{as41} or forty
coefficients of the series~\eqref{as52} in only a~few minutes of computing time.
We should mention, however, that this extrapolation method requires either exact
numbers or else very high precision (often several hundred or even thousand decimal
digits in the calculations we did.)

{\bf Step 4: Recognizing the coefficients exactly.}
Given that the coefficients of \smash{$\Phi^{(K,\s_1)}_{\a}(h)$}
are conjecturally algebraic
numbers in a specific number field, we can then test numerically by using
the standard LLL (Lenstra--Lenstra--Lovasz) algorithm to approximate the
numerically computed coefficients by rational linear combinations of
a basis of this field. If this works to high precision with coefficients
that are not too large and have only small primes in the denominator, then
we have considerable confidence that the approximate equality is an exact one.
The method is self-verifying in the sense that the success at each
stage requires the correctness of the answer guessed at the previous stage.

{\bf Step 5: Computing $\boldsymbol{\Phi_\a^{(K,\s_0)}(h)}$.}
In this step, we explain how to compute the expansion of the element
of the Habiro ring at a root of unity, in linear time. More precisely,
we have:

\begin{Proposition} \label{prop.hab}
The series \smash{$\Phi^{(K,\s_0)}_\a(h) + O(h)^{N+1}$} is computable in $O(N)$ steps.
\end{Proposition}
This follows from the fact that the expansion of the Habiro element at
$q=\z_\a {\rm e}^h$ up to~$O(h)^{N+1}$ requires $cN$ terms of the
cyclotomic polynomial of $K$, which is a linear combination of the
first~$cN$ colored Jones polynomials of $K$. An alternative formula,
inspired by Mahler's ideas on $p$-adic interpolation, gives the
following expansion of the Habiro element:
\be
\label{hab-mahler}
\Phi^{(K,\s_0)}_\a(h) = \sum_{k=1}^{cn} \widehat{J}^K_n\big(\z_\a {\rm e}^h\big) + O\big(h^{n+1}\big),
\ee
where $c$ is the denominator of $\a$, $\z_\a=\e(\a)$ and
\[
\widehat{J}^K_n(q) = \sum_{k=1}^n (-1)^{n-k} \binom{n}{k}_q
q^{\frac{k(k-1)}{2}} J^K_{n-k+1}(q)  q^{-\frac{n(n+1)}{2}}
\]
and \smash{$\binom{n}{k}_q = (q;q)_n/((q;q)_k (q;q)_{n-k})$} is the $q$-binomial
symbol and $(q;q)_n=(1-q)\big(1-q^2\big)\cdots\allowbreak (1-q^n)$ is the $q$-Pochhammer
symbol. The right-hand side of~\eqref{hab-mahler} gives
a well-defined formal power series since $\widehat{J}^K_{cn}\big(\z_\a {\rm e}^h\big)$ lies
in $h^n \BQ[[h]]$.

{\bf Step 6: Computing the remaining power series.}
Once we have the leading term in the original QMC, we can obtain the numerical terms
by successively subtracting the corrections coming from the values of~$\s$ different
from~$\s_1$ as explained in Sections~\ref{sub.4.1} (using ``optimal truncation'')
and~\ref{sub.4.3} (using the more precise ``smooth truncation'' described
in Section~\ref{sub.optimal} and in more detail in~\cite{GZ:optimal}), where numerical
examples were also given.

{\bf Step 7: Alternative methods.}
In Steps 2 and~3, we discussed how to obtain the coefficients of \smash{$\Phi_\a^{(K,\s_1)}(h)$}
from the original QMC together with the high-speed high-precision computation of
Jones polynomials at roots of unity and numerical extrapolation; in Steps~4 and~5
we explained how to get the series \smash{$\Phi_\a^{(K,\s)}(h)$} for $\s$ Galois conjugate
to~$\s_1$ and for $\s=\s_0$, respectively; and in Step~6 we described how to obtain
the remaining series by using the \emph{refined} quantum modularity conjecture
together with optimal truncation. However, there are at least two other ways to get
these other series that are of interest and are sometimes more efficient.

The first way is to use the formal Gaussian integration of Section~\ref{sec.perturbation}
and~\cite{DG, DG2}. This method uses exact arithmetic and allows the computation of
few terms \smash{$A^{(K,\s)}_\a(k)$} (in practice $k \leq 5$) as exact algebraic numbers for
many knots (such as those with ideal triangulations with up to 15 ideal tetrahedra).
See also~\cite{Ga:exact}.

The other, which is completely different, is based on the asymptotics near roots of
unity of~the holomorphic functions in the upper half plane (generalized Nahm sums)
that we study in~\cite{GZ:qseries}. Since Nahm sums converge
quadratically, the values of those~functions at~$\tau=\a+{\rm i}/N$ can be computed in
$ O\big(\sqrt{N}\big)$ steps and after extrapolation give a~nu\-merical~com\-putation of
the algebraic numbers \smash{$A^{(K,\s)}_\a(k)$}. This method, when applicable, is not only
much faster (time $ O\big(\sqrt N\big)$ rather than $ O(N)$), but also has
the major advantage of allowing the simultaneous numerical computation of
\smash{$A^{(K,\s)}_\a(k)$} for all $\a$ of a fixed denominator and for all~$\s$ in a Galois
orbit which (after multiplication by $D_k$ and by a suitable $S$-unit)
reduces the problem of recognizing the list of coefficients \smash{$A^{(K,\s)}_\a(k)$} as
algebraic numbers to the problem of recognizing numerically computed \emph{integers},
albeit of growth $k!^2 C^k$. This allowed us to compute, for instance,
100~coefficients of the series \smash{$\Phi^{(K,\s)}_\a(h)$} for the $5_2$ knot for $\a=0$
and $a=1/2$ and for all three representations $\s$ in the Galois orbit of the
geometric representation,
and it allowed us to compute 37 terms of the series of the $(-2,3,7)$
pretzel knot for $\a=0$ and for both Galois orbits (each of size 3) of
nontrivial representations $\s$. It is remarkable that this method allows
the computation of series for representations not Galois conjugate to the
geometric one, though not for the trivial representation~$\s_0$.

{\bf Orientation conventions.} Finally, we have to discuss an annoying technical
point, namely, the choice of a consistent set of conventions for the two
classical invariants (namely the trace field and the complex volume of
a knot) and the two quantum invariants (namely the colored Jones polynomial
and the Kashaev invariant of a knot). These conventions are especially
important since the tables of knots rarely distinguish a knot from its
mirror, and (for instance) the name~$5_2$ of the unique hyperbolic
5-crossing knot does not convey this distinction.\looseness=1

On the other hand, replacing a knot $K$ by its mirror $K^*$ reverses the
orientation of the knot complement $M_K = S^3 \ssm K$, which implies that
\begin{itemize}\itemsep=0pt
\item
$F_{K^*} = \overline{F_K}$ and $\tV(K^*) =\overline \tV(K)$.
\item
$J_{K^*,N}(q)=J_K\big(q^{-1}\big)$ and $\la K^* \ra_N=\overline{\la K \ra_N}$ and
$\J^{K^*}(x)=\J^K(-x)=\overline{\J^K(x)}$.
\end{itemize}
Thus, a random orientation convention for $K$ might not match the asymptotics
whose coefficients are in a fixed subfield of the complex numbers, and
not on its complex conjugate subfield.

The Jones (hence, also the colored Jones) polynomial
$J_{K}(q) \in \BZ\big[q^{\pm 1}\big]$ of a knot (or a link)~$K$ is uniquely determined
by the following skein-relation~\cite{Jones}
\be
\label{jones}
q  J_{\left(\includegraphics[height=0.1in]{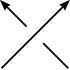} \right)}(q)
-q^{-1} J_{\left(\includegraphics[height=0.1in]{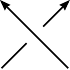} \right)}(q) =
\big(q^{1/2}-q^{-1/2}\big)
J_{\left(\includegraphics[height=0.1in]{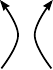} \right)}(q),
\qquad J_{\text{unknot}}(q)=1 .
\ee
In particular, for the left-hand trefoil $3_1$, we have:
$J_{3_1}(q)= -q^{4} + q^{3} + q^{1}$. Moreover, the colored Jones polynomial
$J_{K,N}(q) \in \BZ\big[q^{\pm 1}\big]$ is normalized to be $1$ at the unknot, and
to equal to the Jones polynomial when $N=2$.

Note that the \texttt{SnapPy} program~\cite{snappy} for computing the Jones
polynomial agrees with~\eqref{jones}, whereas the \texttt{Mathematica}
program \texttt{KnotAtlas}~\cite{B-N} polynomial differs by replacing $q$
by $1/q$:

\begin{center}
\begin{tabular}{|l|l|} \hline
\verb|L = Link(braid_closure=[-1,-1,-1])| &
\verb|Jones[BR[2, {-1, -1, -1}]][q]|\\
\verb|L.jones_polynomial()| &
\verb|-q^-4 + q^-3 + q^-1|
\\
\verb|-q^4 + q^3 + q^1| & \\ \hline
\end{tabular}
\end{center}

We will be using the consistent orientation convention for $M_K$ of the
\texttt{SnapPy} program (when~$K$ is given as the closure of a braid, or
via a planar projection, or via an augmented DTcode or Gauss code),
which has the added advantage that it also gives shapes of tetrahedra
corresponding to the hyperbolic structure (exactly or numerically),
as well as the trace field (exactly) and the complex volume $\tV(K)$
(exactly or numerically).

\subsection{Optimal truncation and smoothed optimal truncation}
\label{sub.optimal}

One of the main numerical aspects concerns smoothed optimal truncation, which
was originally an appendix to an earlier draft of this paper but has now been
relegated to a planned independent publication \cite{GZ:optimal} because the methods
are applicable to many problems outside the realm of quantum topology. This is a
method for the numerical summation of factorially divergent series when only a
finite number of coefficients are known and we do not have information about the
possible analytic continuation of the Borel transform, which is the method usually
used.

We already explained in Section~\ref{sub.4.1} the naive optimal truncation,
$\Phi(h)^\opt$ of a factorially divergent series $\Phi(h)=\sum_{n=0}^\infty A_nh^n$,
defined simply as $\sum_{n=0}^NA_nh^n$ where $N$ is the approximate
value of~$n$ at which the term $A_nh^n$ takes on its minimum absolute value,
given explicitly by~${N=|B/h|}$ if $A_n$ grows like $n! B^{-n}$. The idea of
smoothed optimal truncation is very simply to replace $\Phi(h)^\opt$ by a
``smoothed'' version $\Phi(h)^\smo$ which is defined as $\Phi(h)^\opt+\ve_N(h)$
where the exponentially small correction term $\ve_N(h)$ depends on the cutoff
parameter~$N$ in such a way that $\Phi(h)^\smo$ does not jump when one changes $N$
by~1. This means simply that we require~${\ve_{N-1}(h)-\ve_N(h)=A_nh^n}$. Of course, if
we knew how to solve this equation exactly, then the function $\Phi(h)^\smo$
would be completely independent of~$N$ and would give us a canonical way to
lift the power series $\Phi(h)$ to an actual function. This is not the case,
but if $A_n$ has a~known asymptotic expansion, which is true for all of the
series in this paper (see Section~\ref{sub.CoeffA}) then we can define
$\ve_N(h)$ asymptotically as the product of ${\rm e}^{-N}$ and a power series in~$1/n$
chosen in such a way that the desired equality $\ve_{N-1}(h)-\ve_N(h)=A_nh^n$
is true asymptotically to all orders in~$1/n$. The details of how to do this\
if $A_n$ has the asymptotic form $B^{-n} \sum_{\ell=0}^\infty c_\ell \G(n+\k-\ell)$
for some real number~$\k$ and some numerical coefficients~$c_\ell$ (as in
equations~\eqref{AnFirst}, \eqref{BnFirst} or~\eqref{CoeffAsymp}) are given
in~\cite{GZ:optimal} and will not be repeated here. The only thing that is of
importance to us here is that the result of the smoothing gives an evaluation of
$\Phi(h)$ that is independent of all choices (and hence gives a predicted ``right''
definition of the corresponding function) up to an error that is exponentially small
with a better exponent than that given by naive optimal truncation. Specifically,
the new error is ${\rm e}^{-N} (1+|C|)^{-N}$ rather than simply ${\rm e}^{-N}$ as before if, as
is always the case for us, the coefficients $c_\ell$ themselves grow factorially
like $\ell! C^{-\ell}$. Examples of this dramatic numerical improvement were given
in Section~\ref{sub.4.3}, where in one case the error in evaluating a~series
$\Phi(h)$ was of the order of~$10^{-29}$ using naive optimal truncation but of the
order of~$10^{-56}$ using smoothed optimal truncation.

\appendix
\section{Numerical data for five sample knots}\label{sec.num}

In this appendix, we present numerical data that support the quantum modularity
conjecture for a~choice of knots. Initially, we hoped that pairs
of geometrically similar knots (that have identical trace fields and
equal elements in the Bloch group, modulo torsion---henceforth called
``sisters'') might have identical or nearly identical
series~\smash{$\Phi^{(K,\s_1)}_\a$}. With this in mind, and having already performed the
computations for the $4_1$ knot, we were led to consider its sister,
the \texttt{m003} census manifold. The latter is not a knot complement
(it is the complement of a~knot in a~lens space), but its 5-fold cyclic
cover is the complement of the (twisted) $5$-chain link, with a computable
Kashaev invariant, to be compared with the $5$-th power of the Kashaev
invariant of the $4_1$ knot. No relation between the series \smash{$\Phi^{(K,\s_1)}_\a$}
was found for this pair, but the units $\ve(K)_\a$ did match (up to roots
of unity). We then tried the $5_2$ knot, whose sister is the
$(-2,3,7)$ pretzel knot. Here again the series \smash{$\Phi^{(K,\s_1)}_\a$} were
different, but the units $\ve(K)_\a$ matched. The final
example of the $6_1$ knot was chosen because its Bloch group
has rank~$2$ and its $\SL_2(\BC)$-character variety is more complicated,
making the verification of the QMC, the Galois invariance
property~\eqref{cunit} and the match with the unit of~\cite{CGZ}
more subtle, but again all three were verified numerically. For this knot
we did not make any ``sister'' computations.

Recall the coefficients \smash{$A^{(K,\s)}_\a(k)$} of the power series
\smash{$\Phi^{(K,\s)}_\a(h)$} are algebraic numbers. In this appendix,
we present the numerically obtained data for \smash{$A^{(K,\s)}_\a(k)$}
written in the form
\be
\label{PhiA}
A^{(K,\s)}_\a(k) = C^{(K,\s)}_\a \tiA^{(K,\s)}_{\a}(k) ,
\ee
where \smash{$C^{(K,\s)}_\a=\mu_{\s,\a}\cdot(\ve_{\s,\a})^{\frac{1}{c}} \cdot
\delta_\s^{-\frac{1}{2}}$} is given in~\eqref{Phiprec2}.
Note however that \smash{$C^{(K,\s)}_\a$} and \smash{$\tiA^{(K,\s)}_{\a}(k)$} are not
canonically defined numbers, only their product is. (Since we are focusing
on the geometric representation $\s_1$, and we are fixing the knot $K$, we
omit the superscript $(K,\s)$ from the notation in the right-hand side of
the above equation.) We will further specify a choice of an algebraic number
$\lambda_c$ such that \smash{$\lambda_c^k  D_k  \tiA_{\a}(k) \in \calO_{F(\z_\a)}$}
is an algebraic integer, where $c$ is the denominator of~$\a$, $F$ is the
trace field of the knot and $D_n$ is
the universal denominator~\eqref{Dn}. Using a basis of the free abelian group
$\calO_{F(\z_\a)}$, we can represent the above algebraic integers by lists of
integer numbers.

\subsection{The figure eight knot}
\label{sub.41}

In this appendix, we discuss the numerical aspects of the quantum modularity
conjecture for the simplest hyperbolic $4_1$ knot, for which we currently
know how to prove the modularity conjecture.
(The proof was presented in Section~\ref{sec.QMC41}.)
Needless to say, the numerically obtained results agree with the exact
computation of the expansion coefficients given in Section~\ref{sec.QMC41}.
Some information about the numerical aspects of the Kashaev invariant
of the $4_1$ knot were already presented in the introduction, but we give some
additional data (e.g., for the expansion near seventh roots of unity),
since this is the most accessible knot numerically and also to illustrate
the formulas occurring in the proof.

Since the knot is fixed in this section, and so is the
geometric representation $\s_1$, we will suppress them
from the notation. As mentioned in Section~\ref{sec.QMCK}, the
trace field of the $4_1$ knot is~$\BQ\big(\sqrt{-3}\big)$ and the torsion
is $\delta = \sqrt{-3}$. Some terms of the series $\Phi_\a(h)$
when $c=1$ or $c=5$ were given in
equations~\eqref{as41} and~\eqref{at5}, respectively.
The numerical method allows us to compute and identify
the power series $\Phi_{a/c}(h)$ to any desired precision. Although
the coefficients of the series $\Phi_{a/c}(h)$, divided by the
constant term, is an element of the number field $\BQ(\z_{3c})$ (when
$c$ is coprime to $3$) or $\BQ(\z_c)$ when $3$ divides $p$, a judicious
choice of the constant $C_{a/c}$, combined with the Galois invariance
of the coefficients allows us to list the coefficients \smash{$\tiA_{a/p}(k)$} for
$p \neq 3$ prime and for $a=1,\dots,p-1$ by giving a $p-1$ tuple of elements
in the trace field $\BQ(\z_3)=\BQ\big(\sqrt{-3}\big)$ of the knot. Furthermore,
since the knot is amphicheiral, it follows that \smash{$\tiA_{a/p}(k)$} is real or purely
imaginary (for $k$ odd or even, respectively), and combined with the above
discussion, allows to list the vector of coefficients
\smash{$\big(\tiA_{1/p}(k),\dots,\tiA_{(p-1)/p}(k)\big)$} by a $(p-1)$-dimensional vector of rational
numbers. Our numerical extrapolation method allows us to compute this
tuple efficiently, and what is more, our code is self-correcting in
several ways: if a wrong denominator for $\tiA_{\a}(k)$ is guessed
for some $k$, its factorization in primes involves prime larger than $k+1$,
and the precision of the computation drops in the next step by
a factor of two. As a result, we were able to compute 100 terms of the
series $\Phi_0(h)$ when $c=1$, and the results agree with the
computations given in~\cite{DGLZ} as well as computations obtained by
a different method by the first author.

In addition to this, we computed the constant term $\Phi_\a(0)$ for all
$\a$ with denominator a~prime less than 100, and confirmed that its norm
agrees with the predictions of \cite[Section 4.1]{DG2} for~${c \leq 19}$.
We also computed 20 terms of the series $\Phi_\a(h)$ for all $\a$ with
denominator a prime less than 100.

To present a sample of our computations, we start with the special case of
$c=1,2,3,6$, where the $c$-th root of unity is in the trace field
$C_\a$ and $\lambda_\a$
\begin{center}
 \def\arraystretch{1.2}
 \begin{tabular}{|c|c|c|cc|cc|}
\hline
$\a$
& $0$ & $1/2$ & $1/3$ & $2/3$ & $1/6$ & $5/6$ \\
\hline
$C_\a$ & $3^{-{1/4}}$ & $3^{1/4}$ & $2 \cdot 3^{-1/12}$ &
$3^{7/12}$
& $2^2 \cdot 3^{1/12}$ & $3^{17/12}$ \\
\hline
$\lambda_c$ & $72 \sqrt{-3}$ & $72 \sqrt{-3}$ &
\multicolumn{2}{c|}{$24 \sqrt{-3}$}
& \multicolumn{2}{c|}{$36 \sqrt{-3}$} \\
\hline
\end{tabular}
\end{center}
it turns out that $\lambda_c^k D_k \tiA_{\a}(k)$ are integers being given by

\begin{center}
 \def\arraystretch{1.4}\tiny
 \begin{tabular}{|c|c|c|cc|cc|}
\hline
$k$
& $\lambda_1^k D_k \tiA_{0}(k)$ & $\lambda_2^k D_k \tiA_{1/2}(k)$ &
$\lambda_3^k D_k \tiA_{1/3}(k)$ & $\lambda_3^k D_k \tiA_{2/3}(k)$ &
$\lambda_6^k D_k \tiA_{1/6}(k)$ & $\lambda_6^k D_k \tiA_{5/6}(k)$ \\
\hline
$0$ & $1$ & $1$ & $1$ & $1$ & $1$ & $1$ \\
\hline
$1$ & $11$ & $41$ & $37$ & $25$ & $579$ & $201$ \\
\hline
$2$ & $697$ & $12625$ & $7785$ & $6449$ & $1224117$ & $782865$\\
\hline
$3$ & $724351$ & $48022429$ & $21535937$ & $18981677$ & $39903107571$ & $29648832381$ \\
\hline
$4$ & $278392949$ & $72296210981$ & $24220768661$ & $21569737445$ & $535664049856461$ & $412895509718949$
\\
\hline
$5$ & $244284791741$ & $252636824949503$ & $63245072194611$ & $56749680285647$ & $16693882665527364525$
& $13164162601119392223$
\\
\hline
\end{tabular}
\end{center}

When $c=4$ and $a= \pm 1 \bmod 4$, with the choice
$C_{a/4}=\pm \big(3 \big(2 \pm \sqrt{3}\big)\big)^{-1/4}$ and $\lambda_4=6\sqrt{-3}$,
we can write
\[
\lambda_4^k  D_k  \tiA_{\pm 1/4}(k) = \tiB_{1/4}(k) \pm \tiB_{-1/4}(k){\rm i},
\]
where $B_4(k)=\big(\tiB_{1/4}(k),\tiB_{-1/4}(k)\big) \in \BZ^2$ with the first
six values are given by

\begin{center}\renewcommand{\arraystretch}{1.2}
\begin{tabular}{|c|l|}
\hline
$B_4(0)$ & $\la 1, 2 \ra$ \\
\hline
$B_4(1)$ & $\la 365, 370 \ra$ \\
\hline
$B_4(2)$ & $\la 311785, 420386 \ra$ \\
\hline
$B_4(3)$ & $\la 4219048201, 6325027802 \ra$ \\
\hline
$B_4(3)$ & $\la 24805519728725, 38098972914250 \ra$ \\
\hline
$B_4(4)$ & $\la 340419470401244075, 531593492940700894 \ra$ \\
\hline
$B_4(5)$ & $\la 25036998069742932352139, 39557220304220645794918 \ra$ \\
\hline
\end{tabular}
\end{center}

Finally, when $c=p$ is a prime different from $3$, we found
out for that for the primes less than~100, the
constant $C_\a$ of~\eqref{PhiA} can be taken to be
\[
C_{\a} = 3^{(-2 \pm 1)/4}  p^{1/2} (\ve_\a)^{1/p}
\qquad \text{for} \quad p=\pm 1 \bmod 6,
\]
where
\[
\ve_\a = \prod_{|k| \leq \frac{p-1}{2}} (\ve(p'k\a))^k, \qquad
p'=\mp 1/4 \bmod p,
\qquad \ve(x)=2 \cos 2 \pi(x-1/3) .
\]
Note that the unit $\ve_\a$ in $\BQ(\z_{3p})$
that appears in the choice of $C_\a$ agrees, up to $p$-th powers of units,
with the theoretically computed unit from equation~\eqref{tag8} (for $r=0$)
below. With the above choice of $C_\a$, the numbers $A_{\a}(k)$ lie in the
field $\BQ(\z_{3p})$, satisfy the Galois invariance described in detail
in the introduction, and this allows them to be expressed in terms of vectors~${B_p(k)=\big\la \tiB_{1/p}(k),\dots,\tiB_{(p-1)/p}(k) \big\ra \in \BZ^{p-1}}$
as follows:
\[
\lambda_p^k  D_k  \tiA_{a/p}(k) = \sum_{b=1}^{p-1} \eta(ab/p) \tiB_{b/p}(k),
\qquad \eta(x) = 2 \sin(2 \pi (x-1/3)) ,
\]
where $\lambda_p = 3p^2 \sqrt{-3}/2$. The vectors $B_p(k)$ for $k \leq 20$ and
$p$ a prime less than $100$ were numerically obtained and for $p=5$ and
$p=7$ are given by

\begin{center}\renewcommand{\arraystretch}{1.2}\tiny
\begin{tabular}{|c|l|}
\hline
$B_5(0)$ & $\la -1, -4, -4, -6 \ra$ \\
\hline
$B_5(1)$ & $\la -55, -5140, -7660, -9690 \ra$ \\
\hline
$B_5(2)$ & $\la -7586065, -48629140, -58401700, -81382470 \ra$ \\
\hline
$B_5(3)$ & $\la -1066837647875, -5818148628500, -6620399493500, -9407838821250 \ra$ \\
\hline
$B_5(4)$ & $\la -51952598327049125, -274293246490488500, -309180073069692500, -440171876888046750 \ra$ \\
\hline
$B_5(5)$ & $\la -5814113396376116334625, -29960825153926862627500, -33500926926525556664500, -47835527737950677253750 \ra$
\\
\hline
\end{tabular}
\end{center}

\noindent
and

\begin{center}\renewcommand{\arraystretch}{1.2}\tiny
\begin{tabular}{|c|l|}
\hline
$B_7(0)$ & $\la -20, 7, 2, 5, -14, -8 \ra$ \\
\hline
$B_7(1)$ & $\la -98140, 8267, -19670, 27937, -39214, -16576 \ra$ \\
\hline
$B_7(2)$ & $\la -2199415652, 426208447, -172006030, 524259533, -1237405358, -619260152 \ra$ \\
\hline
$B_7(3)$ & $\la -676432728043100, 166452454682479, -15638648253886, 168799271208365, -406506539584838, -215671594628336
\ra$ \\
\hline
$B_7(4)$ & $\la -86350611733284233860, 22591735955847949331, -702673247614974230, 21808440520527403561, $ \\
& $ -52829131820839184902, -28340444966866544008 \ra$ \\
\hline
$B_7(5)$ & $\la -25671367091358132079572196, 6928168872402051353797277, 10873595841062215161670, $ \\
& $ 6492789075493742592974935, -15896921084389159954206466, -8579075179324647599719264 \ra$ \\
\hline
\end{tabular}
\end{center}

\subsection{The sister of the figure eight knot}
\label{sub.s41}

Its quotient by $\BZ/5\BZ$, which
is a knot in the lens space $L(5,1)$ rather than the 3-sphere,
is the sister of the $4_1$ knot, with the same trace field and same
Bloch group invariant. We therefore
expect to find similarities between the asymptotic power series
associated to $\K1$ and to~$\K2$.

Next, we discuss the case of the sister of the $4_1$ knot, the manifold
\texttt{m003} in the hyperbolic knot census~\cite{snappy}, which is not
the complement of a knot in the 3-sphere, but is the complement of a
nullhomologous knot in the Lens space $L(5,1)$. This complicates things
since the sister knot has no Jones (hence, no colored Jones) polynomial,
and although it has a Kashaev invariant, a~formula for it is not available
to us. However, the $5$-fold cyclic cover of the sister of the
$4_1$ knot is the 5-chain link $L$ in $S^3$ \big(denoted by $10^5_3$ and also
by \texttt{L10n113}\big). This is a famous link because virtually every census
manifold is a Dehn filling on it~\cite{Dunfield-Thurston}. The link $L$
has a colored Jones polynomial
$J_{L,N}(q)$ (with all components colored by the same $N$-dimensional
representation) with a~formula available from~\cite{VV}
and a Kashaev invariant. More precisely, we have
\[
J_{L,N}(q) = -\frac{1}{1-q^N}
\sum_{n=0}^{N-1} \big(q^{n+1}-q^{-n}\big)c(N,n)(q)^2 c(N,n)\big(q^{-1}\big)^3,
\]
where
\[
c(N,n)(q)=\frac{q^{-Nn}}{(q;q)_n}\sum_{k=0}^{N-n-1} q^{-Nk}
\prod_{j=k+1}^{n+k} \big(1-q^{N-j}\big)\big(1-q^j\big) .
\]
\noindent
The above formula is $O\big(N^2\big)$ can be rewritten in terms of an $O(N)$
formula that has a recursion relation. However the latter has the
disadvantage that the middle term of the
summand~(${k=N/2}$) now vanishes when evaluated at $\e(1/N)$. To overcome
this, we compute the sum from both sides by differentiation. Having done so,
we tested the QMC and no suprises were found. We computed 10 terms of the
series $\Phi^{L}_\a(h)$ when $\a=0$ (given below)
and 8 terms when $\a=1/2$.

We now give the data for $\a=0$. The trace field of $L$ is $\BQ\big(\sqrt{-3}\big)$,
same as for the $4_1$ knot. The complex volume of $L$ is given by
\[
\V(L)=5  \V(4_1) -3\pi^2 
\]
and its torsion is given by
\[
\delta(L) = 2^7 \sqrt{-3} .
\]
Since $L$ is a link, in~\eqref{K1}, we should replace the exponent
$3/2$ by $5/2$. With these changes, and with the notation of
equation~\eqref{PhiA} we get algebraic integers
$12^k D_k A^{L}_{0}(k)$ in the ring~$\BZ\big[\sqrt{-3}\big]$ and the first
10 are given as follows:
 $$
 \def\arraystretch{1.2}
\begin{array}{|c|l|} \hline
k & 12^{k}  D_k  \tiA^{K_2}_{0}(k) \\ \hline
0 & 1
\\ \hline
1 & -115 \sqrt{-3} + 279
\\ \hline
2 & -49050 \sqrt{-3} + 53286
\\ \hline
3 & -112270440 \sqrt{-3} + 163969920
\\ \hline
4 & -131463532440 \sqrt{-3} + 2948624280
\\ \hline
5 & 4388324675760 \sqrt{-3} - 163377997734672
\\ \hline
6 & -155232475000358400 \sqrt{-3} + 1614884631367642560
\\ \hline
7 & -456051590815208713920 \sqrt{-3} - 409415976078904226880
\\ \hline
8 & 1201424680509251029718400 \sqrt{-3} - 2426468490157451971144320
\\ \hline
9 & 280843674420360230423881689600 \sqrt{-3} + 767958533539384912591107225600
\\ \hline
\end{array}
$$

However, we failed to find any relation between the series for the
$4_1$ knot and for the 5-fold cover of its sister.

\subsection[The 5\_2 knot]{The $\boldsymbol{5_2}$ knot}
\label{sub.52}

The pair of the $4_1$ knot and its sister from the previous section in
unsatisfactory in two ways. For one, the quantum modularity conjecture is
proven for the~$4_1$ knot.
Moreover, the sister of~$4_1$ (and its 5-fold cover) is not a knot.
The next simplest hyperbolic knot after~$4_1$ is the $5_2$ knot,
whose sister is the mirror of the $(-2,3,7)$ pretzel knot. Sister (or
geometrically similar) knots have a decomposition into a finite number of
congruent ideal tetrahedra, hence they have the same trace field
and equal elements in the Bloch group, modulo torsion.

A formula for the Kashaev invariant of the $5_2$ knot
was given in~\cite[equation~(2.3)]{K97},
\be
\label{K52}
\J^{5_2}(x)=\sum_{m=0}^{N-1} \sum_{k=0}^m
q^{-(m+1)k} \frac{(q;q)_m^2}{\big(q^{-1};q^{-1}\big)_k}, \qquad q=\e(x),
\ee
where $N$ is the denominator of $x \in \BQ$. After multiplication of
the above by $\e(x)$, it agrees with the evaluation of the colored
Jones polynomial $J_{5_2,N}(\e(x))$, where the Jones polynomial of $5_2$
is~${J_{5_2}(q)=q - q^2 + 2 q^3 - q^4 + q^5 - q^6}$.
The formula~\eqref{K52} allows a computation of the Kashaev invariant
in~$O\big(N^2\big)$ steps, and a simplification of it was found by one of the
authors~\cite[Section~4.1]{DGLZ}
\[
\J^{5_2}(x)=\sum_{m=0}^{N-1} (q;q)_m^2 
\Biggl(
\big(q^{-1};q^{-1}\big)_m \sum_{k=0}^m \frac{q^{-k^2}}{\big(q^{-1};q^{-1}\big)_k^2} \Biggr),
\qquad q=\e(x)
\]
that allows an $O(N)$-step computation of the Kashaev invariant.
An alternative computation of the latter in $O(N)$-steps can be performed
using the recursion relation for the colored Jones polynomial of
$5_2$~\cite{GS}.

As mentioned in Section~\ref{sub.p}, the trace field of $5_2$ is $F=\BQ(\xi)$, where
\be
\label{a52}
\xi^3-\xi^2+1=0, \qquad \xi = 0.877438833\ldots -
0.74486176661\ldots {\rm i} .
\ee

The trace field has three embeddings labeled by $\s_j$ for $j=1,2,3$
(as discussed in Section~\ref{sub.p}) and their volumes are given
by
\begin{align*}
 \V(\s_1) &= -3 R(\xi_1)+\frac{2 \pi^2}{3}
 = 3.0241283 \ldots + 2.8281220 \ldots  {\rm i},
 \\
 \V(\s_2) &= -3 R(\xi_2)+\frac{2 \pi^2}{3}
 = 3.0241283 \ldots -2.8281220 \ldots  {\rm i},
 \\
 \V(\s_3) &= 3 R(\xi_3/(1-\xi_3))-\frac{\pi^2}{3} = -1.1134545 \ldots,
\end{align*}
where $R(x)$ denotes the Rogers dilogarithm defined by
\[
R(x) = \text{Li}_2(x) + \frac12 \log(x) \log(1-x) - \frac{\pi^2}6
\quad \text{for} \quad
x\in\C\ssm\bigl( (-\infty,-1] \cup [1,\infty) \bigr) .
\]
The torsion of the $5_2$ knot is given by $\delta(5_2)=3\xi-2$.

$\bullet$ Modularity at 0:
We choose $\ve(5_2)_0=1$, and with the notation of equation~\eqref{PhiA}
the first eleven terms are given as follows:
$$
 \def\arraystretch{1.2}
\begin{array}{|c|l|} \hline
k & \big(2^3 \xi^5 (3\xi-2)^3\big)^k  D_k  \tiA^{5_2}_{0}(k) \\ \hline
0 & 1
\\ \hline
1 & -12 \xi^2 + 19 \xi - 86
\\ \hline
2 & -1343 \xi^2 - 12052 \xi + 14620
\\ \hline
3 & 1381097 \xi^2 + 36300408 \xi - 10373787
\\ \hline
4 & -939821147 \xi^2 - 7647561573 \xi - 5587870829
\\ \hline
5 & 114451233224986 \xi^2 - 51239666382079 \xi - 6305751988731
\\ \hline
6 & -2263527400987641127 \xi^2 - 631762147829071739 \xi - 1298875409805289208
\\ \hline
7 & -757944502306007361580 \xi^2 + 1425054483652604079482 \xi + 2654782623273180246011
\\ \hline
8 & 16785033822956024557916646 \xi^2 - 2226340168480665471705515 \xi
\\ &
- 14930684354870794067096358
\\ \hline
9 & -3735848035153601836654158090473 \xi^2 - 3510831690088210470322102227368 \xi
\\ &
- 449224959824265576892987954854
\\ \hline
10 & -34345984964128841574873487072878291 \xi^2 + 25085231887789675521906921078089414 \xi
\\ &
+ 52364404634270110370401111089362065
\\ \hline
\end{array}
$$

$\bullet$ Modularity at 1/2:
We choose $\ve(5_2)_{1/2} =\xi^{-5}$ and with the notation of
equation~\eqref{PhiA} the first six
terms are given as follows:
$$
\def\arraystretch{1.2}
\begin{array}{|c|l|} \hline
k & \big(2  \xi^5 (3\xi-2)^3\big)^k  D_k  \tiA^{5_2}_{1/2}(k) \\ \hline
0 & \xi + 2
\\ \hline
1 & 307 \xi^2 - 138 \xi - 628
\\ \hline
2 & -573109 \xi^2 - 168712 \xi + 457975
\\ \hline
3 & 2096955561 \xi^2 + 5077310601 \xi + 1165885531
\\ \hline
4 & 6470888990010 \xi^2 - 5414463743327 \xi - 10380246225743
\\ \hline
5 & 289484322041800655 \xi^2 - 138373378538474483 \xi - 156775910252412286
\\ \hline
\end{array}
$$

$\bullet$ Modularity at 1/3:
Here, the constant term $\Phi^{5_2}_{1/3}(0)$ was numerically computed to
high precision
\[
\Phi^{5_2}_{1/3}(0)=
-1.3478490468923913068\ldots - 1.5706460265356353326\dots {\rm i}
\]
but was not initially recognized. To identify it,
we used the formula~\eqref{K52} for the Kashaev invariant and performed
a theoretical computation analogous to the constants $S(\a)$ and $E_0(\a)$
(given in~\eqref{tag8} and~\eqref{tag10}) of the $4_1$ knot which produced
the primes
\[
\p_7 = \big(\xi^2-1\big) \z_6 - \xi + 1, \qquad \p_{43} = 2\xi^2-\xi - \z_6
\]
of norm 7 and 43 respectively in the number field $F_3=\BQ(\xi,\z_3)$.
Note that the same primes appear in~\cite[Section~6.2]{DG2}.
In addition, the above constant involves $\delta(5_2)^{-1/2}$ and a number
whose third power is in $F_3$. After some experimentation, we concluded that
\[
\Phi^{5_2}_{1/3}(0)= \e(1/36)
\frac{1}{\sqrt{3\xi-2}}  \p_7^2  \p_{43} .
\]
It follows that a representative of the unit at $\a=1/3$ is given by
\[
\ve(5_2)_{1/3} =\e(1/12) .
\]
It was a bit of a surprise to find that the unit is torsion although the
Bloch group of $F_{6_1}$ has rank~1. On the other hand
$3$ (as well as $2$ and a few other primes) are exceptional ones in the
work~\cite{CGZ}.

Once the constant term was recognized, it turned out that we needed to
separate one factor of $\p_7$ in the constant term \smash{$\Phi^{5_2}_{1/3}(0)$}
from the remaining terms, in order to avoid spurious denominators.
With the choice of \smash{$C_{1/3}=\e(1/36) (3\xi-2)^{-\frac{1}{2}} \p_7  \p_{43}$}
and the notation of~\eqref{PhiA}, the first seven terms were found to be
as follows:
{\small
$$
\def\arraystretch{1.2}
\begin{array}{|@{\,\,}c@{\,\,}|l|l|} \hline
k & \big(\xi^5 (3\xi-2)^3\big)^k  D_k  \tiA^{5_2}_{1/3}(k) \\ \hline
0 & \bigl(-\xi^2 + 2 \xi - 2\bigr) \z_6 + \big(2 \xi^2 - 4 \xi\big)
\\ \hline
1 & \big(717 \xi^2 - 822 \xi + 947\big) \z_6 + \bigl(-2226 \xi^2 + 1856 \xi + 106\bigr)
\\ \hline
2 & \bigl(-680145 \xi^2 + 1283633 \xi - 1844797\bigr) \z_6 + \big(4731470 \xi^2 - 1215426 \xi + 785050\big)
\\ \hline
3 & \bigl(-4879664798 \xi^2 - 15547118437 \xi + 26771206405\bigr) \z_6
\\ &
+ \bigl(-20691193336 \xi^2 - 35194065214 \xi - 73160959238\bigr)
\\ \hline
4 & \big(237593851209955 \xi^2 - 123624865686699 \xi + 65688152000880\big) \z_6
\\ &
+ \bigl(-455730563794746 \xi^2 + 258640669065738 \xi + 244974132213716\bigr)
\\ \hline
5 & \bigl(-8559119253981428654 \xi^2 + 9164193255880569642 \xi - 8506396294603249043\bigr) \z_6
\\ &
+ \bigl(-8914434881967188748 \xi^2 - 7549553228397039176 \xi + 21232362162256499338\bigr)
\\ \hline
6 & \big(1206971041591026374138836 \xi^2 - 1471979903142920023426465 \xi + 1526039068996370402375484\big) \z_6\!
\\ &
+ \big(2034143372251380409655636 \xi^2 + 5390411863643322238842526 \xi - 935392258601663466664696\big)
\\ \hline
\end{array}
$$}%
where $\z_6=\e(1/6)$.

\subsection[The (-2,3,7) pretzel knot]{The $\boldsymbol{(-2,3,7)}$ pretzel knot}
\label{sub.237}

Next, we discuss the case of a sister the $5_2$ knot, namely the
(mirror of) the $(-2,3,7)$ pretzel knot.
Note that the trace fields of $5_2$ and $(-2,3,7)$ coincide,
which allow us to use the notation of~\eqref{a52}.

Unlike the case of the $4_1$ and $5_2$ knots, the Kashaev invariant of
$(-2,3,7)$ can only be computed via the recursion of the colored Jones
polynomial which was guessed in~\cite{GK}, with the convention that the
Jones polynomial of $(-2,3,7)$ is given by
$J^{(-2,3,7)}(q)=q^{-5} + q^{-7} - q^{-11} + q^{-12} - q^{-13}$. The above
inhomogeneous recursion has order 6, maximal degree $(6, 58, 233)$
with respect to the shift variable, the $q^n$ and the $q$ variables,
and contains a total of
90 terms, which can be found in~\cite{Ga:pretzeldata}.
In contrast, the $A$-polynomial of the $(-2,3,7)$ knot has maximal degree~$(6,55)$ with respect to the $(L,M)$ variables and contains 12 terms.
In addition, we multiply the Kashaev invariant of $(-2,3,7)$ by $q^{-4}$.

Since $(-2,3,7)$ is a sister of the $5_2$ knot, they have a common
trace field $\BQ(\xi)$ given in~\eqref{a52}.
The trace field has three embeddings labeled by $\s_j$ for $j=1,2,3$
(as discussed in Section~\ref{sub.p}) and their complex volumes are given
by
\begin{align*}
 \V(\s_1) &= -3 R(\xi_1) + \frac{\pi^2}{3} = 4.6690624 \ldots + 2.8281220 \ldots {\rm i},
 \\
 \V(\s_2) &= -3 R(\xi_2) + \frac{\pi^2}{3} = 4.6690624 \ldots - 2.8281220 \ldots {\rm i},
 \\
 \V(\s_3) &= 3 R(\xi_3/(\xi_3-1)) + \frac{\pi^2}{3} = 0.5314795 \ldots.
\end{align*}
The torsion of $(-2,3,7)$ are given by
$\delta((-2,3,7)) = -2 (3\xi-2) \xi^{-2}$.

$\bullet$ Modularity at 0: Using the notation of~\eqref{PhiA}, we write
\[
\big(\big(2 \xi^2 - 6\big)^3/\xi^5\big)^k  D_k  \tiA^{(-2,3,7),\s_1}_{0}(k) = \big(1,\xi,\xi^2\big)\cdot B^{(-2,3,7),\s_1}_{0}(k),
\]
where $B^{(-2,3,7),\s_1}_{0}(k) \in \BZ^3$ is a vector of integers with the first
11 values given by

{\tiny
$$
\begin{array}{|c|l|} \hline
k & B^{(-2,3,7),\s_1}_{0}(k) \\ \hline
0 & \la 1, 0, 0\ra
\\ \hline
1 & \la -33, 128, -90\ra
\\ \hline
2 & \la 79245, -104172, 50944\ra
\\ \hline
3 & \la 333329999, -597644460, 317584318\ra
\\ \hline
4 & \la -12580573862099, 16160668928488, -9152599685200\ra
\\ \hline
5 & \la 275061075796915969, -366241217321535656, 209464837107544698\ra
\\ \hline
6 & \la -21464059785100413194817, 28432876033981872108244, -16179201892533998639888\ra
\\ \hline
7 & \la 39552725057509518276438631, -52341801268123421363828580, 29838036942620515077356206\ra
\\ \hline
8 & \la 249767901145868199725688538645, -330081248453503483229302323376, 187971265625750854805584690976\ra
\\ \hline
9 & \la -3700925786017810109833640742259950499, 4903075033684898536256604949931358320,
\\ &
-2794204143666309730641613915747239310\ra
\\ \hline
10 & \la 392518725914904741935043787434245408953117, -519977480066306945985500543478969169892188,
\\ &
296298336548750157536627179710807871873120\ra
\\ \hline
\end{array}
$$
}

$\bullet$ Modularity at 1/2:
If we choose $\ve_{1/2}((-2,3,7))=2 \xi^5 $, with the notation
of~\eqref{PhiA}, the first four terms are given by
$$
\def\arraystretch{1.2}
\begin{array}{|c|l|} \hline
k & \big(4  \xi (3\xi-2)^3\big)^k  D_k  \tiA^{((-2,3,7),\s_1)}_{1/2}(k) \\ \hline
0 & 1
\\ \hline
1 & -225 \xi^2 + 404 \xi - 249
\\ \hline
2 & 87535 \xi^2 - 158073 \xi + 123948
\\ \hline
3 & 1981731163 \xi^2 - 3465695160 \xi + 2508787814
\\ \hline
\end{array}
$$

$\bullet$ Modularity at 1/3: Here the constant term and the next two
coefficients of the power series \smash{$\Phi^{(-2,3,7)}_{1/3}(h)$},
\smash{$\Phi^{(-2,3,7)}_{2/3}(h)$} were computed to high precision, and using as a guidance
the appearance of primes of norm 373 (conjectured in~\cite[Section~6.2]{DG}),
we identified the constant terms
\[
\Phi^{(-2,3,7)}_{1/3}(0) = \e(2/9) \sqrt{-\frac{27}{2 (3\xi-2)}} \p_{373},
\qquad
\Phi^{(-2,3,7)}_{2/3}(0) = \e(5/9) \sqrt{-\frac{27}{2 (3\xi-2)}} \p'_{373},
\]
where $\p_{373}=\xi^2+2 \xi \z_6 +1$ and $\p'_{373}=\xi^2+2 \xi (1-\z_6) +1$
are primes in $\BQ(\xi,\z_6)$ of norm~$373$.
It follows that the unit at $\a=1/3$ is given by
\[
\ve((-2,3,7))_{1/3} = \e(2/3) .
\]
The units of $5_2$ and $(-2,3,7)$ at $\a=1/3$ match up to a $24$-th root of
unity.

As mentioned in Section~\ref{sub.p}, the $(-2,3,7)$ pretzel knot
has 6 parabolic nonabelian representations that come in two Galois orbits
of size 3 each: one is defined over the trace field (the cubic field of
discriminant $-23$ discussed above), and another defined over the real
field $\BQ(\eta)$, the abelian field of discriminant 49. At first glance,
the latter three parabolic representations (which are $\SL_2(\BR)$
representations of zero volume) are not seen by the Kashaev invariant. Yet,
one can detect them using the asymptotics of the coefficients of the former
three representations as explained in Section~\ref{sub.optimal}.

In the subsequent paper~\cite{GZ:qseries}, we used the 6 pairs of $q$-series
associated the $(-2,3,7)$ pretzel knot and their asymptotics to compute 37
terms of all six series \smash{$\Phi^{((-2,3,7),\s_j)}_{0}(h)$} for $j=1,\dots,6$.
Below, we give the first 11 terms of the series associated to the abelian
number field $\BQ(\eta)$ given in Section~\ref{sub.p}.
Consider the embeddings $\s_{3+j}$ of the above field for $j=1,2,3$
given in Section~\ref{sub.p} which send $\eta$ to $2\cos(2\pi j/7)$
and let
\smash{$C^{((-2,3,7),\s_{3+j})}_0 = \sqrt{(\eta_j-2)/14}$}. The complex volumes
of~$\s_{3+j}$ are given by
\[
\V(\s_4)= -\frac{1}{21} \pi^2, \qquad \V(\s_5)= \frac{1}{14} \pi^2, \qquad
\V(\s_6)= -\frac{1}{42} \pi^2
\]
and the torsion equals to $\delta((-2,3,7), \s_{3+j})=14/(\eta_j-2)$.
Using the notation of~\eqref{PhiA}, we write
\[
7^k  D_k  \tiA^{((-2,3,7),\s_{3+j})}_{0}(k)
= \big(1,\eta_j,\eta_j^2\big)\cdot B^{((-2,3,7),\s_{3+j})}_{0}(k),
\]
where \smash{$B^{((-2,3,7),\s_{3+j})}_{0}(k) \in \BZ^3$} is a vector of integers with
the first 11 values given by

{\tiny
$$\renewcommand{\arraystretch}{1.2}
\begin{array}{|c|l|} \hline
k & B^{((-2,3,7),\s_{3+j})}_{0}(k) \\ \hline
0 & \la 1, 0, 0\ra
\\ \hline
1 & \la 43, 0, -21\ra
\\ \hline
2 & \la 3928, 63, -1491\ra
\\ \hline
3 & \la -9658210, -2570400, 8759835\ra
\\ \hline
4 & \la -12802855375, 9661452255, 660110430\ra
\\ \hline
5 & \la -42833879089694, 5736063757095, 23026249581258\ra
\\ \hline
6 & \la -360522404258392495, -58094689166990595, 278695629206010765\ra
\\ \hline
7 & \la 108480519886094978165, 114336214602228319050, -161431920455740612440\ra
\\ \hline
8 & \la 420957357301236147078125, -601694281205047856100870, 211820529501946639071105\ra
\\ \hline
9 & \la 276051903390093831791757795950, -105329146895536652560323534375, -93062298372659896456977171525\ra
\\ \hline
10 & \la 3837169849511929903158156720021580, 1712034755788650551262940860512280,
\\ &
-3840647130863172583813306383456135
\ra
\\ \hline
\end{array}
$$
}

\subsection[The 6\_1 knot]{The $\boldsymbol{6_1}$ knot}
\label{sub.61}

In this appendix, we look at one further knot (this time without a sister), the
$6_1$ knot, for two reasons. Firstly, the trace field
is $F_{6_1} = \BQ(\xi)$, a number field of discriminant 257 (a prime) where
\[
\xi^4 + \xi^2 - \xi + 1=0, \qquad \xi=0.5474\ldots + 0.5856\dots {\rm i}.
\]
The trace field has two complex embeddings, so its Bloch group has rank
two, giving a nontrivial test for the unit $\ve(6_1)_\a$. Secondly,
the $\SL_2(\BC)$ character variety (and the corresponding $A$-polynomial)
is a curve whose quotient modulo the involution
$\iota\colon (M,L)\mapsto \big(M^{-1},L^{-1}\big)$ is not a rational curve.
It was observed by Borot that his recent work with Eynard
\cite{BE} suggested a~mechanism (based on the {\em topological recursion})
that could explain at least a weak part of the modularity conjecture, namely
that the asymptotics of $\J^K(\ep)$ (as $\ep$ tends to zero through rational
numbers with bounded denominators), is always given by the same series
$\Phi^K_0(\ep)$ up to a~constant factor, not predicted by their model.
However, Borot could make this argument precise only in the case where the
space of holomorphic differentials of the corresponding spectral curve
was anti-invariant under the involution
$\iota\colon (M,L)\mapsto \big(M^{-1},L^{-1}\big)$. This condition is equivalent to
the statement to the rationality of the quotient of the spectral curve by
$\iota$. This led us to conduct a final experiment for the $6_1$ knot.
The question here was whether the modularity conjecture itself might fail,
or had to be modified in the context where the argument based on the work
of Borot--Eynard no longer applied. Fortunately, however, we found no anomalies.

To fix conventions, the $6_1$ knot is the closure of the braid word
{\tiny $[1, 2, 3, 2, -4, -1, -3, 2, -3, 4, -3, 2]$} where~$j$ (respectively,
$-j$) corresponds to the standard generator \smash{$s_j$} \big(respective, \smash{$s_j^{-1}$}\big)
of the braid group in 4 stands, and
with Jones polynomial $q^{-4} - q^{-3} + q^{-2} - 2 q^{-1} + 2 - q + q^2$.
The complex volume is given by
\[
\V(6_1)= -(2 R_1+R_2+R_3) - \frac{4}{3} \pi^2
= -3.0788629 \ldots + 3.1639632 \ldots {\rm i}
\]
where
\begin{align*}
R_1 &= R\bigl(-\xi^2\bigr) - \frac{1}{2} \pi{\rm i} \log\bigl(-\xi^2\bigr) + \pi{\rm i} \log\big(1+\xi^2\big), \\
R_2 &= R\big(1-\xi^3\big)- \pi{\rm i} \log\big(1-\xi^3\big) + \frac{1}{2} \pi{\rm i} \log\big(\xi^3\big), \qquad
R_3 = R(1-\xi) .
\end{align*}
The torsion, a prime in $F_{6_1}$ of norm 257 is given by
\[
\delta(6_1)= 1+\xi+4 \xi^2+\xi^3 .
\]

Here, the Kashaev invariant of the $6_1$ knot was not computed
using the formula given in~\cite[equation~(24)]{K97}, since the
latter is an $O(N^3)$ computation, but rather in $O(N)$ steps
using the recursion relation of the colored Jones polynomial.
The rather complicated inhomogeneous recursion has order $4$,
has maximal degree $(4,15,31)$ with respect to the shift variable, the
$q^n$ and the $q$ variables, and contains a total of 346 terms, which
can be found in~\cite{Ga:twistdata,GS}.
In contrast, the $A$-polynomial of the $6_1$ knot has maximal degree
$(4,8)$ with respect to the $(L,M)$ variables and contains 21 terms.
Due to the complexity of the recursion, we were forced to use precision~3000 in \texttt{pari} when $N=1000$. The first three coefficients of
$\Phi^{6_1}_0(h)$ were numerically computed at $\a=0$, and using the
prediction of~\cite{DG} and the notation of`\eqref{Phiprec2},
those algebraic numbers were identified as follows:
\begin{align*}
\Phi^{6_1}_0(h) =& {} \frac{1}{\sqrt{\delta}} \left( 1
+ \frac{194 \xi^3 - 331 \xi^2 + 207 \xi - 245}{2^3\cdot 3 \cdot \delta(6_1)^3} h\right.
\\
&\left.{} + \frac{-154734 \xi^3 - 34354 \xi^2 + 127399 \xi - 119864}{2^7 \cdot 3^2
\cdot \delta(6_1)^6}h^2
+ O\bigl(h^3\bigr) \right) .
\end{align*}

\subsection*{Acknowledgments}

The authors would like to thank the anonymous referees for their extraordinarily
careful reading of the manuscript and their detailed suggestions to improve
the exposition.

\newpage

\pdfbookmark[1]{References}{ref}
\LastPageEnding

\end{document}